\documentstyle[12pt,bezier]{report}

\renewcommand{\baselinestretch}{1.0}
\newcommand{\sect}[1]{\setcounter{equation}{0}\section{#1}}

\newcommand{\app}[1]{\setcounter{section}{0} 
\setcounter{equation}{0} \renewcommand{\thesection}{\Alph{section}}
\section{#1}}

\newcommand{\eq}{\begin{equation}}
\newcommand{\eqa}{\begin{eqnarray}}
\newcommand{\en}{\end{equation}}
\newcommand{\ena}{\end{eqnarray}}
\newcommand{\enn}{\nonumber \end{equation}}

\def\dia{{\scriptscriptstyle \Box}}
\def\cvd{{\vskip -0.49cm\rightline{$\Box\!\Box\!\Box$}}\sk}

\newcommand{\spz}{\hspace{0.7cm}}

\def\spinst#1#2{{#1\brack#2}}
\def\sk{\vskip .4cm}
\def\noi{\noindent}
\def\om{\omega}
\def\al{\alpha}
\def\be{\beta}
\def\ga{\gamma}
\def\Ga{\Gamma}
\def\de{\delta}
\def\del{\delta}
\def\la{\lambda}
\def\La{\Lambda}
\def\lam{{1 \over \lambda}}
\def\alb{\bar{\alpha}}
\def\beb{\bar{\beta}}
\def\gab{\bar{\gamma}}
\def\deb{\bar{\delta}}
\def\alp{{\alpha}^{\prime}}
\def\bep{{\beta}^{\prime}}
\def\gap{{\gamma}^{\prime}}
\def\dep{{\delta}^{\prime}}
\def\rhop{{\rho}^{\prime}}
\def\taup{{\tau}^{\prime}}
\def\rhopp{\rho ''}
\def\thetap{{\theta}^{\prime}}
\def\imezzi{{i\over 2}}
\def\unquarto{{1 \over 4}}
\def\unmezzo{{1 \over 2}}
\def\epsi{\varepsilon}
\def\we{\wedge}
\def\th{\theta}
\def\vart{\vartheta}
\def\de{\delta}
\def\cony{i_{\de {\vec y}}}
\def\Liey{l_{\de {\vec y}}}
\def\tv{{\vec t}}
\def\Gt{{\tilde G}}
\def\deyv{\vec {\de y}}
\def\part{\partial}
\def\pdxp{{\partial \over {\partial x^+}}}
\def\pdxm{{\partial \over {\partial x^-}}}
\def\pdxi{{\partial \over {\partial x^i}}}
\def\pdy#1{{\partial \over {\partial y^{#1}}}}
\def\pdx#1{{\partial \over {\partial x^{#1}}}}
\def\pdyx#1{{\partial \over {\partial (yx)^{#1}}}}

\def\qP{q-Poincar\'e~}
\def\R#1#2{ R^{#1}_{~~~#2} }
\def\Rp#1#2{ (R^+)^{#1}_{~~~#2} }
\def\Rm#1#2{ (R^-)^{#1}_{~~~#2} }
\def\Rinv#1#2{ (R^{-1})^{#1}_{~~~#2} }
\def\Rpm#1#2{(R^{\pm})^{#1}_{~~~#2} }
\def\Rpminv#1#2{((R^{\pm})^{-1})^{#1}_{~~~#2} }
\def\RRpm{R^{\pm}}
\def\RRp{R^{+}}
\def\RRm{R^{-}}
\def\Rhat#1#2{ \Rh^{#1}_{~~~#2} }
\def\Rhats#1#2{ \Rhs^{#1}_{~~~#2} }
\def\Rhatinv#1#2{ (\Rh^{-1})^{#1}_{~~~#2} }
\def\Z#1#2{ Z^{#1}_{~~~#2} }
\def\Rt#1{ {\hat R}_{#1} }
\def\Rhs{{\hat R}}
\def\Rh{{\Lambda}}
\def\ff#1#2#3{f_{#1~~~#3}^{~#2}}
\def\MM#1#2#3{M^{#1~~~#3}_{~#2}}
\def\cchi#1#2{\chi^{#1}_{~#2}}
\def\ome#1#2{\om_{#1}^{~#2}}
\def\RRhat#1#2#3#4#5#6#7#8{\Rh^{~#2~#4}_{#1~#3}|^{#5~#7}_{~#6~#8}}
\def\RRhatinv#1#2#3#4#5#6#7#8{(\Rh^{-1})^
{~#2~#4}_{#1~#3}|^{#5~#7}_{~#6~#8}}
\def\U#1#2#3#4#5#6#7#8{U^{~#2~#4}_{#1~#3}|^{#5~#7}_{~#6~#8}}
\def\Cb{{\bf C}}
\def\Rb{{\bf R}}
\def\CC#1#2#3#4#5#6{\Cb_{~#2~#4}^{#1~#3}|_{#5}^{~#6}}
\def\cc#1#2#3#4#5#6{C_{~#2~#4}^{#1~#3}|_{#5}^{~#6}}

\def\C#1#2{ {\bf C}_{#1}^{~~~#2} }
\def\c#1#2{ C_{#1}^{~~~#2} }
\def\q#1{   {{q^{#1} - q^{-#1}} \over {q^{\unmezzo}-q^{-\unmezzo}}}  } 
\def\Dmat#1#2{D^{#1}_{~#2}}
\def\Dmatinv#1#2{(D^{-1})^{#1}_{~#2}}
\def\DR{{}_{\scriptsize \Ga}\Delta}
\def\DL{\Delta_{\scriptstyle \Ga}}

\def\DD{{}_{\scriptsize \Xi}\Delta}
\def\DS{\Delta_{\scriptstyle \Xi}}
\def\FF#1#2{O_{#1}{}^{#2} }
\def\N#1#2{N^{#1}{}_{#2}}
\def\err{\ell^{\cal R}}

\def\f#1#2{ f^{#1}_{~~#2} }
\def\F#1#2{ F^{#1}_{~~#2} }
\def\T#1#2{ T^{#1}_{~~#2} }
\def\Ti#1#2{ (T^{-1})^{#1}_{~~#2} }
\def\Tp#1#2{ (T^{\prime})^{#1}_{~~#2} }
\def\TP{ T^{\prime} }
\def\M#1#2{ M_{#1}^{~#2} }
\def\qm{q^{-1}}
\def\um{u^{-1}}
\def\vm{v^{-1}}
\def\xm{x^{-}}
\def\xp{x^{+}}
\def\fm{f_-}
\def\fp{f_+}
\def\fn{f_0}
\def\D{\Delta}
\def\Mat#1#2#3#4#5#6#7#8#9{\left( \matrix{
     #1 & #2 & #3 \cr
     #4 & #5 & #6 \cr
     #7 & #8 & #9 \cr
   }\right) }   
\def\Ap{A^{\prime}}
\def\Dp{\Delta^{\prime}}
\def\Ip{I^{\prime}}
\def\ep{\epsi^{\prime}}
\def\kp{\kappa^{\prime}}
\def\kpm{\kappa^{\prime -1}}
\def\km{\kappa^{-1}}
\def\gp{g^{\prime}}
\def\qone{$q \rightarrow 1~$}
\def\Fmn{F_{\mu\nu}}
\def\Am{A_{\mu}}
\def\An{A_{\nu}}
\def\dm{\part_{\mu}}
\def\dn{\part_{\nu}}
\def\Ana{A_{\nu]}}
\def\Bna{B_{\nu]}}
\def\Zna{Z_{\nu]}}
\def\dma{\part_{[\mu}}
\def\qsu{$[SU(2) \times U(1)]_q~$}
\def\suq{$SU_q(2)~$}
\def\su{$SU(2) \times U(1)~$}
\def\gij{g_{ij}}
\def\L{{\cal L}}
\def\Un{U}
\def\Uc{{\cal U}}
\def\ad{{a\!d}}
\def\LL{L^*}
\def\ll#1{L^*_{#1}}
\def\RR{R^*}
\def\rr#1{R^*_{#1}}

\def\Lpm#1#2{{L^{\pm}}^{#1}_{~~#2}}
\def\LLpm{L^{\pm}}
\def\LLp{L^{+}}
\def\LLm{L^{-}}
\def\Lp#1#2{{L^{+}}^{#1}_{~~#2}}
\def\Lm#1#2{{L^{-}}^{#1}_{~~#2}}

\def\gu{g_{U(1)}}
\def\gsu{g_{SU(2)}}
\def\tg{ {\rm tg} }
\def\Fun{$Fun(G)~$}
\def\invG{{}_{{\rm inv}}\Ga}
\def\Ginv{{}_{\rm inv}\Ga}
\def\qonelim{\stackrel{q \rightarrow 1}{\longrightarrow}}
\def\viel#1#2{e^{#1}_{~~{#2}}}
\def\invXi{{}_{\rm{inv}}\Xi}

\def\Cbold{\mbox{\scriptsize \bf C}}
\def\Rbold{\mbox{\scriptsize \bf R}}
\def\sma#1{\mbox{\footnotesize #1}}
\def\le{\langle}
\def\re{\rangle}

\begin{document}

\begin{titlepage}
\begin{center}
{\sc{scuola normale superiore di pisa}}
\sk
\sk
\rightline{January 1998}
\sk\sk\sk\sk
\sk
TESI DI PERFEZIONAMENTO

\vskip .5in

{\large \bf On the Geometry of Inhomogeneous Quantum Groups}
\vskip .5in

Paolo Aschieri
\sk
\sk

{\em Scuola Normale Superiore\\
Piazza dei Cavalieri 7, 56100 Pisa, Italy\\
{\it and}\\
Theoretical Physics Group\\
Lawrence Berkeley Laboratory,\\
University of California\\
Berkeley, California 94720, USA}

\vskip .5in

\vskip .5in
\end{center}
\sk
\sk
\sk
\sk
\sk
\sk
\sk\sk
\sk
\sk\sk
\end{titlepage}
\pagebreak

\renewcommand{\thepage}{\roman{page}}
\setcounter{page}{-1}

\clearpage{{}}
\vskip 26cm
\centerline{}
\newpage
\centerline{\bf Abstract}
\sk
We give a pedagogical introduction to the differential calculus 
on quantum groups by stressing at all stages the connection 
with the classical case. 
We further analize the relation between differential calculus and 
quantum Lie algebra of left (right) invariant vectorfields. 
Equivalent definitions of bicovariant differential calculus are studied
and their geometrical interpretation is explained. 
From these data we construct and analize the space of vectorfields,
and naturally introduce a contraction operator and a Lie derivative,
their properties are discussed.

After a review of the geometry of the multiparametric deformation of the
linear  group $GL_{q,r}(N)$ and its real forms,
we then construct the multiparametric linear
inhomogeneous quantum group $IGL_{q,r}(N)$ as a projection from
$GL_{q,r}(N+1)$, or equivalently, as a quotient of
$GL_{q,r}(N+1)$ with respect to a suitable
Hopf algebra ideal. 
The semidirect product structure of   
$IGL_{q,r}(N)$ given by the $GL_{q,r}(N)$ quantum subgroup 
times translations is analized.
A bicovariant 
differential calculus
on $IGL_{q,r}(N)$ is explicitly obtained as a 
projection from the one on $GL_{q,r}(N+1)$. 
The universal enveloping algebra of  $IGL_{q,r}(N)$ and its 
$R$-matrix formulation are 
constructed along the same lines.
This quotient procedure unifies in a single structure the quantum 
plane coordinates and the $q$-group matrix elements $\T{a}{b}$, 
and allows to deduce without effort the differential
calculus on the $q$-plane $IGL_{q,r}(N) / GL_{q,r}(N)$.
The general theory is illustrated on the example
of $IGL_{q,r}(2)$. 

We proceed similarly in the orthogonal and symplectic case.
The inhomogeneous multiparametric $q$-groups
of the $B_n,C_n,D_n$ series are found 
by means of
a projection from $B_{n+1},C_{n+1},D_{n+1}$.
A matrix formulation is given
in terms of
the $R$-matrix of  $B_{n+1},C_{n+1},D_{n+1}$, and
real forms are discussed: in particular we obtain the $q$-groups
$ISO_{q,r}(n+1,n-1)$, including the
quantum Poincar\'e group.
The universal enveloping algebras of the 
multiparametric $B_{n+1},C_{n+1},D_{n+1}$
$q$-groups are studied, they include 
as Hopf subalgebras the 
universal enveloping algebras of  the  inhomogeneous   $B_n,C_n,D_n$.
Bicovariant calculi on the minimal multiparametric deformations (twists) of 
these inhomogeneous
groups
are similarly found by means of a projection
from the bicovariant calculus on
$B_{n+1}$, $C_{n+1}$, $D_{n+1}$.
In particular we obtain the bicovariant calculus
on a dilatation-free minimal deformation of the Poincar\'e group
$ISO_{q}(3,1)$.
\nopagebreak
Then we construct differential calculi on multiparametric quantum
orthogonal planes in any dimension $N$. These calculi are
bicovariant under the action of the full inhomogeneous
multiparametric quantum group $ISO_{q,r}(N)$, and do contain
dilatations. We find a canonical group-geometric 
procedure to restrict these calculi
on the $q$-plane and expressed them in terms of  coordinates $x^a$,
differentials $dx^a$ and partial derivatives $\partial_a$ without
the need of dilatations,  thus generalizing
known results to the multiparametric case.\nopagebreak
Real forms are studied and in particular we obtain the quantum Minkowski space 
$ISO_{q,r}(3,1)/SO_{q,r}(3,1)$.\nopagebreak
The conjugated partial derivatives $\partial_a^*$ can
be expressed as linear combinations of the $\partial_a$.
This allows a deformation of the phase-space with
hermitian operators $x^a$ and $p_a$.  


\renewcommand{\thepage}{\roman{page}}
\setcounter{page}{-1}

\clearpage{{}}
\vskip 26cm
\centerline{}
\newpage

\renewcommand{\thepage}{\roman{page}}
\setcounter{page}{-1}

\vskip 30cm
\vskip 2in
\sk
\sk

\sk
\sk

\renewcommand{\thepage}{\roman{page}}
\setcounter{page}{0}

\newpage

{\bf{Acknowledgements}}
\sk
\sk
\sk
\sk
\sk
\sk
\sk
\sk
\sk
\sk
\sk
\sk
\sk
\sk
\sk
\sk
\sk
\sk
\sk
\sk
I am glad to thank Professor Leonardo Castellani for his guidance, 
for his support and encouragement I have always received, and for
sharing the many moments of research. All this in a
 stimulating and friendly 
atmosphere.

I am grateful to Professor Corrado De Concini for 
his advice and for the many explanations that have enhanced
my mathematical understanding of the subject.
There are many colleagues I am indebted with, in particular Luigi Pilo
and  Peter Schupp.

About this year in Berkeley I would  
like to thank Professor Bruno Zumino for his advice, his  support
and the stimulating environment I experience.
I also thank Professor Nicolai Reshetikhin for his suggestions.
I am pleased to aknowledge the  many  fruitful 
discussions with Bogdan Morariu, Harold Steinacker and
Gaetano Fiore.  

This thesis has been possible also because of the support and affection of
my parents. To them my gratitude.
\sk
\sk
\sk
\sk
\nopagebreak
\noi This work in 1997 has been supported by an exchange fellowship
Scuola  Normale 
Superiore -- University of California and by a research fellowship of
Fondazione Angelo Della Riccia. 
It has been accomplished through
 the Director,
Office of Energy Research, Office af High Energy  and Nuclear  
Physics,
Division of High Energy Physics
of the U.S. Department of Energy under Contract DE-AC03-76SF00098
and by the National Science Foundation under grant PHY-95-14797.

\tableofcontents

\newpage
\renewcommand{\thepage}{\arabic{page}}
\setcounter{page}{1}

\chapter*{Introduction}
\addcontentsline{toc}{part}{Introduction}

Quantum Groups are particular deformations of Lie Groups.
They are algebraic structures $G_q$ depending on one (or more) 
continuous parameter $q$. When $q=1$ we have a standard Lie group.

A method to modify a group could be to continuously  deform 
the structure constants of the relative Lie algebra. This kind
of deformations is uninteresting because it can be shown that for a
semisimple Lie group, with a suitable redefinition of the  generators 
it is always possible to recover the undeformed structure constants.  
At most, for special values of the deformation parameter, we get a
contraction of the original group. (For example the Galilei group is
a contraction of the Lorentz group, and the Poincar\'e 
group is a contraction
of  the de Sitter group).

Quantum deformations, on the contrary, allow to obtain new structures not 
isomorphic to the preceding ones. For example the quantum $SU(2)$ group is
given by 
\eq
[J^+,J^-]={q^{2J^0}-q^{-2J^0}\over q-q^{-1}}~~~,~~~~
[J^0,J^{\pm}]=\pm J^{\pm}\label{uno}
\en
Here we clearly see that the commutator depends continuously
on $q$, moreover we loose the concept of structure constant 
(in sections 2.1 and 2.3, we will see a generalization of this concept).
We however still
have a rich structure (Hopf algebra structure), as rich as the $SU(2)$ one 
\cite{Drinfeld}, \cite{Jimbo}.  In the $q\rightarrow 1$ limit we recover the 
$SU(2)$ algebra $[J^+,J^-]=2J^0,~[J^0,J^{\pm}]=\pm J^{\pm}$.
Relations (\ref{uno}) define the universal enveloping algebra of $SU(2)$,
i.e. the space given by all (formal) polynomials 
$\sum \lambda_{mnp}(J^0)^m(J^+)^n(J^-)^p$ where the 
ordering is always possible thanks to (\ref{uno}).
\sk
We know that a new physical theory can always be seen as a generalization 
of the previous existing one, in the sense that this is recovered as a 
limiting case in which some parameters of the new theory become negligible,
e.g. special relativity where the Galilei symmetry group
is replaced by the Lorentz one. Here the deformation parameter is $c$,
and for $c\rightarrow \infty$ we recover Galiley relativity.
The study of continuous deformation of Lie groups is therefore
physically interesting, the symmetry group underlying many (particle) 
physics theories being Lie Groups.
\sk
The deformations we will deal with  are  in the context of noncommutative 
geometry. Both the  group coordinates and the 
space coordinates on which the group acts, 
consist of non-commuting elements.
It is of interest to apply these  rich mathematical structures to the 
study of spacetime physics at Planck scale. Indeed at this scale 
due to the gravitational forces, Gedanken experiments show that it is not
possible to probe  spacetime 
structure\footnote{For example, in 
relativistic quantum mechanics the position of a particle can
be detected with a precision at most of the order of its 
Compton wave length $\lambda_C=\hbar/mc$.
Probing spacetime at infinitesimal distances implies an extremely 
heavy 
particle that in turn curves spacetime itself. 
When $\lambda_C$ is of the order of the Planck length, the spacetime
curvature 
radius due to the particle has the same order of magnitude and the 
attempt to measure spacetime structure beyond Planck scale fails.}, 
the description of spacetime  as a smooth manifold becomes 
just a mathematical assumption. 
We can relax it and conceive a more general  
noncommutative structure. This is intresting because
in a noncommutative spacetime uncertainty relations 
and discretization naturally
arise. In this way one could
incorporate the impossibility of an operational
definition  of spacetime structure beyond the Planck scale 
(due to gravity) in the  noncommutative geometric structure of
spacetime itself. A dynamic feature of spacetime would be
incorporated at a kinematical level.
This could be a fertile setting
for the study of quantum gravity and field theory at Planck scale.

The easiest example of noncommuting space is given by the quantum plane,
\cite{WZ,Pusz,zumi} where the coordinates $x$ and 
$y$ satisfy $xy=qyx$ where  $q$ is a complex parameter.
The $q$-plane and similar 
higher dimensional $q$-spaces admit a differential 
structure where the derivative
operators are finite difference operators. This resembles
lattice like structures.
At the phase-space level the quantum plane relations imply
expressions of the type $x\partial_x - q\partial_x x= 1$ 
that lead to selfadjoint position and momentum operators 
with discrete heigenvalues: we obtain again a lattice structure \cite{Lorek}. 
In this context $q$ may play the role of short distance 
regularization parameter 
preserving the $q$-symmetries, see also \cite{Majid}. 
In the particular deformations where 
$q$ is a dimensionful parameter, one can also try to relate 
this parameter to the  Planck length  \cite{Camelia}.

Notice that a minimal uncertainty relation in position measurement is also 
in agreement with string theory models \cite{Veneziano}.  
Moreover, non-perturbative attempts to 
describe string 
theories have shown that a noncommutative structure of spacetime emerges
\cite{Banks}, noncommutative geometry can be
the correct geometrical framework for the description of such theories 
\cite{Douglas}.
\sk
In our study of inhomogeneous $q$-groups we will consider examples of 
noncommuting quantum
planes (in 2 and higher dimension) that have quantum symmetry groups. 
This is not the only approach to a 
possible model of noncommutative spacetime.
For example one can consider a
noncommutative Minkowski spacetime which is covariant
under the classical Poincar\'e group, see \cite{Madore}, and 
\cite{Doplicher} where also
first attempts to construct
quantum field theories on noncommutative spaces have shown that in
some cases these theories are equivalent to a non local quantum field
theory on ordinary space.
We also mention a related approach \cite{CFF} that, in the spirit 
of \cite{Lott},
uses noncommutative geometric methods to 
study a Kaluza-Klein gravity theory with a 
discrete two point internal space. 
\sk
In this thesis we generalize to Hopf algebras basic  structures of
differential geometry  on commutative manifolds, we then 
examine examples of inhomogeneous quantum groups and 
consider their differential calculus.

In the first chapters we introduce the concept of quantum group,
and, stressing at all stages the connection with the classical
case ($q\rightarrow 1$ limit), we develop the differential calculus 
on quantum groups, first studied in the seminal paper by Woronowicz
\cite{Wor}.
We then examine in detail the quantum Lie algebra of 
left-invariant vectorfields (see Section 2.3). 
All the properties of the differential calculus can be derived from
intuitive properties of the $q$-Lie algebra, in this way we emphasize the 
space of vectorfields that is more fundamental, for physical applications, 
than the space of $1$-forms. 
For example, this $q$-Lie algebras are a good starting point for the 
formulation of gauge
theories based on deformed Lie groups \cite{gauge}.
Next, a Lie derivative 
and a contraction operator (inner derivative) are found
and, for left-invariant vectorfields, we prove the Cartan identity 
$\ell_{\vec t}=i_{\vec t}d+di_{\vec t}\,$ \cite{AC,SWZ1,SWZ3,PaoloPeter}. 
These are basic tools in differential
geometry, and are of interest in a geometric formulation
of Einstein gravity. For example 
the invariance of Einstein  action under diffeomorphisms can be 
expressed by 
$\ell_{\vec t}\int{\cal L } d^4x=0~,$
this relation 
leads to the covariant conservation of the matter energy-momentum
tensor (if torsion vanishes).
Also, in the soft group manifold approach to gravity theories \cite{CDF},
the Cartan-Maurer 
equation, the Lie derivative and the contraction operator for 
the $q$-Poincar\'e group are fundamental to formulate a geometric 
definition of curvature, covariant 
derivative and Lorentz gauge transformation. There the 
curved general relativity space time (with the Lorentz gauge group) 
is obtained as the coset space Poincar\'e/Lorentz where the rigid
Poincar\'e group structure has been softened allowing for a curvature two
form term in the Cartan-Maurer equations.
For a first example of this construction in the case of 
a minimal $q$-deformation (twist) of the Poincar\'e group see \cite{Cas2}.
\sk
The second part of this thesis  deals with the specific study of
deformations of the inhomogeneous general linear group $GL(N)$ and of 
the inhomogeneous orthogonal and symplectyc groups.
Contrary to the case of semisimple Lie groups [where there is a canonical
Poisson 
(symplectic) structure that can be quantized to give the quantum group]  
there is not a canonical deformation procedure for these groups;
providing examples of inhomogeneous deformation is therefore per se 
interesting. The $GL(N)$ and $SL(N)$ cases are easier to
consider than the orthogonal and symplectic 
ones and the basic structures of inhomogeneous quantum groups are first
found in this cases and later shown for the $B_n,C_n,D_n$
series as well.
There are many studies on
inhomogeneous quantum groups \cite{inhom,Rinhom2,Rinhom1,Rinhom3}
and in particular on deformed Poincar\'e groups \cite{inhompoi}. 
The analysis we
present
\cite{IGLAschieri}, \cite{inson}  is based
on a quotient procedure, we find deformations of the inhomogeneous 
version of the
$A_n,B_n,C_n,D_n$ groups  via a quotient from the $q$-deformed 
homogeneous $A_{n+1},B_{n+1},C_{n+1},D_{n+1}$ groups. 
An $R$-matrix formulation
is thus provided.
The universal enveloping algebra of these groups is also studied as
well as their semidirect product structure of their homogeneous subgroups
times translations.
Then we apply the general theory of the differential calculus on $q$-groups 
to these specific cases. Differential calculi on these inhomogeneous $q$-groups
are found  using again the quotient structure with respect to  
$A_{n+1},B_{n+1},C_{n+1},D_{n+1}$.

In particular we obtain a deformation of the
Poincar\'e group, of its $q$-Lie algebra and differential calculus. 
It turns out that the differential 
calculus on inhomogeneous orthogonal 
and symplectic quantum groups, contrary to the linear case,  
cannot be constructed for general values of the
deformation parameters. It exists only for minimal deformations (twists), these
are the same noncommutative deformations that do not require the presence of
a  dilatation in the fully $q$-deformed inhomogeneous structure. 
Our analysis includes the first example of bicovariant differential calculus 
on  quantum Poincar\'e group with deformed Lorentz subgroup. For bicovariant 
calculi on other deformations of the Poincar\'e group see 
\cite{Maslanka}.
The study of the differential calculi on orthogonal
groups discussed in  Chapter 4,
leads to a good candidate for the 
differential calculus on the $q$-orthogonal planes and in particular 
on the $q$-Minkowski plane, with no restriction on the
deformation parameters.
Using the powerful result:  
$q$-Minkowski =
$q$-Poincar\'e/$q$-Lorentz,
we  derive canonically the $q$-Minkowski geometry
from the $q$-Poincar\'e geometry \cite{ortpplan}. 
In this way 
we rederive
the known results of
\cite{WZ,qplane,Fiore1} using
 the broader setting of the differential calculus on quantum $ISO(N)$.
A detailed analysis of the reality condition on the quantum Minkowski plane
is possible since on the quantum $ISO(3,1)$
differential calculus  there is a canonical definition of
$*$-conjugation. This operation is linear  
and one can canonically obtain real coordinates and momenta
and a  $q$-version of
the Heisenberg $x,p$ commutation relations. 
An interesting issue, in the spirit of  \cite{Lorek} is then  
the analysis of the representations in Hilbert space 
of this algebra in order to study the admissible (discrete) values 
of  momentum and position of particle states.

\chapter{Quantum Groups}

\sect{Hopf structures in ordinary Lie groups and Lie algebras}

Let us begin by considering \Fun, the set of differentiable functions 
from a Lie group $G$ into the complex numbers $\Cb$. \Fun  is an 
algebra with the usual pointwise sum and 
product $(f+h)(g)=f(g)+h(g),~(f \cdot h)=f(g) h(g),
~(\lambda f)(g)=\lambda f(g)$, for $f,h \in Fun(G),~g \in G,~\lambda \in 
\Cb$. The unit of this algebra is $I$, defined by $I(g)=1,~\forall g \in 
G$. 

Using the group structure of $G$, we can introduce on \Fun three other 
linear mappings, the coproduct $\D$, the counit $\epsi$, and 
the coinverse
(or antipode) $\kappa$:
\begin{eqnarray}
\D (f)(g,\gp) &\equiv& f(g\gp),~~~\D:Fun(G) \rightarrow Fun(G)\otimes 
  Fun(G) \label{cop} \\
\epsi(f) &\equiv& f(1_G),~~~~~~\epsi:Fun(G) \rightarrow \Cb \label{cou}\\
(\kappa f)(g) &\equiv& f(g^{-1}),~~~\kappa:Fun(G) \rightarrow Fun(G)
  \label{coi}
\end{eqnarray}
 
\noi where $1_{{}_G}$ is the unit of $G$. It 
is not difficult to verify the following properties of the co-structures:
\begin{eqnarray}
 & & (id \otimes \D)\D=(\D\otimes id)\D~~~({\rm coassociativity~of~\D}) 
\label{prop1}\\
 & & (id\otimes \epsi)\D(a)=(\epsi\otimes id)\D(a)=a \label{prop2}\\
 & & m(\kappa\otimes id)\D(a)=m(id\otimes\kappa)\D(a)=\epsi(a) I 
\label{prop3}
\end{eqnarray}

\noi and
\begin{eqnarray}
 & & \D(ab)=\D(a)\D(b),~~~\D(I)=I\otimes I \label{prop4}\\
 & & \epsi(ab)~=\epsi(a) \epsi(b),~~~~~~\epsi(I)=1 \label{prop5} \\
 & & \kappa(ab) \; =\kappa(b)\kappa(a),~~~~~\kappa(I)=I \label{prop6}
\end{eqnarray}

\noi where $a,b \in A=Fun(G)$ and $m$ is the multiplication map
$m(a\otimes b)\equiv ab$. The product in $\D(a)\D(b)$ is the product in 
$A\otimes A$: $(a\otimes b)(c\otimes d)=ab\otimes cd$.

In general a coproduct can be expanded on
$A \otimes A$ as:
\eq
\D(a)=\sum_i a_1^i \otimes a_2^i \equiv a_1 \otimes a_2, \label{not1}
\en
\noi where $a_1^i, a_2^i \in A$ and $a_1 \otimes a_2$ is a 
shorthand notation we will often use in the sequel. For example for
$A=Fun(G)$ we have:
\eq
\D(f)(g,\gp)=(f_1 \otimes f_2)(g,\gp)=f_1(g)f_2(\gp)=f(g\gp). 
    \label{not2}
\en
\noi Using (\ref{not2}), the proof of (\ref{prop1})-(\ref{prop3}) is 
immediate.
We will also use the following notation: $\D^2(a)\equiv(\D\otimes id)\D(a)=
(id\otimes \D)\D(a)=a_1\otimes a_2\otimes a_3$, more in general 
$\D^n(a)=a_1\otimes a_2\otimes ,... a_{n+1}$.
\sk
An algebra $A$ endowed with the homomorphisms $\D:A \rightarrow A 
\otimes A$ and $\epsi: A \rightarrow \Cb$, and the antimorphism  
$\kappa: A\rightarrow A$ satisfying the properties 
(\ref{prop1})-(\ref{prop6}) 
is a {\sl Hopf algebra}. Thus \Fun is a Hopf algebra.\footnote
{To be precise, \Fun is a Hopf algebra when $Fun(G \times G)$ can be 
identified with $Fun(G) \otimes Fun(G)$, since only then can one define 
a coproduct as in (\ref{cop}).} 
Note that the properties (\ref{prop1})-(\ref{prop6}) 
imply the relations:
\eqa                                    
& &\!\!\!\!\!\!\! \D(\kappa (a))=\kappa (a_2) \otimes \kappa (a_1) 
~~,~~~~\kappa(a){}_1\otimes\kappa(a){}_2\, ,\ldots\kappa(a){}_n=
\kappa(a_n)\otimes\kappa(a_{n-1})\, ,\ldots\kappa(a_1)\nonumber\\
& &\label{Dka}\\
& &\!\!\!\!\!\!\! \epsi (\kappa (a))=\epsi (a). \label{epsika}
\ena
\sk
Consider now the algebra $A$ of polynomials in the matrix
elements $\T{a}{b}$ of the fundamental representation of $G$, i.e. the
algebra $A$  generated by the $\T{a}{b}$. 

It is clear that $A \subset Fun(G)$, since $\T{a}{b} (g)$ are 
functions on $G$. In fact $A$ is dense in $Fun(G)$ 
[the reason is that the matrix elements
of all {\sl finite} irreducible dimensional representations of $G$ can 
be constructed out of appropriate products of $\T{a}{b} (g)$; 
these products span a dense subset in $Fun(G)$ because the matrix
elements of all irreducible representations of $G$ form a complete
basis of $Fun(G)$], therefore, a suitable 
completion $\hat A$ of
$A$ is $Fun(G)$ : $\hat{A}=Fun(G)$.  
In the following we will drop the 
hat and we will not be concerned about these topological aspects.
The group manifold $G$ can be completely characterized by $Fun(G)$, the
co-structures on \Fun carrying the information about the group 
structure of $G$. Thus a classical Lie group can be ``defined" as 
the algebra $A$  generated by the (commuting) matrix 
elements $\T{a}{b}$ of the fundamental representation of $G$, seen 
as functions on $G$.
This definition admits noncommutative generalizations, i.e. the
quantum groups discussed in the next section.
\sk
Using the elements $\T{a}{b}$ we can write an 
explicit formula for the expansion (\ref{not1}) or (\ref{not2}): indeed 
(\ref{cop}) becomes
\eq
\D(\T{a}{b})(g,\gp)=\T{a}{b} (g\gp)=\T{a}{c}(g) \T{c}{b}(\gp), 
\en
\noi since $T$ is a matrix representation of $G$. Therefore:
\eq
\D(\T{a}{b})=\T{a}{c} \otimes \T{c}{b}. \label{copT}
\en
Moreover, using (\ref{cou}) and (\ref{coi}), one finds:
\begin{eqnarray}
 & & \epsi(\T{a}{b})=\de^a_b \label{couT}\\
 & & \kappa(\T{a}{b})=\Ti{a}{b}. \label{coiT}
\end{eqnarray}
Thus the algebra $A=Fun(G)$ of polynomials in the elements $\T{a}{b}$ is 
a Hopf algebra with co-structures defined by 
(\ref{copT})-(\ref{coiT}) and (\ref{prop4})-(\ref{prop6}).
\sk
Another example of Hopf algebra is given by any ordinary Lie algebra 
$\mbox{\sl g}$, or more precisely by the universal enveloping algebra 
$U(\mbox{\sl g})$
of a Lie algebra $\mbox{\sl g}$, i.e. 
(by the Poincar\'e-Birkhoff-Witt theorem) the 
algebra, with unit $I$, of polynomials in the generators $\chi_i$ modulo 
the commutation relations
\eq
[\chi_i,\chi_j]=\c{ij}{k} \chi_k. \label{clcomm}
\en
Here we define the co-structures as:
\begin{eqnarray}
 & & \D'(\chi_i)=\chi_i \otimes I + I \otimes \chi_i~~~\D'(I)=I\otimes I  
  \label{copL}\\
 & & \epsi '(\chi_i)=0~~~~~~~~~~~~~~~~~~~~~~\epsi '(I)=1 \label{couL}\\
 & & \kappa'(\chi_i)=-\chi_i~~~~~~~~~~~~~~~~~~\kappa '(I)=I \label{coiL}
\end{eqnarray}
\noi The reader can check that (\ref{prop1})-(\ref{prop3}) are 
satisfied.

\sect{Quantum groups. The example of $GL_q(2)$}

Quantum groups can be introduced
as noncommutative deformations of the 
algebra $A=Fun(G)$ of the previous section [more precisely as 
noncommutative Hopf algebras obtained by continuous deformations of 
the Hopf algebra $A=Fun(G)$]. The term quantum stems for the fact that
they are obtained quantizing a Poisson (symplectic) structure of the algebra
$Fun(G)$ \cite{Drinfeld}. 
Here, following \cite{FRT} (see also \cite{T}),  we will consider quantum 
groups defined as the associative algebras $A$ freely 
generated by non-commuting matrix entries 
$\T{a}{b}$ satisfying the relation
\eq
\R{ab}{ef} \T{e}{c} \T{f}{d} = \T{b}{f} \T{a}{e} \R{ef}{cd} \label{RTT}
\en
\noi and some other conditions depending on which classical group 
we are deforming (see later). The matrix $R$ controls the 
non-commutativity of the $\T{a}{b}$, and its elements depend 
continuously on a (in general complex) parameter $q$, or even a set of 
parameters. For $q\rightarrow 1$, the so-called ``classical limit", we
have 
\eq 
\R{ab}{cd} \qonelim \de^a_c \de^b_d,  \label{limR}
\en
\noi i.e. the matrix entries $\T{a}{b}$ commute for $q=1$, and one 
recovers the ordinary $Fun(G)$. 

The associativity of $A$ leads to a consistency condition on the $R$
matrix, the quantum Yang--Baxter equation:
\eq
\R{a_1b_1}{a_2b_2} \R{a_2c_1}{a_3c_2} \R{b_2c_2}{b_3c_3}=
\R{b_1c_1}{b_2c_2} \R{a_1c_2}{a_2c_3} \R{a_2b_2}{a_3b_3}. \label{YB}
\en
For simplicity we rewrite the ``RTT" equation (\ref{RTT}) and the 
quantum Yang--Baxter equation as
\eq
R_{12} T_1 T_2 = T_2 T_1 R_{12} \label{rtt}
\en
\eq
R_{12} R_{13} R_{23}=R_{23} R_{13} R_{12}, \label{yb}
\en
\noi where the subscripts 1, 2 and 3 refer to different couples of 
indices. Thus $T_1$ indicates the matrix $\T{a}{b}$, $T_1 T_1$ 
indicates $\T{a}{c} \T{c}{b}$, $R_{12} T_2$ indicates $\R{ab}{cd} 
\T{d}{e}$ and so on, repeated subscripts meaning matrix 
multiplication. The quantum Yang--Baxter equation (\ref{yb}) 
is a condition sufficient 
for the consistency of the RTT equation (\ref{rtt}). Indeed 
the product of three distinct elements $\T{a}{b}$, $\T{c}{d}$ 
and $\T{e}{f}$, indicated by $T_1T_2T_3$, can be reordered as 
$T_3T_2T_1$ via two differents paths:
\eq
 T_1T_2T_3 \begin{array}{c} \nearrow\\ \searrow \end{array}
           \begin{array}{c} T_1T_3T_2 \rightarrow T_3T_1T_2 \\ {} \\
                            T_2T_1T_3 \rightarrow T_2T_3T_1 \end{array}
           \begin{array}{c} \searrow\\ \nearrow \end{array}
 T_3T_2T_1          
\en
\noi by repeated use of the RTT equation. The relation (\ref{yb})
ensures that the two paths lead to the same result. 
\sk
The algebra $A$ (``the quantum group") is a noncommutative Hopf algebra 
whose co-structures are the same of those defined for the commutative 
Hopf algebra \Fun of the previous section, eqs. 
(\ref{copT})-(\ref{coiT}), (\ref{prop4})-(\ref{prop6}).
\sk
Let us give the example of $SL_q(2)\equiv Fun_q(SL(2))$, the algebra freely 
generated by 
the elements $\al,\be,\ga$ and $\de$ of the $2 \times 2$ matrix
\eq
\T{a}{b}= \left( \begin{array}{cc} \al & \be \\ \ga & \de \end{array} 
    \right)   \label{Tmatrix}
\en
\noi satisfying the commutations
\eqa
\al\be=q\be\al,~~\al\ga=q\ga\al,~~\be\de=q\de\be,~~\ga\de=q\de\ga
\nonumber \\
\be\ga=\ga\be,~~\al\de-\de\al=(q-\qm)\be\ga,~~~~q\in\Cb
\label{sucomm} 
\ena
\noi and
\eq
{\det}_q T\equiv\al\de-q\be\ga=I. \label{sudet}
\en
The commutations (\ref{sucomm}) can be obtained from (\ref{RTT}) via the 
$R$ matrix
\eq
\R{ab}{cd}=\left( \begin{array}{cccc} q & 0 & 0 & 0 \\
                                      0 & 1 & 0 & 0 \\
                                      0 & q-\qm & 1 & 0 \\
                                      0 & 0 & 0 & q  \end{array} \right)
  \label{Rsu}
\en
\noi where the rows and columns are numbered in the order 
11, 12, 21, 22. 

It is easy to verify that the ``quantum determinant" defined in 
(\ref{sudet}) commutes with $\al,\be,\ga$ and $\de$, so that the 
requirement ${\det}_q T=I$ is consistent. The matrix inverse of $\T{a}{b}$
is 
\eq
\Ti{a}{b}= ({\det}_q T)^{-1} \left( \begin{array}{cc} \de & -\qm\be \\ 
          -q\ga & \al 
         \end{array} \right)   \label{Timatrix}
\en
\sk
The coproduct, counit and coinverse of $\al,\be,\ga$ and $\de$ are 
determined via formulas (\ref{copT})-(\ref{coiT}) to be:
\eqa
\D(\al)=\al\otimes\al+\be\otimes\ga,~~~\D(\be)=\al\otimes\be+
\be\otimes\de \nonumber\\
\D(\ga)=\ga\otimes\al+\de\otimes\ga,~~~\D(\de)=\ga\otimes\be+
\de\otimes\de \label{copsu}\\
\epsi(\al)=\epsi(\de)=1,~~~\epsi(\be)=\epsi(\ga)=0 \label{cousu}\\
\kappa(\al)=\de,~~\kappa(\be)=-\qm \be,~~\kappa(\ga)=-q\ga,~~
\kappa(\de)=\al . \label{coisu}
\ena
\sk

\noi{\bf Note} 1.2.1 $~$ In general $\kappa^2 \not= 1$, as can be seen from 
(\ref{coisu}). The following useful relation holds \cite{FRT}:
\eq
\kappa^2 (\T{a}{b})=\Dmat{a}{c} \T{c}{d} \Dmatinv{d}{b}=d^a d^{-1}_b
\T{a}{b},  \label{k2}
\en
\noi where $D$ is a diagonal matrix, $\Dmat{a}{b}=d^a \de^a_{b}$, given 
by $d^a=q^{2a-1}$ for the $q$-groups $A_{n-1}$ and $GL_q(n)$. 
\sk
\noi{\bf Note} 1.2.2 $~$ The commutations (\ref{sucomm}) are compatible with 
the 
coproduct $\D$, in the sense that $\D(\al\be)=q\D(\be\al)$ and so on.
In general we must have 
\eq
\D(R_{12}T_1T_2)=\D(T_2T_1R_{12}), \label{DRTT} 
\en
\noi which is easily verified using $\D(R_{12}T_1T_2)=R_{12}\D(T_1)
\D(T_2)$ and $\D(T_1)=T_1 \otimes T_1$. This is equivalent to proving 
that the matrix elements of the matrix product $T_1\TP_1$, where
$\TP$ is a matrix [satisfying (\ref{RTT})] 
whose elements {\sl commute} with those of $\T{a}{b}$,
still obey the commutations (\ref{rtt}).
\sk
\noi{\bf Note} 1.2.3  $~\D({\det}_qT)={\det}_qT\otimes {\det}_qT~$ so that the coproduct 
property $\D(I)=I\otimes I$ is compatible with ${\det}_qT=I$.
\sk
\noi{\bf Note} 1.2.4 $~$ The condition (\ref{sudet}) can be relaxed. Then we have 
to include the central element $\zeta=({\det}_q T)^{-1}$ in $A$, so 
as to be able 
to define the inverse of the $q$-matrix $\T{a}{b}$ 
as in (\ref{Timatrix}),
and the coinverse of the element $\T{a}{b}$ as in (\ref{coiT}). The
$q$-group is then $GL_q(2)$. The reader can deduce the co-structures on 
$\zeta$: $\D(\zeta)
=\zeta \otimes \zeta,~\epsi(\zeta)=1,~\kappa(\zeta)={\det}_q T$.
\sk
\noi{\bf  Note} 1.2.5 $~$ More generally, the quantum determinant of 
$n \times n$ 
$q$-matrices is defined by
${\det}_qT=\sum_{\sigma} (-q)^{l(\sigma)} \T{1}{\sigma (1)} \cdots 
\T{n}{\sigma (n)}$, where $l(\sigma)$ is the minimal number of inversions 
in the permutation $\sigma$. Then ${\det}_q T=1$ restricts $GL_q(n)$ to
$SL_q(n)$. 
\sk
\noi{\bf Note} 1.2.6 $~$ The explicit expression of the $R$-matrices for the 
A,B,C,D 
$q$-groups 
will be given later. Here we recall the important relations \cite{FRT} 
for the $\Rhs$ matrix 
defined by $\Rhats{ab}{cd} \equiv \R{ba}{cd}$, whose $q=1$ limit is the 
permutation operator $\delta^a_d \delta^b_c$:
\eq
\Rhs^2=(q-q^{-1}) \Rhs + I,~~{\rm for}~A_{n-1} ~~~{\rm(Hecke~condition)}
  \label{Hecke}
\en
\eq
(\Rhs -qI)(\Rhs +q^{-1} I)(\Rhs -q^{1-N}I)=0,~~{\rm for}~B_n,C_n,D_n,
  \label{R3}
\en
\noi with $N=2n+1$ for the series $B_n$ and $N=2n$ for $C_n$ and $D_n$.
 Moreover for all $A,B,C,D$ $q$-groups 
the $R$ matrix is lower triangular 
($\R{ab}{cd}=0$ if  [$\sma{$a=c$}$ and $\sma{$b<d$}$] or
$\sma{$a<c$}$)
and satisfies:
\eq
\Rinv{ab}{cd} (q) = \R{ab}{cd} (q^{-1})\label{qqmeno1}
\en
\eq
\R{ab}{cd}=\R{dc}{ba}.
 \en
\sect{Duality and $*$-Structure}
\noi{\bf Duality}
\sk
Consider a finite dimensional Hopf algebra $A$, the vector 
space $A'$ dual to $A$ is also a Hopf algebra with the following  
product, unit and costructures
[we use the notation $\psi(a)=\le\psi,a\re$ in order to stress the
duality between $A'$ and $A$]: $\forall\psi,\phi\in A'$, $\forall a,b\in A$
\eq
\le \psi\phi,a\re=\le\psi\otimes\phi,\D a\re~,~~\le I,a\re=\epsi(a)
\label{pair1}
\en
\eq
\le\D'(\psi), a\otimes b\re=\le\psi,ab\re~,~~\epsi'(\psi)=\le\psi,I\re
\label{pair2}
\en
\eq \le \kappa' (\psi),a\re=\le\psi,\kappa (a)\re\label{pair3}
\en 
where 
$\le \psi\otimes \phi~,~a\otimes b\re\equiv\le\psi, a\re\,\le\phi, b\re~.$ 
Obviously $(A')'=A$ and $A$ and $A'$ are dual Hopf algebras.
\sk
In the infinite dimensional case the definition of duality between 
Hopf algebras is more delicate because the coproduct on $A'$ 
might not take values  in the subspace  $A'\otimes A'$ of
$(A\otimes A)'$ and therefore is ill defined.
However, the space $A^0$ spanned by the matrix 
elements of all finite-dimensional representations of $A$ is a
subalgebra of $A'$ and obviously $\D'(A^0)\subset A^0\otimes A^0$, 
$\kappa(A^0)\subset A^0$
[indeed $\D'(M^i{}_j)=\sum_{k=1}^{dim(M)} M^i{}_k\otimes M^k{}_j$, 
$\kappa(M^i{}_j)=(M^{-1})^i{}_j$]. 
Then  $A^0$ is a Hopf algebra: the Hopf dual of $A$. 
In general $(A^0)^0\not=A$, for example if $\mbox{\sl g}$ is semisimple
$U(\mbox{\sl g})^0=Fun(G)$ while the vector space underlying $Fun(G)^0$
is $U(\mbox{\sl g})\otimes \Cb(G)$ where $\Cb(G)$ is the vector space
on $\Cb$ freely generated by the elments of $G$ \cite{Chari}.
Later on we will consider the quantum groups  of the series $A_n,B_n,C_n,D_n$
and their quantized universal enveloping algebras (as the algebras of regualr
functionals on the deformed $A_n,B_n,C_n,D_n$ \cite{FRT}); 
disregarding the topological aspects we will call these algebras
dually paired or dual,  the quantum group $SL_q(N)$ can indeed be
considered as the Hopf dual of the deformed universal enveloping algebra   
of the $A_n$ series.
\sk
For generic Hopf algebras we will use the notion of non-degenerate pairing:
two Hopf algebras $A$ and $U$ are paired if there exists a bilinear map
$\le~,~\re~:U\otimes A\rightarrow \Cb$ satisfying (\ref{pair1}) and 
(\ref{pair2}), the pairing is non-degenerate if we also have
\eq
\forall~ \psi\in U~~\le\psi,a\re=0\Rightarrow a=0\label{nondeg1} 
\en
and
\eq
\forall~ a\in A~~\le\psi,a\re=0\Rightarrow \psi=0~.\label{nondeg2}
\en
Condition (\ref{nondeg1}) states that $U$ separates the points (elements)
of $A$ and viceversa for (\ref{nondeg2}). 
If $U$ and $A$ are finite dimensional then
(\ref{nondeg1}) and (\ref{nondeg2})
are equivalent to $A'=U$; indeed (\ref{nondeg1}) induces the injection 
$a\rightarrow \le~,a\re$ of $A$ in $U'$, similarly, by 
(\ref{nondeg2}) $U\subseteq A'$ and therefore $A'= U$.
\sk
It is easy to prove that the Hopf algebras $Fun(G)$ and $U(\mbox{\sl g})$ 
described in Section 1.1 are paired when $\mbox{\sl g}$ is the Lie algebra 
of $G$.
Indeed we realize $\mbox{\sl g}$ as left invariant vectorfields $t$ on the 
group manifold and $U(\mbox{\sl g})$ as the algebra generated by composition
of the  operators $t$. Then the pairing is defined by
\[\forall t\in \mbox{\sl g}, \forall f\in Fun(G),~~
\le t,f\re=t(f)|_{{}_{1_{{}_G}}}
\] where $1_{{}_G}$ is the unit of $G$. Notice that        
$t$ is left invariant if $TL_g(t|_{{}_{1_{{}_G}}})=
t|_g$, where $TL_g$ is the tangent map induced by the left multiplication
of the group on itself: $L_gg'=gg'$.  We then have  
\eq
t(f)|_g=\left( TL_g t|_{{}{1_{{}_G}}}\right)(f)=
t[f(g{\tilde g})]|_{{}_{{\tilde g}=1_{{}_G}}}=  
t[f_1(g)f_2({ \tilde g})]|_{{}_{{\tilde g}=1_{{}_G}}}=
f_1(g)\,t(f_2)|_{{}_{1_{{}_G}}}\label{tislinv}
\en
and therefore
\[\le \tilde{t}{\scriptstyle{{}^{{}_{\circ}}}} t ,f\re=\tilde{t} 
(t(f))|_{{}_{1_{{}_G}}}=
\tilde{t}f_1|_{{}_{1_{{}_G}}}tf_2|_{{}_{1_{{}_G}}}
=\le \tilde{t}\otimes t,\D f\re
\]
and, in agreement with (\ref{copL}) and (\ref{coiL}):
\[\le t, f h\re=t(f)|_{{}_{1_{{}_G}}}h|_{{}_{1_{{}_G}}}+f|_{{}_{1_{{}_G}}}
t(h)|_{{}_{1_{{}_G}}}=\le \D'(t),f\otimes h\re~,
\]
\[\le t, \kappa (f)\re=t[f(g^{-1})]|_{{}_{g=1_{{}_G}}}=
-t[f(g)]|_{{}_{g=1_{{}_G}}}
=\le \kappa'(t) , f\re ~.
\]
\sk
\sk
\noi{\bf $*$-Structure}
\sk
The Hopf algebra $SL_q(2)$ we have considered in the previous 
section can be interpreted as the deformation of the algebra
of functions on a group manifold only introducing a 
$*$-structure on $SL_q(2)$ (the analogue of complex conjugation). 
This procedure leads to the quantum groups
$SL_q(2, \Rb)=Fun_q(SL(2,\Rb))$ and $SU_q(2)=Fun_q(SU(2))$.
We need a $*$-structure because the algebra of regular functions on 
$SU(2)$ is isomorphic to
the algebra of regular functions on $SL(2,\Rb)$ and to the
algebra of analytic functions on the {\sl complex} manifold $SL(2,\Cb)$; 
indeed any $f\in Fun(SU(2))$ can be analytically continued in a unique 
function $\hat f\in Fun(SL(2,\Cb))$, then the restriction 
of $\hat{f}$ to the $SL(2,\Rb)$ sub-manifold of $SL(2,\Cb)$
belongs to $Fun(SL(2,\Rb))$.
Therefore, without a $*$-structure we cannot understand if
the polynomials in the symbols $T^i{}_j$ with the relation det$T=1$ 
generate functions on $SU(2)$ or on $SL(2,\Rb)$ or analytic functions on
$SL(2,\Cb)$.
\sk
In general a $*$-algebra over the complex numbers is an algebra with an 
anti-linear map $* : A\rightarrow A$
that is involutive, $*^2=id$ and anti-multiplicative,
$(ab)^*=b^*a^*~\forall a,b\in A$.
A Hopf $*$-algebra $A$ is a Hopf algebra $A$ over $\Cb$ equipped with a
$*$-algebra structure which is compatible with the costructures of $A$:
\eq
\D(a^*)=a_1^*\otimes a_2^*~~;~~~~~~\epsi(a^*)=\epsi(a)~~.
\en
These two conditions imply, forall $a\in A$ 
\eq [\kappa(a)]^*=\kappa^{-1}(a^*)~;
\label{*K*}
\en
indeed the operator $*{\scriptstyle{{}^{{}_{\circ}}}}\kappa^{-1}
{\scriptstyle{{}^{{}_{\circ}}}} *$ 
satisfies all the properties of the
antipode and since (as the inverse of an element in a group) 
the antipode is unique, we have (\ref{*K*}). 
\sk
Let us clarify the interrelation 
between real forms of groups
and $*$-structures on Hopf algebras. 
   
Let $A={\cal{F}}(G_{\Cbold})$ 
be the algebra of analytic functions on a complex group $G_{\Cbold}$.
A $*$-structure on $A$ determines the following real form $G_{\Rbold}$ 
of $G_{\Cbold}$: $G_{\Rbold}=\{g\in G / ~f^*(g)=\overline{f(g)}\}$; viceversa
$G_{\Rbold}$ induces on $Fun(G_{\Rbold})={\cal{F}}(G_{\Cbold})$ 
the following $*$-structure:
$f^*=h~\Leftrightarrow \overline{f(g)}=h(g) ~\forall g\in G_{\Rbold}$. Moreover
a real form of $G_{\Cbold}$ 
determines a real form $\mbox{\sl g}_{\Rbold}$ of its 
Lie algebra, i.e. $\mbox{\sl g}=\mbox{\sl g}_{\Rbold}
\oplus {\sqrt{-1}}\,\mbox{\sl g}_{\Rbold}$. The $*$-operation that acts as 
minus the identity on $\mbox{\sl g}_{\Rbold}$ satisfies $[\chi,\chi']^*=
[{\chi'}^*,\chi^*]\;$ $\forall \chi,\chi'\in\mbox{\sl g}_{\Rbold}$
and is uniquely extended as an anti-linear, anti-multiplicative 
and involutive map on the Hopf algebra 
$U(\mbox{\sl g})$, the  universal enveloping algebra 
of $\mbox{\sl g}$. 
We have seen that a $*$-structure on $A={\cal{F}}(G_{\Cbold})=Fun(G)$ 
determines a $*$-structure
on $U(\mbox{\sl g})$, the viceversa is also obviously true.
The explicit relation is $\forall\psi\in U(\mbox{\sl g})\,,~
\forall a\in Fun(G)$,
\eq
\le\psi^*,a\re={\overline{\le\psi ,[\kappa(a)]^*\re}} ~~~\mbox{ i.e.}~~~~
\le\psi,a^*\re={\overline{\le[\kappa'(\psi)]^*{},a\re}}
\label{U*A}
\en
where we have used  the pairing $\le~,\re$ between the two Hopf algebras  
$Fun(G)$ and $U(\mbox{\sl g})$, and 
${}^{{}^{\overline{~~}}}$ denotes complex conjugation.

More in general 
two Hopf $*$-algebras $A$ and $U$ are paired if, in addition to (\ref{pair1}) 
and (\ref{pair2}), relation (\ref{U*A}) holds
$\forall\psi\in U$ and $\forall a\in A$.

\sk
In a functional-analytic context, the $*$-operation becomes the
hermitian conjugation. For example,
Hopf $*$-algebras that are deformations of compact matrix groups
can be canonically realized (using the Gelfand-Naimark-Segal (GNS) 
construction, since they are dense in a $C^*$-algebra and 
since they have a Haar measure)
as bounded operators on a Hilbert space. Then the
$*$-operation is realized as the usual adjoint map ${}^{\dag}$
on  operators in Hilbert space. 
\sk
We end this section listing  
the $*$-structures that define $SU_q(2)$ and $SL_q(2,R)$;
these $*$-involutions are well defined because they are  
compatible with the $RTT$ relations 
($RT_1T_2=T_2T_1R\Leftrightarrow \overline{R}T^*_2T^*_1=T^*_1T^*_2
\overline{R}$)
and with the determinat condition.
\sk
i) $T^*={T^{-1}}^t 
\Rightarrow \al^*=\de,~\be^*=-q\ga
,~\ga^*=-\qm\be,~\de^*=\al$, where $q$ is a real number.
Gives the Hopf $*$-algebra $SU_q(2)$.

ii)${T^*}=T \Rightarrow \al^*=\al,~\be^*=\be,~\ga^*=
\ga,~\de^*=\de$, $|q|=1$. 
Gives the Hopf $*$-algebra $SL_q(2,\Rb)$.

\chapter{Differential Geometry on Quantum Groups}

In this chapter we 
study basic notions of differential geometry on 
a  Hopf algebra. We first give an introductory review of the 
$q$-differential calculus studied by \cite{Wor}. In the first section, 
following \cite{Wor}, \cite{AC}, 
we define the conditions an exterior differential $d$
has to satisfy and we explain them as natural generalizations of the
conditions on a classical Lie group. The space of $1$-forms and then 
that of $n$-forms is fully 
characterized. Then we introduce the left invariant vectorfields and deduce
their $q$-Lie algebra properties from the properties of the exterior 
differential $d$. The theory is exemplified on the quantum group $GL_{q}(N)$
in Section 2.2. In Section 2.3  we reconsider the basic postulates 
of a differential calculus and describe equivalent ones. 
This section is complementary to Section 2.1 because  emphasizes the space
of vectorfields rather than the space of one forms and the exterior 
differential. We start from a $q$-Lie algebra that closes under
a quadratic $q$-Lie bracket and then we derive the properties of the exterior
differential. 

The last section of this chapter is a deeper study of the geometric
structures on Hopf algebras. We analize the duality between
the space of vectorfields and the space of $1$-forms 
(tangent and the cotangent bundle) and generalize the construction to 
tensorfields. An inner derivative or contraction operator naturally arises
from the above duality. We then introduce a Lie derivative 
and analize its properties. Finally we 
show how these operators and the exterior differential 
form a graded quantum Lie algebra of differantial operators.

\sect{Bicovariant differential calculus}

In this section we review the bicovariant differential 
calculus on $q$-groups as developed by Woronowicz \cite{Wor}. The 
\qone limit will constantly appear in our discussion, so as to make 
clear which classical structure is being $q$-generalized.
\sk
The calculus can be developed on a generic Hopf algebra $A$ with invertible 
antipode. Since we will constantly compare the construction with the 
classical one on Lie groups, we will think of $A$ as the algebra
of the preceding section, i.e. the algebra 
freely generated by the matrix entries $\T{a}{b}$, modulo the relations 
(\ref{RTT}) and possibly some reality or orthogonality conditions. 
However, as said, the construction can be applied to more general cases,
for example it gives differential calculi on finite groups, 
where one cannot apply the  usual techniques of differential geometry. 
\sk
A {\bf {\sl first-order differential calculus}} on $A$ is defined by
\sk
i) a linear map $d$: $A \rightarrow \Gamma$, satisfying the Leibniz rule
\eq
d(ab)=(da)b+a(db),~~\forall a,b\in A; \label{Leibniz}
\en
$\Gamma$ is an 
appropriate bimodule on $A$, which 
essentially means that its elements can be 
multiplied on the left and on the right by elements of $A$.
[More precisely $A$ is a left-module if
$\forall a, b \in A\forall \rho,\rho'\in\Ga$ we have:
$a(\rho+\rho')=a\rho + a\rho'$, $(a+b)\rho=a\rho+b\rho$,
$a(b\rho)=(ab)\rho$, $I\rho=\rho$. Similarly one defines a right-module. 
A left- and right-module is a bimodule if we also have $a(\rho b)=(a\rho)b$].
The space $\Ga$ $q$-generalizes the space of 1-forms on a Lie group.
\sk
ii) the possibility of expressing any $\rho \in \Ga$ as 
\eq
\rho=\sum_k a_k db_k \label{adb}
\en
\noi for some $a_k,b_k$ belonging to $A$.
\sk
\noi {\sl{\bf Bicovariance}}
\sk
The first-order differential calculus $(\Ga,d)$ is said to be 
{\sl bicovariant} if it is both 
left- and right-covariant, i.e.  if we can consistently 
define a left and right action of 
the $q$-group on $\Ga$ as follows
\eqa
& &\!\!\!\!\! \DL(adb)=\D(a)(id \otimes d)\D(b),~~~\DL:\Ga\rightarrow 
A\otimes \Ga
~~~{\rm (left~ covariance)} \label{leftco}\\
& &\!\!\!\!\! \DR(adb)=\D(a)(d \otimes id)\D(b),~~~\DR:\Ga\rightarrow 
\Ga\otimes A
~~~{\rm (right~ covariance)} \label{rightco}
\ena
 
How can we understand these left and right actions on $\Ga$ in the \qone 
limit ? The first observation is that the coproduct $\D$ on $A$ is 
directly related, for $q=1$, to the pullback induced by left 
multiplication of the group on itself
\eq
L_x y \equiv xy,~~~\forall x,y \in G. \label{leftmu}
\en
\noi This induces the left action (pullback) 
$\ll{x}$ on the functions on $G$:
\eq
\ll{x} f(y)\equiv f(xy)|_y,~~~~\ll{x} : Fun(G) \rightarrow Fun(G) 
\label{pullback}
\en
\noi where $f(xy)|_y$ means $f(xy)$ seen as a function of $y$. Let us 
introduce the mapping $\LL$ defined by
\eqa
&(\LL f) (x,y) \equiv (\ll{x} f)(y) = f(xy)|_y&  \nonumber\\
&\LL : Fun(G) \rightarrow Fun(G \times G) \approx Fun(G) \otimes Fun(G).
& 
\label{LL}
\ena
The coproduct $\D$ on $A$, when $q=1$, reduces to the mapping $\LL$. 
Indeed, considering $\T{a}{b}(y)$ as a function on $G$, we have:
\eq
\LL (\T{a}{b}) (x,y) = \ll{x} \T{a}{b}(y) = \T{a}{b} (xy)=\T{a}{c}(x)
\T{c}{b}(y),  
\en
\noi since $\T{a}{b}$ is a representation of $G$. Therefore
\eq
\LL (\T{a}{b}) = \T{a}{c} \otimes \T{c}{b} 
\en
\noi and $\LL$ is seen to coincide with $\D$, cf. (\ref{copT}).
\sk
The pullback $\ll{x}$ can also be defined on 1-forms $\rho$ as
\eq
(\ll{x} \rho) (y) \equiv \rho(xy)|_y \label{410}
\en
\noi and here too we can define $\LL$ as
\eq
(\LL \rho)(x,y) \equiv (\ll{x} \rho)(y) = \rho(xy)|_y .\label{411}
\en
In the $q=1$ case we are now discussing, the left action $\DL$ coincides 
with this mapping $\LL$ for 1-forms. Indeed for $q=1$
\eqa
\lefteqn{\DL (adb) (x,y)=[\D (a) (id \otimes d) \D (b)] (x,y)=
[(a_1 \otimes a_2)
(id\otimes d)(b_1 \otimes b_2)](x,y)} \nonumber\\
& & =[a_1b_1\otimes  a_2db_2](x,y)=a_1(x)b_1(x) a_2(y)db_2(y)=
a_1(x)a_2(y) 
d_y[b_1(x)b_2(y)]\nonumber\\
& & =\LL(a)(x,y) d_y [\LL(b) (x,y)]=a(xy) db(xy)|_y .
\ena
\noi On the other hand:
\eq
\LL(adb)(x,y)=a(xy) db(xy)|_y,
\en
\noi so that $\DL \rightarrow \LL$ when $q \rightarrow 1$. In 
the last equation 
we have 
used the well-known property $\ll{x}(adb)=\ll{x}(a)\ll{x}(db)=
\ll{x}(a) d\ll{x}(b)$ of the 
classical pullback. A similar discussion holds for $\DR$, and we have 
$\DR \rightarrow \RR$ when \qone, where $\RR$ is defined via the 
pullback $\rr{x}$ on 
functions (0-forms) or on 1-forms induced by the right multiplication:
\eq
R_x y=yx,~~~\forall x,y \in G
\en
\eq
(\rr{x} \rho) (y)= \rho (yx)|_y
\en
\eq
(\RR \rho)(y,x) \equiv (\rr{x} \rho)(y).\label{R*rho}
\en
These observations explain why $\DL$ and $\DR$ are called left and right 
actions of the quantum group on $\Ga$ when $q\not= 1$.
\sk
From the definitions (\ref{leftco}) and (\ref{rightco}) one deduces
the following properties \cite{Wor}:
\eq
(\epsi \otimes id) \DL (\rho)=\rho,~~~(id \otimes \epsi) \DR (\rho)=
\rho \label{Dprop1}
\en
\eq
(\D \otimes id)\DL=(id\otimes\DL)\DL,~~~(id\otimes\D)\DR=(\DR\otimes
 id)\DR \label{Dprop2}
\en
\eq
(id \otimes \DR) \DL = (\DL \otimes id) \DR, 
\label{bicovariance} 
\en
this last condition is the $q$-analogue of the fact that  
left and right actions commute for $q=1$
($L^*_x R^*_y = R^*_y L^*_x$). 
\sk\sk
\noi{\sl{\bf Left- and right-invariant $\om$}}
\sk
An element $\om$ of 
$\Ga$ is said to be {\sl left-invariant} if
\eq
\DL (\om) = I \otimes \om \label{linvom}
\en
\noi and {\sl right-invariant} if
\eq
\DR (\om) = \om \otimes I. \label{rinvom}
\en
This terminology is easily understood: in the classical limit,
\eqa
\LL \om = I \otimes \om \\
\RR \om = \om \otimes I
\ena
\noi indeed define respectively left- and right-invariant 1-forms.

\noi{\sl Proof:} the classical definition of left-invariance is
\eq
(\ll{x} \om)(y) = \om(y)
\en
\noi or, in terms of $\LL$,
\eq
(\LL\om)(x,y) = \ll{x}\om(y)=\om(y).
\en
\noi But
\eq
(I\otimes \om)(x,y)=I(x) \om(y)=\om(y),
\en
\noi so that 
\eq
\LL \om = I \otimes \om \label{llinvom}
\en
\noi for left-invariant $\om$. A similar argument holds for 
right-invariant $\om$.
\sk\sk
\noi {\sl{\bf Consequences}}
\sk
For any bicovariant first-order calculus one can prove the following
\cite{Wor} [statements {\bf i)}, {\bf ii)}, {\bf iii)} and 
formulae (\ref{defchi}) and 
(\ref{copchi})-(\ref{coichi}) holds also for a calculus that is only 
left covariant)]:
\sk
{\bf{i) }} Any $\rho \in \Ga$ can be uniquely written in the form:
\eqa
\rho=a_i \om^i \label{rhoaom}\\
\rho=\om^i b_i \label{rhoomb}
\ena
\noi with $a_i$, $b_i \in A$, and {$\om^i$} a basis of 
$\invG$, the linear 
subspace of all left-invariant elements of $\Ga$. Thus, as in the 
classical case, the whole of $\Ga$ is generated by a basis of left
invariant $\om^i$. An analogous theorem holds with a basis of right
invariant elements
$\eta^i \in \Ginv$. Note that in the quantum case we have 
$a \om^i \not= \om^i a$, the bimodule structure of 
$\Ga$ being non-trivial for $q \not= 1$. [This is a consequence of
associativity: $\forall\rho\in\Ga,\forall a,a'\in A,~\rho a=a\rho$
$\Rightarrow$ $(aa')\rho=\rho(aa')=(\rho a)a'=a'(\rho a)=a'a\rho$
$\Rightarrow$ $aa'=a'a$]
\sk
{\bf{ii)}} There exist linear functionals $\f{i}{j}$ on $A$ such that, 
for any $a,b \in A$:
\eqa
& & \om^i b= (\f{i}{j} * b) \om^j \equiv (id \otimes \f{i}{j}) 
\Delta (b)\om^j \label{omb}\\
& & a\om^i=\om^j [(\f{i}{j} \circ \kappa^{-1})* a] \label{aom}
\ena
\sk
\noi Once we have the functionals $\f{i}{j}$, we 
know how to commute elements of $A$ through elements of $\Ga$. The 
$\f{i}{j}$ are uniquely determined by (\ref{omb}) and for consistency 
must satisfy the conditions:
\eqa
& & \f{i}{j} (ab)= \f{i}{k} (a) \f{k}{j} (b) \label{propf1}\\
& & \f{i}{j} (I) = \del^i_j \label{propf2}\\
& & (\f{k}{j} \circ \kappa ) \f{j}{i} = \del^k_i ~\epsi;~~~ 
   \f{k}{j} (\f{j}{i} \circ \kappa)  = \del^k_i ~\epsi,~~~ 
   \label{propf3}
\ena
\noi so that their coproduct, counit and coinverse are given by:
\eqa
& & \Dp (\f{i}{j})=\f{i}{k} \otimes \f{k}{j}   \label{copf}\\
& & \ep (\f{i}{j}) = \del^i_j  \label{couf}\\
& & \kp (\f{k}{j}) \f{j}{i}= \de^k_i \epsi = \f{k}{j} \kp (\f{j}{i})
\label{coif}
\ena
\noi cf. (\ref{pair1})-(\ref{pair3}). Note 
that in the $q=1$ limit $\f{i}{j} 
\rightarrow \de^i_j \epsi$, i.e. $\f{i}{j}$ becomes proportional to the 
identity functional $\epsi(a)=a(1_G)$, and formulas (\ref{omb}), 
(\ref{aom}) become trivial, e.g. $\om^i b = b\om^i$ [use $\epsi * a=a$ 
from (\ref{prop2})].
\sk  
{\bf{iii)}} There exists an {\sl adjoint representation} $\M{j}{i}$ of the 
quantum group, defined by the right action on the (left invariant)
$\om^i$:
\eq
\DR (\om^i) = \om^j \otimes \M{j}{i},~~~\M{j}{i} \in A. \label{adjoint}
\en
It is easy to show that $\DR (\om^i)$ belongs to 
$\invG \otimes A$, which proves the existence of $\M{j}{i}$.
In the classical case, $\M{j}{i}$ is indeed the adjoint representation
of the group. We 
recall that in this limit the left invariant 1-form $\om^i$ can be 
constructed as
\eq
\om^i(y)T_i = (y^{-1} dy)^i T_i ,~~~~~~y\in G.
\en
Under right multiplication by a (constant) element $x \in G: y 
\rightarrow yx$ we have, \footnote{Recall the $q=1$ definition of the 
adjoint 
representation $x^{-1} T_j x \equiv \M{j}{i} (x)T_i$.}
\eqa
&\om^i(yx)T_i=[x^{-1} y^{-1} d(yx)]^i T_i=[x^{-1}(y^{-1}dy)x]^iT_i\\
&=[x^{-1}T_j x]^i(y^{-1}dy)^j T_i=\M{j}{i}(x)\om^j(y) T_i,
\ena
\noi so that
\eq
\om^i(yx)=\om^j(y) \M{j}{i}(x)
\en
\noi or 
\eq
\RR \om^i(y,x)=\om^j \otimes \M{j}{i}(y,x),
\en
\noi which reproduces (\ref{adjoint}) for $q=1$.

The co-structures on the $\M{j}{i}$ can be deduced \cite{Wor}:
\eqa
& & \Delta (\M{j}{i}) = \M{j}{l} \otimes \M{l}{i} \label{copM}\\
& & \epsi (\M{j}{i}) = \delta^i_j \label{couM}\\
& & \kappa (\M{i}{l}) \M{l}{j}=\delta^j_i=\M{i}{l} \kappa (\M{l}{j}).
\label{coiM}\ena
\noi For example, in order to find 
the coproduct (\ref{copM}) it is sufficient to 
apply $(id \otimes \D)$ to both members of (\ref{adjoint}) and use
the second of eqs.(\ref{Dprop2}). 

The elements $\M{j}{i}$ can be used to 
build a right-invariant basis of $\Ga$. Indeed the $\eta^i$ defined by
\eq
\eta^i \equiv \kappa^{-1}(\M{j}{i})\om^j  \label{eta}
\en
are right invariant
[use $\kappa^{-1}(a_2)a_1=\epsi(a)$]:
\eqa
\lefteqn{
\DR (\eta^i)= \D [\kappa^{-1} (\M{j}{i})]\DR (\om^j)  =\nonumber}\\
& &\!\! [\kappa^{-1} (\M{s}{i}) \otimes \kappa^{-1} (\M{j}{s})]
[\om^k \otimes \M{k}{j} ] =\kappa^{-1} (\M{s}{i} )\om^k  \otimes \de^s_k I = 
\eta^i \otimes I
\ena
moreover every element of $\Ga$ can be  written 
as $\rho = a_i\eta^i $ or $\rho = \eta^i b_i$ where $a_i$ and $b_i$
are uniquely determined. 

\noi It can be shown that the functionals $\f{i}{j}$ previously defined 
satisfy:
\eqa
& & \eta^i b = (b * \f{i}{j}) \eta^j\label{etab}\\
& & a \eta^i = \eta^j [a*(\f{i}{j} \circ \kappa)]\label{etabb},
\ena
\noi where $a*f \equiv (f \otimes id) \D(a),~f\in \Ap$.  

{}From (\ref{omb}),  using (\ref{eta}) i.e. $\omega^i=\M{l}{i}\eta^l$
and from (\ref{etab})  one immediately prove the relation 
\eq
\M{i}{j} (a * \f{i}{k})=(\f{j}{i} * a) \M{k}{i}, \label{propM} 
\en
\noi with $a* \f{i}{j} \equiv (\f{i}{k} \otimes id) \D(a)$. 
\sk
\noi{\bf Note} 2.1.1 $~$ Given a first order differential calculus,
the space $\Ga$ is a bicovariant bimodule i.e. $\Ga$ is a bimodule
with left and right actions $\DL$ and $\DR$ that satisfy 
 (\ref{Dprop1}), (\ref{Dprop2}), (\ref{bicovariance})
and that are
compatible with the 
bimodule structure:
\eq
\DL(a\rho b)=\D(a)\DL(\rho)\D(b)~~~~;~~~
\DR(a\rho b)=\D(a)\DR(\rho)\D(b)~. \label{DRrhob}
\en
Points {\bf i), ii), iii)}, and {\bf iv)} below,  hold not only for
a first order differential calculus but also for a generic bicovariant 
bimodule
(where $\DL$ and $\DR$ are not defined via $d$).
Any bicovariant bimodule is uniquely characterized by the functionals 
$\f{i}{j}$ and the elements $\M{i}{j}$ 
satisfying (\ref{propf1}), (\ref{propf2}), (\ref{copM}), (\ref{couM})
and the fundamental condition (\ref{propM}).
Indeed,  cf. Theorem 2.4.3,  
$\DL$ is well defined via (\ref{linvom})
while $\DR$, defined by (\ref{adjoint}),  is compatible 
with the right product in $\Ga$:
\eq
\DR(\om^ia)=\DR(\om^i)\D(a) \label{compRDR}  
\en
if and only if (\ref{propM}) holds.  
{\sl Proof. }  
The projection $P~:~\Ga\rightarrow \Ginv$ defined by $\forall\rho\in\Ga\;,
P(\rho)=m(\kappa \otimes id)\DL(\rho)$ [where $m$ is the multiplication in 
the bimodule $\Ga$] maps $a\om^i$ in $P(a\om^i)=\epsi(a)\om^i$
and $\om^i a$ in $P(\om^i a)= \f{i}{j}(a)\om^j$ and is an epimorphism between 
the two bimodules $\Ga$ and $\Ginv$.
We then have:
\eqa
(P\otimes id)\DR(\om^i a)&=&(P\otimes id)\DR[(\f{i}{j}*a)\om^j]=
\om^k\otimes (\f{i}{j}*a)\M{k}{j}~~;~\nonumber\\
(P\otimes id)\DR(\om^i a)&=&
(P\otimes id)(\om^k a_1\otimes\M{k}{j}a_2)=
(P\otimes id)[(\f{i}{k}*a_1)\om^k\otimes\M{k}{j}a_2]\nonumber\\
&=&\om^k\otimes \M{k}{j}(a*\f{i}{k})~.\nonumber
\ena
This proves the implication (\ref{compRDR}) $\Rightarrow$ (\ref{propM});
the viceversa is also true since  
$(a_1\otimes id)(P\otimes id)\DR(\om^i a_2)=\DR(\om^i a)$.

\cvd
\sk

{\bf{iv)}} An {\sl exterior product}, compatible with the left 
and right actions of the $q$-group, can be defined by a bimodule 
automorphism 
$\Rh$ in $\Gamma \otimes \Gamma$ that generalizes the 
ordinary permutation operator:
\eq
\Rh (\om^i \otimes \eta^j )= \eta^j \otimes \om^i, \label{Rhat} 
\en
\noi where $\om^i$ and $\eta^j$ are respectively left and right 
invariant elements of $\Ga$. Bimodule automorphism means that
\eq
\Rh(a\tau)=a\Rh (\tau) \label{bimauto1}
\en
\eq
\Rh(\tau b)=\Rh(\tau)b \label{bimauto2}
\en
\noi for any $\tau \in \Ga \otimes \Ga$ and $a,b \in A$. 
The tensor product between elements $\rho,\rhop \in \Ga$ 
is defined to 
have the properties $\rho a\otimes \rhop=\rho \otimes a \rhop$, $a(\rho 
\otimes \rhop)=(a\rho) \otimes \rhop$ and $(\rho \otimes \rhop)a=\rho 
\otimes (\rhop a)$. Left and right actions on $\Ga \otimes \Ga$ are 
defined by:
\eq
\DL (\rho \otimes \rhop)\equiv \rho_1 \rhop_1 \otimes \rho_2 \otimes 
\rhop_2,~~~\DL: \Ga \otimes \Ga \rightarrow A\otimes\Ga\otimes\Ga
\label{DLGaGa}
\en
\eq
\DR (\rho \otimes \rhop)\equiv \rho_1 \otimes \rhop_1 \otimes \rho_2 
\rhop_2,~~~\DR: \Ga \otimes \Ga \rightarrow \Ga\otimes\Ga\otimes A
\label{DRGaGa}
\en
\noi where as usual $\rho_1$, $\rho_2$, etc., are defined by
\eq
\DL (\rho) = \rho_1 \otimes \rho_2,~~~\rho_1\in A,~\rho_2\in \Ga
\en
\eq
\DR (\rho) = \rho_1 \otimes \rho_2,~~~\rho_1\in \Ga,~\rho_2\in A.
\en
\noi More generally, we can define the action of $\DL$ 
on $\Ga^{\otimes n}\equiv\underbrace{\Ga \otimes \Ga \otimes \cdots 
\otimes \Ga}_{\mbox{$n$-times}}$ as
\[
\DL (\rho \otimes \rhop \otimes \cdots \otimes \rhopp)\equiv 
\rho_1 \rhop_1 \cdots \rhopp_1 \otimes \rho_2 \otimes 
\rhop_2\otimes \cdots \otimes \rhopp_2 
\]
\eq
\DL: \Ga^{\otimes n}\rightarrow A\otimes  \Ga^{\otimes n}
\label{DLGaGaGa}
\en
\[
\DR (\rho \otimes \rhop \otimes \cdots \otimes \rhopp)\equiv 
\rho_1 \otimes \rhop_1 \cdots \otimes \rhopp_1 \otimes \rho_2  
\rhop_2 \cdots \rhopp_2 
\]
\eq
\DR:  \Ga^{\otimes n}\rightarrow  \Ga^{\otimes n}\otimes A~.
\label{DRGaGaGa}
\en

\noi Left invariance on $\Ga\otimes\Ga$ is naturally defined as
$\DL (\rho \otimes \rhop) = I \otimes \rho \otimes \rhop$ (similar 
definition for right-invariance), so that, for example, $\om^i \otimes
\om^j$ is left invariant, and is in fact a left invariant basis for $\Ga 
\otimes \Ga$.
 
{\bf{---}} In general 
$\Rh^2 \not= 1$, since $\Rh (\eta^j \otimes \om^i )$ is not 
necessarily equal to $ \om^i \otimes \eta^j $. By linearity, $\Rh$ can
be extended to the whole of $\Gamma \otimes \Gamma$.

{\bf{---}} $\Rh$ is invertible and commutes with the left and right action of 
the $q$-group, i.e. $\DL \Rh (\rho \otimes \rhop)=(id \otimes \Rh) \DL 
(\rho \otimes \rhop)
= \rho_1\rhop_1 \otimes \Rh(\rho_2 \otimes \rhop_2)$, and 
similar for $\DR$. Then we see that $\Rh (\om^i \otimes \om^j)$ is 
left invariant, and therefore can be expanded on the left invariant 
basis $\om^k \otimes \om^l$:
\eq  
\Rh (\om^i \otimes \om^j)= \Rhat{ij}{kl} \om^k \otimes \om^l. \label{2.1.61} 
\en
\sk
{\bf{---}} From the definition (\ref{Rhat}) one can prove that \cite{Wor}: 
\eq
\Rhat{ij}{kl} = \f{i}{l} (\M{k}{j}); \label{RfM}
\en
\noi thus the functionals $\f{i}{l}$ and the elements $\M{k}{j} \in A$ 
characterizing the bimodule $\Gamma$ are dual in 
the sense of eq. (\ref{RfM}) 
and determine the exterior product:
\eqa
& & \rho \we \rho ' \equiv W(\rho\otimes\rho')\equiv\rho \otimes \rho ' - 
\Rh (\rho \otimes \rho ')\\ 
& & \om^i \we \om^j \equiv 
W^{ij}_{~kl} \om^k \otimes \om^l\equiv
\om^i \otimes \om^j - 
\Rhat{ij}{kl} \om^k \otimes \om^l. \label{exom}  
\ena
\noi Notice that, given the tensor $\Rhat{ij}{kl}$, we can compute the 
exterior product of any $\rho,\rhop \in \Ga$, since 
any $\rho \in \Ga$ is 
expressible in terms of $\om^i$ [cf. (\ref{rhoaom}), (\ref{rhoomb})]. 
The classical limit of $\Rhat{ij}{kl}$ is
\eq
\Rhat{ij}{kl} \qonelim \de^i_l \de^j_k \label{limRhat}
\en
\noi since $\f{i}{j} \qonelim \de^i_l \epsi$ and $\epsi(\M{j}{k})=\de^
j_k$. Thus in the $q=1$ limit the product defined in 
(\ref{exom}) coincides 
with the usual exterior product. 

From the property (\ref{bimauto1}) and (\ref{bimauto2})
applied to the case $\tau=\om^i \otimes \om^j$, one can derive the 
relation
\eq
\Rhat{nm}{ij} \f{i}{p} \f{j}{q} = \f{n}{i} \f{m}{j} \Rhat{ij}{pq}. 
\label{Rff}
\en
\noi Applying both members of this equation to the element $\M{r}{s}$
yields the braid relation for $\Rh$:
\eq
\Rhat{nm}{ij} \Rhat{ik}{rp} \Rhat{js}{kq}=\Rhat{nk}{ri} \Rhat{ms}{kj}
\Rhat{ij}{pq},~~i.e.~~~~\La_{12}\La_{23}\La_{12}=\La_{23}\La_{12}\La_{23}
\label{QYB}
\en
\noi which is sufficient for the consistency of (\ref{Rff}).
Taking $a=
\M{p}{q}$ in (\ref{propM}), and using (\ref{RfM}), we find the relation
dual to (\ref{QYB}):
\eq
\M{i}{j} \M{r}{q} \Rhat{ir}{pk} = \Rhat{jq}{ri} \M{p}{r} \M{k}{i} 
\label{RMM}~.
\en
This last formula explicitly shows that $W$  commutes with the coaction
$\Delta_{\scriptsize\Ga}$ [the commutation of $W$ with  
${}_{\scriptsize\Ga}\Delta$ is implicit in (\ref{2.1.61})];
the new tensor $\om^i \we \om^j \equiv 
W^{ij}_{~kl} \om^k \otimes \om^l$ transforms covariantly according
to its index structure:
\begin{equation} {}_{\scriptsize\Ga}\Delta (\om^i \wedge \om^j)
\equiv {}_{\scriptsize\Ga}\Delta(\om^k \otimes \om^l ) W_{kl}{}^{ij}
= \om^k \wedge \om^l \otimes M_k{}^i M_l{}^j~.
\end{equation} 
Using also (\ref{Rff}) we conclude that the action of 
$\Delta_{\scriptsize\Ga}$
and ${}_{\scriptsize\Ga}\Delta$ on the tensor $\rho\wedge \rho'$ has the same 
expression as in (\ref{DLGaGa}) and (\ref{DRGaGa}) with the tensor product
 replaced by the wedge product.
\sk
\noi Defining 
\eq
\R{ji}{kl} \equiv \Rhat{ij}{kl}, \label{defR}
\en
\noi we see that  $\R{ij}{kl}$ 
satisfies the quantum Yang--Baxter equation (\ref{YB}),
sufficient for the consistency of (\ref{RMM}). 
Notice that the quantum Yang--Baxter equation is 
typically associated to a quasitriangular Hopf algebra with 
universal ${\cal R}$-matrix. 
In this case, as shown in \cite{Podles} and \cite{Bonechi}, 
$\Rhat{ij}{pq}$ is a representation of the universal 
${\cal R}$-matrix of the quantum double associated to the 
{\sl generic} Hopf algebra A. In other words, (\ref{RMM})
does not rely on the existence of ``RTT'' equations (\ref{RTT}) 
used to define specific examples of Hopf algebras.

\sk
Generalizing  equation (\ref{exom}),
wedge products of $n$ forms are again expressed in terms
of tensorfields:
\begin{equation}
\om_1 \wedge \ldots \wedge \om_n = W_{1\ldots n}\,\om_1 \otimes \ldots 
\otimes \om_n
.\label{wedge1}
\end{equation}
The numerical coefficients $W_{1\ldots n}$ are given through a recursion
relation
\begin{equation}
W_{1\ldots n} = {\cal I}_{1\ldots n} W_{1\ldots n-1} , \label{wedge2}
\end{equation}
where
\begin{equation}
{\cal I}_{1\ldots n} = 1 - \Lambda_{n-1,n} +
\Lambda_{n-2,n-1} \Lambda_{n-1,n}   \ldots
-(-1)^n \Lambda_{12} \Lambda_{23} \cdots \Lambda_{n-1,n}
\label{wedge3}
\end{equation}
and $W^i{}_j = {\cal I}^i{}_j = \delta^i_j$. The space of $n$-forms 
 $\Gamma^{\we n}$ is therefore 
defined as in the 
classical case but with the quantum permutation 
operator $\Rh$.

\noindent As is easily seen writing
\[
{\cal I}_{s\ldots n}= 1-\Lambda_{n-1,n}+\Lambda_{n-2,n-1}\Lambda_{n-1,n}
\ldots -(-1)^{n-s+1}
\Lambda_{s,s+1}\Lambda_{s+1,s+2}\ldots \Lambda_{n-1,n}~,
\]
${\cal I}$ has the following
decomposition property that we will use
later on:
\begin{equation}
{\cal I}_{1\ldots n} = {\cal I}_{s\ldots n}
+(-1)^{n-s+1}{\cal I}_{1\ldots s-1} \Lambda_{s-1,s} \Lambda_{s,s+1} 
\cdots \Lambda_{n-1,n}~. \label{wedge4}
\end{equation}
Due to (\ref{RMM}) and (\ref{Rff}), the action of $\Delta_{\scriptsize\Ga}$
and ${}_{\scriptsize\Ga}\Delta$ on the tensor $\vart\in\Ga^{\wedge n}
\subset\Ga^{\otimes n}$
has the same expression as in (\ref{DLGaGaGa}) and (\ref{DRGaGaGa})
with the tensor product replaced by the wedge product.
Following Note 2.1.1, $\Ga^{\wedge n}$ is a bicovariant bimodule
with left and right coactions $\DL$ and $\DR$. We call the algebra 
$\Ga^{\wedge}$ of exterior forms, 
$\Ga^{\wedge}\equiv\Ga^0\oplus \Ga\oplus\Ga^{\wedge 2}\oplus
\Ga^{\wedge 3}\oplus\ldots$ (with $A\equiv\Ga^0$), 
a bicovariant graded algebra because
it is a graded algebra with $\DL$ and $\DR$ that are grade preserving.
\sk
{\bf{v)}} Having the exterior product we can define the 
{\sl exterior differential}
\eq
d~:~\Gamma \rightarrow \Gamma \we \Gamma
\en
\eq
d (a_k db_k) = da_k \we db_k,
\en
\noi which can easily be extended to 
$\Gamma^{\we n}$:
\eq
d: \Gamma^{\we n} 
\rightarrow \Gamma^{\we (n+1)}
\en
\eq 
d (a_{k_1k_2 ...k_n} db_{k_1}\wedge db_{k_2}\wedge ...db_{k_n})=
da_{k_1k_2, ...k_n}\we db_{k_1}\wedge db_{k_2}\wedge ...db_{k_n}
\en
and has the following properties:
\eq
d(\theta \we \thetap)=d\theta \we \thetap + (-1)^k \theta \we d\thetap 
\label{propd1}
\en
\eq
d(d\theta)=0\label{propd2}
\en
\eq
\DL (d\theta)=(id\otimes d)\DL(\theta)\label{propd3}
\en
\eq
\DR (d\theta)=(d\otimes id)\DR(\theta),\label{propd4}
\en
\noi where $\theta \in \Ga^{\we k}$, $\thetap \in \Ga^{\we n}$. 
The last two properties express the fact that $d$ commutes with the left 
and right action of the quantum group, as in the classical case. 
\sk
{\bf{vi)}} The  $q$-tangent space $T$, 
dual to the left invariant subspace $\invG$ can be 
introduced as a linear subspace of $\Ap$, whose basis elements $\chi_i
\in T\subset \Ap$ are defined by
\eq
da=(\chi_i * a)\om^i,~~~\forall a\in A. \label{defchi}
\en
\noi In the commutative case we write
\eq
da=\pdy{\mu} a(y) dy^{\mu} = \left(\pdy{\mu}a\right) 
\viel{\mu}{i}(y) \viel{i}{\nu}
(y) dy^{\nu}=\left(\pdy{\mu}a\right) \viel{\mu}{i} (y) 
\om^i(y)=t_i|_{{}_y}\om^i
 \label{daclass}
\en
\noi where $\viel{i}{\nu} (y)$ is the vielbein of the group manifold 
and $\viel{\mu}{i}$ is its inverse.
Now recalling (\ref{tislinv}) we write 
\eq
da=(t_i|_{{}_{1_{{}_G}}}*a)\om^i
\en
where $t_i$ are left invariant vectorfields.
Therefore $\chi_i *$ is the $q$-analogue of left invariant vectorfields, 
while $\chi_i$ is the $q$-analogue of the tangent vector at the 
origin $1_G$ of $G$:
\eq
\chi_i * a \qonelim \pdy{\mu} a(y) \viel{\mu}{i} \equiv \part_i a(y)~,~~~~
\chi_i(a) \qonelim \pdx{i} a(x) |_{x=1_G}\,  \label{chiclass}
\en
i.e.  $T\qonelim \mbox{\sl g}$ where $\mbox{\sl g}$ is the Lie algebra of $G$
(and here the Hopf algebra $A$ is the $q$-deformation of $Fun(G)\,$).
\sk
{\bf{vii)}} Given the basis $\{\chi_i\}$ of the $q$-tangent space $T$, 
we can introduce the ``coordinates'' $\{x^i\}$ via the 
following definition: 
consider the linear space $R$ (Woronowicz right ideal) given by
$R=\{a\in A~/~\epsi(a)=0~ \mbox{ and }~ T(a)=0\}$, 
define the linear space $X$ by the relation
\eq
A=X\oplus R\oplus \{I\} \label{axri}
\en
i.e. 
$X$ is maximal in the (ordered) set of all linear subspaces of $A$ disjoint
from $R\oplus\{I\}$. From (\ref{axri}) it follows that the dual vector space
$X'$ is isomorphic to $T$ and therefore
there are $n$ elements
$x^i\in X\subset {\mbox{ker}}\epsi$ uniquely defined by the duality
\eq
\le\chi_i\,,x^j\re=\delta^j_i ~.\label{chix}
\en
Note that $\epsi(x^i)=0$ since $X\subset\mbox{ker}\epsi$ because
$A={\mbox{ker}}\epsi \oplus\{I\}.$ In the classical limit,
in a neighbourhood of the identity $1_G\in G$, we have
$\forall g\in G, ~g=\prod_i e^{x^i(g)\chi_i}$, see for ex. \cite{Fronsdal},
the $x^i$ are called canonical coordinates of the second kind on the group $G$
(while those of the first kind
are given by $g=e^{\scriptsize{\sum_i}y^i(g)\chi_i}$).
\sk

{\bf{viii)}} The $\chi_i$ functionals close on the $q$-Lie algebra:
\eq
\chi_i \chi_j - \Rhat{kl}{ij} \chi_k \chi_l = \C{ij}{k} \chi_k,
\label{qLie}
\en
\noi with $\Rhat{kl}{ij}$ as given in (\ref{RfM}). The product 
$\chi_i\chi_j$ is defined by [cf. (\ref{pair1})]
\eq
\chi_i\chi_j \equiv (\chi_i \otimes \chi_j) \D
\en
\noi and sometimes indicated by $\chi_i * \chi_j$. Note that this $*$ 
product (called also convolution product) is associative:
\eq
\chi_i * (\chi_j * \chi_k)= (\chi_i * \chi_j)*\chi_k
\en
\eq
\chi_i * (\chi_j * a) = (\chi_i * \chi_j) * a,~~a\in A.
\en
\noi We leave the easy proof to the reader. The 
$q$-structure constants $\C{ij}{k}$ are given by
\eq
\C{ij}{k} = \chi_j(\M{i}{k}). \label{Cijk}
\en
\noi This last equation is easily seen to hold in the $q=1$ limit, since 
the $(\chi_j)_i^{~k} \equiv \C{ij}{k}$ are indeed in this case 
the infinitesimal generators of the adjoint representation:
\eq
\M{i}{k}=\de^k_i + \C{ij}{k} x^j +  O(x^2).
\en
\noi Using $\chi_j \qonelim \pdx{j} |_{x=1_G}$ indeed yields (\ref{Cijk}).

By applying both sides of (\ref{qLie}) to $\M{r}{s} \in A$, we find the 
$q$-Jacobi identities:
\eq
\C{ri}{n} \C{nj}{s}-\Rhat{kl}{ij} \C{rk}{n} \C{nl}{s} = 
\C{ij}{k} \C{rk}{s},
\label{Jacobi}
\en
\noi which give an explicit matrix realization (the adjoint 
representation) of the generators $\chi_i$:
\eq
(\chi_i)_k^{~l}=\chi_i(\M{k}{l})=\C{ki}{l}. \label{chiadj}
\en
\noi Note that the $q$-Jacobi identities (\ref{Jacobi}) can also be given 
in terms of the $q$-Lie algebra generators $\chi_i$ as :
\eq
[[\chi_r,\chi_i],\chi_j]-\Rhat{kl}{ij} [[\chi_r,\chi_k],\chi_l]=
[\chi_r,[\chi_i,\chi_j]],\label{qJacobi}
\en
\noi where 
\eq
[\chi_i,\chi_j] \equiv \chi_i \chi_j - \Rhat{kl}{ij} \chi_k\chi_l
\en
\noi is the deformed commutator of eq. (\ref{qLie}). 
\sk
{\bf{ix)}} The left invariant $\om^i$ satisfy the $q$-analogue of the 
{\sl Cartan-Maurer equations}:
\eq
d\om^i+\unmezzo\c{jk}{i} \om^j \we \om^k=0, \label{CM}
\en
\noi where
\eq
\c{jk}{i} \equiv 2\chi_j \chi_k (x^i) \label{cijk}
\en
The structure constants $C$ satisfy the Jacobi 
identities obtained by taking 
the exterior derivative of (\ref{CM}):
\eq
(\c{jk}{i} \c{rs}{j} - \c{rj}{i} \c{sk}{j} )\om^r\we\om^s\we\om^k=0.
\label{cJacobi}
\en
In the $q=1$ limit, $\om^j \we \om^k$ becomes antisymmetric 
in {\sl j} and {\sl k}, and we have
\eq
\c{jk}{i} \qonelim =(\chi_j \chi_k - \chi_k\chi_j)(x^i)
=\C{jk}{l}~ \chi_l (x^i) =  \C{jk}{i},
\en
\noi where $\C{jk}{l}$ are now the classical structure constants. Thus
when $q=1$ we have $\c{jk}{i}=  \C{jk}{i}$ and (\ref{CM}) 
reproduces the classical Cartan-Maurer equations.  

For $q \not= 1$, we find the following relation:
\eq
\C{jk}{i} = \unmezzo\c{jk}{i} - \unmezzo\Rhat{rs}{jk} \c{rs}{i}  
\label{Ccrelation}
\en 
\noi after applying both members of eq. (\ref{qLie}) to $x^i$. Note
that, using (\ref{Ccrelation}), the Cartan-Maurer equations (\ref{CM})
can also be written as:
\eq
d\om^i+\C{jk}{i} \om^j \otimes \om^k=0. 
\en 
\sk
{\bf{x)}} Finally, we derive two operatorial identities that become trivial in 
the limit $q \rightarrow 1$. From the formula
\eq
d(h * \theta)=h * d\theta,~~~h \in \Ap,~\theta \in \Ga^
{\we n}\label{dhthe}
\en
\noi [a direct consequence of (\ref{propd4})] with $h=\f{n}{l}$, we find
\eq
\chi_k  \f{n}{l}=\Rhat{ij}{kl} \f{n}{i} \chi_j  \label{bic4}~.
\en

By requiring consistency between the external derivative and the 
bimodule structure of $\Ga$, i.e. requiring that
\eq
d(\om^i a)=d[(\f{i}{j} * a) \om^j],
\en
\noi one finds the identity 
\eq
\C{mn}{i} \f{m}{j} \f{n}{k} + \f{i}{j} \chi_k= \Rhat{pq}{jk} \chi_p 
\f{i}{q} + \C{jk}{l} \f{i}{l}. \label{bic3}
\en
\noi See Appendix A for the derivation of (\ref{bic4}) and 
(\ref{bic3}); see Subsection 2.3.5 for an alternative derivation. 
\sk
In summary, a bicovariant calculus on a Hopf algebra $A$
(``the algebra of functions on the quantum group")  is characterized 
by functionals $\chi_i$ and $\f{i}{j}$ on $A$ 
satisfying, cf. 
\cite{Bernard},
\eq 
\chi_i \chi_j - \Rhat{kl}{ij} \chi_k \chi_l = \C{ij}{k} \chi_k
\label{bico1}
\en
\eq
\Rhat{nm}{ij} \f{i}{p} \f{j}{q} = \f{n}{i} \f{m}{j} \Rhat{ij}{pq} 
\label{bico2}
\en
\eq
\C{mn}{i} \f{m}{j} \f{n}{k} + \f{i}{j} \chi_k= \Rhat{pq}{jk} \chi_p 
\f{i}{q} + \C{jk}{l} \f{i}{l} \label{bico3}
\en
\eq
\chi_k  \f{n}{l}=\Rhat{ij}{kl} \f{n}{i} \chi_j,  \label{bico4}
\en
\noi where the $q$-structure constants are given by $\C{jk}{i}=
\chi_k(\M{j}{i})$ and the braiding matrix by $\Rhat{ij}{kl}=
 \f{i}{l} (\M{k}{j})$.  

The co-structures on the 
quantum Lie algebra generators $\chi_i$  are:
\eqa
& & \Dp (\chi_i)=\chi_j 
     \otimes \f{j}{i} + \epsi \otimes \chi_i \label{copchi}\\
& & \ep(\chi_i)=0 \label{couchi}\\
& & \kp(\chi_i)=-\chi_j \kp(\f{j}{i}), \label{coichi}
\ena
\noi which $q$-generalize the ones given in (\ref{copL})-(\ref{coiL}).
These co-structures derive from the duality relations 
(\ref{pair2}) and (\ref{pair3}). For
example, using the Leibniz rule for the exterior differential and 
(\ref{omb}), we have $\chi_i(ab)=\chi_j(a)\f{j}{i}(b)+\epsi(a)\chi_i(b)$ i.e. 
(\ref{copchi}); a straighforward way to obtain (\ref{coichi}) is
to apply $m(id\otimes \kappa)$ to (\ref{copchi}).
The co-structures on the functionals 
$\f{i}{j}$ have been given in (\ref{copf})-(\ref{coif}) and can be 
easily derived from (\ref{copchi}) and (\ref{couchi}) using the 
coassociativity of the coproduct [eq.  (\ref{prop1})].
These costructures are consistent with the bicovariance
conditions (\ref{bico1})-(\ref{bico4}).  

Relations (\ref{bico1}), (\ref{bico4}) and (\ref{copchi}) are also
{\sl sufficient} to construct a bicovariant calculus on $A$:
\sk
\noi{\bf Proposition } 2.1.1 $~$ Consider a set $\{\chi_i,\f{i}{j}\}$
of functionals on $A$ that satisfies:
\eqa
&&  \chi_i\chi_j -\La^{ef}_{~ij}\chi_e\chi_f\in T\\
&&  \chi_k\f{n}{l}=\Rhat{ij}{kl} \f{n}{i} \chi_j\\
&&  \chi_i(ab)=\chi_j(a)\f{j}{i}(b)+\epsi(a)\chi_i(b)~~\forall a,b\in A
\label{copchiii}\ena
where $T$ is the vector space spanned by the $\chi_i$, 
and assume that the algebra $U$ of polynomials in $\chi_i$ and $\f{i}{j}$
separates the points of $A$. Then these data determine a bicovariant 
differential calculus on $A$.
\sk
\noi{\sl Proof } This proposition is easily proven in the framework of 
Section 2.3.
Formula (\ref{bico4}) can also be written
\eq
a\!d_{\f{k}{j}}\chi_i=\Lambda^{kh}_{~ij}\chi_h
\en
where $a\!d$ is the adjoint action: 
$\forall \psi, \phi\in U,~ a\!d_{\psi}\phi\equiv\kappa'(\psi_1)\phi\psi_2$. 
We similarly have [see the first three terms in (\ref{trick})]:
\eq
a\!d_{\chi_j}\chi_i\equiv\kappa'({\chi_j}_{{}_1}) \chi_i{\chi_j}_{{}_2}
= \chi_i\chi_j -\La^{ef}_{~ij}\chi_e\chi_f\in T
\label{adbraket}        
\en
Notice that the $a\!d$ action is a right representation of $U$ on $U$:
$
a\!d_{\psi\zeta}\varphi=
a\!d_{\zeta}(a\!d_{\psi}\varphi)
$
and therefore we conclude that, $\forall \psi\in U$,  
$ad_{\psi}\chi_k$  is a linear combination
of $\chi_i$ elements.
This last condition and (\ref{copchiii}) are formulae (\ref{adbico}) 
and (\ref{leibdef}). 
In Section 2.3 the differential calculus is explicitly 
constructed out of these two conditions.

\cvd
\sk
By applying (\ref{bico1})-(\ref{bico4})
to the element $\M{r}{s}$ we express these
relations in the adjoint representation, thus obtaining a set of 
numerical equations necessary for the existence of a 
bicovariant calculus:
\eqa
& &\!\!\!\!\!\! \C{ri}{n} \C{nj}{s}-\Rhat{kl}{ij} \C{rk}{n} \C{nl}{s} = 
\C{ij}{k} \C{rk}{s}
~~\mbox{({\sl q}-Jacobi identities)} \label{bicov1}\\
& &\!\!\!\!\!\! \Rhat{nm}{ij} \Rhat{ik}{rp} \Rhat{js}{kq}=\Rhat{nk}{ri}  
\Rhat{ms}{kj}
\Rhat{ij}{pq}~~~~~~~~~\mbox{(Yang--Baxter)} \label{bicov2}\\
& &\!\!\!\!\!\! \C{mn}{i} \Rhat{ml}{rj} \Rhat{ns}{lk} + \Rhat{il}{rj} 
\C{lk}{s} =
\Rhat{pq}{jk} \Rhat{il}{rq} \C{lp}{s} + \C{jk}{m} \Rhat{is}{rm}
\label{bicov3}\\
& &\!\!\!\!\!\! \C{rk}{m} \Rhat{ns}{ml} = \Rhat{ij}{kl} \Rhat{nm}{ri} 
\C{mj}{s}
\label{bicov4}
\ena
\sk
In the next section, we describe a constructive procedure due to 
Jur\v co 
\cite{Jurco} 
for a bicovariant differential calculus on any $q$-group of 
the $A,B,C,D$ series considered in \cite{FRT}. The procedure is 
illustrated on the example of $GL_q(2)$, for which all the objects 
$\f{i}{j},~\M{r}{s},~\Rhat{ij}{kl},~\c{jk}{i}$ and $\C{jk}{i}$ 
are explicitly computed. 

\sect{Constructive procedure and the example of $GL_q(2)$}

The  $q$-groups discussed in Section 1.2 are characterized by the 
matrix $\R{ab}{cd}$. In terms of this matrix, it is possible to 
construct a bicovariant differential calculus on these $q$-groups
 \cite{Jurco},  see also 
\cite{Zumino1}, \cite{Watamura}. The 
general procedure is described in this section, and the results for the
specific case of $GL_q(2)$ are collected in the table. For a detailed study
of the $GL_q(3)$ case see \cite{AscCas}. 
\sk
{\bf The $L^{\pm}$ functionals}
\sk
We start by introducing the linear functionals $\Lpm{a}{b}$, defined
by their value on the elements $\T{a}{b}$:
\eq
\Lpm{a}{b} (\T{c}{d})=\Rpm{ac}{bd}, \label{defL}
\en
\noi where
\eq
\Rp{ac}{bd} \equiv c^+ \R{ca}{db} \label{Rplus} 
\en
\eq
\Rm{ac}{bd} \equiv c^- \Rinv{ac}{bd}, \label{Rminus}
\en
\noi where $c^+$, $c^-$ are free parameters (see later). The 
inverse matrix $R^{-1}$ is defined by
\eq
\Rinv{ab}{cd}\R{cd}{ef} \equiv \de^
a_e \de^b_f \equiv \R{ab}{cd}\Rinv{cd}{ef}. 
\en
\noi We see that the $\Lpm{a}{b}$ functionals are dual to 
the $\T{a}{b}$ elements (fundamental representation) in the same way 
the $\f{i}{j}$ functionals are dual to the $\M{i}{j}$ elements
of the adjoint representation. To extend the definition (\ref{defL}) 
to the whole algebra $A$, we set:
\eq
\Lpm{a}{b} (ab)=\Lpm{a}{g} (a) \Lpm{g}{b} (b),~~~\forall a,b\in A 
\label{Lab}
\en
\noi so that, for example,
\eq
\Lpm{a}{b} (\T{c}{d} \T{e}{f}) = \Rpm{ac}{gd}\Rpm{ge}{bf}.\label{perQYB}
\en
\noi In general, using the compact notation introduced in Section 1.2,
\eq
\LLpm_1(T_2T_3...T_n)=\RRpm_{12} \RRpm_{13} ... \RRpm_{1n}. \label{LTT}  
\en
\noi As it is esily seen from (\ref{perQYB}), the quantum Yang-Baxter equation 
(\ref{QYB}) is a {necessary} and {sufficient} condition for the 
compatibility of (\ref{defL}) and (\ref{Lab}) with the $RTT$ relations:
$\LLpm_1(R_{23}T_2T_3 - T_3T_2R_{23})=0$
\sk
Finally, the value of $\LLpm$ on the unit $I$ is defined by
\eq
\Lpm{a}{b} (I)=\de^a_b. \label{LI}
\en
Thus the functionals $\Lpm{a}{b}$ have the same properties as their
adjoint counterpart $\f{i}{j}$, and not surprisingly the latter will be
constructed in terms of the former.    

From (\ref{LTT}) we can also find the action of $\Lpm{a}{b}$ on $a\in A
$, i.e. $\Lpm{a}{b} * a$. Indeed
\eqa
\lefteqn{\Lpm{a}{b} * (\T{c_1}{d_1} \T{c_2}{d_2}\cdots\T{c_n}{d_n})=
 [id \otimes \Lpm{a}{b}] \D(\T{c_1}{d_1} \T{c_2}{d_2} \cdots
\T{c_n}{d_n})=}
 \nonumber\\
& & [id \otimes \Lpm{a}{b}]\D(\T{c_1}{d_1})\cdots\D(\T{c_n}{d_n})=
 \nonumber\\
& & [id\otimes \Lpm{a}{b}] (\T{c_1}{e_1}\cdots\T{c_n}{e_n}\otimes
  \T{e_1}{d_1}\cdots\T{e_n}{d_n})\nonumber\\
& & \T{c_1}{e_1}\cdots\T{c_n}{e_n}\Lpm{a}{b}(\T{e_1}{d_1}\cdots
\T{e_n}{d_n})=\nonumber \\
& & \T{c_1}{e_1}\cdots\T{c_n}{e_n} \Rpm{ae_1}{g_1d_1} 
\Rpm{g_1e_2}{g_2d_2}\cdots\Rpm{g_{n-1}e_n}{~~~bd_n}
\ena
\noi or, more compactly,
\eq
\LLpm_1 * T_2...T_n=T_2...T_n \RRpm_{12} \RRpm_{13} ...\RRpm_{1n},
\en
\noi which can also be written as the cross-commutation relation
\eq
\LLpm_1 T_2=T_2 \RRpm_{12}\LLpm_1.   
\en

It is not difficult to find the commutations between $\Lpm{a}{b}$ 
and $\Lpm{c}{d}$:
\eq
R_{12} \LLpm_2 \LLpm_1=\LLpm_1 \LLpm_2 R_{12} \label{RLL}
\en
\eq 
R_{12} \LLp_2 \LLm_1=\LLm_1 \LLp_2 R_{12}, \label{RLpLm}
\en
\noi where as usual the product $\LLpm_2 \LLpm_1$ is the convolution 
product: $\LLpm_2 \LLpm_1 (a)\equiv (\LLpm_2 \otimes \LLpm_1)\D(a)$ 
$\forall a\in A$. 
Consider
\eq
R_{12} (\LLp_2 \LLp_1)(T_3)=R_{12}(\LLp_2 \otimes\LLp_1)\D (T_3)=
R_{12}(\LLp_2 \otimes\LLp_1)(T_3 \otimes T_3)=(c^+)^2~R_{12}R_{32}R_{31}
\nonumber
\en
\noi and
\eq
\LLp_1 \LLp_2(T_3) R_{12}=(c^+)^2 ~R_{31}R_{32}R_{12}
\nonumber
\en
\noi so that the equation (\ref{RLL}) is proven for $\LLp$ 
by virtue of the 
quantum Yang--Baxter equation (\ref{YB}), where the indices have been 
renamed $2\rightarrow 1,3 \rightarrow 2,1\rightarrow 3$. Similarly,
one proves the remaining ``RLL" relations.
\sk
\noi{\bf Note} 2.2.1 $~$ As mentioned in \cite{FRT}, $L^+$ is upper 
triangular, $L^-$ is lower triangular (this is due to the upper 
and lower 
triangularity of $R^+$ and $R^-$, respectively).
{}From (\ref{RLL}) and (\ref{RLpLm}) we have
\eq
\Lpm{A}{A} \Lpm{B}{B}=\Lpm{B}{B} \Lpm{A}{A}~;~~ 
\Lp{A}{A} \Lm{B}{B}=\Lm{B}{B} \Lp{A}{A} 
\en
\sk
\noi {\bf Note} 2.2.2 $~$
A determinant can be defined for the matrix 
$\Lpm{A}{B}$ 
as in Note 1.2.5,
with $q \rightarrow \qm$.
Indeed the ``$RLL$"
relations are identical to the  "$RTT$" with 
$R \rightarrow R^{-1}$ (which means
 $q \rightarrow \qm,~r \rightarrow r^{-1}$,
cf. eq. (\ref{Rinv})). Then , because
of the upper or lower triangularity of $\LLp$ and $\LLm$
respectively, we have
\eq
{\det}_q \LLpm=\Lpm{1}{1} \Lpm{2}{2} \cdots \Lpm{N}{N} 
\en
\sk
\noi{\bf Note} 2.2.3 $~$ 
From (\ref{defL}) we deduce:
\eq
\Lpm{A}{B} ({\det}_q T)= \de^A_B (c^{\pm})^N r^{\pm 1}  
\label{LondetT0}
\en
{\sl Proof:} observe that $\Lpm{A}{B} 
({\det}_q T)=\Lpm{A}{B} (\T{1}{1}  
\T{2}{2}
\cdots \T{N}{N})$ since all the other permutations 
do not contribute,  
due to the structure of the $R^{\pm}$ matrix. Then it 
is easy to see that
\eq
\!\!\!\!\!\!\!\!\!\!\!\!\Lpm{A}{B} (\T{1}{1} \T{2}{2} \cdots \T{N}{N})=
\nonumber
\en
\eq
{}~~~~\de^A_B (c^{\pm})^N 
\Rpm{A1}{A1} \Rpm{A2}{A2} \cdots \Rpm{AN}{AN}=
\de^A_B (c^{\pm})^N
r^{\pm 1}~.\label{cmeno} 
\en
If we set  ${\det}_q T =I$,
then $\Lpm{A}{B} ({\det}_q T)= \de^A_B (c^{\pm})^N r^{\pm 1}$
must be equal to $\de^A_B$,
or $c^{\pm}=r^{\mp {1\over N}} \al^{\pm}$
with $(\al^{\pm})^{N}=1$.
In this case $[{\det}_q \LLpm](\T{A}{B})=\delta^A_B$
so that $\det \LLpm=\epsi$ [Proof:
${\det}_q \LLpm (\T{A}{B})$ $=$  $\de^A_B (c^{\pm})^N 
\Rpm{1A}{1A} \cdots 
\Rpm{NA}{NA} $ $=$ $
\de^A_B (c^{\pm})^N r^{\pm 1}$ ].
Thus for 
$(c^{\pm})^{N} r^{\pm 1}=1 $, the functionals $\LLpm$and 
$\epsi$ generate
the Hopf algebra $U(sl_{q}(N))$.
In the case of 
$GL_q(n)$, $c^{\pm}$ are extra free parameters.
In fact, they appear only in the combination $s=(c^+)^{-1} c^-$.
They do not enter in the $\Rh$ matrix, nor in the structure
constants or the Cartan-Maurer 
equations, they 
however enter the $\omega~-~T$ commutation relations  (see the table), 
so that different values 
of $s$ give different bimodules of $1$-forms and different 
bicovariant differential calculi on $GL_q(n)$. (This accounts for the 
one parameter family of differential calculi found in the classification of 
$GL_q(n)$ calculi \cite{SchmuedgenGL}).
\sk 
The co-structures are defined by the duality (\ref{defL}):
\eq
\Dp(\Lpm{a}{b})(\T{c}{d} \otimes \T{e}{f}) \equiv \Lpm{a}{b}
(\T{c}{d}\T{e}{f})=\Lpm{a}{g}(\T{c}{d}) \Lpm{g}{b} (\T{e}{f})
\en
\eqa
& & \ep (\Lpm{a}{b})\equiv \Lpm{a}{b} (I)\\
& & \kp (\Lpm{a}{b})(\T{c}{d})\equiv \Lpm{a}{b} (\kappa (\T{c}{d})) 
\ena
\noi cf. [(\ref{pair3}), (\ref{pair3})], so that
\eqa
& & \Dp (\Lpm{a}{b})=\Lpm{a}{g} \otimes \Lpm{g}{b}\\
& & \ep (\Lpm{a}{b})=\de^a_b \\
& & \kp (\Lpm{a}{b})=\Lpm{a}{b} \circ \kappa 
\ena
\noi The matrix $\kp(\LLpm)=(\LLpm)^{-1}$ is a polynomial in the 
$\Lpm{a}{b}$ elements and therefore  the $\Lpm{a}{b}$ 
generate a Hopf algebra, the Hopf algebra $U_q(gl(n))$ 
paired to the quantum group $GL_q(n)$\footnote{The pairing between these 
two Hopf algebras is nondegenerate, indeed $U_q(gl((n))$ as a Hopf algebra is 
isomorphic \cite{Frenkel, FRT} to $\Uc_q(gl(n))$ 
(as defined by Jimbo
\cite{JimboGL}) then the Hopf isomorphism $\Uc^0_h(sl(n))\cong SU_q(n)$
(where $\Uc^0_h(sl(n))$ 
is the Hopf dual of Drinfeld universal enveloping algebra $\Uc_h(sl(n))$ 
\cite{Drinfeld}) allows to conclude that the pairing between $GL_q(n)$ and 
$U_q(gl(n))$ is nondegenerate.}.
Note that 
\eq
\Lpm{a}{b} (\kappa (\T{c}{d}))=\Rpminv{ac}{bd}, \label{LkT}
\en
\noi since 
\eq
\Lpm{a}{b} (\kappa (\T{c}{d}) \T{d}{e} )=\de^c_d \Lpm{a}{b} (I)=
\de^c_e \de^a_b
\en
\noi and
\eqa
\Lpm{a}{b} (\kappa (\T{c}{d}) \T{d}{e})&=& \Lpm{a}{f} (\kappa(\T{c}{d}))
\Lpm{f}{b} (\T{d}{e}) \nonumber\\
&=&\Lpm{a}{f} (\kappa (\T{c}{d})) \Rpm{fd}{be}.
\ena
\sk
{\bf The space of quantum 1-forms}
\sk
The bimodule $\Ga$ (``space of quantum 1-forms") can be constructed as 
follows. First we define $\ome{a}{b}$ to be a basis of left invariant 
quantum 1-forms. The index pairs ${}_a^{~b}$ or ${}^a_{~b}$ 
will replace in the 
sequel the indices ${}^i$ or ${}_i$ of the previous section. The 
dimension of $\invG$ is therefore $N^2$ at this stage.
Since the $\ome{a}{b}$ are left invariant, we have:
\eq
\DL (\ome{a}{b})=I\otimes \ome{a}{b},~~~a,b=1,...,N. \label{DLome}
\en
\noi The left action $\DL$ on the whole of $\Ga$ is then defined by
(\ref{DLome}), since $\ome{a}{b}$ is a basis for $\Ga$. The bimodule 
$\Ga$ is further characterized by the commutations 
between $\ome{a}{b}$ and $a \in A$ [cf. eq. (\ref{omb})]:
\eq
\ome{a_1}{a_2} b=(\ff{a_1}{a_2b_1}{b_2} * b)\ome{b_1}{b_2}, \label{omeb}
\en
\noi where
\eq
\ff{a_1}{a_2b_1}{b_2} \equiv \kp (\Lp{b_1}{a_1}) \Lm{a_2}{b_2}. 
\label{defff}
\en
\noi Finally, the right action $\DR$ on $\Ga$ is defined by
\eq
\DR (\ome{a_1}{a_2}) = \ome{b_1}{b_2} \otimes \MM{b_1}{b_2a_1}{a_2}, 
\en
\noi where $\MM{b_1}{b_2a_1}{a_2}$, the adjoint representation, 
is given by
\eq
\MM{b_1}{b_2a_1}{a_2} \equiv \T{b_1}{a_1} \kappa (\T{a_2}{b_2}). 
\label{defMM}
\en
It is easy to check that $\ff{a_1}{a_2b_1}{b_2}$ fulfill the 
consistency conditions (\ref{propf1})-(\ref{propf3}), where 
the {\small{\sl i,j,...}} indices stand for pairs of {\small{\sl a,b,...}
} indices. Also, the 
co-structures of $\MM{b_1}{b_2a_1}{a_2}$ 
are as given in (\ref{copM})-(\ref{coiM}).
The last compatibility condition between the bimodule $\Ga$
and the action $\DR$, as explained in Note 2.1.1, 
is  (\ref{propM}). This relation is easily checked for 
$a=\T{A}{B}$ since in  
this case it is implied by the $RTT$ relations;
it holds for a generic $a$ because of property
(\ref{propf1}).

\sk
{\bf The $\Rh$ tensor and the exterior product}
\sk
The $\Rh$ tensor defined in (\ref{defR}) can now be computed:
\eqa
\lefteqn{\RRhat{a_1}{a_2}{d_1}{d_2}{c_1}{c_2}{b_1}{b_2} 
\equiv \ff{a_1}{a_2b_1}{b_2} (\MM{c_1}{c_2d_1}{d_2})
=\kp(\Lp{b_1}{a_1}) \Lm{a_2}{b_2} (\T{c_1}{d_1} \kappa(\T{d_2}{c_2}))
}\nonumber \\
=& & [\kp (\Lp{b_1}{a_1}) \otimes \Lm{a_2}{b_2}]\D (\T{c_1}{d_1} \kappa
(\T{d_2}{c_2}))\nonumber\\
=& & [\kp (\Lp{b_1}{a_1}) \otimes \Lm{a_2}{b_2}] 
(\T{c_1}{e_1} \otimes \T{e_1}{d_1})(\kappa(\T{f_2}{c_2})\otimes \kappa 
(\T{d_2}{f_2})) \nonumber \\
=& & [\kp (\Lp{b_1}{a_1}) \otimes \Lm{a_2}{b_2}]
[\T{c_1}{e_1}\kappa(\T{f_2}{c_2}) \otimes \T{e_1}{d_1}\kappa 
(\T{d_2}{f_2})] \nonumber\\
=& & \Lp{b_1}{a_1} (\kappa^2 (\T{f_2}{c_2}) \kappa (\T{c_1}{e_1}))
~\Lm{a_2}{b_2} (\T{e_1}{d_1} \kappa (\T{d_2}{f_2}))\nonumber \\
=& & d^{f_2} d^{-1}_{c_2} \Lp{b_1}{a_1} (\T{f_2}{c_2} \kappa (\T{c_1}{e_1
}))~\Lm{a_2}{b_2} (\T{e_1}{d_1} \kappa(\T{d_2}{f_2})\nonumber \\
=& & d^{f_2} d^{-1}_{c_2} \Lp{b_1}{g_1} (\T{f_2}{c_2}) ~\Lp{g_1}{a_1}
(\kappa (\T{c_1}{e_1}))~ \Lm{a_2}{g_2} (\T{e_1}{d_1}) ~\Lm{g_2}{b_2}
(\kappa(\T{d_2}{f_2}))\nonumber \\
=& & d^{f_2} d^{-1}_{c_2} \R{f_2b_1}{c_2g_1} \Rinv{c_1g_1}{e_1a_1}
    \Rinv{a_2e_1}{g_2d_1} \R{g_2d_2}{b_2f_2} \label{RRffMM}
\ena
\noi where we made use of relations (\ref{Dka}), (\ref{pair3}), 
(\ref{k2}), (\ref{defL}) and (\ref{LkT}). The 
$\Rh$ tensor allows the definition of the exterior product as in 
(\ref{exom}). For future use we give here also the inverse $\Rh^{-1}$
of the $\Rh$ tensor, defined by:
\eq
\RRhatinv{a_1}{a_2}{d_1}{d_2}{b_1}{b_2}{c_1}{c_2} \RRhat{b_1}{b_2}{c_1}
{c_2}{e_1}{e_2}{f_1}{f_2}= \de^{a_2}_{e_2} \de^{e_1}_{a_1} \de^{f_1}_
{d_1} \de^{d_2}_{f_2}. \label{defRRhatinv}
\en
It is not difficult to see that
\eqa
\lefteqn{\RRhatinv{a_1}{a_2}{d_1}{d_2}{b_1}{b_2}{c_1}{c_2}=
\ff{d_1}{d_2b_1}{b_2}
(\T{a_2}{c_2} \km (\T{c_1}{a_1}))=}\nonumber\\
& & \R{f_1b_1}{a_1g_1} \Rinv{a_2g_1}{e_2d_1} \Rinv{d_2e_2}{g_2
c_2} \R{g_2c_1}{b_2f_1} (d^{-1})^{c_1} d_{f_1}
\label{RRffMMinv}
\ena
\noi does the trick. Another useful relation gives a particular
trace of the $\Rh$ matrix:
\eq
\RRhat{c_1}{c_2}{b}{b}{a_1}{a_2}{b_1}{b_2} = \de^{a_1}_{a_2}
\de^{b_1}_{c_1} \de^{c_2}_{b_2}. \label{RReqdeltas}
\en
\noi This identity is simply proven. Indeed:
\eqa
\lefteqn{\RRhat{c_1}{c_2}{b}{b}{a_1}{a_2}{b_1}{b_2} 
\equiv \ff{c_1}{c_2b_1}{b_2} (\MM{a_1}{a_2 b}{b})=}\nonumber\\
& & \kp(\Lp{b_1}{c_1}) \Lm{c_2}{b_2} (\T{a_1}{b} \kappa(\T{b}{a_2}))=
\kp(\Lp{b_1}{c_1}) \Lm{c_2}{b_2} (\de^{a_1}_{a_2} I)=\nonumber\\
& & \de^{a_1}_{a_2} [\kp(\Lp{b_1}{c_1}) \otimes \Lm{c_2}{b_2}] 
    (I\otimes I)=\de^{a_1}_{a_2}\de^{b_1}_{c_1} \de^{c_2}_{b_2}.
\label{Lambda}
\ena
\sk
The relations (\ref{Hecke}), (\ref{R3}) for the $R$ matrix 
reflect themselves in relations
for the $\Rh$ matrix (\ref{RRffMM}). For example, the Hecke condition
(\ref{Hecke}) implies:
\eq
(\Rh + q^2)(\Rh + q^{-2})(\Rh-I)=0 \label{RRHecke}
\en
\noi for the $A_{n-1}$ $q$-groups, and replaces the 
classical relation $(\Rh -1)(\Rh+1)=0$, $\Rh$ being for $q=1$ the 
ordinary permutation operator, cf. (\ref{limRhat}). 

With the help of 
(\ref{RRHecke}) we can give explicitly the commutations of the
left invariant forms $\om$. Indeed, reverting to the 
{\small{\sl i,j...}} indices,
relation (\ref{RRHecke}) implies:
\eqa
\lefteqn{(\Rhat{ij}{kl} + q^2 \de^{i}_{k} \de^{j}_{l})
(\Rhat{kl}{mn} + q^{-2} \de^{k}_{m} \de^{l}_{n})
(\Rhat{mn}{rs} -  \de^{m}_{r} \de^{n}_{s})\om^r \otimes \om^s=}
\nonumber\\
& & (\Rhat{ij}{kl} + q^2 \de^{i}_{k} \de^{j}_{l})
(\Rhat{kl}{mn} + q^{-2} \de^{k}_{m} \de^{l}_{n})
\om^m \we \om^n = 0
\ena
\noi and it is easy to see that the last equality can be rewritten as
\eq
\om^i \we \om^j = - \Z{ij}{kl} \om^k \we \om^l \label{commom}
\en
\eq
\Z{ij}{kl} \equiv {1\over {q^2 + q^{-2}}} [\Rhat{ij}{kl} + \Rhatinv
{ij}{kl}]. \label{defZ}
\en
\sk
{\bf The exterior differential}
\sk 
The exterior differential on $\Ga^{\we k}$ is defined by means of the 
bi-invariant (i.e. left  and right invariant) element $\tau=\sum_a 
\ome{a}{a} \in \Ga$ as follows:
\eq
d\theta \equiv \lam [\tau \we \theta - (-1)^k \theta \we \tau], 
\label{defd} 
\en
\noi where $\theta \in \Ga^{\we k}$, and $\la$ is a 
normalization factor depending on $q$, necessary in order to obtain the 
correct classical limit. It will be later determined to be 
$\la = q-q^{-1}$. Here we can only see that it has to vanish for $q=1$, 
since otherwise $d\theta$ would vanish in the classical limit. For 
$a \in A$ we have 
\eq
da=\lam [\tau a - a \tau]. \label{defd2}
\en
\noi This linear map satisfies the Leibniz rule (\ref{Leibniz}), and 
properties 
(\ref{propd1})-(\ref{propd4}), as the reader can easily check (use the 
definition of exterior product and the bi-invariance of $\tau$). A proof 
that also the property (\ref{adb}) holds can be 
obtained by considering the 
exterior differential of the adjoint representation:
\eq
d\M{j}{i}=(\chi_k * \M{j}{i} )\om^k=\M{j}{l} \C{kl}{i}\om^k
\en
\noi or
\eq
\kappa (\M{l}{j})d\M{j}{i}=\C{kl}{i} \om^k.
\en
\noi Multiplying by $\C{ni}{l}$, we have:
\eq
\C{ni}{l} \kappa (\M{l}{j}) d\M{j}{i}=\C{kl}{i} \C{ni}{l}\om^k\equiv
g_{nk} \om^k,
\en
\noi where $g_{nk}$ is the $q$-Killing metric. 
The explicit example of this section  being $GL_q (2)$, one may wonder
what happens to the invertibility of the $q$-Killing metric, since its 
classical limit is no more invertible [$GL(2)$ being nonsemisimple].
The answer is that for $q \not= 1$ the $q$-Killing metric of $GL_q(2)$
{\sl is} invertible, as can be checked explicitly from the values of the 
structure constants given in the table. Therefore $GL_q(2)$ 
could be said to 
be ``$q$-semisimple". With an analogous procedure (using $\T{a}{b}$ 
instead of $\M{j}{i}$) we have derived in the table the 
explicit expression of the 
$\om^i$ in terms of the $d\T{a}{b}$ for $GL_q(2)$.  
\sk
{\bf The $q$-Lie algebra}
\sk
The ``quantum generators" $\cchi{a_1}{a_2}$ are introduced as in 
(\ref{defchi}):
\eq
da=\lam[\tau a - a\tau] =(\cchi{a_1}{a_2} * a) 
\ome{a_1}{a_2}. \label{defcchi}
\en
Using (\ref{omeb}) we can find an explicit expression for the
$\cchi{a_1}{a_2}$ in terms of the $\LLpm$
functionals. Indeed
\eq
\tau a= \ome{b}{b} a= (\ff{b}{bc_1}{c_2} * a)\ome{c_1}{c_2}=
([\kp (\Lp{c_1}{b}) \Lm{b}{c_2}] * a) \ome{c_1}{c_2}.
\en
\noi Therefore
\eq
da=\lam [(\kp (\Lp{c_1}{b}) \Lm{b}{c_2}- \de^{c_1}_{c_2} \epsi) * a] 
\ome{c_1}{c_2} \label{explicitda}
\en
\noi (recall $\epsi * a= a$), so that the $q$-generators take the 
explicit form
\eq
\cchi{c_1}{c_2}=\lam [\kp (\Lp{c_1}{b})\Lm{b}{c_2}-\de^{c_1}_{c_2} 
\epsi ]=\lam(\ff{b}{bc_1}{c_2} -\de^{c_1}_{c_2} \epsi)~.
\label{defchi2}
\en
\noi The commutations between the $\chi$'s can now be obtained by taking 
the exterior derivative of eq. (\ref{explicitda}). We find
\eqa
\lefteqn{d^2(a)=0=d[(\cchi{c_1}{c_2} * a) \ome{c_1}{c_2}]=
(\cchi{d_1}{d_2}*\cchi{c_1}{c_2} * a) \ome{d_1}{d_2} \we \ome{c_1}{c_2}
 + (\cchi{c_1}{c_2} * a)d\ome{c_1}{c_2} }\nonumber\\
& & =(\cchi{d_1}{d_2}*\cchi{c_1}{c_2} * a) (\ome{d_1}{d_2} \otimes 
\ome{c_1}{c_2}- \RRhat{d_1}{d_2}{c_1}{c_2}{e_1}{e_2}{f_1}{f_2} 
\ome{e_1}{e_2} \otimes \ome{f_1}{f_2})  \nonumber\\
& & +\lam (\cchi{c_1}{c_2} * a)
(\ome{b}{b} \we \ome{c_1}{c_2} + \ome{c_1}{c_2} \we \ome{b}{b}).
\label{d2}
\ena
\noi Now we use the fact that $\tau=\ome{b}{b}$ is bi-invariant, and 
therefore also right-invariant, so that we can write
\eqa
\lefteqn{\ome{b}{b} \we \ome{c_1}{c_2} + \ome{c_1}{c_2} \we \ome{b}{b}
\equiv}\nonumber\\
& & \ome{b}{b} \otimes \ome{c_1}{c_2} - \Rh (\ome{b}{b} \otimes 
\ome{c_1}{c_2})+\ome{c_1}{c_2} \otimes \ome{b}{b} - 
\Rh (\ome{c_1}{c_2} \otimes \ome{b}{b})=\nonumber\\
& & \ome{c_1}{c_2} \otimes \ome{b}{b} - \Rh (\ome{b}{b} \otimes 
\ome{c_1}{c_2})=\nonumber\\
& & \ome{c_1}{c_2} \otimes \ome{b}{b} - \RRhat{b}{b}{c_1}{c_2}{e_1}{e_2
}{f_1}{f_2} \ome{e_1}{e_2} \otimes \ome{f_1}{f_2}, 
\label{om2}
\ena
\noi where we have used $\Rh(\ome{c_1}{c_2} \otimes \tau) = 
\tau \otimes \ome{c_1}{c_2}$, cf. (\ref{Rhat}). After substituting
(\ref{om2}) in (\ref{d2}), and factorizing $\ome{d_1}{d_2} \otimes 
\ome{c_1}{c_2}$, we arrive at the $q$-Lie algebra relations:
\eq
\cchi{d_1}{d_2} \cchi{c_1}{c_2} - \RRhat{e_1}{e_2}{f_1}{f_2}
{d_1}{d_2}{c_1}{c_2} ~\cchi{e_1}{e_2} \cchi{f_1}{f_2} =
\lam [-\de^{c_1}_{c_2} \cchi{d_1}{d_2} + \RRhat{b}{b}{e_1}{e_2}{d_1}{d_2}
{c_1}{c_2} ~\cchi{e_1}{e_2}].
\label{qLieexplicit}
\en
\noi The structure constants are then explicitly given by:
\eq
\CC{a_1}{a_2}{b_1}{b_2}{c_1}{c_2} =\lam [- \de^{b_1}_{b_2} 
\de^{a_1}_{c_1} 
\de^{c_2}_{a_2} + \RRhat{b}{b}{c_1}{c_2}{a_1}{a_2}{b_1}{b_2}]. \label{CC}
\en
Here we determine $\la$. Indeed we first observe that
\eq
\RRhat{a_1}{a_2}{d_1}{d_2}{c_1}{c_2}{b_1}{b_2}=
\de^{b_1}_{a_1} \de^{a_2}_{b_2} \de^{c_1}_{d_1} \de^{d_2}_{c_2}
+ (q-q^{-1}) \U{a_1}{a_2}{d_1}{d_2}{c_1}{c_2}{b_1}{b_2},
  \label{Rlam}
\en
\noi where the matrix $U$ is finite and different from zero in the limit 
$q=1$. This can be proven by considering the explicit form of the
$R$ and $R^{-1}$ matrices. In the case of the $A_{n-1}$ 
$q$-groups, for example, these matrices have the form \cite{FRT}:
\eq
\R{ab}{cd}= \de^{a}_{c} \de^b_d + (q-q^{-1}) \left[{{q-1}\over{q-q^{-1}}}
\de^a_c \de^b_d \de^{ab} + \de^b_c \de^a_d \theta (a-b) \right]
\label{explicitR}
\en
\eq
\Rinv{ab}{cd}= \de^{a}_{c} \de^b_d - (q-q^{-1}) \left[{{1-q^{-1}}
\over{q-q^{-1}}}
\de^a_c \de^b_d \de^{ab} + \de^b_c \de^a_d \theta (a-b) \right],
\label{explicitRinv}
\en
\noi where $\theta (x) = 1$ for $x > 0$ and vanishes for $x \leq 0$.  
Substituting these expressions in the formula for $\Rh$ (\ref{RRffMM}) 
we find (\ref{Rlam}). Using (\ref{Rlam}) in the expression
(\ref{CC}) for the $q$-structure constants $\Cb$, we find that the
terms proportional to $\lam$ do cancel, and we are left with:
\eq
\CC{a_1}{a_2}{b_1}{b_2}{c_1}{c_2} =- \lam (q-q^{-1}) 
\U{b}{b}{c_1}{c_2}{a_1}{a_2}{b_1}{b_2}. \label{CClam}
\en
\noi A simple choice for $\la$ is therefore $\la=q-q^{-1}$, ensuring
that $\Cb$ remains finite in the limit \qone; moreover
with this normalization the differential $d$ reduces for $q\rightarrow 1$
to the classical differential, cf. Section 4.6.
\sk
{\bf The Cartan-Maurer equations}
\sk
The Cartan-Maurer equations are found as follows:
\eq
d\ome{c_1}{c_2}=\lam (\ome{b}{b} \we \ome{c_1}{c_2} + \ome{c_1}{c_2} \we 
\ome{b}{b}) \equiv -\unmezzo 
\cc{a_1}{a_2}{b_1}{b_2}{c_1}{c_2} ~\ome{a_1}{a_2} \we
\ome{b_1}{b_2}. \label{CartanMaurer}
\en
In order to obtain an explicit and, for $q\rightarrow 1$, well defined 
expression for the $C$ structure 
constants in (\ref{CartanMaurer}), we must use the relation 
(\ref{commom}) for the commutations of $\ome{a_1}{a_2}$ with
$\ome{b_1}{b_2}$. Then the term $\ome{c_1}{c_2} \we \ome{b}{b}$
in (\ref{CartanMaurer}) can be written as 
$-Z\om\om$ via formula (\ref{commom}), and we find the 
$C$-structure constants to be:
\eqa
\cc{a_1}{a_2}{b_1}{b_2}{c_1}{c_2}&=&
-{2\over \la}(\de^{a_1}_{a_2} \de^{b_1}_{c_1}
\de^{c_2}_{b_2} - {1 \over {q^2 + q^{-2}}}[ \RRhat{c_1}{c_2}{b}{b}{a_1}
{a_2}{b_1}{b_2}+\RRhatinv{c_1}{c_2}{b}{b}{a_1}{a_2}{b_1}{b_2}])
\nonumber\\
&=& -{2\over \la}(\de^{a_1}_{a_2} \de^{b_1}_{c_1}
\de^{c_2}_{b_2} - {1 \over {q^2 + q^{-2}}}[ \de^{a_1}_{a_2} 
\de^{b_1}_{c_1}\de^{c_2}_{b_2}+ 
\RRhatinv{c_1}{c_2}{b}{b}{a_1}{a_2}{b_1}{b_2}]),
\label{explicitcc}
\ena
\noi where we have also used eq. (\ref{RReqdeltas}). By considering 
the analogue of (\ref{Rlam}) for $\Rh^{-1}$,
it is not difficult to see that the terms proportional to $\lam$
cancel, and the \qone limit of (\ref{explicitcc}) is well defined.
For a similar result on the $B_n,C_n$ and
$D_n$ $q$-groups see (\ref{cc}) and ref. \cite{Monteiro}. 
\sk
In the table we summarize the results of this section for the 
case of $GL_q(2)$. The composite indices ${}_a^{~b}$ are translated into 
the corresponding indices ${}^i$, $i=1,+,-,2$, according to 
the convention:
\eq
{}_1^{~1} \rightarrow {}^1,~~{}_1^{~2} \rightarrow {}^+,~~{}_2^{~1} 
\rightarrow {}^-,~~{}_2^{~2} \rightarrow {}^2
.\en
\noi A similar convention holds for ${}^a_{~b} \rightarrow {}_i$.
\sk


\vfill\eject

\subsection{Table of $GL_q(2)$ }

\centerline{The bicovariant $GL_q(2)$ algebra}
\sk
\sk
\noi {\sl $R$ and $D$-matrices:}
\[
\R{ab}{cd}=\left( \begin{array}{cccc} q & 0 & 0 & 0 \\
                                      0 & 1 & 0 & 0 \\
                                      0 & q-\qm & 1 & 0 \\
                                      0 & 0 & 0 & q  \end{array} \right)
\]
\[
\Rm{ab}{cd}\equiv c^- \Rinv{ab}{cd}= c^- \left( 
                   \begin{array}{cccc} \qm & 0 & 0 & 0 \\
                                      0 & 1 & 0 & 0 \\
                                      0 & -(q-\qm) & 1 & 0 \\
                                     0 & 0 & 0 & \qm  \end{array} \right)
\]
\[
\Rp{ab}{cd}\equiv c^+\R{ba}{dc}=c^+ 
                                      \left( \begin{array}{cccc} 
                                      q & 0 & 0 & 0 \\
                                      0 & 1 & q-\qm & 0 \\
                                      0 & 0 & 1 & 0 \\
                                      0 & 0 & 0 & q  \end{array} \right)
,~~D^a_{~b}=\left( \begin{array}{cc} q & 0  \\ 0 & q^3 \end{array} 
\right)
\]
\sk

\noi {\sl Non-vanishing components of the $\Rh$ matrix:}

\[ 
\begin{array}{llll}
\Rhat{11}{11}=1 &\Rhat{1+}{+1}=q^{-2} &\Rhat{1-}{-1}=q^2
&\Rhat{12}{21}=1\\
\Rhat{+1}{1+}=1 &\Rhat{+1}{+1}=1-q^{-2}
&\Rhat{++}{++}=1 &\Rhat{+-}{11}=1-q^2 \\
\Rhat{+-}{-+}=1 &\Rhat{+-}{21}=1-q^{-2}
&\Rhat{+2}{+1}=-1+q^{-2} &\Rhat{+2}{2+}=1\\
\Rhat{-1}{1-}=1 &\Rhat{-1}{-1}=1-q^{2}
&\Rhat{-+}{11}=-1+q^2 &\Rhat{-+}{+-}=1\\
\Rhat{-+}{21}=-1+q^{-2} &\Rhat{--}{--}=1
&\Rhat{-2}{-1}=-1+q^2 &\Rhat{-2}{2-}=1 \\
\Rhat{21}{11}=(q^2-1)^2 &\Rhat{21}{12}=1
&\Rhat{21}{+-}=q^2-1 &\Rhat{21}{-+}=1-q^2 \\
\Rhat{21}{21}=2-q^2-q^{-2} &\Rhat{2+}{1+}=-q^2+q^4
&\Rhat{2+}{+2}=q^2 &\Rhat{2+}{2+}=1-q^2\\
\Rhat{2-}{1-}=1-q^2 &\Rhat{2-}{-1}=q^{-2}-1-q^2+q^4
&\Rhat{2-}{-2}=q^{-2} &\Rhat{2-}{2-}=1-q^{-2}\\
\Rhat{22}{11}=-(q^2-1)^2 &\Rhat{22}{+-}=1-q^2
&\Rhat{22}{-+}=q^2-1 &\Rhat{22}{21}=(q^{-1}-q)^2\\
\Rhat{22}{22}=1& & &  \end{array}
\]
\sk
\noi {\sl Non-vanishing components of the $\Cb$ structure constants:}

\[
\begin{array}{llll}
\C{11}{1}=q(q^2-1) &\C{11}{2}=-q(q^2-1) &\C{1+}{+}=q^3
 &\C{1-}{-}=-q\\
\C{21}{1}=\qm-q &\C{21}{2}=q-\qm &\C{2+}{+}=-q
 &\C{2-}{-}=\qm\\
\C{+1}{+}=-\qm &\C{+2}{+}=q &\C{+-}{1}=q
  &\C{+-}{2}=-q\\
\C{-1}{-}=q(q^2+1)-\qm &\C{-2}{-}=-\qm &\C{-+}{1}=-q
  &\C{-+}{2}=q \end{array} 
\]
\sk
\renewcommand{\baselinestretch}{1.2}
\noi {\sl Non-vanishing components of the $C$ structure constants:}
\[
\begin{array}{llll}
\c{11}{1}=\frac{q(q^2-1)^2}{1+q^4} &\c{11}{2}=\frac{q^3 (1-q^2)}
{1+q^4} &\c{1+}{+}=\frac{q^5}{1+q^4}
 &\c{1-}{-}=\frac{-q^3}{1+q^4}\\
\c{12}{1}=\frac{q(1-q^2)}{1+q^4}  &\c{+1}{+}=\frac{-q^3}{1+q^4} 
&\c{+-}{1}=\frac{q^3}{1+q^4} &\c{+-}{2}= \frac{-q^3}{1+q^4}\\
\c{+2}{+}=\frac{q}{1+q^4} &\c{-1}{-}=\frac{q^5}{1+q^4} 
&\c{-+}{1}=\frac{-q^3}{1+q^4}  &\c{-+}{2}=\frac{q^3}{1+q^4}\\
\c{-2}{-}=\frac{-q^3}{1+q^4} &\c{21}{1}=\frac{q(1-q^2)}{1+q^4}
&\c{2+}{+}=\frac{-q^3}{1+q^4} &\c{2-}{-}=\frac{q}{1+q^4}\\
\c{22}{2}=\frac{q(1-q^2)}{1+q^4} & & & \end{array} 
\]
\sk
\renewcommand{\baselinestretch}{1.0}
\noi {\sl Cartan-Maurer equations:}
\[ d\om^1+q\om^+ \we \om^-=0 \]
\[ d\om^+ + q \om^+(-q^2 \om^1 + \om^2)=0 \]
\[ d\om^- + q (-q^2 \om^1 + \om^2)\om^-=0 \]
\[ d\om^2 - q \om^+ \we \om^-=0 \]
\sk

\noi {\sl The q-Lie algebra:}
\sk
\[ \chi_1 \chi_+ - \chi_+ \chi_1 -(q^4-q^2)\chi_2 \chi_+ = q^3 \chi_+\]
\[ \chi_1 \chi_- - \chi_- \chi_1 +(q^2-1)\chi_2 \chi_- = -q \chi_-\]
\[ \chi_1 \chi_2 - \chi_2 \chi_1 =0\]
\[ \chi_+ \chi_- - \chi_- \chi_+ + (1-q^2) \chi_2 \chi_1-(1-q^2) 
\chi_2 \chi_2 = q (\chi_1-\chi_2)
\]
\[ \chi_+ \chi_2-q^2 \chi_2 \chi_+ = q\chi_+\]
\[ \chi_- \chi_2-q^{-2} \chi_2 \chi_- = -\qm \chi_-\]
\sk
\noi {\sl Commutation 
relations between left invariant $\om^i$ and $\om^j$:}
\[\om^1 \we \om^+ + \om^+ \we \om^1 = 0\]
\[\om^1 \we \om^- + \om^- \we \om^1 = 0\]
\[\om^1 \we \om^2 + \om^2 \we \om^1 = (1-q^2)\om^+ \we \om^-
\]
\[\om^+ \we \om^- + \om^- \we \om^+ = 0\]
\[\om^2 \we \om^+ + q^2 \om^+ \we \om^2 = q^2 (q^2 - 1)\om^+ \we \om^1
\]
\[\om^2 \we \om^- + q^{-2} \om^- \we \om^2 = (1-q^2)\om^- \we \om^1\]
\[\om^2 \we \om^2 =(q^2 - 1)\om^+ \we \om^-
\]
\[\om^1 \we \om^1 = \om^+ \we \om^+ = \om^- \we \om^-=0\]
\sk

\noi {\sl Commutation relations 
between $\om^i$ and the basic elements of $A$ ($s=(c^+)^{-1} c^-$):}

\[ 
\begin{array}{ll}
\om^1 \al= sq^{-2} \al \om^1
& \om^+ \al= s\qm\al \om^+\\
\om^1\be=s\be\om^1 & \om^+\be=s\qm\be\om^+ + s(q^{-2} -1)\al\om^1\\
\om^1\ga=sq^{-2}\ga\om^1 & \om^+\ga=s\qm\ga\om^+\\
\om^1\de=s\de\om^1 & \om^+\de=s\qm\de\om^+ + s(q^{-2} -1)\ga\om^1
\end{array}
\]
\[
\begin{array}{ll}
\om^-\al=s\qm\al\om^- + s(q^{-2}-1)\be\om^1 & \om^2\al=s\al\om^2+s
(\qm-q)\be\om^+\\
\om^-\be=s\qm\be\om^- & \om^2\be=sq^{-2}\be\om^2+s(\qm-q)\al\om^-+s
(q^{-1}-q)^2\be\om^1\\
\om^-\ga=s\qm\ga\om^-+s(q^{-2}-1)\de\om^1 & \om^2\ga=s\ga\om^2+s(\qm-q)
\de\om^+\\
\om^-\de=s\qm\de\om^- & \om^2\de=sq^{-2}\de\om^2+s(\qm-q)\ga\om^-+s(q^
{-1}-q)^2\de\om^1
\end{array}
\]
\sk
\renewcommand{\baselinestretch}{1.2}
\noi {\sl Values and action of the generators on the q-group elements:}
\[ 
 \begin{array}{llll}
   \chi_1 (\al)={{s-q^2}\over {q^3-q}} &\chi_+ (\al)=0 &\chi_-(\al)=0 
   &\chi_2(\al)={{s-1}\over{q-\qm}}\\
   \chi_1 (\be)=0 &\chi_+ (\be)=0 &\chi_-(\be)=-s &\chi_2(\be)=0\\
   \chi_1 (\ga)=0 &\chi_+ (\ga)=-s &\chi_-(\ga)=0 &\chi_2(\ga)=0\\
   \chi_1 (\de)={{-q^2+s(1-q^2+q^4)}\over{q^3-q}} &\chi_+ (\de)=0 
    &\chi_-(\de)=0 &\chi_2(\de)={{s-q^2}\over{q^3-q}}\\
  \end{array}
\]
\[
  \begin{array}{llll}
   \chi_1*\al={{s-q^2}\over{q^3-q}}~\al &\chi_+*\al=-s\be &\chi_-*\al=0 
   &\chi_2*\al={{s-1}\over{q-\qm}}~\al\\
   \chi_1*\be={{-q^2+s(1-q^2+q^4)}\over{q^3-q}}~\be &\chi_+*\be=0 
   &\chi_-*\be=-s\al &\chi_2*\be={(s-q^2)\over{q^3-q}}~\be\\
   \chi_1*\ga={{s-q^2}\over{q^3-q}}~\ga &\chi_+*\ga=-s\de &\chi_-*\ga=0 
   &\chi_2*\ga={{s-1}\over{q-\qm}}~\ga\\
   \chi_1*\de={{-q^2+s(1-q^2+q^4)}\over{q^3-q}}~\de &\chi_+*\de=0 
   &\chi_-*\de=-s\ga &\chi_2*\de={{s-q^2}\over{q^3-q}}~\de\\
  \end{array}
\]
\sk
\noi {\sl Exterior derivatives of the basic elements of $A$:}
\[
\begin{array}{l}
 d\alpha={{s-q^2}\over{q^3-q}}\al\om^1-s\beta\om^++{{s-1}\over{q-
\qm}}\al\om^2\\
 d\beta={{-q^2+s(1-q^2+q^4)}\over{q^3-q}}\beta\om^1-s\alpha\om^-
+{{s-q^2}\over{q^3-q}}\be\om^2\\
d\gamma={{s-q^2}\over{q^3-q}}\ga\om^1-s\de\om^++{{s-1}\over{q-\qm}}\ga
\om^2\\
d\de={{-q^2+s(1-q^2+q^4)}\over{q^3-q}}\de\om^1-s\ga\om^-+{{s-q^2}\over{q
^3-q}}\de\om^2
\end{array}
\]
\noi {\sl The $\om^i$ in 
terms of the exterior derivatives on $\al,\be,\ga,
\de$:}
\[
\begin{array}{l}
 \om^1={q\over{s(-q^2-q^4+s+sq^4)}} [(q^2-s)(\kappa (\al)d\al+\kappa 
(\be)d\ga)+q^2(s-1)(\kappa (\ga)d\be + \kappa (\de) d\de)]\\
 \om^+=-{1\over s} [\kappa (\ga) d\al + \kappa (\de) d\ga]\\
\om^-=-{1\over s} [\kappa (\al) d\be + \kappa (\be) d\de]\\
\om^2={q\over{s(-q^2-q^4+s+sq^4)}}[(s-q^2-sq^2+sq^4)(\kappa (\al)
d\al+\kappa (\be) d\ga)+(q^2-s)(\kappa (\ga) d\be+\kappa (\de) d\de)]
\end{array}
\]
\vfill\eject
\renewcommand{\baselinestretch}{1.0}
\noi {\sl Lie derivative on $\om^i$:} (See Subsection 2.4.5)
\[
 \begin{array}{ll}
  \chi_1*\om^1=q(q^2-1)\om^1 + (q^{-1}-q)\om^2 &\chi_+*\om^1=-q\om^-\\
  \chi_1*\om^+=-q^{-1}\om^+ &\chi_+*\om^+=-q\om^2+q^3\om^1\\
  \chi_1*\om^-=[q(q^2+1)-q^{-1}]\om^- &\chi_+*\om^-=0\\
  \chi_1*\om^2=-q(q^2-1)\om^1 - (q^{-1}-q)\om^2 &\chi_+*\om^2=q\om^-\\
 \end{array}
\]
\[
 \begin{array}{ll}
  \chi_-*\om^1=q\om^+ &\chi_2*\om^1=0\\
  \chi_-*\om^+=0 &\chi_2*\om^+=q\om^+\\
  \chi_-*\om^-=q^{-1}\om^2-q\om^1 &\chi_2*\om^-=-q^{-1}\om^-\\
  \chi_-*\om^2=-q\om^+ &\chi_2*\om^2=0\\
 \end{array}
\]

\vfill\eject

\sect{Differential calculus from the $q$-Lie Algebra.
(A more intuitive presentation of the differential calculus on
$q$-groups)}
In the previous sections we have analized the differential calculus on
$q$-groups starting from the properties of the exterior differential
and the space of  1-forms. The left invariant vectorfields $\chi_i$ 
and their $q$-Lie algebra were then introduced at the very end. Here 
we would like to invert the exposition procedure and
following   \cite{poinAschieri}, \cite{Delius} and in the spirit of
\cite{S1,SWZ3,SWZ2}, \cite{PaoloPeter,Paolo}, \cite{Su2},
we derive differential calculi on $q$-groups from basic properties of 
$q$-Lie algebras.

This clarify the important role played by the adjoint action in the 
$q$-Lie algebra and in the construction of a bicovariant 
differential calculus. In this way we also give an alternative proof
of the Woronowicz theorems that we stated in  Section 2.1. 
 
This approach is also suitable for the study of generalizations
of the Woronowicz theory. As is evident from Subsection 2.3.3
bicovariant calculi that do not satisfy the undeformed Leibniz rule
can be found studying quantum Lie algebras that are closed under the 
adjoint action \cite{Delius}. 

In this section we consider a generic Hopf algebra $A$ and a 
Hopf algebra $U$ paired to $A$, we  also consider the pairing 
non-degenerate. Intuitively $A$ and $U$ are  
the quantum analogue of $Fun(G)$ and of the universal enveloping algebra
$U(\mbox{\sl g})$.
The differential calculus on $A$ and the quantum Lie 
algebra structure can be formulated, as in Section 2.1,  
without the introduction of $U$, however to gain a better geometrical 
understanding of the structures we are condisering, the universal 
enveloping algebra $U$ is helpful. 
We think that this presentation is more intuitive than the one in 
Section 2.1, because it is closer to the classical case, where the 
exterior differential on a group manifold can be introduced via a basis
of left invariant vectorfields and the dual basis of $1$-forms.
In this way we emphasize the role played by left invariant vectorfields 
i.e. the $q$-tangent space, a more intuitive and basic concept
than that of left invariant $1$-forms ($q$-cotangent space).
The $q$-tangent bundle of general vectorfields will be studied 
in the next section.
 
\subsection{Left invariant Vectorfields}

Classically the differential calculus on a group is uniquely
determined by the Lie algebra of the tangent vectors to  
the origin of the group. 
Locally we write a basis as $\{\partial_i|_{1_G}\}$.
Once we have this basis,
using the tangent map (namely
$TL_{g}$) induced by the left multiplication of the group on itself:
$L_gg'=gg'\;,~\,\forall g,g' \in G$ we can construct a basis of left
invariant vectorfields $\{t_i\}$. The action of these vectorfields
on a generic function $a$ on the group manifold
is
\eq t_i(a) = a_1 (\partial_i a_2 |_{1_G}) \equiv \partial_i|_{1_G} * a
\label{leftinvclass}
\en
in compliance with the following picture:\\
\unitlength=1.00mm
\special{em:linewidth 0.4pt}
\linethickness{0.4pt}
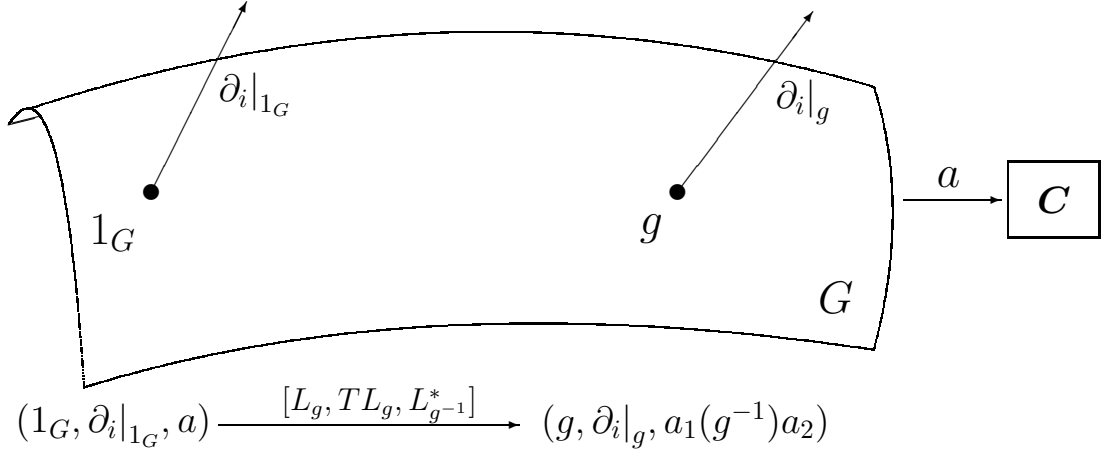
\begin{picture}(144.00,68.00)
\put(18.00,38.00){\circle*{2.00}}
\put(88.00,38.00){\circle*{2.00}}
\put(18.00,38.00){\vector(1,2){12.67}}
\put(88.00,38.00){\vector(3,4){18.00}}
\put(109.00,24.00){\makebox(0,0)[cc]{\Large$G$}}
\put(16.00,35.00){\makebox(0,0)[rt]{\Large$1_G$}}
\put(86.00,35.00){\makebox(0,0)[rt]{\Large$g$}}
\put(27.00,51.00){\makebox(0,0)[lc]{\large${\partial_i|}_{1_G}$}}
\put(101.00,50.00){\makebox(0,0)[lc]{\large${\partial_i|}_g$}}
\put(89.00,7.00){\makebox(0,0)[cc]{\large$
(g,{\partial_i|}_g,a_1(g^{-1})a_2)$}}
\put(26.00,7.00){\makebox(0,0)[rc]{\large$(1_G,{\partial_i|}_{1_G},a)$}}
\put(49.00,9.00){\makebox(0,0)[cb]{[$L_g,TL_{g},L^*_{g^{-1}}$] }}
\put(132.00,32.00){\framebox(12.00,10.00)[cc]{\large\mbox{\boldmath$C$}}}
\put(124.00,40.00){\makebox(0,0)[cc]{\Large$a$}}
\bezier{144}(114.00,52.00)(119.00,36.00)(114.00,17.00)
\bezier{236}(-1.00,47.00)(6.00,58.00)(9.00,12.00)
\bezier{432}(9.00,12.00)(54.00,26.00)(114.00,17.00)
\bezier{472}(1.50,49.00)(59.00,68.00)(114.00,52.00)
\put(27.00,7.00){\vector(1,0){40.00}}
\put(-1.00,47.00){\line(4,1){4.07}}
\put(118.00,37.00){\vector(1,0){13.00}}
\end{picture}\\
where $L^*_g(a)(h) \equiv a(gh) = a_1(g) a_2(h)$ [cf. (\ref{410})].
Explicitly, $t$ is left invariant if $TL_g(t|_{1_{{}_G}})=
t|_g$, then we have  
$$t(a)|_g=\left( TL_g t|_{1_{{}_G}}\right)(a)=
t[a(g{\tilde g})]|_{{\tilde g}=1_{{}_G}}=  
t[a_1(g)a_2({\tilde g})]|_{{\tilde g}=1_{{}_G}}=
a_1(g)t(a_2)|_{1_{{}_G}}
$$
and  \ref{leftinvclass} follows [cf. (\ref{tislinv})].

\noi In the commutative case, since the space $\invXi$ of left invariant 
vectorfields provides a 
trivialization of the tangent bundle of the group manifold the 
space $\Xi$ of 
vectorfields is isomorphic to $C^{\infty}(G)\otimes \invXi\simeq 
\invXi\otimes C^{\infty}(G)$,  i.e., a
generic vectorfield $V$ can be written as $V=b^it_i$ where 
$b^i$ are functions on the group manifold. 
Similarly a generic $1$-form can be written $\rho= b_i\om^i~~[b_i\in
C^{\infty}(G)]$ where $\{\om^i\}$ is the dual basis of
$\{t_i\}$.
Finally, the exterior differential on a generic function $b$ is
\eq
db= t_i(b)\om^i~\label{donf}
\en
and is compatible with the left and right action of the group on the space
of $1$-forms: $\ll{x}(adb)=\ll{x}(a)\ll{x}(db)=
\ll{x}(a) d\ll{x}(b)$ and $\rr{x}(adb)=\rr{x}(a)\rr{x}(db)=
\rr{x}(a) d\rr{x}(b)$. 
\sk
Following this classical construction, in this section we show 
that a differential calculus on a 
$q$-group $A$, with universal 
enveloping algebra ${\Un}$ (${\Un}\equiv U_q(\mbox{g})$), is determined 
by a $q$-Lie algebra $T$: the $q$-deformation of $\mbox{g}$. 
The exterior differential 
is then given by (\ref{donf}) where now $t_i$ are left invariant vectorfield 
on $A$ and $\{\om^i\}$ is the dual basis of $\{t_i\}$.   

It is natural to look for a linear space $T$, 
$T\subset \mbox{ker}\epsi\subset 
{\Un}$
satisfying the following three conditions:
\eq
T ~\mbox{ generates }~{\Un}\label{ZERO}
\en
\eq
\D '(T)\subset T\otimes {\Un} + \epsi\otimes T ,\label{UNO}
\en
\eq
[T,T]\subset T \label{DUE}
\en
where the bracket is the adjoint action defined by
\eq
\forall\, 
\chi, {\tilde\chi} \in {T}~,~~[{\tilde\chi},\chi]\equiv
\kappa'(\chi_{{}_1}){\tilde\chi}\chi_{{}_2}~.
\label{adjact}
\en
Conditions (\ref{ZERO}) and (\ref{DUE}) 
encode the quantum group properties of the $q$-tangent space
because they involve 
the product  and the antipode $\kappa'$ of $U$.
Condition (\ref{UNO}) states that 
the  elements of $T$ are generalized tangent vectors, 
and in fact, if  $\{\chi_i\}$ is a basis of the linear space $T$, 
we have 
\eq
\D(\chi_i)=\chi_j\otimes f^j{}_i + \epsi\otimes \chi_j\label{leibdef}
\en 
that is equivalent to   
\eq
\chi_i(ab)=\chi_j(a)\,f^j{}_i(b)+\epsi(a)\,\chi_j(b)\label{leibdef2}
\en
where  $f^j{}_i\in {\Un}$ and $\epsi'(f^j{}_i)\equiv
f^j{}_i(I)=\delta_i^j$.
[Hint: apply  
$(id \otimes \epsi')$  to (\ref{UNO})]. 
In the commutative limit we expect $f_i{}^j\rightarrow \epsi$.
Notice that we follow the historical convention
which consider derivative
operators acting from the right to the left, as is seen from their
deformed Leibniz rule (\ref{leibdef}). Of course one can consider 
deformed derivative operators ${\check{\chi_i}}$ acting from the left to 
the right, they are for example 
given by  ${\check{\chi_i}}=-{\kappa'}^{-1}(\chi_i)$,
similarly ${\check{f_j{}^i}}={\kappa'}^{-1}(f^i{}_j)$ and then 
${\D(\check{\chi_i})=\check{\chi_i}\otimes \epsi +\check{f_i{}^j} \otimes 
\check{\chi_j}}$.

Following (\ref{leftinvclass}) we can also consider the $q$-deformed 
left invariant vectorfields 
\eq
t_i\equiv\chi_i*\label{qleftinv}
\en 
defined by 
$\chi_i*a\equiv a_1\chi_i(a_2)$, 
then (\ref{UNO}) states that the $\chi_i*$ are generalized 
derivations 
\eq
\chi_i*(ab)=(\chi_j*a)(f^j{}_i*b)+ a\,\chi_j*b~,\label{leibdef3}
\en
where we have also defined $f^j{}_i*b\equiv b_1f^j{}_i(b_2)$,
and $\epsi*a\equiv a_1\epsi(a_2)=a$..

There is a one-to-one correspondence $\chi_i \leftrightarrow t_i=
\chi_i*$. In order to obtain $\chi_i$ from $\chi_i*$ we simply apply 
$\epsi$ :
$$
(\epsi \circ t_i)(a) = \epsi(id \otimes \chi_i)\D(a) = 
\epsi(a_1\chi_i(a_2)) = \epsi(a_1)\chi_i(a_2) = 
\chi_i(\epsi \otimes id)\D(a) =\chi_i(a)\;.\nonumber
$$

\subsection{Adjoint action}
\sk 
In this subsection we examine and study the consequences of
conditions (\ref{ZERO}) and (\ref{DUE}). We see that they define the adjoint 
representation $\M{i}{j}$ that we will later identify with the one 
studied in Section 2.1; we rewrite  
(\ref{ZERO}) and (\ref{DUE}) in a more geometric language using
left and right invariant vectorfields and finally, in Note 2.3.1,  we 
establish the equivalence between (\ref{ZERO})-(\ref{DUE}) and Woronowicz
theory.
\sk 
Condition (\ref{DUE}) is the closure of 
$T$ under the adjoint action,
in the classical case, if $\chi$ is a tangent vector: $\D(\chi)
=\chi\otimes \epsi +\epsi\otimes \chi\,,~ \kappa'(\chi)=-\chi$
and the adjoint action of $\chi$ on ${\tilde\chi}$ is given by the 
commutator ${\tilde\chi}\chi -\chi{\tilde\chi}$. 
Expression (\ref{adjact})
is a natural and simple way to write, both at the quantum and at the
classical level, the adjoint action. We can derive (\ref{adjact}) from 
the classical conjugation on the group elements.
First we notice that the conjugation $\mbox{\sl Conj}^{-1}~:~G \times G
\rightarrow G$, $\mbox{\sl Conj}^{-1}_{\tilde g}(g)={\tilde g}^{-1}(g)
{\tilde g}$ induces the following action on $A=Fun(G)$: 
\eq
a\!d(a)=a_2\otimes \kappa(a_1)a_3~.\label{ada} 
\en
Proof:
$~a\!d(a)(g,{\tilde g})=a_2(g)(a_1)({\tilde g}^{-1})a_3({\tilde g})
=a({\tilde g}^{-1}g{\tilde g})~.$
Expression (\ref{ada}) is independent from the points $g$ of $G$ and 
therefore holds also for a generic Hopf algebra $A$.
Now we use the pairing between $A$ and ${\Un}\subseteq A'$, discussed in
Section 1.3, to deduce the adjoint action of the universal enveloping
algebra $U$ on itself :
$\forall \psi,\varphi\in {\Un},~\forall a\in A$
\eqa
\langle a\!d_{\psi} \varphi , a\rangle
&\equiv&\langle\varphi\otimes \psi, a\!d a\rangle\label{addef}\\
&=& \langle \varphi \otimes\kappa'(\psi_1)\otimes \psi_2\,,\,
 a_2\otimes a_1\otimes a_3\rangle\\
&=&\langle\kappa'(\psi_1)\varphi\psi_2\,,\, a\rangle
\ena
so that 
\eq
a\!d_{\psi}\varphi = \kappa'(\psi_1)\varphi\psi_2\label{adjonU}~.
\en
Notice that $\ad$ is a right action of $\Un$ on  $\Un$, 
$\forall \varphi ,  \psi , \zeta\in \Un$:
\eq
a\!d_{\psi\zeta}\varphi=
 a\!d_{\zeta}(a\!d_{\psi}\varphi)\label{Uaction}~;
\en
in particular, 
 for $\chi,\tilde\chi,\chi',...\chi''\in {T},$
formulae (\ref{adjonU}) and (\ref{Uaction}) read:
\eq
a\!d_{{}{\chi}}{\tilde \chi}=[{\tilde\chi},\chi]~~~~\mbox{ and }~~~~ 
a\!d_{{}_{(\chi\chi'...\chi'')}}{\tilde\chi}=
[...[[{\tilde\chi},\chi],\chi'],\ldots\chi'']~.\label{manybrakets}
\en 
The first expression proves that the bracket in ({\ref{adjact}) is indeed the 
adjoint action; from the second expression,
(\ref{ZERO}) and (\ref{DUE}), we have
\eq
\forall\psi\in {\Un}, \forall\chi\in T~, ~~~~~~~~~~~~~~~~~
a\!d_{\psi}\chi\in T~~~~~~~~\mbox{i.e.}~~~~~~~
a\!d_{\psi}\chi_i=\M{i}{j}(\psi)\chi_j
\label{adbico}
\en
where we have introduced 
the basis $\{\chi_i\}_{i=1,...n}$ of  $T$ and $\M{i}{j}(\psi)$ are complex 
numbers depending on $\psi$.
The linearity of the adjoint map imply that the functionals $\M{i}{j}$
are linear: $\M{i}{j}(\al\psi+\phi)=\al\M{i}{j}(\psi)
+\M{i}{j}(\phi)$ while, due to (\ref{Uaction}), they are a representation
of $U$: $\M{i}{j}(\psi\phi)=\M{i}{k}(\psi)\M{k}{j}(\phi)$.
Since the pairing $A\leftrightarrow \Un$ is nondegenerate and $\psi$ 
is a generic element of $\Un$, the second expression in (\ref{adbico})
uniquely defines the functionals  $\M{i}{j}$ as elements of $A$. They 
are  the adjoint 
representation,
in Hopf algebra notations [see (\ref{coiM})]:
\eq
\varepsilon(\M{i}{j})=\delta_i^j~~~,~~~~~~
\Delta (\M{i}{j}) = \M{i}{l} \otimes \M{l}{j} \label{1copM}.
\en
We can also obtain an explicit expression for the elements $\M{i}{j}\in A$;
since $A$ separates the points of $U$, and therefore of $T$, 
we can consider $n$ elements $y^i\in A$
such that $\langle\chi_i\,, y^j\rangle=\chi_i(y^j)=\delta_i^j.$ Then
\eq
M_i{}^j= (\chi_i\otimes id)a\!d\,y^j= \chi_i(y_2^j)\kappa(y_1^j)y_3^j
\label{defMy}
\en
Proof: apply a generic element $\psi\in \Un$ to (\ref{defMy}) and recall 
(\ref{addef}). $~\Box$\\
In particular, formula (\ref{defMy}) holds if we consider 
the coordinates $x^i$ on the quantum group,
as defined in (\ref{axri}) and (\ref{chix}), we therefore have:
\eq
M_i{}^j= (\chi_i\otimes id)a\!d\,x^j= \chi_i(x_2^j)\kappa(x_1^j)x_3^j~.
\label{defM}
\en
\sk
\sk
There is an equivalent expression for (\ref{adbico}) 
that shows its $q$-group geometric content. 
We have shown that $t=\chi*$ is a
left-invariant vectorfield, similarly 
\eq
h\equiv *\chi~~~~~i.e.~~~\forall a\in A~~~ h(a)\equiv\chi(a_1)a_2
\en
is a right invariant vectorfield.
\sk
\noi{\bf Theorem} 2.3.1 $~$ Relation (\ref{adbico}) is equivalent to,
$~ \forall a\in A~~ t_j(a)M_i{}^j=h_i(a)~$ i.e.
\eq
(\chi_j*a)M_i{}^j=a*\chi_i
{}~~~~~\mbox{i.e.}~~~~a_1\chi_j(a_2)M_i{}^j=\chi_i(a_1)a_2
\label{leftrightrans}
\en 
\noi{\sl Proof } 
Multiplying (\ref{adbico}) i.e. 
$a\!d_{\psi}\chi_i=\kappa'({\psi_1})\chi_i\psi_2=\psi(M_i{}^j)\chi_j$
by $\psi_0$, 
where $(\D'\otimes id)\D'(\psi)=\psi_0\otimes \psi_1\otimes\psi_2$,
we obtain the equivalent expression
\eq
\forall \psi\in \Un~~~~~~
\chi_i\psi=\psi_1\chi_j\;\psi_2(M_i{}^j)\label{Ucommrel}
\en
This relation gives the $q$-commutation relations between 
any $\psi\in \Un$ and
the $\chi_j$ elements.
We apply a generic element $a\in A$ to (\ref{Ucommrel}) and 
rewrite the right hand side as, $\forall \psi\in \Un,\;\forall a\in A$
\eqa
\langle \psi_1\chi_j\psi_2(\M{i}{j})\,,\,a\rangle
&=&\langle\psi_1\chi_i\otimes\psi_2\,,\,a\otimes \M{j}{i}\rangle\\
&=&\langle \psi\otimes\chi_j\,,\,a_1\M{j}{i}\otimes a_2\rangle
\ena
Since $\psi\in \Un$ and $a\in A$ are arbitrary elements and since 
the pairing $\Un\leftrightarrow A$ is nondegenerate
(we actually only need $\Un$ to separate the points of $A$) we conclude
\eq
\langle\chi_i\otimes\psi\,,\,a_1\otimes a_2\rangle=
\langle \psi\otimes\chi_j\,,\,a_1\M{j}{i}\otimes a_2\rangle~
\Leftrightarrow ~ a*\chi_i=(\chi_j*a)M_i{}^j
\en
that proves the theorem.
\cvd

{}Formula (\ref{leftrightrans}) relates the left invariant vectorfields
$t_i=\chi_i*$ to the right invariant ones $h_i=*\chi$ via the 
adjoint representation
$\M{i}{j}$. We can write
$h_i=t_j\dia \M{i}{j}$ where  
$(t_j\dia \M{i}{j})(a)\equiv t_j(a)\M{i}{j}~\forall a\in A$.
This formula is the analogue of (\ref{eta}).
The space of vectorfields  is analized in the next section.

\sk
\noi {\bf Note} 2.3.1 $~$
In \cite{Wor} Woronowicz has shown that bicovariant differential calculi
are in one-to-one correspondence with $a\!d$-invariant right ideals $R$ of 
$A$ :
$Ra\subset R~\forall a\in A$ $a\!d(R)\subset R\otimes A$.
These two conditions are slightly weaker than (\ref{ZERO})-(\ref{DUE}). 
Relation (\ref{UNO}) can also be written
$\D (T\oplus\{\epsi\})\subset {(T\oplus\{\epsi\})\otimes {\Un}} $ ,
where $T\oplus\{\epsi\}$ is the vector space spanned by $\chi_i$ and 
$\epsi$; therefore $(T\oplus\{\epsi\})$ is a right {\sl co-ideal}, 
it is the space orthogonal to the Woronowicz \cite{Wor} 
right {\sl ideal} $R\equiv\{a\in A~/~\epsi(a)=0~ \mbox{ and }~ T(a)=0\}$.
We have seen that relations (\ref{ZERO}) and (\ref{DUE}) imply 
(\ref{adbico}) this condition
is then equivalent to the $a\!d$ invariance of $R$:
$a\!d(R)\subset R\otimes A$. 

\noi Proof: 
$\forall r\in R, \forall \chi\in T, \forall \psi\in {\Un},$
\eq
0=\langle a\!d_{\psi}\chi\,,\,r\rangle=\langle \chi\otimes 
\psi\,,\, a\!d\, r\rangle \Rightarrow a\!d\,r\in (R\oplus \{I\})\otimes A
\Rightarrow a\!d\,r\in R\otimes A 
\en
where the last implication holds because 
$\langle\epsi\otimes id\,,\,a\!d\,r\rangle=0.$ Viceversa
$a\!d\,R\subset R\otimes A\Rightarrow a\!d_{\psi}\chi\in T\oplus \{\epsi\}
\Rightarrow a\!d_{\psi}\chi\in T$ since 
$\langle a\!d_{\psi}\chi\,,\, I\rangle =0~~\Box$
\sk

Notice that a Woronowicz type bicovariant differential calculus 
is given by a set $\{\chi_i\}$ of linear functionals on $A$  satisfying 
(\ref{leftrightrans}) and (\ref{leibdef2}), the full structure of
the dual Hopf algebra $\Un$ and the nondegeneracy of the 
$A\leftrightarrow \Un$ pairing is not needed to formulate the calculus.
In particular 
(\ref{1copM}) can be derived from  
(\ref{leftrightrans}).\footnote{\vskip -0.7em $~~~
\begin{array}{rcl}
\D\M{i}{j}&=&\chi_i(x_2^j)\D(\kappa(x_1^j)x_3^j)=
\kappa(x_1^j)\chi_i(x^j_2)x_3^j\otimes\kappa(x_0^j)x_4^j\nonumber\\
&=&\kappa(x^j_1)x_2^j\chi_n(x_3^j)\M{i}{n}\otimes \kappa(x_0^j)x_4^j
=\M{i}{n}\otimes\chi_n(x_3^j)\kappa(x_1^j)\epsi(x_2^j)x_4^j\nonumber\\
&=&\M{i}{n}\otimes\M{n}{j}~.
\end{array}
$}
\sk
\subsection{The space of $1$-forms and the exterior differential}

We now proceed in the construction of the differential calculus 
introducing the space $\Ga$ of $1$-forms.
The space of left-invariant $1$-forms $\Ginv$ is
defined as the space dual to that of the
tangent vectors $T$, let $\{\omega^i\}$ be the base of $\invG$ dual 
to $\{\chi_i\}$, we use the notation 
\eq
\langle \chi_i~,~\omega^j\rangle=
\delta^j_i~.\label{linvduality}
\en
By definition a generic $1$-form is then uniquely written as 
[see (\ref{rhoaom})] $\rho=a_i\omega^i$ i.e. 
the space of $1$-forms is the left $A$-module freely generated by the elements
$\om^i$. This corresponds to the classical property that the cotangent 
bundle of a group manifold is trivial.
The differential is defined by 
\eq
\forall a\in A  ~~~~da=(\chi_i*a)\;\omega^i\label{defdiff}~.
\en
\sk
\noi{\bf Note} 2.3.2 $~$ We can rewrite the exterior differential using 
right-invariant vectorfields:
\eq
da=(a*\chi_i)\kappa^{-1}(\M{j}{i})\om^j=(a*\chi_i)\eta^i
\en
where we have defined the $1$-forms $\eta^i=\kappa^{-1}(\M{j}{i})\om^j$.
It is easy to check that the $\eta^i$ are right invariant, see (\ref{eta}).
Using (\ref{etabb}) we also have:
\eq
da=-\eta^i(a*\kappa'(\chi_i))\label{etabbb}~.
\en
\sk
Any $1$-form $\rho\in\Ga$ can be written as $\rho=\sum_k a_k\,db_k$ for some 
$a_k,b_k\in A$ because we have the following expression for the left invariant 
$1$-forms:
\eq
\om^i=\kappa(y^i_1)dy^i_2 ~;~~\mbox{ in particular }~
\om^i=\kappa(x^i_1)dx^i_2 \label{kxdx}
\en
where $y^i$ and $x^i$ are the same elements that appear in (\ref{defMy}) and
(\ref{defM}). 

\noi Proof: $~\kappa(y_1^i)dy_2^i=\kappa(y_1^i)y_2^i\chi_j(y_3^i)\om^j=
\chi_j(y^i)\om^j=\om^i~~\Box$
\sk
On the space $\Ga$ of $1$-forms we can introduce a left 
and a right coaction of 
$A$ [the analogue of the left and right pullback on $1$-forms, see 
(\ref{410})-(\ref{411})]  by the following definitions, [see (\ref{linvom}) 
and 
(\ref{adjoint})]:
\eq
\DL (a_i\om^i) = \D(a_i) (I \otimes \om^i) ~~;~~~
\DR (a_i\om^i) = \D(a)(\om^j \otimes \M{j}{i}) \label{bicomodule}
\en
notice that the right $A$ coaction on the left-invariant $1$-forms $\Ginv$ 
corresponds to the (left) adjoint action of $\Un$ on $\Ginv$ : 
$\forall\/\psi\in \Un ,~a\!d_{\psi}\om^i\equiv\om^j\psi(\M{j}{i})$. 
We say that $\Ga$ is a {\sl left and right covariant} left-module because 
properties
(\ref{Dprop1}), (\ref{Dprop2}) and (\ref{bicovariance}) are satisfied. 
[Hint: use (\ref{1copM})].
We are now able to prove that the differential calculus is left and
right covariant i.e. it is  bicovariant on the left module $\Ga$:
\sk
\noi{\bf Proposition} 2.3.1 $~$ The exterior differential defined in 
(\ref{defdiff}) 
is bicovariant on the left module $\Ga$:
\eqa
& &\!\!\!\!\!\!\!\!\!\!\!\!
\DL(adb)=\D(a)(id \otimes d)\D(b),~~~\DL:\Ga\rightarrow A\otimes \Ga
~~~{\rm (left~ covariance)} \label{leftco2}\\
& &\!\!\!\!\!\!\!\!\!\!\!\!
\DR(adb)=\D(a)(d \otimes id)\D(b),~~~\DR:\Ga\rightarrow \Ga\otimes A
~~~{\rm (right~ covariance)} \label{rightco2}
\ena
Proof: since $\DL(a\rho)=\D(a)\DL(\rho)$ and
$\DR(a\rho)=\D(a)\DR(\rho)$  where $\rho$ is a 
generic $1$-form it is sufficient to prove:
\eq
\DL(db)=\D(\chi_i*b)I\otimes \om^i=b_1\otimes b_2\chi_i(b_3)\om^i=
(id\otimes d)\D(b)~;
\label{p1}
\en
\eqa
\!\!\!\!\!\!\!\!\DR(db)&=&\D(\chi_i*b)\om^j\otimes\M{j}{i}=b_1\om^j\otimes 
b_2\chi_i(b_3)\M{j}{i}=b_1\om^j\otimes \chi_j(b_2)b_3\nonumber\\
\!\!\!\!\!\!\!\!&=&(d\otimes id)\D(b)~~\label{p2}
\ena
\cvd
We have seen that from the closure of the $q$-Lie algebra $T$ under the
adjoint action of $\Un$ on $T$, equation (\ref{adbico}) [or from (\ref{ZERO}) 
and (\ref{DUE})] or equivalently from the relation (\ref{leftrightrans})
between left and right invariant vectorfields, one can construct an exterior 
differential $d~:~A\rightarrow\Ga$; where $\Ga$ is the left 
$A$-module of $1$-forms freely generated by the space of 
left-invariant one forms $\Ginv$. 
We have introduced a left and a right coaction of the quantum group 
$A$ on $\Ga$ and proved that the exterior differential is 
compatible with these coactions, 
see (\ref{leftco2})-(\ref{rightco2}).
This clarify the importance of the adjoint action in the construction of a 
differential calculus on a quantum group. 

We now analize the consequences of (\ref{UNO}) 
that, so far, we have never used  in 
this section. We show that (\ref{UNO}) is equivalent
to the Leibniz rule for the exterior differential and that it implies the 
$q$-antisymmetry of the $q$-Lie algebra.

\sk
\subsection{The Leibniz rule and the bicovariant bimodule of \\
\mbox{$1$-forms}}
\sk
\noi{\bf Lemma } The deformed Leibniz rule (\ref{UNO}) and 
(\ref{leftrightrans})
imply [see (\ref{propM})]:
\eq
\M{i}{j} (a * \f{i}{k})=(\f{j}{i} * a) \M{k}{i}, \label{propM22} 
\en
Proof: from  (\ref{leftrightrans}) we have, $ \forall\, b\in A$:
\eqa
\kappa(x_1^l)(\chi_j*x_2^lb)\M{i}{j}=\kappa(x_1^l)(x_2^lb*\chi_i) 
&\Leftrightarrow &\nonumber\\ 
\kappa(x_1^l)x_2^lb_1\chi_j(x_3^lb_2)\M{i}{j}=
\kappa(x_1^l)\chi_i(x_2^lb_1)(x_3^lb_2)
&\Leftrightarrow &\nonumber\\
(\f{l}{j}*b)\M{i}{j}=\kappa(x_1^l)[\chi_n(x_2^l)\f{n}{i}(b_1)+
\epsi(x_2^l)\chi_i(b_1)]x^l_3b_2
&\Leftrightarrow &\nonumber\\
(\f{l}{j}*b)\M{i}{j}=\kappa(x_1^l)\chi_n(x^l_2)
\f{n}{i}(b_1)x_3^lb_2
&\Leftrightarrow &\nonumber\\
(\f{l}{j} * b) \M{i}{j}=\M{n}{l} (b * \f{n}{i})& &\nonumber
\ena
where in the  
left hand side of the second passage
we have used $\f{l}{j}(b)=\chi_j(x^lb)$ that is obtained from 
(\ref{leibdef}) when $x^i=a$.
\cvd

The Leibniz rule on the left $A$-module $\Ga$
can be introduced if we know how to multiply $1$-forms with functions 
from  the right, i.e. if $\Ga$ is also a right module.
Consider the functionals $\f{i}{j}$ given in (\ref{leibdef}),
we define the following right product:
\sk
\noi{\sl Definition }
\eq
\om^ic=(\f{i}{j}*c)\om^j ~~~~,~~~~~(a_i\om^i)c=a_i(\f{i}{j}*c)\om^j
\label{rightpro2}
\en
the definition is well given because
\eq
\Dp (\f{i}{j})=\f{i}{k} \otimes \f{k}{j} ~~,~~~
\ep (\f{i}{j}) = \del^i_j~,
\en
[see (\ref{copf})-(\ref{couf})];
these two properties immediately follow, respectively, from the coassociativity
of the coproduct on the $\chi_i$ elements, and from $\chi_i(x^j)=\delta_i^j$.
\sk
We now prove the compatibility of (\ref{rightpro2}) with the left and 
right coactions $\DL$ and $\DR$; i.e. we prove that $\DL$ and $\DR$ are also, 
respectively,  left and  right coactions on $\Ga$ seen as a right module:
\eq
\forall \rho \in \Ga ,\forall a\in A ,~~~~~
\DL(\rho a)=\DL(\rho)\D(a)~~,~~~~
\DR(\rho a)=\DR(\rho)\D(a)\label{compLRR}~.
\en
Since any $\rho \in \Ga$ is of the form $\rho=a_i\om^i$ and since 
$\DL\om^i=I\otimes\om^i$, the only nontrivial espression in (\ref{compLRR}) is
$\DR(a_i\om^ia)=\DR(a_i\om^i)\D(a)$. As proven in Note 2.1.1, this is 
equivalent to
(\ref{propM22}) and we conclude that $\Ga$ is a bicovariant bimodule, i.e. 
that $\DL$ and $\DR$ are compatible with the bimodule structure of $\Ga$.

\sk
{}From the deformed Leibniz rule for the tangent vectors $\chi_i$, 
see (\ref{leibdef}) or (\ref{leibdef3}), and from (\ref{rightpro2}), the 
Leibniz rule for the 
exterior differential immediately follows. Viceversa, suppose that  
on $\Ga$ a right module structure can be introduced such that $d$ 
satisfies the Leibniz rule and is well defined in the following sense: 
$$
adb=a'db'~\Rightarrow ~(adb)c=(a'db')c ~~~\mbox{i.e.}~~~~ 
a\,d(bc)-ab\,dc=a'd(b'c)-a'b'dc~
$$
($adb$ is a shorthand notation for
$\sum_ka_kdb_k$; $a,b,c$ are generic elements).   
Then we can use the Leibniz rule to express the
right module structure on $\Ga$ as : $(adb)c=a\,d(bc)-ab\,dc$.
In the particular case $adb=\om^i$ we have, see (\ref{kxdx}),
\eqa
\om^i\,c
&=&\kappa(x^i_1)dx^i_2\;c=\kappa(x^i_1)d(x^i_2 c)=
\kappa(x^i_1)x_2^ic_1\chi_j(x_3^ic_2)\om^j\\
&=&c_1\chi_j(x^ic_2)\om^j
\equiv c_1\f{j}{i}(c_2)\om^j=(\f{j}{i}*c)\;\om^j
\label{omc}
\ena
where, in the last but one passage we have defined
$
\forall c\in A~~\f{j}{i}(c)\equiv\chi_j(x^ic)~.
$
Now from $d(ab)=d(a) \;b+adb$ and (\ref{omc}) we obtain:
\eq
\chi_i*ab=(\chi_j*a)(\f{j}{i}*b) + a(\chi_i*b)
\en
that is equivalent to (\ref{leibdef3}).
\sk
\subsection{$q$-antisymmetry of the $q$-Lie algebra bracket}

The coproduct (\ref{UNO}) implies that the espression $[\chi_i,\chi_j]$
is quadratic and $q$-antisym-
metric. We first write
\eq
[\chi_i,\chi_j]=\kappa'({\chi_j}_{_1})
\chi_i{\chi_j}_{_2}=\kappa'(\chi_l)\chi_i\f{l}{j}+\chi_i\chi_j~;
\en
now we apply $m(id\otimes \kappa')$ to 
$\D(\chi_l)=\chi_n\otimes f^n{}_l + \epsi\otimes \chi_l$
to obtain $\kappa'(\chi_l)=-\chi_n\kappa'(\f{n}{l})$ and therefore
\eq
\kappa'(\chi_l)\chi_i\f{l}{j}=-\chi_n\kappa'(\f{n}{l})\chi_i\f{l}{j}=
-\chi_n a\!d_{(\f{n}{j})}\chi_i=-\chi_n\f{n}{j}(\M{i}{l})\chi_l\label{trick}
\en
so that 
\eq
[\chi_i,\chi_j]=\chi_i\chi_j - 
\f{n}{j}(\M{i}{l})\chi_n\chi_l\label{quadratic}~.
\en
In the above framework it is easy to derive the bicovariance conditions
(\ref{bico1})-(\ref{bico4}); indeed
recalling equations
(\ref{RfM}): $\Rhat{ij}{kl} = \f{i}{l} (\M{k}{j})$ and 
(\ref{Cijk}):
$\C{ij}{k} = \chi_j(\M{i}{k})$  we immediately derive from 
(\ref{quadratic}) the bicovariance condition
(\ref{bico1}) and from 
[recall (\ref{adbico})]
$
a\!d_{(\f{n}{j})}\chi_i=\f{n}{j}(\M{i}{l})\chi_l
$
the bicovariance condition  (\ref{bico4}).
Relation (\ref{bico3}) can be derived applying to $(x^i\otimes id)$ the 
coproduct of (\ref{bico1}) and then using
(\ref{bico4}) and (\ref{Ccrelation}).
Finally (\ref{bico2}) can be obtained applying the functionals $\f{l}{n}$
to (\ref{propM22}).
\sk
\subsection{$q$-Jacoby identities}

We end this section briefly commenting on the $q$-Jacoby identities.
We have seen that the  $a\!d$ map given in (\ref{adjonU}) is a right action of
$\Un$ on $\Un$ [see (\ref{Uaction})]: 
$
\forall \varphi ,  \psi , \zeta~$ $ a\!d_{\psi\zeta}\varphi=
 a\!d_{\zeta}(a\!d_{\psi}\varphi)~.
$
We use the identity $\psi\zeta=\zeta_1 ad_{\zeta_2}\psi$  to find
\eq
a\!d_{\zeta}(ad_{\psi}\varphi)=a\!d_{\psi\zeta}\varphi=
a\!d_{\zeta_1 a\!d_{\zeta_2}\psi}\varphi=a\!d_{(a\!d_{\zeta_2}\psi)}
(a\!d_{\zeta_1}\varphi)\label{jacobigen}~.
\en
The above relation, written for elements $ \chi_i,\,\chi_j,\,\chi_l$
of the $q$-Lie algebra $T$, is the $q$-Jacobi identity
(with abuse of notation we define, $\forall\chi\in T$ , 
$[\chi,\f{i}{j}]\equiv a\!d_{(\f{i}{j})}\chi$, 
$[\chi,\epsi]\equiv a\!d_{\epsi}\chi=\chi$):
\eq
[[\chi_i,\chi_j],\chi_l]=[[\chi_i,{\chi_l}_{1}],[\chi_j,{\chi_l}_{2}]]
\label{derjac}~;
\en
it express the property that the bracket operation 
is a derivation of the $q$-Lie algebra $T$
(i.e. $\ad_{\zeta}$ is a generalized derivation with respect to the product
in $\Un$ given by the $\ad$ map).
Using the explicit coproduct expression
$\D(\chi_i)=\chi_j\otimes f^j{}_i + \epsi\otimes \chi_j$
in (\ref{derjac}), we obtain (\ref{qJacobi}).
There is also a second Jacobi identity :
$a\!d_{a\!d_{(\zeta}\psi)}\varphi=a\!d_{\kappa(\zeta_1)\psi\zeta_2}\varphi
=a\!d_{\zeta_2}(a\!d_{\psi}(a\!d_{\kappa(\zeta_1)}\varphi))$. On the $\chi$
elements it reads:
$[\chi_i,[\chi_j,\chi_l] ]=[[[\chi_i,\kappa({\chi_l}_{1})],\chi_j],
{\chi_l}_{2}]$. The two Jacobi identities are not independent
because $a\!d_{\zeta}(a\!d_{\psi}\varphi)=
a\!d_{[a\!d_{\zeta_2}\kappa^2(\zeta_1)\psi]}\varphi$.

Notice also  that
the map $a\!d$ is compatible with the product of $\Un$ in the sense that:
$a\!d_{\zeta}(\psi\varphi)=a\!d_{\zeta_1}(\psi)a\!d_{\zeta_2}(\varphi)$
(i.e. $\ad_{\zeta}$ is a generalized derivation with respect to the product 
of $\Un$); in particular $a\!d_{\chi_l}(\chi_i\chi_j)
\equiv[\chi_i\chi_j,\chi_l]$=$[\chi_i,\chi_s]\f{s}{l}(\M{j}{n})
\chi_n+\chi_i[\chi_j,\chi_l]$.
\sk
\subsection{$*$-Structure}
Given a $*$-Hopf algebra $A$ we have a canonical $*$-structure on the dual 
$\cal U$. This is compatible with the quantum Lie algebra $T$ if
$T^*\subseteq T$, i.e., if the tangent vectors 
$(\chi_i)^*$ are linear combination  
of the $\chi_i$, so that we have a real form of the $q$-Lie algebra.
In this case we can construct a differential calculus that is real:
\eq
 (db)^* = db^* ~\mbox{ and more in general }~  (adb)^*=db^*\,a^*~~.
\en
We now explicitly perform this construction. 
There is a canonical $*$-structure on the space of $1$-forms. We first 
define a 
$*$-involution on $\invG$ via the expression:
\eq
\forall\chi\in T\,,~\forall\om\in\invG\,,~~~ 
\le \omega^*\,,\,\chi\re\equiv-\overline{\le\omega\,,\,\chi^*\re}\label{om*chi}
\en
and generalize it to $\Ga$ as 
\eq
\forall a\in A\,,~\forall\omega\in \invG\,,~~~~~(a\omega)^*=\omega^*a^*~.
\en
We have to check that these definitions are consistent with the bicovariant 
bimodule structure on $\Ga$. 
We recall from Section 1.3 
that the $*$ operation becomes the hermitian conjugation 
${}^{\dag}$ when we realize the 
elements of $A$ and $U$ as operators on Hilbert space (in the \qone limit
the elements of $A$ commute and  correspond to diagonal operators,
so that the $*$-operation becomes complex conjugation).
It is then natural to consider a basis of antihermitian 
$q$-Lie algebra generators $\chi_i$: $(\chi_i)^*=-\chi_i$, 
then the dual basis of
$1$-forms is real: ${\om^{i*}}=\om^i$. 
{}From $\D(\chi_j)=-\D(\chi_j^*)=-(\D\chi_j)^{*\otimes *}$
it follows that  the $\f{i}{j}$ are real and we have
\eq
{\om^{i*}}a=(\f{i}{j}*a){\om^{j*}}\label{omreal}
\en 
that is compatible with ${\om^{i*}}=\om^i$. 
Proof of (\ref{omreal}): $(\omega^i{}^*a)^*=a^*\omega^i=$
$\omega^j\f{i}{j}{\circ}\kappa^{-1}*(a^*)=$
$\omega^ja_1^*\f{i}{j}(\kappa^{-1}(a_2^*))=$
$\omega^ja_1^*\overline{\f{i}{j}(a_2)}=$
$[(\f{i}{j}*a)\omega^j{}^*]^*$,
where we have used (\ref{U*A}).
\sk
\noi Also 
the elements $\M{i}{j}$ are real so that $\DR(\om^i)=\om^j\otimes \M{i}{j}$
is well defined.\\ 
Proof: 
$\M{i}{j}{}^*(\psi)\chi_j=$
$\overline{\M{i}{j}{}(\kappa'^{-1}(\psi^*))}\chi_j=$
$-[\M{i}{j}(\kappa'^{-1}(\psi^*))\chi_j]^*=$
$-[ad_{\kappa'^{-1}(\psi^*)}\chi_i]^*=$
$-[\psi_2^*\chi_i\kappa'^{-1}(\psi_1^*)]^*=$
$\M{i}{j}(\psi)\chi_j$.
\sk
We are now able to show that the 
differential, given by $da=(\chi_i*a)\om^i$ is real:
\eqa
(da^*)^*&=&\om^i(\chi_i*a^*)^*=[\f{i}{j}*(\chi_i*a^*)^*]\om^j
=(\f{i}{j}*a_1)\overline{\chi_i(a_2^*)}\om^j\nonumber\\
&=&(\f{i}{j}*a_1)(\kappa'^{-1}\,\chi_i^*)(a_2)
=-[\f{i}{j}(\kappa'^{-1}\chi_i)*a]\om^j
=(\chi_j*a)\om^j\nonumber\\
&=&da\nonumber
\ena
where we have used $m[\tau(\kappa'^{-1}\otimes id)\D'(\chi)]=\epsi'(\chi)=0$, 
a consequence of 
(\ref{kappadelta}).
\sk
\sk
\noi {\large{\bf {Conclusions}}}
\nopagebreak
\sk
\nopagebreak
We have seen, from     
(\ref{ZERO})--(\ref{DUE}), or more in general from  
(\ref{UNO}) and (\ref{adbico}),
or from (\ref{UNO}) and (\ref{leftrightrans}), that
the construction of the differential calculus associated 
to the $q$-Lie algebra $T$ spanned by the $\chi_i$ elements 
is quite straighforward, the main ingredients are
\sk
{\bf{i) }} the left invariant vectorfields $t_i=\chi_i*~,$ with deformed 
Leibniz rule:
$t_i(ab)=t_j(a)\f{j}{i}(b)+at_i(b)$
 
{\bf{ii) }} the adjoint representation $\M{i}{j}$ defined 
via (\ref{adbico}) [or explicilty via the coordinates $x^i$ (\ref{defM})].
The adjoint representation satisfies $\D({M_i}^j)=M_i{}^k
\otimes M_k{}^j$ and $\epsi(M_i{}^j)=\delta^i_j$.

{\bf{iii) }} the space of left invariant $1$-forms, 
defined as the space dual to that of the
tangent vectors: $\,\langle \chi_i~,~\omega^j\rangle=\delta^j_i$.
A generic $1$-form is then given by
$\rho=a_i\omega^i$.
[The space of $1$-forms is the bicovariant bimodule
freely generated by the $\omega^i$ with 
$\omega^ia=(f^i{}_j*a)\omega^j,~ $ 
$\D{}_L \omega^i\equiv I\otimes 
\omega^i,~ \D{}_R \omega^i\equiv\omega^j \otimes M_j{}^i$].

{\bf{iv)}} The differential, defined by $da=(\chi_i*a)\omega^i$; it
satisfies the undeformed Leibniz rule.
\sk
\sk
\noi {\bf Note} 
2.3.3 $~$ Following \cite{Bonechi} we here briefly characterize in a
cohomological context the bicovariant differential calculus on quantum groups.
{}For any $\om\in \Ginv,~ a\in A$, we define the left and right products
\[
a.\om=\epsi(a)\om~~~~~~;~~~~~~~~~~
 \om.a=\kappa(a_1)\om a_2~.
\]
In particular $\om^i.a=\kappa(a_1)\om^i 
a_2=\kappa(a_1)(\f{i}{j}*a_2)\om^j$ $=$ $\f{i}{j}(a)\om^j$,
this shows that the  right product is well defined and, using the
property (\ref{copf}) and (\ref{couf}) of the $f$ functionals, 
that $\Ginv$ is a left and right $A$-module.
The projection $P~:~\Ga\rightarrow \Ginv$ defined in Note 2.1.1
is an epimorphism between the two bimodules $\Ga$ and $\Ginv$.
\sk
We now characterize the differential $d$ through a $1$-cocycle of
the Hochschild coboundary operator $\delta$ relative to the $A$-bimodule 
$\Ginv$.
Given an algebra $\cal A$ and a bimodule $M$ over $\cal A$, a Hochschild
$k$-cochain $C \in C^k({\cal A}, M)$  is a $k$-multilinear map from
${\cal A}\otimes{\cal A}\otimes ...{\cal A}$ ($k$-times) to $M$, with  
$C^0({\cal A}, M)=M$.
The coboundary operator $\delta~:~ C^k({\cal A}, M)\rightarrow  
C^{k+1}({\cal A}, M)$ is defined by
$$
\begin{array}{rl}
\!\delta C (a_1,\ldots a_{k+1})=&\!\!a_1.C(a_2,\ldots a_{k+1})\\
                  &\!\!+\sum_{i=1}^k (-1)^iC(a_1,\ldots a_ia_{i+1},\ldots 
a_{k+1})+ 
(-1)^{k+1}C(a_1,\ldots a_k).a_{k+1}
\end{array}
$$
and satisfies $\delta^2=0$. [ We have denoted by ``$.$'' 
the multiplication in the bimodule $M$]

To a bicovariant differential calculus with differential $d$, we associate 
the map 
\[
\begin{array}{rl}
c~:~ & A\rightarrow \Ginv\\
     & a \mapsto    P(da)=\kappa(a_1)da_2~.
\end{array}
\]
It is easy to see that $\delta c=0$, i.e., $c$ is a $1$-cocycle :
$c(ab)=c(a)b+ac(b)$ . \mbox{Vice}versa, given a $1$-cocycle $c$ we 
immediately obtain
a left covariant differential calculus defining
\[
da\equiv a_1c(a_2)~.
\]
[Proof of the left covariance: $\DL(da)=a_1\otimes a_2c(a_3)=(id\otimes d)\D 
a$].

The right covariance (\ref{rightco}) of a differential calculus is equivalent 
to
the following property for the cocycle $c$:
\eq
\forall \psi\in U ~~~~~~~(id\otimes \psi_2)\DR[c(a*\psi_1)]=c(\psi*a)
\label{Uinv}
\en
[Hint: (\ref{rightco}) is equivalent to 
$(\kappa(a_1)\otimes id)\DR(da_2)=\kappa(a_1)da_2\otimes a_3$; apply  
$(id\otimes \psi)$ to this last expression].
Therefore $1$-cocycles satisfying (\ref{Uinv}) are in one-to-one
correspondence with bicovariant differential calculi.         
In the notations of Subsection 2.4.5  (\ref{Uinv}) reads 
$\ell_{\psi_2}c(a*\psi_1)=
c(\psi*a)$.
In \cite{Bonechi} it is shown that there is a one-to-one correspondence between
bicovariant $A$-bimodules on $\Ga$ and $\cal D$-bimodules on $\Ginv$ where 
$\cal D$
is the quantum double of $A$; moreover the cocycles satisfying (\ref{Uinv})
correspond to cocycles $c$ in the set of the Hochschild cochains 
$C^1({\cal D},\Ginv)$ that have the simple property $c(U)=0$.

Notice also that 
the $0$-cochains are the left invariant one forms: $C^0(A,\Ginv)=\Ginv$, 
it can be checked  
that the coboundary of any  bi-invariant $1$-form, i.e. of any 
right invariant $0$-cochain satisfies (\ref{Uinv}) and therefore 
defines a bicovariant differential calculus.
The differential studied in (\ref{defd}) and (\ref{defd2}) corresponds to the
$1$-cocycle $\delta({-1\over \la}\tau)$. Indeed  
$\delta({-1\over \la}\tau)\,(a)=
{-1\over\la}(a.\tau -\tau.a) ={-1\over \la}\epsi(a)\tau+
{1\over\la}\kappa(a_1)\tau a_2$
and 
$da=a_1\,\delta({-1\over\la}\tau)\,(a_2)=
{1\over\la}(\tau a-a\tau)$
as in (\ref{defd2}). 
\sk
The differential calculus on classical Lie groups corresponds to a nontrivial
$1$-cocycle since in the commutative case $\delta\om=0$ for any $\om\in\Ginv$. 
All the differential calculi we will examine   
correspond to $1$-cocycles that are coboundaries, the only exception being
those  on the twisted homogeneous and inhomogeneous orthogonal 
groups of sections 4.6 and 4.7.

The existence of a bi-invariant $1$-form trivializes the calculus from the 
Hochschild cohomology viewpoint (it is associated to a coboundary) but  
it is interesting geometrically and for physical speculations because 
introduces a discretized geometry; indeed $d$ given by (\ref{defd2}) is a 
finite difference operator as well as the partial derivatives $\chi_i$ 
in (\ref{defchi2}).

%
%
%
%
%
%
%


\def\spinst#1#2{{#1\brack#2}}
\def\sk{\vskip .4cm}
\def\noi{\noindent}
\def\om{\omega}
\def\al{\alpha}
\def\be{\beta}
\def\ga{\gamma}
\def\Ga{\Gamma}
\def\de{\delta}
\def\del{\delta}
\def\la{\lambda}
\def\lam{{1 \over \lambda}}
\def\alb{\bar{\alpha}}
\def\beb{\bar{\beta}}
\def\gab{\bar{\gamma}}
\def\deb{\bar{\delta}}
\def\alp{{\alpha}^{\prime}}
\def\bep{{\beta}^{\prime}}
\def\gap{{\gamma}^{\prime}}
\def\dep{{\delta}^{\prime}}
\def\rhop{{\rho}^{\prime}}
\def\taup{{\tau}^{\prime}}
\def\rhopp{\rho ''}
\def\thetap{{\theta}^{\prime}}
\def\imezzi{{i\over 2}}
\def\unquarto{{1 \over 4}}
\def\unmezzo{{1 \over 2}}
\def\epsi{\varepsilon}
\def\we{\wedge}
\def\th{\theta}
\def\de{\delta}
\def\cony{i_{\de {\vec y}}}
\def\Liey{l_{\de {\vec y}}}
\def\tv{{\vec t}}
\def\Gt{{\tilde G}}
\def\deyv{\vec {\de y}}
\def\part{\partial}
\def\pdxp{{\partial \over {\partial x^+}}}
\def\pdxm{{\partial \over {\partial x^-}}}
\def\pdxi{{\partial \over {\partial x^i}}}
\def\pdy#1{{\partial \over {\partial y^{#1}}}}
\def\pdx#1{{\partial \over {\partial x^{#1}}}}
\def\pdyx#1{{\partial \over {\partial (yx)^{#1}}}}

\def\qP{q-Poincar\'e~}
\def\R#1#2{ R^{#1}_{~~~#2} }
\def\Rp#1#2{ (R^+)^{#1}_{~~~#2} }
\def\Rm#1#2{ (R^-)^{#1}_{~~~#2} }
\def\Rinv#1#2{ (R^{-1})^{#1}_{~~~#2} }
\def\Rpm#1#2{(R^{\pm})^{#1}_{~~~#2} }
\def\Rpminv#1#2{((R^{\pm})^{-1})^{#1}_{~~~#2} }
\def\RRpm{R^{\pm}}
\def\RRp{R^{+}}
\def\RRm{R^{-}}
\def\Rhat#1#2{ \Rh^{#1}_{~~~#2} }
\def\Rhatinv#1#2{ (\Rh^{-1})^{#1}_{~~~#2} }
\def\Z#1#2{ Z^{#1}_{~~~#2} }
\def\Rt#1{ {\hat R}_{#1} }
\def\Rh{\hat R}
\def\ff#1#2#3{f_{#1~~~#3}^{~#2}}
\def\MM#1#2#3{M^{#1~~~#3}_{~#2}}
\def\cchi#1#2{\chi^{#1}_{~#2}}
\def\ome#1#2{\omega_{#1}^{~#2}}
\def\RRhat#1#2#3#4#5#6#7#8{\Rh^{~#2~#4}_{#1~#3}|^{#5~#7}_{~#6~#8}}
\def\RRhatinv#1#2#3#4#5#6#7#8{(\Rh^{-1})^
{~#2~#4}_{#1~#3}|^{#5~#7}_{~#6~#8}}
\def\U#1#2#3#4#5#6#7#8{U^{~#2~#4}_{#1~#3}|^{#5~#7}_{~#6~#8}}
\def\CC#1#2#3#4#5#6{\Cb_{~#2~#4}^{#1~#3}|_{#5}^{~#6}}
\def\cc#1#2#3#4#5#6{C_{~#2~#4}^{#1~#3}|_{#5}^{~#6}}

\def\C#1#2{ {\bf C}_{#1}^{~~~#2} }
\def\c#1#2{ C_{#1}^{~~~#2} }
\def\q#1{   {{q^{#1} - q^{-#1}} \over {q^{\unmezzo}-q^{-\unmezzo}}}  } 
\def\Dmat#1#2{D^{#1}_{~#2}}
\def\Dmatinv#1#2{(D^{-1})^{#1}_{~#2}}
\def\f#1#2{ f^{#1}_{~~#2} }
\def\F#1#2{ F^{#1}{\!}_{#2} }
\def\T#1#2{ T^{#1}_{~~#2} }
\def\Ti#1#2{ (T^{-1})^{#1}_{~~#2} }
\def\Tp#1#2{ (T^{\prime})^{#1}_{~~#2} }
\def\TP{ T^{\prime} }
\def\M#1#2{ M_{#1}^{~#2} }
\def\qm{q^{-1}}
\def\um{u^{-1}}
\def\vm{v^{-1}}
\def\xm{x^{-}}
\def\xp{x^{+}}
\def\fm{f_-}
\def\fp{f_+}
\def\fn{f_0}
\def\D{\Delta}
\def\Mat#1#2#3#4#5#6#7#8#9{\left( \matrix{
     #1 & #2 & #3 \cr
     #4 & #5 & #6 \cr
     #7 & #8 & #9 \cr
   }\right) }   
\def\Ap{A^{\prime}}
\def\Dp{\Delta^{\prime}}
\def\Ip{I^{\prime}}
\def\ep{\epsi^{\prime}}
\def\kp{S^{\prime}}
\def\kpm{S^{\prime -1}}
\def\km{S^{-1}}
\def\gp{g^{\prime}}
\def\Cb{{\bf C}}
\def\qone{$q \rightarrow 1~$}
\def\Fmn{F_{\mu\nu}}
\def\Am{A_{\mu}}
\def\An{A_{\nu}}
\def\dm{\part_{\mu}}
\def\dn{\part_{\nu}}
\def\Ana{A_{\nu]}}
\def\Bna{B_{\nu]}}
\def\Zna{Z_{\nu]}}
\def\dma{\part_{[\mu}}
\def\qsu{$[SU(2) \times U(1)]_q~$}
\def\suq{$SU_q(2)~$}
\def\su{$SU(2) \times U(1)~$}
\def\gij{g_{ij}}
\def\L{{\cal L}}

\def\LL{L^*}
\def\ll#1{L^*_{#1}}
\def\RR{R^*}
\def\rr#1{R^*_{#1}}

\def\Lpm#1#2{(L^{\pm})^{#1}_{~~#2}}
\def\LLpm{L^{\pm}}
\def\LLp{L^{+}}
\def\LLm{L^{-}}
\def\Lp#1#2{(L^{+})^{#1}_{~~#2}}
\def\Lm#1#2{(L^{-})^{#1}_{~~#2}}

\def\gu{g_{U(1)}}
\def\gsu{g_{SU(2)}}
\def\tg{ {\rm tg} }
\def\Fun{$Fun(G)~$}
\def\qonelim{\stackrel{q \rightarrow 1}{\longrightarrow}}
\def\viel#1#2{e^{#1}_{~~{#2}}}

\sect{More $q$-Geometry: vectorfields, inner derivative and Lie derivative}

In the previous section we have studied the space of left invariant 
vectorfields i.e. the $q$-tangent space, we now
construct, for a generic quantum group, the space of vectorfields.
Its elements are products of elements of the quantum group
itself with left invariant vectorfields.
We study the duality between vectorfields and 1-forms and generalize
the construction to tensorfields.
As in the classical case, using the duality between covariant and 
contravariant tensorfields, we can introduce the contraction
operator. This is defined on the space of covariant tensorfields,
and therefore acts in particular on forms; indeed the algebra of forms,
as defined in (\ref{wedge1}), is a subalgebra of the algebra of covariant
tensorfields. We then prove that the contraction operator is a (inner) 
derivation in the space of forms. On the other hand the right action of the 
$q$-group on the space of $1$-forms naturally define the Lie 
derivative along left invariant vectorfields. The Cartan identity 
$\ell_{t_i}= i_{t_i}  d +d  i_{t_i}$ is proven and 
the Cartan calculus of inner derivatives, Lie derivatives and
the exterior derivative generalized to $q$-group geometry.
Not all properties of the classical Cartan Calculus can however
be generalized, while the contraction operator is defined for general 
vectorfields, there is no completely satisfactory 
expression for the Lie derivative along general vectorfields $V$. We
propose the  definition  $\ell_{V}\equiv i_{V}  d +d  i_{V}$
and analize and  discuss its properties.
The topics discussed in this sections have been studied in
\cite{SWZ1}, \cite{AC} and more extensively in 
\cite{Paolo},
\cite{SWZ3}, 
\cite{SW}, \cite{S2},
\cite{PaoloPeter}, \cite{Vladimirov}.
We follow \cite{PaoloPeter} and  \cite{AC}.
Here we give a  self-contained exposition, and all 
the theorems are proved starting from only one data:
a bicovariant differential calculus on a generic Hopf algebra.

\subsection{From Left invariant Vectorfields to general Vectorfields}
In this subsection we  study the space $\Xi$ of vectorfields 
over the generic Hopf algebra $A$ 
defining a
right product between left invariant 
vectorfields and  elements of $A$.
\sk
In the commutative case a generic vectorfield can be written in the
form $f^it_i$ where $\{t_i\}~\,i=1,\,\ldots,n$ is a basis of left invariant
vectorfields and $f^i$ are $n$ smooth functions on the group
manifold.
In the commutative case $f^it_i=t_i\dia f^i$ i.e. left and right
products (that we have denoted with $\dia$) are  the same, indeed 
$ (t_i\dia f^i)(h)\equiv t_i(h)f^i=f^it_i(h)$.
These considerations lead to the following definition.

Let ${t_i} = {\chi_i *}$ be  a basis in  $\invXi$, the space of left
invariant vectorfields, and 
let $a^i$, $\,i=1,\,\ldots,n$ be  generic
elements of $A$:
\sk
\noi {\bf Definition} 
\eq
\Xi \equiv \{V~ /~~ V: A \longrightarrow A ~;~ V = t_i\dia a^i\}~,
\label{defffX}
\en
where the definition of the right product $\dia$ is given below:
\sk
\noi {\bf Definition}
\eq
\forall a,b \in A  ,\forall t \in {}_{inv}\Xi
~~~~~~(t \dia a) b \equiv t(b)a = (\chi * b)a ~.
\en
The product $\dia$ has a natural generalization to the whole $\Xi$ :
\eq \begin{array}{rcl} 
\dia~~:~~& \Xi \times A \longrightarrow & \Xi \nonumber\\
         & ~ (V,a) ~        \longmapsto     &  V\dia a 
\end{array} 
{}~~~ \mbox{ where }~~~ \forall b \in A~~~ (V \dia a)(b) \equiv V(b)a ~.
\en

\sk
\noi It is easy to prove that $( \Xi ,\dia )$ is a right $A$-module:
\eq
V\dia (a + b) = V\dia a + V\dia b~;~~~ 
  V\dia (ab) = (V \dia a) \dia b~;~~~
  V\dia(a+b) = V\dia a + V\dia b
\label{amodule}
\en
(we have also $V\dia \lambda a = \lambda V \dia a$ with $\lambda\in    
\mbox{\boldmath$C$}\,$).

\noi For example $V\dia (ab) = (V \dia a) \dia b$ because
\[\forall c\in A~~~~ [(V\dia a)\dia b]c = [(V\dia a)(c)]b = (V(c)a)b =
V(c)ab = [V\dia ab]c
.\]
\noi Notice that to distinguish the elements $V\dia (ab) \in \Xi $ and $V(
ab)\in A$ we have not omitted the simbol $\dia$ representing the
right product.
\sk
$\Xi$ is the analogue of the space of derivations on the ring
$C^{\infty}(G)$ of the
smooth functions on the group $G$. 
Indeed we have:

\eq V(a+b)=V(a)+ V(b)~~,~~~~~ V(\lambda a) = \lambda
V(a)~~~~\mbox{Linearity}\label{Linearity}
\en
\eq
{}\:~~~V(ab) \equiv (t_i\dia c^i)(ab) = t_j(a)(f^j{}_i*b)c^i + aV(b)
~~~~~~\mbox{ Leibniz rule}\label{Leibnizrule}
\en

\noi in the classical case $t_j(a)(f^j{}_i *b)c^i =
V(a)b${}$\,$ (recall $\f{j}{i}=\delta^j_i\epsi\; ;~\:\epsi *b=b).$

\sk
We have seen the duality between $\invG$ and $\invXi$. We now extend it to
$\Ga$ and $\Xi$, where $\Ga$ is seen as a left $A$-module (not
necessarily a bimodule) and $\Xi$ is our right $A$-module.

\sk

\noi{\bf Theorem} 2.4.1 $~$ 
There exists a unique map
$$
\langle~~,~~\rangle ~:~~~\Ga \times \Xi \longrightarrow A
$$
\indent such that:
\sk
\noi 1) $ \forall\: V \in \Xi$; the application
\[
 \langle ~~,V \rangle  ~:~~\Ga \longrightarrow A
\]
is a left $A$-module morphism, i.e. is linear and $\langle
a\rho, V\rangle=a\langle\rho, V\rangle$.

\noi 2) $\forall \rho\in \Ga $; the application
\[ \langle \rho,~~\rangle ~:~~\Xi \longrightarrow A\]
is a right $A$-module morphism, i.e. is linear and $\langle
\rho, Vb\rangle=\langle\rho, V\rangle b$.

\noi 3) Given  $\rho\in\Ga$
\eq
\langle\rho,~~\rangle=0
  ~\Rightarrow ~ \rho=0~, \label{Duality1}
\en
\sk
\noi where $\langle\rho,~~\rangle=0$ means 
$\langle\rho,V\rangle = 0 ~~\forall V \in \Xi  $.

\noi 4) Given  $V \in \Xi$
\eq
\langle~~,V\rangle=0 
~\Rightarrow ~ V=0~, \label{Duality2}
\en
\sk
\noi where $\langle~~,V\rangle=0 $ 
means $\langle\rho,V\rangle = 0 ~~\forall \rho \in \Ga  .$

\noi 5) On $\invG \times \invXi$ the bracket $ \langle~~,~~\rangle$ 
acts as the one introduced in the previous section.

\sk
\noi {\sl Remark } Properties 3) and 4) state that 
$\Ga$ and $\Xi$ are dual $A$-moduli, in the sense that
they are dual with respect to $A$.

\sk
\noi {\sl Proof}

\noi Properties 1), 2) and 5) uniquely characterize this map . To prove 
the existence of such a map we show that the following bracket
\sk
\noi {\bf Definition}       
\eq
\langle\rho,V\rangle =\langle a_{\alpha}db_{\beta}\,,\,V\rangle\equiv 
a_{\al}V(b_{\al})~,
\en
 where $a_{\al},b_{\al}$ are elements of $A$ such that 
$\rho =a_{\al}db_{\al}$,
satisfies 1),2) and 5).
\sk
\noi We first verify that the above definition is well given, that is:
\[
\mbox{Let } \rho =a_{\al}db_{\al}=a'_{\beta}db'_{\beta}~~\mbox{ then }  
~~a_{\al}V(b_{\al})= a'_{\beta}V(b'_{\beta})
~.\]
\noi Indeed, since \[a_{\al}db_{\al} = 
a'_{\beta}db'_{\beta}  ~\Leftrightarrow ~ 
       a_{\al}t_i(b_{\al})\om^i = a'_{\beta}t_i(b'_{\beta})\om^i
~\Leftrightarrow ~ 
       a_{\al}t_i(b_{\al}) = a'_{\beta}t_i(b'_{\beta})\]\\
{[we used  the uniqueness of the decomposition (\ref{rhoaom})]}
\noi the definition is consistent
 because $$a_{\al}V(b_{\al})=a'_{\beta}V(b'_{\beta})
~\Leftrightarrow ~ a_{\al}t_i(b_{\al})c^i=a'_{\beta}t_i(b'_{\beta})c^i $$
where $V=t_i\dia c^i.$
\sk
\noi {\sl Property} 1) is trivial since $a\rho =a(a_{\al}db_{\al}) 
= (aa_{\al})db_{\al}.$ \\
\noi {\sl Property} 2) holds since
\[\langle \rho,V\dia c\rangle = a_{\al}(V\dia c)(b_{\al}) =
a_{\al}V(b_{\al})c=\langle \rho,V\rangle c
~.\]

\noi {\sl Property} 5).
Let $\{\om^i\}$ and $\{t_i\}$ be dual bases in $\invG$ and $\invXi$. 
Since  $\om^i \in
\Ga ~,~ \om^i = a_{\al}db_{\al}$ for some $a_{\al}$ and $b_{\al}$ in $A$.
We can also write $\om^i = a_{\al}db_{\al}= a_{\al}t_k(b_{\al})\om^k$ ,
so that, due to the uniqueness of the decomposition (\ref{rhoaom}),  
we have 
\[~~~~~~~~~~~~~~~ a_{\al}t_k (b_{\al})=\delta_k^i I~~~~~~~~~~~~(I 
\mbox{ unit of} A);\]

\noi we then obtain  \[\langle\om^i,t_j\rangle
=a_{\al}t_j(b_{\al})=\delta^i{}_j I~.\]

\noi{\sl Property} 3).
Let  $\rho=a_i \om^i \in\Ga~.$ \\
If $\langle\rho , V\rangle =0 ~~\forall V\in \Xi$, in particular
$\langle\rho , t_j\rangle = 0 ~~\forall j=1,...,n$; then 
$\mbox{$a_i\langle\om^i , t_j\rangle=0 \Leftrightarrow $} \\ a_j=0$ , 
and therefore $\rho =0~.$
\sk
\noi Property 4).
Let $V = t_i\dia a^i \in \Xi~.$ \\
If $\langle \rho , V\rangle =0 ~~\forall\rho\in\Ga$, in particular
$\langle\om^j  , V\rangle = 0 ~~\forall j=1,...,n$; then 
$\mbox{$\langle\om^j , t_i\rangle a^i=0 \Leftrightarrow $} \\ a^j=0$ , 
and therefore $V =0~.$

\cvd
By construction every $V$ is of the form 
\[ V=t_i\dia a^i.\]
\indent We can now show the unicity of such a decomposition.
\sk
\noi {\bf Theorem} 2.4.2 $~$ Any $V \in \Xi$ can be uniquely written in 
the form
\[V=t_i\dia a^i
\]
\noi {\sl Proof}\\
Let $ V=t_i\dia a^i =t_i\dia a'^i$ then 
\[\!\!\!\!\!\!\!\!\!\!\!\!\!\!\!\!\!\!\!
\forall i=1,\,\ldots,n~~~~~~ a^i=  \langle\om^i , t_j\rangle a^j = 
 \langle\om^i , V\rangle = \langle\om^i , t_j\rangle a'^j = a'^i~. 
\]
\cvd
\noi Notice that once we know the decomposition of $\rho$ and $V$ in terms of 
$\om^i$ and $t_i$, the evaluation of $ \langle~~,~~\rangle$ is trivial:
\[
\langle\rho,V\rangle = \langle a_i\om^i,t_j\dia b^j\rangle = a_i\langle
\om^i,t_j\rangle b^j = a_ib^i~.
\]
Viceversa from the previous theorem $V=t_i\dia\langle\om^i,V\rangle$ and
$\rho = \langle\rho,t_i\rangle \om^i~.$
\sk
\sk
\sk
$\!{}$ We conclude 
this section by remarking the three different $\mbox{ways 
of looking at $\Xi.$}$

\noi (I)${}~{}~{}$   \spz $\Xi$ as the 
set of all deformed derivations over $A$ [see (\ref{defffX}),
(\ref{Linearity}) and 
(\ref{Leibnizrule})].

\noi (II)${}~{}$  \spz $\Xi$ as the right $A$-module freely  generated by
the {\sl elements} $t_i\,, ~~i=1,\,\ldots,n$. The latter is  the 
set of all the {\sl formal} products and sums of the type $t_ia^i$,  
where $a^i$ are generic elements of $A$. 
Indeed, by virtue of Theorem 2.4.2, the map that associates to each 
$V=t_i\dia a^i$ in $\Xi$ the corresponding element $t_ia^i $
 is an isomorphism between right $A$-moduli.

\noi (III) $\!$\spz  $\Xi$ as $ \Xi '=
\left\{ U ~:~\Ga \longrightarrow A ,
~U \mbox{ linear and } U(a\rho)=aU(\rho) ~\forall a\in A\right\}$, i.e.
$\Xi$ as the dual (with respect to $A$) of the space of 1-forms $\Ga$. 
The space $\Xi '$ has a trivial right $A$-module structure: $(Ua)(\rho)
\equiv U(\rho)a~.$
$\Xi$ and $\Xi '$ are isomorphic right $A$-moduli because of property 
(\ref{Duality2}) which states that to each $\langle~~,V\rangle \: :~\Ga
\rightarrow A$ 
 there 
corresponds one and only one
$V$.$[\langle~~,V\rangle=\langle~~,V'\rangle\Rightarrow V=V'].$ 
Every $U\in \Xi '$ is of the form $U = \langle~~, V\rangle$; more precisely,
if $a^i$ is such that 
$ U(\om^i)=a^i$ then $U = \langle~~~, t_i\dia a^i\rangle~.$
\sk
These three ways of looking at $\Xi$ will correspond to different
aspects of the Cartan Calculus: the Lie derivatives $\ell_V$ will
generalize (I), inner derivations $i_V$ will correspond to (III), while
the transformation properties of $\ell_V$ and $i_V$ are governed by (II).
\sk
\subsection{Bicovariant Bimodule Structure}

In Section 2.1 we have studied the space $\Ga$ of $1$-forms, we have seen
that $\Ga$ is a bimodule 
over $A$ because there is a right and a left product
between elements of $\Ga$ and of $A$.
The left and the right  product are related by 
$\om^ia=(f^i{}_j*a)\om^j$.
Since the coactions $\DL$ and $\DR$ are compatible with the bimodule 
structure and since they commute:
$$
(id \otimes \DR) \DL = (\DL \otimes id) \DR 
$$
the bimodule $\Ga$ is a  bicovariant bimodule (cf. Note 2.1.1). 
\sk
In the previous subsection 
we have studied the right product $\dia$ and we have seen 
that $\Xi$ is a right module over $A$ [see (\ref{amodule})].
Here we introduce a left product and a left and right coaction
of the Hopf algebra $A$ on $\Xi$. The left and right coactions $\DS$ and 
$\DD$ are the $q$-analogue of the push-forward of tensorfields on a 
group manifold. Similarly to $\Ga$ also $\Xi$ is a bicovariant bimodule.
\sk
The construction of the left product on $\Xi$, of the right coaction $\DD$
and of the left coaction $\DS$ will be effected  along the lines of 
Woronowicz' 
Theorem 2.5 in \cite{Wor}, whose  statement can be explained  in the 
following steps (cf. Note 2.1.1):
\sk
\noi {\bf Theorem} 2.4.3   
Consider the {\sl symbols} ${t_i}~ (i=1,\,\ldots,n)$
 and let $\Xi$ be the right $A$-module freely generated by them:

\[
\Xi \equiv \{t_i a^i ~/~~a^i \in A\}
\]
Consider functionals $\FF{i}{j}:~A\longrightarrow \mbox{\boldmath$C$}~$ 
satisfying [see (\ref{propf1}) and (\ref{propf2})]
\begin{eqnarray}
& & \FF{i}{j} (ab)= \FF{i}{k} (a) \FF{k}{j} (b) \label{propF1}\\
& & \FF{i}{j} (I) = \del^i_j~~~~~~~~~~~~~~~. \label{propF2}
\end{eqnarray}
\indent Introduce a left product via the definition [see (\ref{aom})]
\sk
\noi {\sl Definition}  
\eq 
b(t_ia^i) \equiv t_j{[}(\FF{i}{j}\circ {\kappa}^{-1})* b{]}a^i~. \label{bta}
\en
It is easy to prove that

\noi i) \spz $\Xi$ is  a bimodule over $A$. (A proof of this first 
statement as well as of the following ones is contained in \cite{Wor}).

\rightline{$\Box$}
\sk
Introduce an action (push-forward) of the Hopf algebra $A$ on $\Xi$
\sk
\noi {\sl Definition}
\eq
\DS(t_ia^i)\equiv(I\otimes t_i)\D(a^i)\label{DDta}~.
\en
It follows that

\noi ii) \spz $(\Xi, \DS)$ is a left covariant
 bimodule over $A$, that is
\[
 \DS(a V b)=\D (a)\DS(V)\D(b)~;~~~
 (\epsi \otimes id) \DS (V)=V~;~~~
 (\D \otimes id)\DS=(id\otimes\DS)\DS ~. 
\]
\rightline{$\Box$}
\sk
Introduce $n^2$ elements $\N{i}{j}\in A$ satisfying 
[see (\ref{propM}),(\ref{copM}) and 
(\ref{couM})]
\begin{eqnarray}
& \N{i}{j} (a * \FF{i}{k})=(\FF{j}{i} * a) \N{k}{i} & \label{propN} \\
 & \Delta (\N{j}{i}) = \N{j}{l} \otimes \N{l}{i} &  \label{copN}\\
 & \epsi (\N{j}{i}) = \delta^i_j~~~~~~, &  \label{couN}
\end{eqnarray}
and introduce $\DD $ such that [see (\ref{adjoint})] 
\sk
\noi {\sl Definition}
\eq
\DD(a^it_i)\equiv\D(a^i)t_j\otimes \N{j}{i}\label{DS}~.
\en
Then it can be proven that 

\noi iii) \spz The elements  [see (\ref{eta})]
\eq
h_i \equiv t_j{\kappa}(\N{j}{i})\label{heich}
\en
are right invariant: $\DD(h_i) = h_i\otimes I$.
Moreover any $V \in \Xi$ can be expressed in a unique 
way respectively as $V=h_ia^i$ and as
$V=b^ih_i$, where $a^i, b^i \in A$.

\rightline{${}~~~~~~~~~~~~~~~~~~~~~~~~~\Box$}

\noi iv) \spz $ (\Xi, \DD)$ is a right covariant bimodule over $A$,
that is 
\[
 \DD(a V b)=\D (a)\DD(V)\D(b)~;~~~ (id \otimes \epsi) \DD (V)=V~;~~~
 (id \otimes \D)\DD =(\DD\otimes id)\DD~.
\]{}
\rightline{$\Box$}

\noi v) \spz The left and right covariant bimodule $(\Xi,\DS,\DD) $ 
is a bicovariant 
bimodule, that is left and
right coactions are compatible:
\[
(id\otimes\DD)\DS=(\DS\otimes id)\DD~.
\]
\cvd

In the previous section we have seen [remark (II)] that the space of 
vectorfields
$\Xi$ is the free right $A$-module generated by the symbols $t_i$, so
that the above theorem applies to our case.

There are many bimodule structures (i.e. choices of $\FF{i}{j}$) $\Xi$ can
be endowed with.
Using the fact that $\Xi$ is dual to $\Ga$ we request compatibility with
the $\Ga$ bimodule.\\
In the commutative case $\langle f\om^i,t_j\rangle=\langle \om^if,t_j\rangle=
\langle\om^i,ft_j\rangle=\langle\om^i,t_jf\rangle.$\\
In the quantum case we know that $\langle a\om^i,t_j\rangle = 
\langle\om^i,t_j\dia a\rangle$ and we require 
\eq
\langle\om^i a,t_j\rangle =
 \langle\om^i,at_j\rangle~; \label {wat}
\en
this condition uniquely determines the bimodule structure of $\Xi$.
Indeed we have 
 \eqa
 \langle\om^i,at_j\rangle &=& \langle\om^i a,t_j\rangle =
\langle (\f{i}{k} *a )\om^k,t_j \rangle = (\f{i}{k} *a)\langle\om^k,t_j\rangle=
\f{i}{k} * a\,\delta^k_j = \delta^i_l\f{l}{j} * a \nonumber \\
& = &\langle\om^i,t_l \dia
(\f{l}{j} * a)\rangle
\ena
so that
\eq
at_i = t_j \dia (\f{j}{i} * a~). \label{ffS}                    
\en
We then define
\eq
\FF{i}{j} \equiv \f{j}{i}\circ {\kappa}
\en
it follows that $\FF{i}{j} \circ {\kappa}^{-1} = \f{j}{i}$  and
(\ref{ffS}) can be rewritten [see (\ref{aom}) and (\ref{bta})]
\eq
at_j=t_j\dia {[}(\FF{i}{j}\circ {\kappa}^{-1})*a{]}~.
\en
\noi{\bf Theorem} 2.4.4 $~$ The functionals $ \FF{i}{j} $ 
satisfy conditions (\ref{propF1}) and
(\ref{propF2}).
\sk
\noi{\sl Proof }
The first condition $\FF{i}{j}(I)=\delta^j_i $ 
holds trivially.\\
The second one is also easily checked:
\sk
\noi
\[\begin{array}{rcl}
\!\! \FF{i}{j}(ab) &\!\!\! =\!\!\! & (\f{j}{i}\circ {\kappa})(ab) 
= \f{j}{i}{[}{\kappa}(b){\kappa}(a){]}  
 =\f{j}{k}{[}{\kappa}(b){]}\f{k}{i}{[}{\kappa}(a){]}
 =\f{k}{i}{[}{\kappa}(a){]}\f{j}{k}{[}({\kappa}(b){]} \\
               &\!\!\! =\!\!\! &\FF{i}{k}(a) \FF{k}{j}(b)
\end{array}
\]
{} \cvd

So far $\Xi$ has a bimodule structure.
We now define a left coaction $\DS$ so that 
$\Xi$ becomes a left covariant bimodule.
The left invariant vectorfields were characterized in
(\ref{qleftinv}) and (\ref{leftinvclass})
through their  action $t_i(a) = \chi_i * a$ on functions;
following the same 
derivation as in (\ref{linvom})--(\ref{llinvom})
the left invariant property of the  $t_i$ can also be expressed
via the coaction $\DS$ as defined in 
(\ref{DDta}): 
\eq
\DS(t_ia^i)\equiv(I\otimes t_i)\D(a^i)\,.
\en
  
Similarly the right coaction $\DR$ is defined to act trivially on the 
right invariant vectorfields $h_i=*\chi_i$. This uniquely defines
the elements $\N{l}{k} \in A$; indeed we want relation (\ref{heich})
and (\ref{leftrightrans}) to coincide and therefore:
\eq
\N{l}{k}={\kappa}^{-1}(\M{k}{l})~.\label{NeqSM}
\en

\noi Notice that (\ref{NeqSM}) implies
\eq
\langle t_i , \om^j \rangle = \delta_i^j = \langle h_i , \eta^j \rangle,
\en
where $t_i$ and $\om^j$ are left-invariant and $h_i$ and $\eta^j$ are
the canonically associated right-invariant objects; see
(\ref{heich}) and (\ref{eta}).
Notice also that  $\M{k}{l}$ and $\f{i}{j}$ are dual, and likewise 
$\N{l}{k}$ and $\FF{j}{i}$, 
in the sense that 
$\f{i}{j}(\M{k}{l}) = \FF{j}{i}(\N{l}{k})=\Lambda^{il}{}_{kj}$ with
$\Lambda^{il}{}_{kj} =\delta^i_j\delta^l_k $ when $q = 1$. 
\sk
\noi {\bf Theorem} 2.4.5 $~$  The $\N{l}{k}$ elements defined above 
satisfy relations
(\ref{couN}),
(\ref{copN}) and
(\ref{propN}):
$$
1)~ \epsi (\N{j}{i}) = \delta^i_j
{}~~~~~2)~ \Delta (\N{j}{i}) = \N{j}{l} \otimes \N{l}{i}
{}~~~~~3)~ \N{i}{k} (a * \FF{i}{j})=(\FF{k}{i} * a) \N{j}{i}
$$
\sk
\noi {\sl Proof}
\sk
\noi $1)~$ This expression is trivial.
\sk
\noi $2)~$ Use $N^i{}_j = {\kappa}^{-1}M_j{}^i$ and $
\D\circ {\kappa}^{-1} = \sigma \circ ({\kappa}^{-1} \otimes 
{\kappa}^{-1})\circ \D$,
where $\sigma_A$ is the flip map in $A\otimes A$.
\sk
\noi $3)~$ We know that [see (\ref{propM})] 
\[
\!\!\!\!\!\!\forall\/ a\in A~~~~~~~\M{i}{j} (a * \f{i}{k})=
(\f{j}{i} * a) \M{k}{i} 
\]
or equivalently,
\[
~~~ ~~~~~~\M{i}{j} [{\kappa}(a) * \f{i}{k}]=[\f{j}{i} * 
{\kappa}(a)] \M{k}{i} ~.
\]
Now
$$\begin{array}{lll}
\!\!\!{[}{\kappa}(a)*\f{i}{k}{]}&\!\!=\!\!&(\f{i}{k}\otimes
id)\D{[}{\kappa}(a){]}=(id\otimes\f{i}{k})({\kappa}\otimes 
{\kappa})\D(a)={\kappa}(id\otimes
\f{i}{k}\circ {\kappa})\D(a)\\
&\!\!=\!\!&{\kappa}(\FF{k}{i}*a)~.
\end{array}$$
Similarly,
\[
{[}\f{j}{i}*{\kappa}(a){]}={\kappa}(a*\FF{i}{j})~.
\]
So we can write
\[
{\kappa}(a*\FF{i}{j})\M{k}{i} = \M{i}{j}{\kappa}(\FF{k}{i}*a) 
\]
for all $a\in A.$ 
Applying $\:{\kappa}^{-1}$ to both members of this last expression we obtain
relation $3).$
\sk \cvd
Following Theorem 2.4.3 the construction of the bicovariant bimodule
$\Xi$ is now easy and straightforward, and
we can conclude that $(\Xi ,\DS ,\/\DD)$ is a bicovariant
bimodule.
\sk
We end this subsection observing that in the expression 
(\ref{defdiff}) for the exterior differential,
elements of $\Xi$ and $\Gamma$ make a
joint appearance. To be still able to talk about transformation
properties of such expressions we need to combine the previously
introduced coactions into one object, $\Delta_A$, simply by putting
$\Delta_A \equiv {}_\Xi\Delta$ on $\Xi$ and $\Delta_A \equiv
{}_\Gamma\Delta$ on $\Gamma$ and requiring $\Delta_A$ to be an
algebra homomorphism. From this definition we
get the following important corollary:
\sk
\noindent {\bf Corollary.} The expression $\omega^i t_i$ in (\ref{donf})
is invariant in the sense that
$$
\Delta_A( t_i\,\omega^i) = {}_\Xi\Delta(t_i){}_\Gamma\Delta(\omega^i) =
 t_k \omega^j \otimes  N^k{}_i M_j{}^i =t_i \omega^i  \otimes 1.
$$
Similar statements apply to ${}_A\Delta$.

\sk
Notice that, since  Theorem 2.4.3  completely characterizes a bicovariant
bimodule all the formulas containing the symbols $\f{i}{j}$ or
$\M{k}{l}$ or elements of $\Ga$ are still valid under the
substitutions $\f{i}{j}\rightarrow \FF{i}{j},~~\M{k}{l}
\rightarrow\N{k}{l} $ and $ \Ga\rightarrow \Xi$.
\sk

\subsection{Tensorfields}

The construction completed for vectorfields is readily generalized 
to $p$-times contravariant tensorfields.
We proceed as in (\ref{DLGaGa})--(\ref{DRGaGaGa}) and
define $\Xi\otimes\Xi$ to be the space of all elements that
can be written as finite sums of the kind $\sum_i V_i\otimes V'_i $
with $V_i,V'_i\in\Xi$. The tensor product (in the algebra $A{)}$ 
between $V_i$ and $V'_i$ has the following properties:\\
 ${}~~~~~~V\dia a\otimes
V'=V\otimes aV' ~, a(V\otimes V') =(aV)\otimes V' $ and $ (V\otimes
V')\dia a=V\otimes(V'\dia a) $\\
 so that $\Xi\otimes\Xi $ is naturally a bimodule
over $A$.\\
Left and right coactions on $\Xi\otimes\Xi$ are defined by:
\eq
\DS (V \otimes V')\equiv   V_1   {V'} _1 \otimes   V_2 \otimes 
  {V'}_2,~~~\DS: \Xi \otimes \Xi \rightarrow A\otimes\Xi\otimes\Xi
\label{DSXiXi}
\en
\eq
\DD (  V \otimes   {V'})\equiv   V_1 \otimes   {V'}_1 \otimes   V_2 
  {V'}_2,~~~\DD: \Xi \otimes \Xi \rightarrow \Xi\otimes\Xi\otimes A
\label{DDXiXi}
\en
\noi where as usual $  V_1$, $  V_2$ etc. are defined by
\eq
\DS (  V) =   V_1 \otimes   V_2,~~~  V_1\in A,~  V_2\in \Xi
\en
\eq
\DD (  V) =   V_1 \otimes   V_2,~~~  V_1\in \Xi,~  V_2\in A~.
\en
\noi More generally, we can introduce the coaction of $\DS$ 
on $\Xi^{\otimes p}\equiv\underbrace{\Xi \otimes \Xi \otimes \cdots
\otimes \Xi}_{\mbox{$p$-times}}$ as 
\[
\DS (  V \otimes   {V'} \otimes \cdots \otimes   {V''})\equiv 
  V_1   {V'}_1 \cdots   {V''}_1 \otimes   V_2 \otimes 
  {V'}_2\otimes \cdots \otimes   {V''}_2 
\]
\eq
\DS~:~~ \Xi^{\otimes p} \longrightarrow 
A\otimes\Xi^{\otimes p}~;
\label{DSXiXiXi}
\en
\[
\DD (  V \otimes   {V'} \otimes \cdots \otimes   {V''})\equiv 
  V_1 \otimes   {V'}_1 \otimes \cdots \otimes   {V''}_1 \otimes   V_2  
  {V'}_2 \cdots   {V''}_2 
\]
\eq
\DD~:~~ \Xi^{\otimes p} \longrightarrow 
\Xi^{\otimes p}\otimes A~~.
\label{DDXiXiXi}
\en

\noi Left invariance on $\Xi\otimes\Xi$ is naturally defined as
$\DS (  V \otimes   {V'}) = I \otimes   V \otimes   {V'}$ (similar 
definition for right invariance), so that for example $t_i \otimes
t_j$ is left invariant, and is in fact a left invariant basis for $\Xi 
\otimes \Xi$: each element can be written as $t_i\otimes t_j\dia  a^{ij}$ in a 
unique way.

 It is not difficult to show that $\Xi \otimes \Xi$ is a bicovariant bimodule. 
In the same way also
$(\Xi^{\otimes p},\DS,\DD)$ is a bicovariant bimodule. 
\sk
Any element $v\in \Xi^{\otimes p}$ can be written as
$v=t_{i_1}\otimes\ldots t_{1_p}\dia b^{i_1...i_p}$ in a unique way,
similarly any element $\tau\in\Ga^{\otimes n}$, 
the $n$-times tensor product of $1$-forms,
can be written as
$\tau=a_{i_n...i_1}\om^{i_n}\otimes\ldots \om^{i_1}$ in a unique way.
\sk
It is now possible to generalize the previous bracket 
$\langle~~,~~\rangle\,
:\;\Ga\times\Xi\rightarrow A\:$ to $\Ga^{\otimes n}$ 
and $\Xi^{\otimes p}\;{}$:

\eq
\begin{array}{rccll}
\!\!\!\!\!\!\!\langle~~,~~\rangle ~:~~&\Ga^{\otimes n} \times \Xi^{\otimes p}&
\longrightarrow & A &{} \\
\!\!\!\!\!\!\!&(\tau,v) &\longmapsto
&\langle\tau,v\rangle &=a_{i_n...i_1}\langle\om^{i_n}\otimes
...\om^{i_1}\, ,\,t_{j_1}\otimes ...t_{j_p}\rangle b^{j_1...j_p}\\
\!\!\!\!\!\!\!&{}&{}&{}&=a_{{i_n}...i_1}\om^{i_n}\otimes ...\om^{i_{p+1}} 
b^{i_1...i_p}
\end{array}\label{uniprod}
\en
where $\Ga^{\otimes 0}
\equiv A$, $\Ga^{\otimes 1} \equiv \Ga$ and we have defined
\eqa
\!\!\!\!\!\!\!\!\!\!\langle\om^{i_n}\otimes
...\om^{i_1}\, ,\,t_{j_1}\otimes ...t_{j_p}\rangle&\equiv &
\om^{i_n}\otimes ...\om^{i_{p+1}}\langle\om^{i_p}\otimes ...\om^{i_1}\, 
,\,t_{j_1}\otimes ...t_{j_p}\rangle\label{mirrorcoup}\\
&\equiv & \om^{i_n}\otimes ...\om^{i_{p+1}}
\langle \om^{i_1},t_{j_1}\rangle ...\langle \om^{i_p},t_{j_p}\rangle 
\nonumber\\
&=&\delta^{i_1}_{j_1} ...\delta^{i_p}_{j_p}\om^{i_n}
\otimes ...\om^{i_{p+1}} \; ~.
\nonumber
\ena
Using definition (\ref{mirrorcoup}) it is easy to prove that 
\eq
\langle\tau a,v\rangle=\langle\tau,av\rangle \label{pass}~,
\en
namely
\[\begin{array}{rcl}
\langle\om^{i_n}\otimes\ldots\om^{i_{p+1}}\otimes\om^{i_p}
\otimes\ldots &&\!\!\!\!\!\!\!\!\!\!\!\!\!\!\om^{i_1}a \,,
t_{j_1}\otimes\ldots t_{j_p}\rangle=\\
&=&\om^{i_n}\otimes ...\om^{i_{p+1}}
(\f{i_p}{k_p}*\ldots 
\f{i_1}{k_1}*a)
\langle\om^{k_p}\otimes\ldots\om^{k_1} , t_{j_1
}\otimes\ldots t_{j_p}\rangle
\end{array}
\]
\[
\begin{array}{rcl}
\langle\om^{i_n}\otimes\ldots\om^{i_{p+1}}\otimes
\om^{i_p}\otimes\ldots &&\!\!\!\!\!\!\!\!\!\!\!\!\!\!\om^{i_1} \,, a
t_{j_1}\otimes\ldots t_{j_p}\rangle=\\
&=&\om^{i_n}\otimes ...\om^{i_{p+1}}\langle\om^{i_p}\otimes\ldots\om^{i_1} , 
t_{l_1}\otimes\ldots t_{l_p}\rangle 
(\f{l_p}{j_p}*\ldots \f{l_1}{j_1}*a)
\end{array}
\]
and these last two expressions are equal if and {\sl only if}
(\ref{mirrorcoup}) holds.

Therefore we have also shown that definition (\ref{mirrorcoup}) is the
only one compatible with property (\ref{pass}), i.e. property
(\ref{pass}) uniquely determines the coupling between $\Xi^{\otimes}$
and $\Ga^{\otimes}.$

It is easy to prove that
the bracket $\langle~~,~~\rangle$ extends
to  $\Ga^{\otimes p}$ and $\Xi^{\otimes p}$ the duality between $\Ga$
and $\Xi$.
\sk
More generally we can define $\Xi^{\otimes}\equiv
A\oplus\Xi\oplus\Xi^{\otimes 2}\oplus\Xi^{\otimes 3}...$
 to be the algebra  
of contravariant tensorfields.
The coactions $\DS$ and $\DD$ have a natural generalization to
$\Xi^{\otimes}$ so that we can conclude that $(\Xi^{\otimes},\DS,\DD)$
is a bicovariant graded algebra, the graded algebra of tensorfields
over the ring ``of functions on the group'' $A$, with the left and
right ``push-forward'' $\DS$ and $\DD$.
Similarly 
$\Ga^{\otimes}$ is the bicovariant graded algebra of covariant tensorfields
on $A$.

\subsection{Contraction operator}

In this subsection we study the contraction operator $i_V$  
along a generic vectorfield $V\in\Xi$
and we prove that it acts as a (deformed) derivative operator on the 
space of $1$-forms.
The definition of  the contraction operator 
$i_V $ with
$V\in\Xi$ is based on equation (\ref{uniprod}).
For a generic vectorfield  $V=b^j t_j$ we define:
\sk
\noi {\bf Definition of right inner derivative}
$$
(\vartheta)
\!\stackrel{\leftarrow}{\imath_V}
 \equiv \langle \vartheta , V \rangle
{}~~~~~~~~~~~~~~~~~\forall \vart\in\Ga^{\otimes}~.
$$
this definition  applies when $\vart$ is a 
generic covariant tensorfield and 
in particular when $\vart$ is  a generic form.
\sk
\noi {\bf Theorem} 2.4.6 $~$ 
The contraction operator $\!\stackrel{\leftarrow}{\imath_V}$ satisfies the
following properties:\\
$a, a_{i_1\ldots i_n} \in A$;
$V=t_i\dia b^i$; $\lambda \in$ {\bf C}, property d) holds only if $\vartheta, 
\vartheta'$ are forms;
\begin{enumerate}
\item[a)] $(\vartheta)\!\stackrel{\leftarrow}{\imath_V} = 
(\vartheta)\,{\!\stackrel{\leftarrow}{\imath_{t_j}}}{}_{\!\dia b^j} =
(\vartheta)\!\stackrel{\leftarrow}{\imath_j}\,b^j$
\item[b)] $(a)\!\stackrel{\leftarrow}{\imath_V}=0$
\item[c)] $( \omega^{j} )\!\stackrel{\leftarrow}{\imath_V}= b^j$
\item[d)] $(a_{i_1\ldots i_n}\om_1 \wedge \ldots \wedge \om_n)
\!\stackrel{\leftarrow}{\imath_V}=
      \begin{array}[t]{l}
      (a_{i_1\ldots i_n}\om_1\wedge \ldots \wedge 
      \om_{s-1})\wedge(\om_{s}\wedge\ldots
      \wedge\om_n)\!\stackrel{\leftarrow}{\imath_{t_i}}\\
      +(-1)^{n-s+1}(a_{i_1\ldots i_n}\om_1 \wedge \ldots \wedge \om_{s-1})
      \!\stackrel{\leftarrow}{\imath_{t_j}}
      \wedge\; f^j{}_i * (\om_{s}\wedge\om_{s+1}\ldots\wedge\om_n)
      b^i
      \end{array}$

\item[e)] $(a\vartheta   +\vartheta')\!\stackrel{\leftarrow}{\imath_V} =
               a\, ( \vartheta )\!\stackrel{\leftarrow}{\imath_V} + 
(\vartheta')\!\stackrel{\leftarrow}{\imath_V}$
\item[f)]  $( \vartheta a)\!\stackrel{\leftarrow}{\imath_V} = 
(\vartheta)\!\stackrel{\leftarrow}{\imath_{t_j}}(f^j{}_i * a)\,b^i$
\item[g)] $\!\stackrel{\leftarrow}{\imath_{\lambda V}} = 
\lambda \!\stackrel{\leftarrow}{\imath_V}$
\end{enumerate}
\sk
\noi {\sl Proof }\\
Properties a), b), c), f), g) are direct consequences of (\ref{uniprod}),
e) follows from (\ref{pass}). To proof d) we have to use the definition
of the wedge product (\ref{wedge1})--(\ref{wedge3}): first we note that
(in tensor product notation)
$$\begin{array}{rcl}
(\om_1 \wedge \ldots \wedge \om_n)\!\stackrel{\leftarrow}{\imath_V}
& = &  W_{1\ldots n}(\om_1 \otimes \ldots \otimes \om_n)
\!\stackrel{\leftarrow}{\imath_V}\\
& = &   W_{1 \ldots n}\;\om_1 \otimes \ldots \otimes \om_{n-1} 
(\om_n)\!\stackrel{\leftarrow}{\imath_V} \\
& = & {\cal I}_{1 \ldots n}\;\om_1 \wedge \ldots \wedge\om_{n-1} 
(\om_n)\!\stackrel{\leftarrow}{\imath_V} ,
\end{array}$$
where we have used (\ref{wedge2}) in the last step --- in index notation:
$$(\om^{i_1}\wedge \ldots \wedge \om^{i_n})\!
\stackrel{\leftarrow}{\imath_{t_i}}
= {\cal I}_{j_1 \ldots j_{n-1}i}^{i_1 \ldots i_n}\,\om^{j_1}\wedge 
\ldots \wedge \om^{j_{n-1}} .$$
Next we can show
$$\begin{array}{rcl}
\!\!\!\f{i_{s-1}}{i} * (\omega^{i_{s} } \wedge \ldots\wedge\omega^{i_{n}} )
&\!\!\! =\!\!\! & \omega^{k_{s-1} } \wedge \ldots\wedge\omega^{k_{n-1}}
      \f{i_{s-1}}{i}(M_{k_{s-1}}{}^{i_{s}} \cdots M_{k_{n-1}}{}^{i_n})\\
& = &      \Lambda^{i_{s-1} i_{s}}{}_{k_{s-1}l_{s}}{} 
\Lambda^{l_{s} i_{s+1}}{}_{k_{s} l_{s+1}}{} \ldots
      \Lambda^{l_{n} i_n}{}_{k_{n-1} i}{}
 \omega^{k_{s-1} } \wedge  \omega^{k_{s} }\wedge\ldots\omega^{k_{n-1}}
\end{array}$$
that in tensor product notation can be written:
$$
(\om_{s-1})\!\stackrel{\leftarrow}{\imath_{t_j}}
f^j{}_i *(\om_s\wedge\om_{s+1}\wedge\ldots\om_n)
=\Lambda_{s-1,s}\Lambda_{s,s+1}\ldots\Lambda_{n-1,n}
(\om_{s-1}\wedge
\om_s\wedge\ldots\om_{n-1})\;(\om_n)\!\stackrel{\leftarrow}{\imath_{t_i}}
$$ 
Finally we utilize the decomposition
property (\ref{wedge4}) and associativity of the wedge product
(in tensor product notation)
$$\begin{array}{rcl}
(\om_1 \wedge &&\!\!\!\!\!\!\!\!\! 
\ldots \wedge \om_n)\!\stackrel{\leftarrow}{\imath_V} ~=\\
\!\!\!= &&\!\!\!\!\!\!\! {\cal I}_{1 \ldots n}
      \om_1 \wedge  \ldots \om_{n-1}(\om_n) \!\stackrel{\leftarrow}{\imath_V}\\
\!\!\!\!\!= &&\!\!\!\!\!\!\! [{\cal I}_{s \ldots n} +
      (-1)^{n-s+1} {\cal I}_{1 \ldots s-1}\Lambda_{s-1,s}\cdots\Lambda_{n-1,n}]
       (\om_1 \wedge \ldots \om_{s-1})
       \wedge (\om_{s}\wedge\ldots \om_{n-1})\,(\om_n)\!
\stackrel{\leftarrow}{\imath_V}
       \\
\!\!\!\!= &&\!\!\!\!\!\!\! \begin{array}[t]{l}
      (\om_1\wedge \ldots \wedge \om_{s-1})\wedge(\om_{s}\wedge\ldots
      \wedge\om_n)\!\stackrel{\leftarrow}{\imath_V}\\
      +(-1)^{n-s+1} {\cal I}_{1 \ldots s-1}(\om_1 \wedge \ldots \om_{s-2})
\wedge  \Lambda_{s-1,s}\Lambda_{s,s+1}\cdots\Lambda_{n-1,n}
      (\om_{s-1}\wedge\ldots \om_{n-1})
      (\om_n)\!\stackrel{\leftarrow}{\imath_{t_i}}
      b^i
      \end{array}\\

\!\!\!\!= &&\!\!\!\!\!\!\!\begin{array}[t]{l}
      (\om_1\wedge \ldots \wedge \om_{s-1})\wedge(\om_{s}\wedge\ldots
      \wedge\om_n)\!\stackrel{\leftarrow}{\imath_{t_i}}\\
      +(-1)^{n-s+1}{\cal I}_{1 \ldots s-1}(\om_1 \wedge \ldots \wedge 
\om_{s-2})
      \,(\om_{s-1})\!\stackrel{\leftarrow}{\imath_{t_j}}
      \wedge\; f^j{}_i * (\om_{s}\wedge\om_{s+1}\ldots\wedge\om_n)
      b^i
      \end{array}\\

\!\!\!\!= &&\!\!\!\!\!\!\!\begin{array}[t]{l}
      (\om_1\wedge \ldots \wedge \om_{s-1})\wedge(\om_{s}\wedge\ldots
      \wedge\om_n)\!\stackrel{\leftarrow}{\imath_{t_i}}\\
      +(-1)^{n-s+1}(\om_1 \wedge \ldots \wedge \om_{s-1})
      \!\stackrel{\leftarrow}{\imath_{t_j}}
      \wedge\; f^j{}_i * (\om_{s}\wedge\om_{s+1}\ldots\wedge\om_n)
      b^i
      \end{array}
\end{array}$$
With property e) this proves d).
\cvd
\sk
\noi {\sl Remark } A slight generalization of property d) for two
generic forms $\vartheta$ and $\vartheta'$ is also true [use f)]:
\eq
(\vartheta \wedge \vartheta')\!\stackrel{\leftarrow}{\imath_V} =
\vartheta\wedge (\vartheta')\!\stackrel{\leftarrow}{\imath_V} 
+ (-1)^{deg}(\vart') \!\stackrel{\leftarrow}{\imath_{t_j}}\wedge
(\f{j}{i}* \vartheta') b^i~.\label{itwotheta}
\en
\sk
We have defined the exterior differential as an operator 
acting from the left to the right, indeed we have the following 
behaviour under grading, as opposed to the one in (\ref{itwotheta}):
\eq
d(\vartheta\wedge\vartheta')=d\vartheta\wedge\vartheta' +(-1)^{deg(\vartheta)}
\vartheta\wedge d\vartheta'\label{dgraduato}~.
\en
In order to find the Cartan expression for the Lie derivative :
$\ell_{V}= i_{V} d +d  i_{V}$, we therefore have to introduce an inner 
derivation
$i_{V}$ that has the same behaviour as in (\ref{dgraduato}). 
This motivates the following 
\sk
\noi {\bf Definition of inner derivative}
\eq
i_V(\vartheta)\equiv (-1)^{deg(\vartheta-1)}
(\vartheta)\!\stackrel{\leftarrow}{\imath_V}~~~~~~~~~~~~~~~
\forall\vart\in\Ga^{\otimes}~.\label{contraction}
\en
this definition applies when $\vart$ is a 
generic covariant tensorfield and 
in particular when  $\vart$  generic form.
We immediately have:
\sk
\noi {\bf Theorem} 2.4.7 $~$ 
The ${i_V}$ contraction operator satisfies the
following properties:\\
$a, a_{i_1\ldots i_n} \in A$;
$V=t_i\dia b^i$; $\lambda \in$ {\bf C}, property d) holds only if $\vartheta, 
\vartheta'$ are forms
\begin{enumerate}
\item[a$'$)]  $i_V(\vartheta)=i_{t_j\dia{b^j}}(\vartheta)
=i_{t_j}(\vartheta)b^j$
\item[b$'$)] $
        i_{V}(a)=0 $
\item[c$'$)] $
        i_{V} ( \omega^{j} )= b^j  $
\item[d$'$)] $  
             i_{V} (\vartheta\wedge\vartheta')=
             i_{t_j}(\vartheta)\wedge(\f{j}{i}*{\vartheta'})b^i
             + (-1)^{deg(\vartheta)}\vartheta\wedge i_V(\vartheta')$
\item[e$'$)] $
        i_{V} ( a \vartheta +\vartheta') =
        a i_{V} ( \vartheta ) + i_{V}(\vartheta')
        $
\item[f$'$)] $  i_{V} (\vartheta a) = i_{t_j}(\vartheta)(\f{j}{i}*a)b^i$
\item[g$'$)] $  i_{\lambda V} =
                  \lambda i_{V}$
\end{enumerate}

\cvd
\sk
Notice that properties a$'$) e$'$) and f$'$) reduce in
the commutative case to the familiar formulae:
\[
i_{fV}\vart=f\,i_V\vart ~~ \mbox{ and } ~~ i_V(f\vart) = fi_V\vart ~.
\]
It is also straightforward to see that 
\eq
(id \otimes i_t) \DL(\vart) = \DL i_t(\vart)~~~~~~~~~~~
\forall\vart\in \Ga^{\otimes}~.
\en
This formula $q$-generalizes the classical 
commutativity of $i_t$ with the left coaction $\DL$, when $t$ is a left
invariant vectorfield.
\sk
\subsection{Lie Derivative and Cartan identity}

In  (\ref{leibdef2}), or in (\ref{chiclass}),  we have seen 
that the $\chi_i$ 
are the quantum analogues of the tangent 
vectors at the origin of the group :
\eq 
\chi_i \stackrel{q\rightarrow 1}{\longrightarrow} \pdxi\mid_{x=0} 
\en
\noi and that the left-invariant vectorfields $t_i$ 
constructed from the
$\chi_i$ are :
\eq 
t_i = \chi_i* = (id \otimes \chi_i) \D 
\en

In the commutative case, the Lie derivative along a generic vectorfield 
$V$ is given by:
\eq
\ell_{V}\tau =\lim_{\epsi\rightarrow 1}{1\over \epsi}
[{\varphi_{\epsi}^V}^*(\tau)-\tau]\label{liecomm}
{}~~~~~~~~~~~~~~~~~~~~~\forall\tau\in\Ga^{\otimes}
\en
where $\varphi_{\epsi}^V$ is the flow of the vectorfield $V$ and
${\varphi_{\epsi}^V}^*$ the pullback. If $t$ is a left invariant
vectorfield then 
\eq
\varphi_{\epsi}^t=R_{e^{\epsi t}}~~~~~~~~~\mbox{i.e.}~~~~~~~~~~~~~
\varphi_{\epsi}^t(g)=g{e^{\epsi t}}~~~\forall~g\in G~.
\en
We have $\ell_{t}\vartheta=\lim_{\epsi\rightarrow 0}{1\over \epsi}
[{R_{e^{\epsi t}}}^*(\vartheta)-\vartheta]$
and therefore the Lie derivative $\ell_t$ is given by the right action 
${R_{e^{\epsi t}}}^*$ of 
the group on covariant tensorfields.
At the quantum level, recalling (\ref{R*rho}) and that 
$\DR \rightarrow \RR$ when \qone, it is natural to define:
\sk
\noi{\bf Definition } The quantum Lie derivative 
along the left-invariant vectorfield
$t = (id \otimes \chi)\D$ is the operator:
\eq
\ell_t \equiv (id \otimes \chi)\DR~  \label{defqlie}
\en
\noi that is
$$ 
\forall\tau\in\Ga^{\otimes}~~~~~~~~~~~~
\ell_t(\tau) \equiv (id \otimes \chi)\DR (\tau)~ \;\;
\;\;\;\;\; \ell_t\; : \;\; \Gamma^{\otimes n} \longrightarrow 
\Gamma^{\otimes n} 
\;.  $$
\noi For example we have :
\eq
\ell_{t}(a) = t(a),~~~a \in A, \label{la}
\en
\eq
\ell_{t_i}(\om^j)=(id\otimes\chi_i)
\DR(\om^j)=\om^k\chi_i(M_k{}^j)=\C{ki}{j}\om^k, \label{lom}
\en
\noi the classical limit being evident.
\sk 
It is useful to define the $*$ product of a functional with any $\tau \in 
\Ga^{\otimes n}$ as
\eq
\chi * \tau \equiv (id \otimes \chi)\DR(\tau),
\en
where the $\DR$ acts on a generic element 
$\tau = \rho^1 \otimes \rho^2 \otimes\cdots \rho^n \;\; \in \; 
\Gamma^{\otimes n}$ as in (\ref{DRGaGaGa}).
In these notations we then have 
\eq
\ell_t =\chi *~~
\en

The quantum Lie derivative has properties analogous to that of the 
ordinary Lie derivative:

i) it is linear in $\tau$: 
\eq
\ell_t(\lambda\tau+\tau')=\lambda\ell_t(\tau)+\ell_t(\tau');
\en

ii) it is linear in $t$:
\eq
\ell_{\lambda t+t'} = \lambda\ell_t + \ell_{t'},~~\la \in \Cb.
\en
By virtue of this last property we can just study $\ell_{t_i}$, where 
$\{t_i\}$ is a basis of $\invXi$.\sk
\noi{\bf Theorem} 2.4.8 $~$
The following relation holds:
\eq
\ell_{t_i}(\tau\otimes \tau')=\ell_{t_j}(\tau)\otimes \;f^j{}_i\ast\tau'+
     \tau\otimes \ell_{t_i}(\tau')\label{lieprod}
\en

\noi{\sl Proof}
\begin{eqnarray*}
\ell_{t_i}(\tau\otimes \tau')
&=&(id\otimes\chi_i)
  \DR(\tau\otimes \tau')\\
&=&(id\otimes\chi_i)(\tau_1\otimes \tau'_1\;\otimes\tau_2\tau'_2)\\
&=&(\tau_1\otimes \tau'_1)
\chi_i(\tau_2\tau'_2)
=(\tau_1\otimes \tau'_1)
 [\chi_j(\tau_2) f^j_i(\tau'_2) +
 \varepsilon(\tau_2)\chi_i(\tau'_2)]\\
&=&\tau_1\chi_j(\tau_2)
 \otimes  \tau'_1 f^j{}_i(\tau'_2) + \tau_1\varepsilon(\tau_2)\otimes 
\tau'_1 \chi_i (\tau'_2)\\
&=&\ell_{t_j}(\tau)\otimes (id\otimes f^j{}_i)*
                           \tau' +\tau\otimes \ell{_{t_i}}(\tau')
\end{eqnarray*}

\noi [remember that $\chi_j(a)$ and $f{^j}_i(a)$ 
are $\Cb$ numbers].
The same argument leads to:
\eq \ell_{t_i}(a\om^j) =  \ell_{t_k}(a)
(f^k{}_i*\om^j) + a\ell_{t_i}(\om^j) \en
\eq \ell_{t_i}(\om^j a) =  \ell_{t_k}(\om^j)
(f^k{}_i*a) + \om^j\ell_{t_i}(a). \en
\noi The classical limit of (\ref{lieprod}) is easy to recover if we 
remember that $\epsi * \tau = \tau$. Formulas (\ref{lieprod}), 
(\ref{la}) and (\ref{lom}) uniquely define the quantum $\ell_t$, which 
reduces, for \qone, to the classical Lie derivative.
\sk
\noi{\bf Theorem}  2.4.9 $~$ 
The Lie derivative commutes 
with the exterior derivative:
\eq
\ell_{t_{i}} (d \vartheta ) = d ( \ell_{t_{i}} \vartheta ), 
\;\;\;\;\;\;\;\;\;\;\;\;\;\;\vartheta\in\Ga^{\wedge}\subset\Ga^{\otimes} 
: \;\;\mbox{generic form}.
\en

\noi{\sl Proof:}
$$\ell_{t_i}(d\vartheta)=(id\otimes\chi_i)\DR(d\vartheta)
        =(id\otimes\chi_i)(d\otimes id)\DR(\vartheta)=$$
$$(d\otimes\chi_i)\DR(\vartheta)=
d\vartheta_1\underbrace{\chi_i(\vartheta_2)}
      _{\in\Cb}=d[\vartheta_1
                       \chi_i(\vartheta_2)]=d(\ell_{t_i}\vartheta),$$

\noi where in the second equality we have used property (\ref{propd4}).
\sk

\noi{\bf Theorem} 2.4.10 $~$
The Lie derivative commutes 
with the left and right coactions 
$\DL$ and $\DR$, $\forall\tau\in\Ga^{\otimes}$:
\eq
(id \otimes \ell_t)\DL(\tau)=\DL(\ell_t \tau) \label{liec}
\en
\eq
(id \otimes \ell_t)\DR(\tau)=\DR(\ell_t \tau)~.
\en
\noi The proof is easy and relies on the fact that left and right 
coactions commute, cf. eq. (\ref{bicovariance}). In 
the classical limit, eq. (\ref{liec}) becomes :
\eq
\ell_t (\ll{g} \theta) = \ll{g} ( \ell_t \theta).
\en

\noi{\bf Note} 2.4.1 $~$ It is not difficult to prove the associativity of the 
generalized $*$ product, for example that 
$(\chi *\chi')*\tau = \chi *(\chi'*\tau)$. From 
this property it follows that the $q$-Lie derivative is a 
representation of the $q$-Lie algebra:
$$[\ell_{t} , \ell_{t'}](\tau) = \ell_{[t ,t']}(\tau),$$
where the left hand side is defined via the adjoint action :
\[[\ell_{t} , \ell_{t'}](\tau) \equiv
\ell_{\kappa'{(\chi'_1)*}}{\scriptstyle{{}^{{}_{\circ}}}}
\ell_{t}{\scriptstyle{{}^{{}_{\circ}}}}\ell_{t'_2}~.\]
In the $\{t_i\}$ basis: $
[\ell_{t_i},\ell_{t_k}]
=\ell_{t_i}{\scriptstyle{{}^{{}_{\circ}}}}\ell_{t_k}-\La_{~ik}^{ef}
\ell_{t_e}{\scriptstyle{{}^{{}_{\circ}}}}\ell_{t_f}~.$
\sk
We can now prove the Cartan identity:
\sk
\noi{\bf Theorem} 2.4.11 $~$
The contraction operator $i_t$ defined in
(\ref{contraction}) , the Lie 
derivative and the exterior differential satisfy
(we omit the composition product ${\scriptstyle{{}^{{}_{\circ}}}}$) :
\eq
\ell_{t_i}= i_{t_i}  d +d  i_{t_i}.
\en
\noi A proof of this theorem is given in Appendix B. 
\cvd

Led by Theorem 2.4.11, 
it is natural to introduce the Lie derivative along a generic
vectorfield $V$ through the following 
\sk
\noi {\bf Definition}
\eq
\ell_{V}= i_{V}  d +d  i_{V}   \label{ellV} ~\label{iright}
\en
\noi {\bf Theorem} 2.4.12 $~$ The Lie derivative satisfies the following 
properties:
\sk
\noi 1) \spz $\ell_Va=V(a)$
\sk
\noi 2) \spz $\ell_Vd\vart=d\ell_V\vart$
\sk
\noi 3) \spz
$\ell_V(\lambda\vart+\vart')=\lambda\ell_V(\vart)+\ell_V(\vart')$
\sk
\noi 4) \spz $\ell_{V\dia
b}(\vart)=(\ell_V\vart)b-(-1)^pi_V(\vart)\wedge db$ 
\sk
\noi 5) \spz \( \ell_V(\vartheta \wedge\vartheta')=
\vartheta\wedge\ell_V(\vartheta')+
\ell_{t_k}(\vartheta)\wedge(\f{k}{j}*\vartheta')
b^j+(-1)^{deg(\vartheta')}i_{t_k}(\vartheta)\wedge(\f{k}{j}*\vartheta')\wedge
db^j\)
\sk
\noi where $\vart$ and $\vartheta'$ are generic forms and
$V=t_j\dia b^j$.
\sk
\noi {\sl Proof}
\sk
\noi Properties 1), 2), 3) and 4) follow directly from the definition
(\ref{ellV}).\\
Property 5) is also a consequence of definition (\ref{ellV}); the
proof  uses relation d$'$) and the identity
$d(\f{k}{j}*\vartheta)=\f{k}{j}*d\vartheta$ [see (\ref{dhthe})]. 
\cvd

\noi{\bf Note} 2.4.2 $~$ 
It is natural to define a Lie derivative 
$\err_{h_i}$ of a generic covariant 
tensorfield $\tau\in \Ga^{\otimes}$ along a
{\em right}-invariant vectorfield $h_i$
in terms of the {\em left} coaction $\DL ~:$
$$
\err_h(\tau) \equiv (\chi \otimes i\!d)\DL(\tau)
 \equiv \tau * \chi~,
$$
just like it was natural that we used the {\em right} coaction,
when we defined $\ell_{t_i}$ in (\ref{defqlie}).
In this note we compare the two definitions.

{}From the above definition we find
\eq
\begin{array}{rcl}
\err_{h_i}(\vart \wedge\vart') & = & \chi_i(\vart_1 {\vart'}_1) \vart_2 
\wedge{\vart'}_2\\
& = & \chi_j(\vart_1) f^j{}_i(\vart'_1)\vart_2\wedge \vart'_2 + 
\epsi(\vart_1) \chi_i({\vart'}_1) \vart_2\wedge {\vart'}_2\\
& = & \err_{h_j}(\vart) \wedge(\vart'* f^j{}_i)+ \vart\wedge \err_{h_i}(\vart')
\end{array}
\en
where $\vart'*\f{j}{i}\equiv (\f{j}{i}\otimes id)\DL(\vart')$.
In particular:
\eq
\begin{array}{rcl}
\err_{h_i}(adb) & = & h_i(a) d(b*f^j{}_i) + a  d(h_j(b)) ~.
\end{array} \label{errab}
\en
where we have used (\ref{propd3}).
On the other hand, since   $h_i=t_j\dia M_i{}^j$, 
we can give an alternative expression for the Lie derivative along the
right invariant vectorfield $h_i$:
\eq
\begin{array}{lcl}
\ell_{h_i}(adb)=\ell_{t_j\dia M_i{}^j}(adb)& = &
a\ell_{t_j\dia M_i{}^j} (db)+\ell_{t_k}(a)\wedge(f^k{}_j*db)M_i{}^j\\
          & = &a\,d(h_i(b))+t_k(a)d(f^k{}_j*b)\,M_i{}^j~. \label{ellab}
\end{array}
\en
The difference between expressions (\ref{errab}) and (\ref{ellab}) is
a good index for the ``defect'' between left and right transports on
a quantum group:
\eq
\begin{array}{rcl}
(\ell_{h_i} - \err_{h_i})(adb)
& = &t_k(a) [(f^k{}_j * db)M_i{}^j  -  M_j{}^k(db * f^j{}_i)]\\
& = & -t_k(a)D\!I^k{}_i(b);
\end{array}
\en
where
\eq
D\!I^k{}_i(b) \equiv
[(f^k{}_j* b)d(M_i{}^j) - d(M_j{}^k)(b * f^j{}_i) ]
\spz {\mbox{ (Defect Index)}}.
\en
In the last passage
we have used  the Leibniz rule for $d$ combined with the bicovariance condition
(\ref{propM}).
The term in the square brackets is
always zero in the classical (undeformed) case. Note that
$(\ell_{h_i} - \err_{h_i})$
vanishes on $a$ and $db$ separately but not necessarily on $adb$.
The case of ``$a$'' confirms (\ref{leftrightrans}):

\eq
\ell_{h_i}(a) = t_j(a)M_i{}^j = h_i(a) = \err_{h_i}(a)
\en
and shows that we will not encounter any ambiguities or inconsistencies
as long as we deal with general vectorfields and functions alone.
Problems can occur however when we start to introduce forms.
For example in the $GL_q(2)$ differential calculus of Section 2.2 we have,
sum over $\sma{$a$}$ understood,
$\kappa({x_a{}^a}_1)(\ell_{h^i{}_j}-\err_{h^i{}_j})({x_a{}^a}_2db)=
-(\ff{k}{ke}{f}*b)d\MM{i}{je}{f}$
$=\lambda t^e{}_f(b)d\MM{i}{je}{f}$ for any $b$ such that $\epsi(b)=0$. 
This expression is clearly 
$\not= 0$ in general.

\subsection{Algebra of Differential Operators}
In the previous sections, given a Woronowicz differential calculus
on a generic Hopf algebra $A$, we have defined the quantum analogue of Lie 
derivative and of inner derivative by a natural generalization of their
defining classical formulae.
The Lie derivative and  contraction operators act on 
the space $\Ga^{\otimes}$ of covariant tensorfields, we have in 
particular studied their properties 
on the space $\Ga^{\wedge}\subset\Ga^{\otimes}$ of forms where
the exterior differential is also present.

These operators form a  graded quantum Lie algebra
\begin{eqnarray}
\{ d , d \} & = & 0\\
\left[ d , \ell_V \right] & = & 0\\
\{ d , i_V \} & = & \ell_V
\end{eqnarray}
which is supplemented by two more
relations
\begin{eqnarray}
& &\left[ \ell_{t_i} ,
\ell_{t_k} \right] =
\ell_{[t_i,t_k]} =\C{ik}{l}  \ell_{t_l}\label{ellinner0} \\
& &\left[i_{t_i},\ell_{t_k} \right] = i_{[t_i,t_k]}  
= \C{ik}{l} i_{t_l}\label{ellinner}
\end{eqnarray}
where the definition of the brackets in the left hand side of 
(\ref{ellinner0}) and (\ref{ellinner}) is  the generalization of the adjoint 
action:
$$
[\ell_{t_i},\ell_{t_k}]\equiv
\ell_{\kappa'{({\chi_k}_{1}) *}}{\scriptstyle{{}^{{}_{\circ}}}}
\ell_{t_i}{\scriptstyle{{}^{{}_{\circ}}}}\ell_{{\chi_k}_2*}=
\ell_{t_i}{\scriptstyle{{}^{{}_{\circ}}}}\ell_{t_k}+
\ell_{\kappa'{({\chi_e) *}}}
{\scriptstyle{{}^{{}_{\circ}}}}\ell_{t_i}{\scriptstyle{{}^{{}_{\circ}}}}
f^e{}_k*
=\ell_{t_i}{\scriptstyle{{}^{{}_{\circ}}}}\ell_{t_k}-\La_{~ik}^{ef}
\ell_{t_e}{\scriptstyle{{}^{{}_{\circ}}}}\ell_{t_f}
$$
(this last equality is explained in Note 2.4.1)
\eq
\left[i_{t_i}, \ell_{t_k}, \right]\equiv
\ell_{\kappa'{({\chi_k}_{1})*}}{\scriptstyle{{}^{{}_{\circ}}}}
i_{t_i}{\scriptstyle{{}^{{}_{\circ}}}}\ell_{{\chi_k}_2 *}=
i_{t_i}{\scriptstyle{{}^{{}_{\circ}}}} \ell_{t_k}-\La_{~ik}^{ef}
\ell_{t_e}{\scriptstyle{{}^{{}_{\circ}}}} i_{t_f}\label{adieell}
\en
The proof of (\ref{ellinner}) and of the last equality in 
(\ref{adieell}), similarly to the proof of the Cartan identity,
is by induction. It is given in Appendix B.
\sk
The cross-commutation relations between forms, exterior derivative, Lie
derivative and inner derivative, that we have derived from the actions
of $i_V$ and $\ell_t$ on generic tensors $\tau\in\Ga^{\otimes}$ and essentially
(see the definition of $i_V$) from
the $\Ga^{\otimes}\leftrightarrow\Xi^{\otimes}$ duality --i.e. the 
bicovariant bimodule structure of
$\Ga^{\otimes}$ and $\Xi^{\otimes}$-- can be formally   derived also
from the cross product algebra
$\Ga^{\wedge}{\mbox{$\times \!\rule{0.3pt}{1.1ex}\,$}}{\Ga^{\wedge}}^*$
\cite{S2, Vladimirov}. 
Here $\Ga^{\wedge}$ is seen as a graded Hopf algebra: the product
in $\Ga^{\wedge}\otimes \Ga^{\wedge}$ is given by $(I\otimes\mu)(\nu\otimes I) 
=(-1)^{deg(\mu)deg(\nu)}(\nu\otimes\mu)$, the costructures generalize those 
of $A$ and are:
$\D(\om^i)=({}_{\scriptsize\Ga}\D +\D_{\scriptsize\Ga})(\om^i)=
\om^i\otimes M_j{}^i+I\otimes\om^i$, $\epsi(\om^i)=0$,
$\kappa(\om^i)=-\om^j\kappa(M_j{}^i)$
\cite{Brzezinski}.  
${\Ga^{\wedge}}^*$
is the graded Hopf algebra dual to ${\Ga^{\wedge}}$, 
${\Ga^{\wedge}}^*=U\oplus\Ga^*\oplus{\Ga^{\wedge 2}}^*\oplus\ldots~.$ 
For example $\chi_i\vart=\vart_1\le{\chi_i}_1,\vart_2
\re{\chi_i}_2=(\chi_j*\vart)f^j{}_i+\vart{\chi_i}$ corresponds to
$\ell(\vart\wedge\vart')=\ell_{\chi_i}(\vart)\wedge(f^j{}_i*\vart')
+\vart\wedge\ell_{\chi_i}(\vart')$.
As shown in \cite{PaoloPeter} the graded $q$-Lie algebra of the operators 
$i_V,d,\ell$ can also be interpreted as a braided tensor algebra.

\chapter{Geometry of the quantum Inhomogeneous Linear Groups $IGL_{q,r}(N)$}

\def\spinst#1#2{{#1\brack#2}}
\def\noi{\noindent}
\def\om{\omega}
\def\al{\alpha}
\def\be{\beta}
\def\ga{\gamma}
\def\Ga{\Gamma}
\def\del{\delta}
\def\linv{{1 \over \lambda}}
\def\rinv{{1\over {r-r^{-1}}}}
\def\alb{\bar{\alpha}}
\def\beb{\bar{\beta}}
\def\gab{\bar{\gamma}}
\def\deb{\bar{\delta}}
\def\ab{\bar{a}}
\def\Ab{\bar{A}}
\def\Bb{\bar{B}}
\def\Cb{\bar{C}}
\def\Db{\bar{D}}
\def\ab{\bar{a}}
\def\cb{\bar{c}}
\def\db{\bar{d}}
\def\bb{\bar{b}}
\def\eb{\bar{e}}
\def\fb{\bar{f}}
\def\gb{\bar{g}}
\def\xih{\hat\xi}
\def\Xih{\hat\Xi}
\def\uh{\hat u}
\def\vh{\hat v}
\def\ub{\bar u}
\def\vb{\bar v}
\def\xib{\bar \xi}

\def\alp{{\alpha}^{\prime}}
\def\bep{{\beta}^{\prime}}
\def\gap{{\gamma}^{\prime}}
\def\dep{{\delta}^{\prime}}
\def\rhop{{\rho}^{\prime}}
\def\taup{{\tau}^{\prime}}
\def\rhopp{\rho ''}
\def\thetap{{\theta}^{\prime}}
\def\imezzi{{i\over 2}}
\def\unquarto{{1 \over 4}}
\def\onehalf{{1 \over 2}}
\def\unmezzo{{1 \over 2}}
\def\epsi{\varepsilon}
\def\we{\wedge}
\def\th{\theta}
\def\de{\delta}
\def\cony{i_{\de {\vec y}}}
\def\Liey{l_{\de {\vec y}}}
\def\tv{{\vec t}}
\def\Gt{{\tilde G}}
\def\deyv{\vec {\de y}}
\def\part{\partial}
\def\pdxp{{\partial \over {\partial x^+}}}
\def\pdxm{{\partial \over {\partial x^-}}}
\def\pdxi{{\partial \over {\partial x^i}}}
\def\pdy#1{{\partial \over {\partial y^{#1}}}}
\def\pdx#1{{\partial \over {\partial x^{#1}}}}
\def\pdyx#1{{\partial \over {\partial (yx)^{#1}}}}

\def\qP{q-Poincar\'e~}
\def\A#1#2{ A^{#1}_{~~~#2} }

\def\R#1#2{ R^{#1}_{~~~#2} }
\def\Rp#1#2{ (R^+)^{#1}_{~~~#2} }
\def\Rpinv#1#2{ [(R^+)^{-1}]^{#1}_{~~~#2} }
\def\Rm#1#2{ (R^-)^{#1}_{~~~#2} }
\def\Rinv#1#2{ (R^{-1})^{#1}_{~~~#2} }
\def\Rsecondinv#1#2{ (R^{\sim 1})^{#1}_{~~~#2} }
\def\Rinvsecondinv#1#2{ ((R^{-1})^{\sim 1})^{#1}_{~~~#2} }

\def\Rpm#1#2{(R^{\pm})^{#1}_{~~~#2} }
\def\Rpminv#1#2{((R^{\pm})^{-1})^{#1}_{~~~#2} }

\def\Rbo{{\cal R}}
\def\Rb#1#2{{ \Rbo^{#1}_{~~~#2} }}
\def\Rbp#1#2{{ (\Rbo^+)^{#1}_{~~~#2} }}
\def\Rbm#1#2{ (\Rbo^-)^{#1}_{~~~#2} }
\def\Rbinv#1#2{ (\Rbo^{-1})^{#1}_{~~~#2} }
\def\Rbpm#1#2{(\Rbo^{\pm})^{#1}_{~~~#2} }
\def\Rbpminv#1#2{((\Rbo^{\pm})^{-1})^{#1}_{~~~#2} }

\def\RRpm{R^{\pm}}
\def\RRp{R^{+}}
\def\RRm{R^{-}}

\def\Rh{{\hat R}}
\def\Rbh{{\hat {\Rbo}}}
\def\Rhat#1#2{ \Rh^{#1}_{~~~#2} }
\def\Rbar#1#2{ {\bar R}^{#1}_{~~~#2} }
\def\L#1#2{ \La^{#1}_{~~~#2} }
\def\Linv#1#2{ \La^{-1~#1}_{~~~~~#2} }
\def\Rbhat#1#2{ \Rbh^{#1}_{~~~#2} }
\def\Rhatinv#1#2{ (\Rh^{-1})^{#1}_{~~~#2} }
\def\Rbhatinv#1#2{ (\Rbh^{-1})^{#1}_{~~~#2} }
\def\Z#1#2{ Z^{#1}_{~~~#2} }
\def\Rt#1{ {\hat R}_{#1} }
\def\La{\Lambda}
\def\Rha{{\hat R}}
\def\ff#1#2#3{f_{#1~~~#3}^{~#2}}
\def\MM#1#2#3{M^{#1~~~#3}_{~#2}}
\def\cchi#1#2{\chi^{#1}_{~#2}}
\def\ome#1#2{\om_{#1}^{~#2}}
\def\RRhat#1#2#3#4#5#6#7#8{\La^{~#2~#4}_{#1~#3}|^{#5~#7}_{~#6~#8}}
\def\RRhatinv#1#2#3#4#5#6#7#8{(\La^{-1})^
{~#2~#4}_{#1~#3}|^{#5~#7}_{~#6~#8}}
\def\LL#1#2#3#4#5#6#7#8{\La^{~#2~#4}_{#1~#3}|^{#5~#7}_{~#6~#8}}
\def\LLinv#1#2#3#4#5#6#7#8{(\La^{-1})^
{~#2~#4}_{#1~#3}|^{#5~#7}_{~#6~#8}}
\def\U#1#2#3#4#5#6#7#8{U^{~#2~#4}_{#1~#3}|^{#5~#7}_{~#6~#8}}
\def\Cb{{\bf C}}
\def\CC#1#2#3#4#5#6{\Cb_{~#2~#4}^{#1~#3}|_{#5}^{~#6}}
\def\cc#1#2#3#4#5#6{C_{~#2~#4}^{#1~#3}|_{#5}^{~#6}}

\def\C#1#2{ {\bf C}_{#1}^{~~~#2} }
\def\c#1#2{ C_{#1}^{~~~#2} }
\def\q#1{   {{q^{#1} - q^{-#1}} \over {q^{\unmezzo}-q^{-\unmezzo}}}   
}
\def\Dmat#1#2{D^{#1}_{~#2}}
\def\Dmatinv#1#2{(D^{-1})^{#1}_{~#2}}
\def\DR{\Delta_R}
\def\DL{\Delta_L}
\def\f#1#2{ f^{#1}_{~~#2} }
\def\F#1#2{ F^{#1}_{~~#2} }
\def\T#1#2{ T^{#1}_{~~#2} }
\def\Ti#1#2{ (T^{-1})^{#1}_{~~#2} }
\def\Tp#1#2{ (T^{\prime})^{#1}_{~~#2} }
\def\Th#1#2{ {\hat T}^{#1}_{~~#2} }
\def\TP{ T^{\prime} }
\def\M#1#2{ M_{#1}^{~#2} }
\def\qm{q^{-1}}
\def\rminus{r^{-1}}
\def\um{u^{-1}}
\def\vm{v^{-1}}
\def\xm{x^{-}}
\def\xp{x^{+}}
\def\fm{f_-}
\def\fp{f_+}
\def\fn{f_0}
\def\D{\Delta}
\def\DN{\Delta_{N+1}}
\def\kN{\kappa_{N+1}}
\def\eN{\epsi_{N+1}}
\def\Mat#1#2#3#4#5#6#7#8#9{\left( \matrix{
     #1 & #2 & #3 \cr
     #4 & #5 & #6 \cr
     #7 & #8 & #9 \cr
   }\right) }
\def\Ap{A^{\prime}}
\def\Dp{\Delta^{\prime}}
\def\Ip{I^{\prime}}
\def\ep{\epsi^{\prime}}
\def\kp{\kappa^{\prime}}
\def\kpm{\kappa^{\prime -1}}
\def\kpsq{\kappa^{\prime 2}}
\def\km{\kappa^{-1}}
\def\gp{g^{\prime}}
\def\qone{q \rightarrow 1}
\def\rone{r \rightarrow 1}
\def\qrone{q,r \rightarrow 1}
\def\Fmn{F_{\mu\nu}}
\def\Am{A_{\mu}}
\def\An{A_{\nu}}
\def\dm{\part_{\mu}}
\def\dn{\part_{\nu}}
\def\Ana{A_{\nu]}}
\def\Bna{B_{\nu]}}
\def\Zna{Z_{\nu]}}
\def\dma{\part_{[\mu}}
\def\qsu{$[SU(2) \times U(1)]_q~$}
\def\suq{$SU_q(2)~$}
\def\su{$SU(2) \times U(1)~$}
\def\gij{g_{ij}}
\def\qL{SL_q(2,{\bf C})}
\def\GLqrN{GL_{q,r}(N)}
\def\IGLqrN{IGL_{q,r}(N)}
\def\IGLqrtwo{IGL_{q,r}(2)}
\def\GLqrNo{GL_{q,r}(N+1)}
\def\SLqrN{SL_{q,r}(N)}
\def\UglqrN{U_{q,r}(gl(N))}
\def\UglqrNo{U_{q,r}(gl(N+1))}
\def\UiglqrN{U_{q,r}(igl(N))}

\def\RR{R^*}
\def\rr#1{R^*_{#1}}

\def\Lpm#1#2{L^{\pm #1}_{~~~#2}}
\def\Lmp#1#2{L^{\mp#1}_{~~~#2}}
\def\LLpm{L^{\pm}}
\def\LLmp{L^{\mp}}
\def\LLp{L^{+}}
\def\LLm{L^{-}}
\def\Lp#1#2{L^{+ #1}_{~~~#2}}
\def\Lm#1#2{L^{- #1}_{~~~#2}}

\def\gu{g_{U(1)}}
\def\gsu{g_{SU(2)}}
\def\tg{ {\rm tg} }
\def\Fun{$Fun(G)~$}
\def\invG{{}_{{\rm inv}}\Ga}
\def\Ginv{\Ga_{{\rm inv}}}
\def\qonelim{\stackrel{q \rightarrow 1}{\longrightarrow}}
\def\ronelim{\stackrel{r \rightarrow 1}{\longrightarrow}}
\def\qronelim{\stackrel{q=r \rightarrow 1}{\longrightarrow}}
\def\viel#1#2{e^{#1}_{~~{#2}}}
\def\ra{\rightarrow}
\def\detq{{\det}}
\def\detqr{{\det}}
\def\detqrm{{\det} }
\def\detqrTAB{{\det} \T{A}{B}}
\def\detqrTab{{\det} \T{a}{b}}
\def\P{P}
\def\Qt{Q}
\def\chit{\stackrel{\leftarrow}{\partial}}

\def\pp#1#2{\Pi_{#1}^{(#2)}}

\def\square{{\,\lower0.9pt\vbox{\hrule \hbox{\vrule height 0.2 cm
\hskip 0.2 cm \vrule height 0.2 cm}\hrule}\,}}

\def\sma#1{\mbox{\footnotesize #1}}

\def\H{H^\bot}
\def\Mat2#1#2#3#4{\left( \matrix{
     #1 & #2 &  \cr
     #3 & #4 &  \cr
       }\right) }
\def\Uigl{U_{q,r}(igl(N))}
\def\Lcm{{{\cal L}^-}}
\def\le{\langle}
\def\re{\rangle}


In this chapter, following \cite{IGLAschieri},
we analize the geometry of the inhomogeneous quantum 
linear groups $\IGLqrN$.
Quantum deformations of inhomogeneous Lie groups
have been studied in  \cite{inhom}  \cite{Rinhom1} \cite{Rinhom2}. 
An $R$-matrix 
approach has been independently proposed for $IGL_q(N)$
in ref.s \cite{Rinhom1} and \cite{Rinhom2}.

We construct the {\sl multiparametric}
$\IGLqrN$ $q$-groups, 
their universal enveloping algebra and their bicovariant differential 
calculus using a projection $P\;:~ GL_{q,r}(N$ $+1)\rightarrow \IGLqrN$;
this projection procedure was first introduced in \cite{Rinhom3}. 
 
All the quantities relevant to the $\IGLqrN$
(bicovariant)  differential calculus are given explicitly:
exterior derivatives,  left-invariant $1$-forms, Cartan-Maurer
equations, tangent vectors and their $q$-Lie algebra and so on.
The method is illustrated in the case of $\IGLqrtwo$:  the
general formulas are applied and tested on this example.

In this framework we  construct the differential
geometry of the (multiparametric) quantum plane in a novel
and easy way. 
\sk
Deformations of inhomogeneous Lie groups and Lie algebras usually 
include a dilatation generator, moreover the determinant of the 
fundamental representation of the $q$-group is in general not central. 
It is studying the most general (multiparametric) deformation that we 
understand the interplay between the absence or presence of the
dilatation and the properties of the determinant.
This also clarify the relation between the non-commutativity of
the quantum plane coordinates $x^a$ discussed in Section 3.7 (due to the 
auxiliary deformation parameters $q_i$), the non-comutativity of the
generators $\T{a}{b}$ of the homogeneous linear subgroup and  the finite 
difference structure of the differential calculus, that is due to the main 
deformation parameter $r$ (called $q$ in the previous chapters),
cf.  \cite{Zumino00}. 
\sk
In Section 3.1 we recall the basics of  the 
linear quantum groups and in Section 3.2 we discuss their duals
in some detail.  In fact, Sections 3.1 and 3.2 are a short review of
the multiparametric deformations of
$GL_{q,r}(N)$, where $q$ indicates a set of parameters $q_i$,
and of  their universal enveloping algebras. 
The usual uniparametric case is recovered for $r=q_i=q$. 
For references on multiparametric deformations, see 
\cite{Schirrmacher,multiparam1, multiparam2}. 
\sk
In Section 3.3, we first present 
the quantum group  $IGL_{q,r}(N)$ as a Hopf algebra
with given generators, commutation relations and co-structures.
We then reobtain it as the image of a
projection $P$ from $GL_{q,r}(N+1)$, and show how
the ``mother" Hopf algebra $GL_{q,r}(N+1)$ 
determines the Hopf algebra structure on $IGL_{q,r}(N)$. 
In the language of Hopf algebra
ideals  $IGL_{q,r}(N)$ 
is seen as the quotient
of $GL_{q,r}(N+1)$
with respect to a suitable Hopf ideal.

The fundamental representation
of $IGL_{q,r}(N)$ contains the $GL_{q,r}(N)$ elements $\T{a}{b}$
and the ``coordinates" $x^a$ as in the classical case,
in addition, there is also an element $u$ playing
the role of a dilatation.
By fixing some of the parameters $q$, we find that 
this element $u$ can be made central, and hence consistently set  
equal to
the identity $I$. 

A quantum determinant can be defined, and is central
only in a subclass of the multiparametric
deformations. In this subclass, however, the element $u$ is not  
central.
We end the section analizing the semidirect product structure
of  $\IGLqrN$ given by $GL_{q,r}(N)$ and the quantum plane:
this construction is based on the observation that $GL_{q,r}(N)$ 
is both a Hopf subalgebra in $IGL_{q,r}(N)$ and a quotient of $IGL_{q,r}(N)$
obtained projecting to zero the quantum plane coordinates $x^a$.
\sk
The explicit
construction of the bicovariant differential calculus for
$\GLqrN$, in terms of the dual algebra,  is given in Section 3.4.
In Section 3.5 we project the
bicovariant differential calculus of $GL_{q,r}(N+1)$ 
to $IGL_{q,r}(N)$ and
study the bicovariant bimodules of $1$-forms and 
tangent vectors on 
$IGL_{q,r}(N)$. In particular, the $q$-Lie algebra is given
explicitly. We also study in detail
the exterior algebra and the exterior derivative, and  find the 
Cartan-Maurer equations. 
In Section 3.6 we then study  the universal enveloping algebra 
$U_{q,r}(igl(N))$ its semidirect product structure [given by
$U_{q,r}(gl(N))$ and the translation generators] 
and the duality with $\IGLqrN$. 
The Universal enveloping algebra $U_{q,r}(igl(N))$ is the natural setting 
where 
to study  $q$-Lie algebras and therefore differential calculi. 
Using the general theory of Section 2.3, we easily obtain another differential 
calculus
on $IGL_{q,r}(N)$ that differs from the previous one by the presence of a 
dilatation 
generator corresponding to the dilatation $u\in IGL_{q,r}(N)$.
\sk
In Section 3.7 we discuss the multiparametric quantum plane,
i.e. the quantum coset space  $\IGLqrN / \GLqrN$
spanned by the  
coordinates $x^a$, and find a generalization
of the differential geometry of the $q$-plane 
of \cite{WZ}, \cite{zumi}, see also Schirrmacher in \cite{multiparam2}.
\sk
In the Table at the
end of the chapter we specialize
our general treatment to  $\IGLqrtwo$
and collect 
all the relevant  formulas for  its
bicovariant differential calculus.
\sk

\sect{$GL_{q,r}(N)$ and its real forms}

We here introduce the  multiparametric $q$-group 
$GL_{q,r}(N)$, where now the index $q \equiv q_{ab}$
represents a set of parameters, and $r$ is the parameter
we called $q$ in the previous chapters.
$GL_{q,r}(N)$  is the algebra
(over the complex field)
freely generated by the 
non-commuting matrix elements $\T{A}{B}$, 
({\small \sl  A,B=1,..N}),
 the identity $I$ and  the inverse $\Xi$ of
the $q$-determinant of $T$ defined in (\ref{qrdet}), modulo 
the ``$RTT$" relations:
\eq
\R{AB}{EF} \T{E}{C} \T{F}{D} = \T{B}{F} \T{A}{E} \R{EF}{CD}
\label{RTTGL}
\en
\noi where the $R$-matrix is given by \cite{Schirrmacher, multiparam1}:
\eq
\R{AB}{CD}=\de^A_C \de^B_D [{r\over q_{AB}}+(r-1) \de^{AB}]
+(r-r^{-1})~ \de^A_D  
\de^B_C
 \theta^{AB} \label{Rmp}
\en
\noi with $\theta^{AB}=1$ 
for {\footnotesize $A > B$ } and zero otherwise, and
\eq 
q_{AB}={r^2\over q_{BA}},~~q_{AA}=r  \label{qabqba}
\en
 
It is useful to list the nonzero
complex
components of the $R$ matrix (no sum on repeated indices):
\eqa
& &\R{AA}{AA}=r \cr
& &\R{AB}{AB}={r \over q_{AB}} , ~~~~~~~\mbox{\footnotesize 
 $A \not= B$ }\cr
& &\R{BA}{AB}=r-r^{-1},~~~\mbox{\footnotesize $B>A $}
\label{Rnonzero} 
\ena

\noi The $R$ matrix in (\ref{Rmp}) satisfies the 
quantum Yang-Baxter 
equation.
\sk
The standard uniparametric $R$ matrix \cite{FRT} 
is obtained from
(\ref{Rmp}) by setting all deformation 
parameters $q_{AB},r$ equal to a
single parameter $q$. For a further insight about the relationship between the
multiparametric and the uniparametric $R$-matrix see Section 4.1.
\sk
The quantum determinant of $T$ and 
its inverse $\Xi$ are defined by:
\eq 
\Xi~ {\det} T={\det} T~\Xi=I\label{Xi}
\en
\eq 
{\det} T \equiv \sum_{\sigma}  \left [
\prod_{A<B,\sigma (A) > \sigma(B)}
 (-{r^2\over q_{\sigma(B) \sigma(A)}})  \right ]~
\T{1}{\sigma(1)} \cdots
\T{N}{\sigma(N)} \label{qrdet} 
\en

\noi{\bf Note } 3.1.1 $~$ In the uniparametric case $r=q_{AB}=q$ we recover
 the usual formula
\eq 
\detq T \equiv \sum_{\sigma} (-q)^{l(\sigma)}
\T{1}{\sigma(1)} \cdots
\T{N}{\sigma(N)} \label{qdet}
\en
\noi where $l(\sigma)$ is the minimum number of
transpositions in the permutation
$\sigma$.
\sk
\noi {\bf Note } 3.1.2 $~$
In more mathematical terms, the algebra $GL_{q,r}(N)$ is the quotient
of the non-commuting algebra ${\bf C}\langle T^A{}_B, I, \Xi\rangle$ freely
generated by the elements $T^A{}_B, I, \Xi$ with respect to 
the two-sided ideal in 
${\bf C}\langle T^A{}_B, I, \Xi\rangle$ 
generated by the $RTT$ relations (\ref{RTTGL}).
\sk
\noi{\bf Note } 3.1.3 $~$ The  
inverse matrix $R^{-1}$, defined as
\eq
\Rinv{AB}{CD}\R{CD}{EF} \equiv \de^
A_E \de^B_F \equiv \R{AB}{CD}\Rinv{CD}{EF}
\en
\noi is given by 
\eq
R^{-1}_{q,r}=R_{\qm,\rminus} \label{Rinv}
\en 
\sk
\noi {\bf Note } 3.1.4 $~$ The $\Rh$ matrix defined by
$\Rhat{AB}{CD} \equiv \R{BA}{CD}$ 
satisfies the spectral decomposition
(Hecke condition):
\eq
(\Rh-rI)~ (\Rh + r^{-1} I)=0 
\en
\sk
\noi{\bf Note } 3.1.5 $~$ The determinant in (\ref{qrdet}) is 
central if and only if the following conditions on the
 parameters are satisfied
(see ref. \cite{Schirrmacher}):
\eq
q_{1,A} q_{2,A} \cdots q_{A-1,A} {r^2\over q_{A,A+1}} 
{r^2\over  
q_{A,A+2}}
\cdots {r^2\over q_{A,N}}=const. \label{centralitycondition}
\en
\noi for all {\small A=1,...N}. This results in 
{\small N-1} conditions among
the $q_{AB}$ and determines $const=r^{N-1}$. Using
 (\ref{qabqba}), and defining
\eq
Q_A \equiv 
\prod_{C=1}^N ({q_{CA}\over r})
\label{Qdefinition}
\en
the centrality conditions
(\ref{centralitycondition}) become:
\eq
Q_A=1 \label{centralityconditionQ} 
\en
We have used also $const=r^{N-1}$, so that
only $N-1$ of the above conditions are independent. Indeed 
the $Q_A$ satisfy the relation
\eq
Q_1 Q_2 \cdots Q_N=1
\en
In  general we have:
\eq
(\detqr T ) \T{A}{B} = {Q_A \over Q_B} \T{A}{B} (\detqr T),
~~\Xi \T{A}{B} = {Q_B \over Q_A} \T{A}{B} \Xi
\label{detcomm}
\en
When (\ref{centralityconditionQ})
holds\footnote{We disregard the solutions $\forall A\in a,~Q_A=\sqrt[N]{1}$
because we want a continuous deformation of the classical limit.}, we can 
consistently set $\detqr \T{A}{B} = I = \Xi$,
and obtain the multiparametric deformations $SL_{q,r}(N)$.
\sk
The algebra $GL_{q,r}(N)$ becomes a Hopf algebra with
the following coproduct
$\D$, counit $\epsi$ and coinverse $\kappa$:
\eqa
& & \D(\T{A}{B})=\T{A}{B} \otimes \T{B}{C}  \label{cos1} \\
& & \epsi (\T{A}{B})=\delta^A_B\\
& & \kappa(\T{A}{B})=\Ti{A}{B} \label{coinverse}\\
& &  \D (\detqr T)=\detqr T \otimes \detqr T, ~~
\D (\Xi)=\Xi \otimes  
\Xi,
~~\D(I)=I\otimes I\\
& &  \epsi (\detqr T)=1,~~\epsi (\Xi)=1,~~\epsi (I)=1\\
& &  \kappa (\detqr T)=\Xi,~~\kappa (\Xi)=
\detqr T,~~\kappa (I)=I
\label{cos2}
\ena
\noi The quantum inverse of $\T{A}{B}$ in 
(\ref{coinverse}) is given  
by:
\eq
\Ti{A}{B}=\Xi~\pp{AB}{1,N}~ t_B^{~A}  \label{Tinverse}
\en
\noi where $t_B^{~A}$ is the quantum minor, i.e. the quantum  
determinant
of the submatrix of $T$ obtained by removing the 
{\small $B$}-th row  
and
the {\small $A$}-th column, and $\pp{AB}{1,N}$ is a 
function of the
parameters $q$:
\eq
\pp{AB}{1,N} \equiv  { {\prod_{C=B+1}^N (-q_{BC})}  \over
 {\prod_{D=A+1}^N (-q_{AD})}} \label{defPi}
\en
\noi The superscript (1,N) reminds the range 
of the indices {\small A,B,C,..}.
 In the uniparametric case, the quantum  
inverse
has the simpler expression: 
\eq
\Ti{A}{B}=\Xi~ (-q)^{A-B} t_B^{~A}
\en
\noi{\bf Note } 3.1.6 $~$ 
As in Note 1.2.1, we recall that
in general $\kappa^2 \not= 1$ and
\eq
\kappa^2 (\T{A}{B})=\Dmat{A}{C} \T{C}{D}
 \Dmatinv{D}{B}=d^A d^{-1}_B
\T{A}{B},  
\en
\noi where $D$ is a diagonal matrix, 
$\Dmat{A}{B}=d^A \de^A_{B}$, given 
by $d^A=r^{2A-1}$ for $\GLqrN$. This matrix satisfies:
\eq
d^A d^{-1}_C \Rinv{BA}{DC} \R{EC}{BF}=
\de^A_F \de^E_D,~~~
d^A d^{-1}_C \R{AB}{CD} \Rinv{CE}{FB}=
\de^A_F \de^E_D  \label{Rsecinv1}
\en
\eq
d^B d^{-1}_D \Rinv{AB}{CD} \R{CE}{FB}=
\de^A_F \de^E_D,~~~
d^B d^{-1}_D \R{BA}{DC} \Rinv{EC}{BF}=
\de^A_F \de^E_D \label{Rsecinv2}
\en
\eq
\R{AC}{CB} d^{-1}_C=\de^A_B=\Rinv{AC}{CB} d_{C}
\label{RD}
\en
This last condition fixes the normalization of $D$.
\noi Relations (\ref{Rsecinv1}) and (\ref{Rsecinv2})
define a second inverse
$R^{\sim 1}$ of
the $R$ matrix  and a second inverse $(R^{-1})^{\sim 1}$
of the $R^{-1}$ matrix as:
\eq
\Rsecondinv{AB}{CD} \equiv
d^B d^{-1}_D \Rinv{AB}{CD}
\en
\eq
\Rinvsecondinv{AB}{CD} \equiv
d^A d^{-1}_C \R{AB}{CD}
\en
\noi Using (\ref{RD})
we can relate the $D$ matrix to this second inverse:
 \eq
\Dmatinv{A}{B}=\Rsecondinv{AC}{CB},~~~\Dmat{A}{B}=
\Rinvsecondinv{AC}{CB}
\en
This generalizes the analogous discussion for 
the uniparametric $D$ matrix
given in \cite{FRT}.

\sk
We turn now to the real forms of $\GLqrN$,
that are defined by 
$*$-conjugations of the $\GLqrN$ Hopf algebra;
see Section 1.3.
These  conjugations
must be compatible with the $RTT$ relations:
this restricts the range of the parameters $q,r$.
Three such conjugations are known (cf. \cite{Schirrmacher}):
\sk
i)   $~~T^*=T$, i.e. the elements $\T{A}{B}$ are ``real".  
Applying
the $*$-conjugation to the $RTT$ equations (\ref{RTTGL}) 
yields again the $RTT$ relations if the $R$ matrix satisfies
${\bar R}=R^{-1}$. This happens for 
$|q_{AB}|=|r|=1$, i.e. for
deformation parameters lying on the unit circle in ${\bf C}$
(cf. eq. (\ref{Rinv})).  The quantum group is then 
denoted by
$GL_{q,r}(N;{\bf R})$.
\sk
ii)   $~~(\T{A}{B})^*=\T{A'}{B'}$ with primed indices 
defined as {\small $A'=N+1-A$}.  Here compatibility
with the $RTT$ relations requires 
$\Rbar{AB}{CD}=\R{B'A'}{D'C'}$,
satisfied when ${\bar q}_{AB}=q_{B'A'},~r \in {\bf R}$. 
\sk
iii) $~~(\T{A}{B})^*=\kappa (\T{B}{A})$, the 
generalization of the  
unitarity
condition for the matrix $T$. In this case (left as 
an exercise in 
\cite{Schirrmacher}) the restriction on the $R$
matrix is $\Rbar{AB}{CD}=\R{DC}{BA}$, leading to 
the conditions
${\bar q}_{AB}=q_{BA},~r \in {\bf R}$. The 
corresponding quantum
groups are denoted by $U_{q,r}(N)$. 
\sk
Imposing  also $\detqr T=I$ yields the quantum groups
$SL_{q,r}(N;{\bf R})$ or
$SU_{q,r}(N)$.


\sect{The universal enveloping algebra of $\GLqrN$}


We construct the universal enveloping algebra of $\GLqrN$ as the
algebra of regular functionals \cite{FRT} on $\GLqrN$: 
it is generated
by the functionals $\LLpm , \epsi$ and $\Phi$ defined below. 
\sk
The  $\LLpm$ functionals are defined as in Section 2.2 where the uniparametric
$R$-matrix is now replaced by the multiparametric one.

A determinant can be defined for the matrix 
$\Lpm{A}{B}$, as in Note 2.2.1, this is given by:
\eq
\detqrm \LLpm=\Lpm{1}{1} \Lpm{2}{2} \cdots \Lpm{N}{N} \label{detL}~.
\en
A quantum inverse for $\Lpm{A}{B}$ can be found, using 
an expression analogous to (\ref{Tinverse}) 
with $q_{AB} \rightarrow 
\qm_{AB}$. For this we need to introduce the 
element $\Phi$ defined by:
\eq
\Phi\detqrm L^+\detqrm L^- =  \detqrm L^+\detqrm L^- \Phi = \epsi ~.
\en
Then the quantum inverse of
 $\Lpm{A}{B}$ is given by:
\eq
(\Lpm{A}{B})^{-1}=\Phi ~\detqrm L^{\mp} ~ \pp{BA}{1,N}
~ \ell_B^{~A}  \label{Linverse}
\en
\noi where $\ell_B^{~A}$ is the quantum minor and $\pp{BA}{1,N}$
is given in (\ref{defPi}). Notice that $\Phi ~\detqrm L^{\mp}$ is the
inverse of $\detqrm L^{\pm}$ because of property (\ref{LL3}) below.
\sk
\sk
The co-structures of the algebra generated by the 
functionals $L^{\pm}$, $\epsi$ and $\Phi$
are as in Section 2.2 :\eqa
& & \Dp (\Lpm{A}{B})=\Lpm{A}{G} \otimes \Lpm{G}{B}\\
& & \ep (\Lpm{A}{B})=\de^A_B \label{couLpm}\\
& & \kp (\Lpm{A}{B})=\Lpm{A}{B} \circ \kappa \label{coiLpm}\\
& &\Dp(\detqrm \LLpm)=\detqrm \LLpm \otimes \detqrm\LLpm,\\
 & &\Dp(\Phi)=\Phi\otimes\Phi, \Dp(\epsi)=\epsi\otimes \epsi\\
& &\ep(\detqrm\LLpm)=1,~\ep(\Phi)=1,~\ep(\epsi)=1\\
& &\kp(\detqrm\LLpm)=\Phi~\detqrm\LLmp,~\\
& &\kp(\Phi)=\detqrm \LLp \detqrm \LLm,~\kp(\epsi)=\epsi
\ena
\sk
\noi{\bf Note} 3.2.1 $~$ In (\ref{coiLpm}) we have defined 
$\kp$ using $\kappa$, 
we now prove that $\kp(\Lpm{A}{B})=(\Lpm{A}{B})^{-1}$
as defined in (\ref{Linverse}).
 This shows that \(\kp(\Lpm{A}{B})\) is expressible by
polynomials in
\(\Lpm{A}{B}, \Phi\).
\sk
\noi {\sl Proof :} From  \((\LLpm)^{-1}\LLpm=\epsi\) we have
\({ 1}=[(\LLpm_1)^{-1}\LLpm_1](T)=(\LLpm_1)^{-1}(T_2)\LLpm_1(T_2)=
(\LLpm_1)^{-1}(T_2)R_{12}^{\pm}\) so that 
\((\LLpm_1)^{-1}(T_2)={R_{12}^{\pm}}^{-1}\).

\noi {}From $\kappa(T)T= 1$ 
we similarly have
\([\kp(\LLpm_1)](T_2)={R_{12}^{\pm}}^{-1}\) and therefore
$\kp(\Lpm{A}{B})=(\Lpm{A}{B})^{-1}$.
\sk

Since $\kp$ is an inner operation in the algebra generated by the
functionals \(\Lpm{A}{B},~\epsi\) and \(\phi\) we 
conclude that these 
elements generate the Hopf algebra $\UglqrN$  of the regular 
functionals on the quantum
group $\GLqrN$. 
\sk

In the following we list some useful properties of the $\LLpm$
functionals.
\sk
\noi {\bf Properties of $\LLpm$}
\sk
\noi {\bf i)} Similarly to the uniparametric case, cf. Note 2.2.1 -- 
Note 2.2.3,
we have
\eq
\Lpm{A}{A} \Lpm{B}{B}=\Lpm{B}{B} \Lpm{A}{A}~;~~ 
\Lp{A}{A} \Lm{B}{B}=\Lm{B}{B} \Lp{A}{A} \label{LL2}
\en
\noi As a consequence:
\eq
\detqrm \LLp \detqrm \LLm = \detqrm \LLm \detqrm \LLp.
\label{LL3}
\en
We also have
\eq
\Lpm{A}{B} (\detqr T)= \de^A_B (c^{\pm})^N r^{\pm 1}  
Q_A^{-1}\label{LondetT}
\en
\eq
\detqrm \LLpm (\T{A}{B})=\de^A_B (c^{\pm})^N
 r^{\pm 1}Q_A \label{LL4}
\en
\eq
\detqrm \LLpm (\T{A}{B})= \Lpm{A}{B} (\detqr T) Q_A^2.
\en
\noi {}From (\ref{detL}) it is easy to see that 
$\detqrm \LLpm (I)=1 \label{LL5}.$
\sk
\noi {\bf ii)}  Since the $RLL$ relations are 
the same as the $RTT$ relations with 
$q_{AB}\rightarrow(q_{AB})^{-1}$,  $r \rightarrow r^{-1}$, 
we obtain a formula analogous to (\ref{detcomm}):
\eq
(\detqrm \LLpm)\Lpm{A}{B}={Q_B\over Q_A}\Lpm{A}{B}
(\detqrm \LLpm)
\label{LpmdetqrmLpm},
\en
moreover
\eq
(\detqrm \LLmp)\Lpm{A}{B}={Q_B\over Q_A}
\Lpm{A}{B}(\detqrm \LLmp)
\label{LpmdetqrmLmp}~.
\en
\sk
\noi {\bf iii)} From (\ref{LpmdetqrmLpm}) and 
(\ref{LpmdetqrmLmp}) the following 
element:
\eq
\detqrm \LLp (\detqrm \LLm)^{-1}= 
(\detqrm \LLm)^{-1}\detqrm \LLp
\en
is seen to be central.  Notice that it is 
also group-like since
\eq
\Dp (\detqrm L^{\pm})=\detqrm L^{\pm} \otimes 
\detqrm L^{\pm}.
\en
In general even if $\detqrm \LLp (\detqrm \LLm)^{-1}$
is central and group-like it is not equal to $\epsi$ because
\eq
\detqrm \LLp (\detqrm \LLm)^{-1}(\T{A}{B})=
(c^+)^N (c^-)^{-N} ~r^{2}\delta^{A}_{B}.
\en
\sk
\noi {\bf iv)} The elements $\Lp{A}{A} \Lm{A}{A}$ 
(no sum on {\small $A$})
 play a special role for particular
values of the deformation parameters $q_{AB}, r$;  if we set
\eq
\Lp{A}{A} \Lm{A}{A} \equiv \epsi_A
\en
\noi we leave as an exercise to deduce 
that (no sum on repeated indices):
\eq
\epsi_A(\T{B}{C})\equiv c^+ c^- \de^B_C ~{q^2_{AB} \over r^2},
{}~~\epsi_A(I)=1, ~~
\epsi_A(\Xi)=[\epsi_A (\detqr T)]^{-1}
\en
\eq
\epsi_A (ab)=\epsi_A (a) \epsi_A (b), ~~~\mbox{$a,b ~\in \GLqrN$}
\en
\eq
\kp (\Lpm{A}{A})=\Lmp{A}{A}\epsi_A^{-1}
\en
\eq
\epsi_A \epsi_B=\epsi_B \epsi_A,
~~~~\epsi_A \Lpm{B}{B}=\Lpm{B}{B} \epsi_A
\en
\eq
\detqrm \LLp \detqrm \LLm = 
\epsi_1 \cdots \epsi_N~;
\en
\eq
\kp (\detqrm L^{\pm})=\detqrm L^{\mp}
(\epsi_1 \cdots \epsi_N)^{-1}=(\epsi_1 \cdots \epsi_N)^{-1}
\detqrm L^{\mp}
\en
\sk

\noi{\bf Note} 3.2.2 $~$ When $\detqr T$ is central  
($Q_A=1$) we also have that
$\det \LLpm$ is central (cf.  (\ref{LpmdetqrmLpm}) 
and (\ref{LpmdetqrmLmp}) ). 
As in Note 2.2.2, for $Q_A=1$ and
$(c^{\pm})^{N} r^{\pm 1}=1 $, the functionals $\LLpm$and 
$\epsi$ generate
the Hopf algebra $U(sl_{q,r}(N))$, and we have 
the simplified relations:
\eqa
& & \detqrm \LLp (\detqrm \LLm)^{-1}=\epsi\\
& & [\Lpm{A}{B}](\detqr T)=\delta^A_B~~ 
{\mbox{ no sum on {\small$A$}}}\\
& & [\detqrm \LLpm](\T{A}{B})=\delta^A_B\\
& & [\detqrm \LLpm](\detqr T)=1
\ena
\sk
\noi{\bf Note} 3.2.3 $~$ When $q_{AB}=r $  we recover 
the standard uniparametric $R$ matrix,
we have also  $Q_A=1$ and, for $c^+c^-=1$, 
\eq
\forall\, A ~~\epsi_A=\epsi  
~~\mbox{ i.e. }  
\Lp{A}{A}\Lm{A}{A}=\Lm{A}{A}\Lp{A}{A}
=\epsi~. \label{epsiaepsi}
\en
In this case the Hopf algebra of functionals 
$\UglqrN$ is equivalent to 
the algebra generated by  the symbols $\LLpm , \Phi$ 
and $\epsi$ 
modulo relations 
(\ref{RLL}),(\ref{RLpLm}) and
(\ref{epsiaepsi}) \cite{FRT}.
\sk
\noi{\bf Note} 3.2.4  $\GLqrN$ and $U_{q,r}(gl(N))$ 
are graded Hopf algebras: 
$T^A{}_B$ has grade $+1$,
$\kappa(T^A{}_B)$ has grade $-1$, 
$I$ has grade $0$, 
$\det T$ has grade $+N$ etc., and similarly for $L^{\pm}$.
\sk

\noi {\bf Conjugation}
\sk

The canonical  $*$-conjugation on 
$\UglqrN$ induced by the  $*$-conjugation on $\GLqrN$ 
is given by:
\eq
\psi^{*} (a)\equiv {\overline {\psi(\kappa^{-1}(a^*))}} 
\label{sharpconjugation}
\en
\noi where $\psi \in \UglqrN$, $a \in \GLqrN$, and  the overline denotes
 the usual complex conjugation.
It is not difficult 
to determine the action
on the basis elements $\Lpm{A}{B}$.
The three $\GLqrN$ $*$-conjugations  i), ii), iii) of
the previous section induce
respectively the following conjugations on the 
$\Lpm{A}{B}$:
\eqa
& &i) ~~~~(\Lpm{A}{B})^*={\kappa'}^2(\Lpm{A}{B}) \cr
& &ii) ~~~(\Lpm{A}{B})^*={\kappa'}^2(\Lmp{A'}{B'}) \cr
& &iii)~~(\Lpm{A}{B})^*=\kappa'  (\Lmp{B}{A}) ~.\label{Lstar}
\ena


\sect{The quantum group $IGL_{q,r} (N)$}


The $q$-inhomogeneous group $\IGLqrN$  is freely generated by the
non-commuting matrix elements $T^A{}_B~ [A=(0,a) ; a: 1,..N]$, the
identity $I$ and the inverse $\xi$ of the $q$-determinant of $T$ as
defined in (\ref{qrdet}), modulo the relations:
\eq
T^0{}_a=0 \label{proietto}
\en
and the relations:  

\eqa
& &\R{ab}{ef} \T{e}{c} \T{f}{d}=\T{b}{f} \T{a}{e} \R{ef}{cd}
\label{RTTIGL1}\\
& &\R{ab}{ef} \T{e}{c} x^f= {q_{0c} \over r}
 x^b \T{a}{c} \label{RTTIGL2} \\
& &\R{ab}{ef} x^e x^f=r x^b x^a \label{RTTIGL3}\\
& &q_{0a}\T{a}{c} u=q_{0c} u\T{a}{c} \label{RTTIGL4}\\
& &q_{0a} x^a u = u x^a \label{RTTIGL5}
\ena
\noi where $x^a \equiv \T{a}{0}$ and $u\equiv T^0{}_0$.

It is not difficult to check that this algebra, 
endowed with the coproduct $\D$, the counit 
$\epsi$ and the coinverse
$\kappa$ defined by :
\eqa
\D(\T{A}{B})=\T{A}{C} \otimes \T{C}{A} ;
~ & \epsi(\T{A}{B})=\de^A_B
;~ & \kappa(T)=T^{-1} \label{TperT}\\
\D(\xi)=\xi\otimes \xi ;~ & \epsi(\xi)=1 ;
~ & \kappa(\xi)=\det T
\label{xiperxi} \\
\D(I)=I\otimes I ;~ & \epsi(I)=1 ;~ & \kappa(I)=I 
\label{IperI}
\ena

\noi where the quantum inverse of $T^A{}_B$ is given
by $\Ti{A}{B}=\xi~\pp{AB}{0,N}~ t_B^{~A} $ 
 [see eq. (\ref{defPi}):  $t_B{}^A$ is the 
quantum minor ], is a Hopf algebra. The proof goes
as in uniparametric case (see the second ref. of
\cite{inhom}).
\sk
In the commutative limit it is the algebra of 
functions on $IGL(N)$
plus the dilatation $T^0{}_0$..
\sk
Relations (\ref{TperT})-(\ref{IperI}) explicity read:
\eqa
& &\D(\T{a}{b})=\T{a}{c}
\otimes \T{c}{b},~~\D (I)=I\otimes I,\label{DT}\\
& & \D(x^a)=\T{a}{b} \otimes x^b + x^a \otimes u \label{Dx}\\
& & \D(u)=u\otimes u,~~\D(\xi)=\xi\otimes \xi \label{Du}\\
& & \D(\detqrTab)=\detqrTab \otimes \detqrTab
\label{Ddet}
\ena
\eqa
& & \epsi(\T{a}{b})=\de^a_b,~~\epsi (I)=1,\label{epsiT}\\
& & \epsi(x^a)=0 \label{epsix}\\
& & \epsi(u)=\epsi(\xi)=1 \label{epsiu}\\
& & \epsi(\detqrTab)=1 \label{epsidet}
\ena
\eqa
& & \kappa(\T{a}{b})=\Ti{a}{b}=\xi u ~\pp{ab}{1,N} t_b^{~a}\\
& & \kappa(I)=I, \label{kappaT}\\
& & \kappa(x^a)=-\kappa(\T{a}{b}) x^b \kappa (u) \label{kappax}\\
& & \kappa(u)= \detqrTab ~\xi \label{kappau}\\
& & \kappa(\xi)=u ~\detqrTab,~~\kappa(\detqrTab)=\xi u  
\label{kappadet}
\ena
\noi where for completeness we have included the 
expressions for the
$q$-determinant of $T$. Note that 
$\kappa (u) u = I = u \kappa(u)$.

\sk
This procedure is very similar to that 
discussed for $\GLqrNo$ 
in Section 3.1: indeed both these Hopf 
algebrae are obtained from 
the algebra freely generated by 
$T^A{}_B, I,\Xi\mbox{ or } \xi$ through the
introduction of moduli relations i.e. 
as quotients of suitable
two-sided ideals: the one generated by the 
$RTT$ relations in the
$\GLqrNo$ case, and the one generated by the
(\ref{proietto})-(\ref{RTTIGL5}) relations 
in the $\IGLqrN$ case. 
\sk
We now rederive the quantum group $\IGLqrN$ as a
quotient of $\GLqrNo$: all Hopf algebra properties
of $\IGLqrN$ will descend from 
those of $\GLqrNo$. The formalism
employed will be useful in the next 
section to deduce the differential
calculus on $\IGLqrN$ from the one on $\GLqrNo$.
\sk
We start from the observation that the $R$-matrix of $\GLqrNo$ can
be written as  ({\small A=(0,a)}):
\eq
\R{AB}{CD}=\left( \begin{array}{cccc}
                  r&0&0&0\\0&{r\over q_{0b}}\de^b_d&0&0\\
                 0&(r-r^{-1})\de^a_d &{q_{0a}\over r} \de^a_c&0\\
                  0&0&0&\R{ab}{cd}\\
                \end{array} \right) \label{RGLqrNo}
\en
\noi where $\R{ab}{cd}$ is the $R$-matrix of $GL_{q,r}(N)$, and
the indices {\small AB} are
ordered as {\small $00,0b,a0,ab$}.

It is apparent that the $\GLqrNo$ $R$ matrix 
 contains the information on $\GLqrN$.  
We will show that it also contains the information
about the quantum group $\IGLqrN$.

In the index notation $A=(0,a)$ the $RTT$ relations explicity 
read :
\eqa
& &\R{ab}{ef} \T{e}{c} \T{f}{d}=\T{b}{f} \T{a}{e} \R{ef}{cd}
\label{exp0}\\
& &T^{a}{}_{c}T^{b}{}_{0} = \frac{q_{ab}}
{q_{c0}}T^{b}{}_{0}T^{a}{}_{c}
\label{exp1} \\
& &T^{a}{}_{0}T^{b}{}_{d} =
\frac{q_{ab}}{q_{0d}}T^{b}{}_{d}T^{a}{}_{0}
+\frac{r}{q_{0d}}(r-r^{-1})T^{a}{}_{d}T^{b}{}_{0} 
\label{exp2} \\
& &T^{a}{}_{c}T^{0}{}_{d} = \frac{q_{a0}}
{q_{cd}}T^{0}{}_{d}T^{a}{}_{c}
\label{exp3} \\
& &T^{0}{}_{c}T^{b}{}_{d} =
\frac{q_{0b}}{q_{cd}}T^{b}{}_{d}T^{0}{}_{c}
+\frac{r}{q_{cd}}(r-r^{-1})T^{0}{}_{d}T^{b}{}_{c} 
\label{exp4} \\
& &\T{a}{0} \T{b}{0}=q_{ab} \T{b}{0} \T{a}{0} \\
& &T^{0}{}_{c}T^{b}{}_{0} = \frac{q_{0b}}{q_{c0}}T^{b}
{}_{0}T^{0}{}_{c}
\label{exp5} \\
& &\T{0}{c} \T{0}{d}=q_{dc} \T{0}{d} \T{0}{c}\\
& &T^{0}{}_{0}T^{b}{}_{d} =
\frac{q_{0b}}{q_{0d}}T^{b}{}_{d}T^{0}{}_{0}
+\frac{r}{q_{0d}}(r-r^{-1})T^{0}{}_{d}T^{b}{}_{0} 
\label{exp6} \\
& &T^{0}{}_{0}T^{b}{}_{0} = {q_{0b}}T^{b}{}_{0}T^{0}{}_{0}
\label{exp7} \\
& &T^{0}{}_{0}T^{0}{}_{d} = {q_{d0}}T^{0}{}_{d}T^{0}{}_{0}
\label{exp8} 
\ena  
where $a<b$ and $c<d$.
\sk
Consider now in $\GLqrN$ the space $H$ of all sums of
monomials containing at least an element 
of the kind $T^0{}_a$
(i.e. $H$ is the ideal in $\GLqrNo$ generated 
by the elements $T^0{}_a$
as we will see). 
Notice that 
$T^{0}{}_{0}T^{b}{}_{d} -\frac{q_{0b}}
{q_{0d}}T^{b}{}_{d}T^{0}{}_{0}$
is an element of $H$ because of relation (\ref{exp6}).
\sk
We now prove that $H$ is a Hopf ideal, i.e. an ideal 
in the $\GLqrNo$
algebra that is also compatible with the 
co-structures of $\GLqrNo$;
this allows to structure ${\GLqrNo}/{H}$ as a Hopf algebra
\cite{Sweedler}. We denote by $\D_{N+1}$,  $\epsi_{N+1}$ and
$\kappa_{N+1}$ the co-structures of  $\GLqrNo$.
\sk
\noi {\bf Theorem} 3.3.1 $~$ The space $H$ is a Hopf ideal 
in $\GLqrNo$, that is:
\sk
\noi ${}~~{}$i) ${}~~{}$ $H$ is a two-sided ideal in $\GLqrNo$  
\sk
\noi ${}~{}$ii) ${}~~{}$ $H$ is a co-ideal i.e. 
\eq
{}~~~~~\DN(H)\subseteq H \otimes \GLqrNo + 
\GLqrNo\otimes H ~;~~ \eN (H) =0 
\label{coideal}
\en
\noi ${}{}$iii) ${}~~{}$ $H$ is compatible with
$~\kN$ :
\eq
\kN(H)\subseteq H~.\label{Hideal}
\en

\noi {\sl Proof:}
\sk
\noi ${}~~{}$i) $H$ is trivially a subalgebra of $\GLqrNo$.
It is a right and left ideal since $\forall h\in H,\/ \forall a\in
\GLqrNo ~ha\in H \mbox{ and }
 ah\in H.$
This follows immediately from the definition of $H$ as sums of
monomials containing at least a factor $T^0{}_a$.
$H$ is the
ideal in $\GLqrNo$ generated by the elements $T^0{}_a$.
\sk
\noi ${}~{}$ii) First notice that 
$\DN(\T{0}{b})\in H \otimes \GLqrNo + \GLqrNo \otimes H .$  Now by definition
of  $H$ we have
\eq
\forall ~h\in H ,~~~h=a~\T{0}{b} c,~~~~~a,c \in \GLqrNo . 
\en
where $a~ \T{0}{b} c$ represents a sum of monomials.
Then we find
\eq
\DN(h)=\DN(a)\DN(\T{0}{b})\DN(c) \in  H \otimes \GLqrNo + \GLqrNo \otimes H ~.
\en
Moreover, since $\eN$ vanishes on $\T{0}{b}$ we
have:
\eq
\eN (h)= 0,~~~~\forall h \in H.
\en
These relations ensure that (\ref{coideal}) hold.
\sk
\noi ${}{}$iii) 
\eq
\kN(T^0{}_b) =\Xi~\pp{ob}{0,N}~ t_b^{~0} 
\en
where $\pp{ob}{0,N}$ is defined in (\ref{defPi})   
and  it is easy to see that the quantum minor 
$t_b{}^0 \in H$ since it is
the determinant of a matrix that has elements  
$T^0{}_a$ in the first row. Then
\eq
\kN (h)=\kN(a~\T{0}{b} c)=\kN(c)\kN(\T{0}{b})\kN(a) 
\in H
\en
and  Theorem 3.3.1 is proved. \cvd
\sk
We now consider the quotient 
\eq
\frac{\GLqrNo}{H}~, 
\label{quotient}
\en
and the canonical
projection 
\eq
P ~:~~ \GLqrNo\longrightarrow  \GLqrNo/{H}  
\en
Any element of ${\GLqrNo/H}$ is of the form $P(a)$. Also, $P(H)=0$,
i.e. $H=Ker(P)$.
\sk
Since $H$ is a two-sided ideal, ${\GLqrNo/H}$ is an algebra with
the following sum and products:
\eq
P(a)+P(b)\equiv P(a+b) ~;~~ P(a)P(b)\equiv P(ab) ~;~~
\mu P(a)\equiv P(\mu a),~~~\mu \in \mbox{\bf C} \label{iglalgebra}
\en
We will use the following notation: 
\eq
x^a\equiv P(T^a{}_0) ~;~~ u\equiv P(T^0{}_0) ~;~~ \xi\equiv P(\Xi)
\en
and with abuse of symbols: 
\eq
T^a{}_b\equiv P(T^a{}_b) ~;~~ I\equiv P(I) ~;~~ 0\equiv P(0)
\en
notice that $P(T^0{}_a)=P(0)=0$.
Using (\ref{iglalgebra}) it is easy to show that 
$T^a{}_b , ~x^a , ~u, ~\xi$ and $I$ 
generate the algebra ${\GLqrNo/H}$. From 
the $RTT$  relations
$R_{12}T_1T_2=T_2T_1R_{12}$ in $\GLqrNo$  we find the ``$P(RTT)$"
relations in ${\GLqrNo/H}$:
\eq
P(R_{12}T_1T_2)=P(T_2T_1R_{12}) ~~~ i.e. ~~~
R_{12}P(T_1)P(T_2)=P(T_2)P(T_1)R_{12}  \label{PRTT}
\en
that are explicity given in  (\ref{RTTIGL1})-(\ref{RTTIGL5}).
\sk
Since H is a Hopf ideal then ${\GLqrNo/H}$ is also a Hopf
algebra with co-structures:
\eq
\D(P(a))\equiv (P\otimes P)\DN(a) ~;~~ \epsi(P(a))
\equiv\eN(a) ~;~~
\kappa(P(a))\equiv P(\kN(a)) \label{co-igl}
\en
Indeed (\ref{coideal}) and (\ref{Hideal}) ensure  
that $\D ,~\epsi ,$
and $\kappa$ are well defined. For example
\eq
(P\otimes P)\DN(a) = (P\otimes P)\DN(b)~~ 
\mbox{ if }~~ P(a)=P(b) ~. \label{welldefined}
\en
In order to prove the Hopf algebra axioms of Appendix A for
$\D,~\epsi,~\kappa$ we just have to project 
those for $\DN,~\eN,~\kN~.$
For example, the first axiom is proved by applying
$P\otimes P \otimes P$ to $(\DN \otimes id)\DN(a) =
 (id \otimes \DN)\DN(a)$. The other axioms are proved
in a similar way.
\sk  
Notice that on the generators $T^a{}_b , ~x^a , ~u, 
~\xi$ and $I$ the 
co-structures (\ref{co-igl}) act as in  
(\ref{TperT})-(\ref{IperI}).
\sk
In conclusion: the elements $T^a{}_b , ~x^a , ~u, 
~\xi$ and $I$ 
generate the Hopf algebra ${\GLqrNo/H}$ and satisfy the
``$P(RTT)$'' commutation rules 
(\ref{RTTIGL1})-(\ref{RTTIGL5}). The 
co-structures act on them exactly as the co-structures defined in 
(\ref{TperT})-(\ref{IperI}). Therefore the quotient
${\GLqrNo/H}$ is the $q$-inhomogeneous group defined at the
beginning of this section:
\eq
\IGLqrN = \frac{\GLqrNo}{H}  ~.
\en
The canonical projection $P ~:~~ \GLqrNo \rightarrow 
\IGLqrN$ is an
epimorphism between these two Hopf algebrae.
\sk
\noi{\bf Note } 3.3.1 $~$  The consistency of the $P(RTT)$ 
relations with the co-structures $\D,\epsi$ and $\kappa$
is easily proved. For example, 
\eq
\D (P(R_{12}T_1T_2)-P(T_2T_1R_{12}))=0
\en
\noi is a particular case of eq.  (\ref{welldefined}).
Similarly for 
$\epsi$ and $\kappa$. 
\sk
We have thus obtained a $R$ matrix formulation of the inhomogeneous
$\IGLqrN$ quantum groups. Indeed the results of this section 
can be summarized in the following theorem:
\sk
\noi{\bf Theorem } 3.3.2 $~$The quantum inhomogeneous groups
$\IGLqrN$ are freely generated by the
non-commuting matrix elements $\T{A}{B}$
 [{\small A}=$(0, a)$, with $a=1,...N$)] and the identity  
$I$,
modulo the relations:
\eq
\T{0}{b}=0  \label{Tprojected}
\en
\noi and the $RTT$ relations
\eq
\R{AB}{EF} \T{E}{C} \T{F}{D} = \T{B}{F} \T{A}{E} \R{EF}{CD}
\label{RTTbig}
\en
The co-structures of $\IGLqrN$ 
are simply given by:
\eqa
& &\D (\T{A}{B})=\T{A}{C} \otimes \T{C}{B}\\
& &\kappa (\T{A}{B})= {T^{-1}}^D_C\\
& &\epsi (\T{A}{B})=\de^A_B \label{costructuresbig}
\ena
\cvd

\noi{\bf Note }  3.3.2 $~$ From the commutations 
(\ref{RTTIGL4}) - (\ref{RTTIGL5})
 we see that
one can set $u=I$ only when $q_{0a}=1$ for all $a$.
\sk
\noi{\bf Note } 3.3.3 $~$  $P(\detqrTAB)=u~\detqrTab$ 
is central in $\IGLqrN$
only when $Q_A=1$, {\small {\it A=0,1,..N}}
(apply the projection $P$ to eq. (\ref{detcomm})). Note
that here we have $Q_A \equiv 
\prod_{C=0}^N ({q_{CA}\over r})$.
\sk
\noi{\bf Note } 3.3.4 $~$ 
It is not difficult to see how the real forms of $\GLqrNo$  
are inherited by $\IGLqrN$. In fact, only the 
conjugation i) of $\GLqrNo$, discussed in Section 3.1, is
compatible with the coset structure of
$\IGLqrN$. More precisely, $H$ is a $*$-Hopf ideal, i.e.
$(H)^* \subseteq H$, only for
$T^*=T$. Then we can define a $*$-structure on $\IGLqrN$
as $[P(a)]^* \equiv P(a^*)$. 

\sk
\noi{\bf Theorem } 3.3.3 $~$ The centrality of $u$ is incompatible with
the centrality of $\detqrTab$.
\sk
\noi{\sl Proof:}  Suppose that
$q_{0a}=1$ so that $u$ is central.
Then the centrality of $\detqrTab$
is equivalent to the centrality
of $P(\detqrTAB)$ and requires $Q_A=1$ (Note 3.3.3);
in particular $Q_0 \equiv 
\prod_{c=1}^N {r\over q_{0c}}=1$, which cannot be since
for $q_{0a}=1$ we find $Q_0=r^N$. 
\sk
 \cvd
\sk
The commutations 
of $\det \T{a}{b}$ and $\xi$ with all the generators
are given by:
\eq
(\det \T{c}{d}) \T{a}{b}={\Qt_a \over \Qt_b} \T{a}{b}
 (\det \T{c}{d}),~~~~~\zeta \T{a}{b}={\Qt_b \over \Qt_a} \T{a}{b}
 \zeta
\en
\eq
(\det \T{c}{d}) x^a ={\Qt_a \over Q_0} x^a (\det \T{c}{d}),
~~~~~\zeta x^a ={Q_0 \over \Qt_a} x^a \zeta
\en
\eq
(\det \T{c}{d}) u = u (\det \T{c}{d}),
~~~~~\zeta u =u \zeta
\en
\noi where here $Q_a \equiv 
\prod_{c=1}^N ({q_{ca}\over r})$ and $\zeta$ is the
inverse of $\det \T{c}{d}$, i.e. $\zeta \equiv u \xi$.
We see that the commutations of $\det \T{c}{d}$ with
$\T{a}{b}$ are the correct ones for $\GLqrN$ (i.e. are
identical to the ones deduced in Section 3.1).
In the standard uniparametric case $Q_a=1$, and 
the $q$-determinant $\det \T{c}{d}$ becomes central
(and likewise $\zeta$),
provided that also $Q_0=1$.
\sk
We have derived the properties of the quantum group 
$\IGLqrN$ from those of $IGL_{q,r}(N+1)$, we now study the 
structure of $\IGLqrN$ with respect to its Hopf subalgebra
$\GLqrN$; this is explicitly done in Theorem 3.3.4, while in Theorem 
3.3.5 the same construction is seen in a more general and abstract setting.
\sk
We first notice that
the $x^0\equiv u$ and $x^a$ elements generate a subalgebra of
$\IGLqrN$ because their commutation relations  
do not involve the $\T{a}{b}$ elements.
Moreover these elements can be ordered using 
(\ref{RTTIGL3}) and (\ref{RTTIGL5}), and the Poincar\'e series of this 
subalgebra is the same as that 
of the commutative algebra in $N+1$ 
indeterminates, indeed (\ref{RTTIGL3}) and (\ref{RTTIGL5}) read
\eq
{}~~~~~~~~~~x^Ax^B=q_{AB}x^Bx^A~~~~~~~~\forall~\sma{$ A<B$}~.
\en 
A linear basis 
of this subalgebra is  therefore given by the ordered monomials:
$\zeta^i\equiv u^{i_{0}} (x^1)^{i_{1}}...\,
(x^N)^{i_{N}}$.
Then, using (\ref{RTTIGL2}) and (\ref{RTTIGL4}),
a generic element of $\IGLqrN$
can be written as $\zeta^ia_i$ where $a_i\in GL_{q,r}(N)$ and 
we conclude that $\IGLqrN$ is a right $GL_{q,r}(N)$--module 
generated by the ordered 
monomials $\zeta^i.$ Since the $RTT$ relations
of $\IGLqrN$ (\ref{RTTIGL1})--(\ref{RTTIGL5})
are homogeneous both in the $x^a$ and in $x^0$ we can naturally 
introduce a ({\bf{Z}},{\bf{N}}) grading : the generators 
$x^a$ have grade $(0,1)$,
$x^0$ has degree $(1,0)$, $(x^0)^{-1}$ has degree $(-1,0)$,
the elements of $\GLqrN$ have degree $(0,0)$.
Then 
\eq
\IGLqrN ~=\sum_{(h,k)\in
({\mbox{\scriptsize \bf Z}},{\mbox{\scriptsize \bf N}})}
{}^{{\hskip -0.64 cm}\oplus} ~~\Ga^{(h,k)}\label{grad1}
\en
where $\Ga^{(0,0)}=\GLqrN$, 
\eqa
&&\!\!\!\Ga^{(0,1)}=\{x^{a}b_{a}\;~/~~b_a\in GL_{q,r}(N)\}~~,~~~~
\Ga^{(\pm 1,0)}=\{({x^0})^{\pm 1}b\;~/~~b\in GL_{q,r}(N)\}\nonumber\\
&&\!\!\!\Ga^{(h,k)}=\{{(x^0)}^{h}
x^{a_1}x^{a_2}\ldots x^{a_k}b_{a_1a_2...a_k}~~/~~b_{a_1a_2...a_k}
\in GL_{q,r}(N)\}~~~~~~\sma{$\forall~h\in \mbox{\bf{Z}},~k\in
\mbox{\bf{N}}$}~.\nonumber
\ena
Therefore $\IGLqrN$ is a direct sum of right
$\GLqrN$-modules; it is also a graded algebra with the product
$\zeta^ib_i\cdot{\zeta'}^jb'_j$ trivially inherited from the 
$\IGLqrN$ algebra structure (in the sequel we will
omit the ``$\cdot$'').

We now show that 
each right module $\Ga^{(h,k)}$
is a bicovariant bimodule on $\GLqrN$, also 
$\IGLqrN=\sum_{(h,k)\in
({\mbox{\scriptsize \bf Z}},{\mbox{\scriptsize \bf N}})}^\oplus \Ga^{(h,k)}
$ is a bicovariant bimodule
with left and right coactions $\de_L$ and $\de_R$ that
are multiplicative: forall $\mbox{\sl a},\,\mbox{\sl b}\in\IGLqrN$, 
$\de_L(\mbox{{\sl a}{\sl b}})$ $=$ $\de_L(\mbox{\sl a})\de_L(\mbox{\sl b})$,
$\de_R(\mbox{{\sl a}{\sl b}})=\de_R(\mbox{\sl a})\de_R(\mbox{\sl b})$.
This shows that the structure
of a inhomogeneous quantum group is similar to that of the exterior
algebra of a generic Hopf algebra 
[as discussed at the end of point {\bf{iv)}}, Section 2.1];
also recall that, as noticed in the end of Subsection 2.4.6,  
the exterior algebra of forms is a Hopf algebra. 

\sk
\noi{\bf Theorem} 3.3.4  $\IGLqrN$, when  $q_{a0}=const~\forall a$,  
is a bicovariant graded algebra, i.e. it is a graded algebra 
with left and  right coactions
\eqa
&&\de_L~:~~\IGLqrN\rightarrow \GLqrN\otimes\IGLqrN\nonumber\\
&&\de_R~:~~\IGLqrN\rightarrow \IGLqrN\otimes\GLqrN\nonumber
\ena
that commute, see  (\ref{Bicco}), are multiplicative: 
$\de_L(\mbox{{\sl a}{\sl b}})=\de_L(\mbox{\sl a})\de_L(\mbox{\sl b})$,
$\de_R(\mbox{{\sl a}{\sl b}})=\de_R(\mbox{\sl a})\de_R(\mbox{\sl b})$,
$\forall ~\mbox{\sl a},\,\mbox{\sl b}\in\IGLqrN$, and preserve the grading .
\sk
\noi{\sl Proof }$~\;$
Consider the linear map 
$\de_R~:~~\IGLqrN\rightarrow \IGLqrN\otimes\GLqrN$
defined by
\eq
\de_R(x^A)=x^A\otimes I~;~~\de_R(a)=\D(a)~~~~
\forall\,a\in \GLqrN.
\label{derIGL}
\en
and extended multiplicatively on all $\IGLqrN$. 
This grade preserving map is obviously well defined on 
$\GLqrN$ because it coincides with the coproduct on $\GLqrN$
[$\GLqrN$ is the Hopf subalgebra of  $\IGLqrN$ with degree zero]; it is 
also well defined on 
all $\IGLqrN$ since it is multiplicative and compatible with 
(\ref{RTTIGL1})-(\ref{RTTIGL5}). We check for example 
(\ref{RTTIGL2}) with 
$q_{a 0}=const\equiv q_{0}\;~\forall a$:
$$
\de_R(x^{a}\T{b}{d})=
x^{a}\T{b}{c}\otimes\T{c}{d}=
{q_0\over r}R^{ba}{}_{\!ef}\T{e}{c}x^{f}\otimes
\T{c}{d}=\de_R\!\left({q_0\over r}R^{ba}{}_{\!ef}
\T{e}{d}x^{f}\right) .
$$
This shows that $\de_R\;:~\IGLqrN\rightarrow \IGLqrN\otimes \GLqrN$ 
is well defined.

To show that $\de_R$ is a right coaction notice that 
\eq
\forall \;\zeta^ia_i\, ,~~~(\de_R\otimes id)\de_R(\zeta^ia_i)=
(id\otimes \D )\de_R(\zeta^ia_i)~;~~(id\otimes\epsi )\de_R(\zeta^ia_i)=
\zeta^ia_i~.
\en

{}For the left coaction we proceed as in the previous case, defining the 
linear map  $\de_L:\IGLqrN\rightarrow \GLqrN\otimes \IGLqrN$,
\eq
\de_L(x^{a})=\T{a}{b}\otimes x^{b}~;~~
\de_L(x^0)=I\otimes x^0~;~~\de_L(a)=\D (a)~~~
\forall\,a\in \GLqrN
\label{delIGL}
\en
which is  extended multiplicatively on all $\IGLqrN$. 
As was the case for $\de_R$, it is well defined on 
$\GLqrN$  and it is also well defined on all $\IGLqrN$  because it 
is multiplicative and compatible with 
(\ref{RTTIGL1})-(\ref{RTTIGL5}).

To prove that  $\de_L$ is a left coaction  notice that
\eq
(\epsi\otimes id)\de_L(x^{a})=x^{a}\;,~
(\D \otimes id)\de_L(x^{a})=\T{a}{d}\otimes\T{d}{b}\otimes x^{b}
=(id\otimes\de_L)\de_L(x^{a})~
\en
and similarly for $x^0$. Now since $\de_L(a)=\D {(a)}$ if 
$a\in \GLqrN$, and since $\de_L$ is
multiplicative, we have on all $\IGLqrN$:
\eq
(\epsi\otimes id)\de_L=id~;~~(\D \otimes id)\de_L=(id\otimes \de_L)\de_L~.
\en
Finally, the compatibility of $\de_L$ and $\de_R$:
\eq
(id\otimes\de_R)\de_L=(\de_L\otimes id)\de_R~
\label{Bicco}
\en
follows directly
from:
\eqa
&(id\otimes \de_R)\de_L(x^{a})=\T{a}{b}\otimes x^{b}\otimes I=
(\de_L\otimes id)\de_R(x^{a})&\nonumber\\
&(id\otimes \de_R)\de_L(x^0)= I\otimes x^0
\otimes I=(\de_L\otimes id)\de_R(x^0)&\nonumber
\ena
\cvd
\sk
\noi{\bf Corollary }3.3.5 $~$
$\IGLqrN$, for $q_{a0}=const~\forall a$, is a bicovariant bimodule over 
$\GLqrN$ freely generated, as a right module, 
by the elements $\zeta^i$; also any submodule $\Ga^{(h,k)}$ is a 
bicovariant bimodule freely generated by the elements $\zeta^i$ with degree
\sma{$(h,k)$}.
\sk
$\!${\sl Proof }$\,$ 
We immediately have that $\IGLqrN$ 
and $\Ga^{(h,k)}$ are bimodules with the left module 
structure trivially inherited from the algebra $\IGLqrN$.
$\IGLqrN$ is a bicovariant bimodule
because, since the left and right coactions $\de_L$ and $\de_R$
are multiplicative, they are compatible with the left and right product of 
$\GLqrN$ on $\IGLqrN$; moreover they satisfy (\ref{Bicco}). 
Also the submodules $\Ga^{(h,k)}$ 
are bicovariant bimodules since the coactions $\de_L$ and $\de_R$ are grade 
preserving.

We now recall that a bicovariant bimodule is always
freely generated by a basis of 
right invariant elements, [cf. the text after (\ref{eta})].
We also know that the $\zeta^i$
are right invariant. Now, since they generate $\IGLqrN$, they linearly span 
the space of right invariant elements [$\IGLqrN$]${}_{\rm{inv}}$, 
and since they are linearly independent,
they form a basis of [$\IGLqrN$]${}_{\rm{inv}}$.
We conclude that $\IGLqrN$ is freely generated by the $\zeta^i$:
$\zeta^ia_i=0\Rightarrow a_i=0~\forall \sma{$i$}$.
The same arguments
apply also to each submodule $\Ga^{(h,k)}$. 
\cvd

In conclusion, the Hopf algebra $\IGLqrN$ is very rich because
it is both a bicovariant and a graded algebra on $\GLqrN$.
The bicovariant structure of $\IGLqrN$ can be seen as an example of 
a general 
theory by Radford \cite{Radford} on the properties of Hopf algebras 
$A$ with 
a Hopf subalgebra $H$ that is also a quotient of $A$. 
On the structure of inhomogeneous quantum groups see also the last reference 
in \cite{inhom}.  We summarize some results of 
\cite{Radford} in the following 
\sk
\noi{\bf Theorem } 3.3.6 $~$ Let $A$ and $H$ be Hopf algebras and 
suppose 
there exist Hopf algebras homomorphisms $\pi\;:~A\rightarrow H$ and
$i\;:~H\hookrightarrow A$ such that 
$\pi{\scriptstyle{{}^{{}_{\circ}}}}i=id_H$. Consider the projection $\Pi$ on 
$A$ defined by:
\eq
\Pi(a)=a_1i[\kappa\pi (a_2)] ~~a\in A
\en
and let
\eq
B\equiv \Pi(A)~.
\en
Then:
\begin{itemize}
\item[a)]{$B$ is a subalgebra of $A$, $\D(B)\subseteq A\otimes B$ and
\eq
B\equiv\Pi(A)=\{b~/~~b_1\otimes \pi(b_2)=b\otimes I\}.\label{matemqplane}
\en}
\item[b)]{$B$ is also a coalgebra with  counit $\underline{\epsi}$ that is the 
restriction to $B$ of the counit $\epsi$ of $A$, and with coproduct 
${\underline{\D}}$
given by 
\eq
{\underline{\D}}(\Pi(a))=\Pi(a_1)\otimes \Pi(a_2)~.\label{underlineD}
\en [Notice that ${\underline{\D}}$ is in 
general not compatible with the algebra structure of $B$, only
(\ref{prop1}) and (\ref{prop2}) hold].}
\item[c)]{$B$ is an $H$-bicovariant algebra with trivial right action 
and coaction and with left action given by the adjoint map 
$a\!d_h b=i(h_1)b i(\kappa(h_2))~~ \forall h\in H, ~\forall b\in B$ ,
and left coaction
given by $\de_L(b)=\pi(b_1)\otimes b_2\,\in H\otimes B\,,~\forall b\in B$.}
\item[d)]{As a  coalgebra $B$ is compatible with the left action
$a\!d$ and with the left coaction $\de_L$ (we use the notation
${\underline{\D}}(b)=b_{\underline{1}}\otimes b_{\underline{2}}$):
\eqa
&&{\underline{\D}}(a\!d_h b)=a\!d_{h_1}b_{\underline{1}}\otimes 
a\!d_{h_2}b_{\underline{2}}~~,~~~
{\underline{\epsi}}(ad_h b)=\epsi(h){\underline{\epsi}}(b)~,\nonumber\\ 
&&(id\otimes {\underline{\D}})\de_L(b)=
(m_H\otimes id)(\de_L\otimes\de_L){\underline{\D}}(b)~,\nonumber\\
&& (id \otimes {\underline{\epsi}})\de_L(b)=
{\underline{\epsi}}(b)I_H~.\nonumber\\
\ena
Moreover
\eq
{\underline{\D}}(bb')=b_{\underline{1}}\,a\!d_{b_{\underline{2}}{}^{(1)}}
b'_{\underline{1}}
\otimes b_{\underline{2}}{}^{(2)}b'_{\underline{2}}\label{braiddelta}
\en
where we have used the notations $\de_L(b)=b^{(1)}\otimes b^{(2)}$.}
\item[e)]{$B\otimes H$ has a canonical Hopf algebra structure 
(cross-product and cross-coproduct construction)  denoted
$B{\mbox{$\times \!\rule{0.3pt}{1.1ex}\;\!\!\!\cdot\,$}}H$.
The product is given by:
\eq
(b\otimes h)(b'\otimes h')=b(a\!d_{h_1}b')\otimes h_2h'~~~ \forall h\in H,~
b\in B
\en
the counit 
is given by ${\underline{\epsi}}\otimes \epsi$ and the coproduct is given by
\eq
\D(b\otimes h)=(b_{\underline{1}}\otimes b_{\underline{2}}{}^{(1)}h_1)\,
\otimes\,(b_{\underline{2}}{}^{(2)}\otimes h_2)~.
\en }
\item[d)]{$B{\mbox{$\times \!\rule{0.3pt}{1.1ex}\;\!\!\!\cdot\,$}}H$
and $A$
are isomorphic Hopf algebras via the isomorphism $\vartheta$:
\eq
\vartheta(b\otimes h)=bi(h)~~,~~~\vartheta^{-1}(a)=\Pi(a_1)\otimes \pi(a_2)
\label{ordino}
\en}
\end{itemize}
\cvd
\noi{\bf Note } 3.3.5 $~$
The algebra $B$ has a   braided Hopf algebra structure, and the
quantum group $A$ has a natural interpretation
via Majid bosonization procedure \cite{Majidbraidcmp}, \cite{Majidbraidalg}.
Since $B$ is a bicovariant bimodule over $H$, with trivial right action and 
right coaction [i.e. $B$ is a bimodule on the
quantum double ${\cal{D}}(H)$, cf. Note 2.3.3] the braiding $\Psi$ is 
given via the left action and the  left coaction
of $H$ on $B$:
\eq
\Psi(b\otimes b')=ad_{b^{(1)}}b'\otimes b^{(2)}
\en
then (\ref{braiddelta}) states that ${\underline{\D}}$ as defined in 
(\ref{underlineD})
is braided 
multiplicative: 
\eq
{\underline{\D}}(bb')=b_{\underline{1}}
\Psi(b_{\underline{2}}\otimes b'_{\underline{1}})b'_{\underline{2}}~.
\en
Finally $B$ is a braided Hopf algebra with antipode
${\underline{\kappa}}{(b)}=i[\pi(b_1)]\kappa(b_2)$.
\sk
In our case $A=\IGLqrN$, the projection $\pi\,:~\IGLqrN\rightarrow \GLqrN$,
that is well defined only if $q_{a 0}=const\equiv q_{0}\;~\forall a$,
is given by
\eq
\pi(\T{a}{b})=\T{a}{b}~,~~\pi(u)=I~,~~\pi(x^a)=0~,
\en
the subalgebra $B$ is generated by 
the $x^0\equiv u$ and $x^a$ elements 
that satisfy (\ref{RTTIGL3}) and (\ref{RTTIGL5}) 
i.e.
$
x^Ax^B=r{R^{-1}}^{AB}_{~CD}x^Dx^C
$.
The braiding is:
$\Psi(x^a\otimes x^b)
= {r\over q_0}{R^{-1}}^{AB}_{~CD}x^Dx^C$, 
$\Psi(x\otimes u) = u\otimes x $,  $\Psi(u\otimes x)=x\otimes u$,
the coproduct (coaddition) is  
${\underline{\D}}(u)=u\otimes u$, 
${\underline{\D}}(x^a)=x^a\otimes I +I\otimes x^a$
and the braided antipode and counit are 
$\underline{\kappa} (u)=u^{-1}$, 
$\underline{\kappa} (x^a)=-x^a$
and $\underline{\epsi}(u)=1$, 
$\underline{\epsi}(x^a)=0$.
\sk
\noi{\bf Note} 3.3.6 $~$ We also have $\IGLqrN = 
M{\mbox{$\times \!\rule{0.3pt}{1.1ex}\;\!\!\!\cdot\,$}}GL_{q,r}^u(N)$
where $M$ is the subalgebra of $\IGLqrN$ generated by the elements $x^au^{-1}$
and $GL_{q,r}^u(N)$ is the quantum group generated by $\GLqrN$ and 
the dilatation $u$.


\sect{Differential calculus on $\GLqrN$}


The bicovariant differential calculus on the uniparametric 
$q$-groups of the $A,B,C,D$ series can be formulated in terms
of the corresponding $R$-matrix and the associated
$\LLpm$ functionals.  This holds also for the multiparametric
case.  In fact all formulas are the same, modulo substituting
the $q$ parameter with $r$  when
it appears explicitly (typically as ${1 \over {q-q^{-1}}}$).
\sk
We list here some relevant formulae that do not appear in Section 2.2:
\eq
\ome{A_1}{A_2} \T{R}{S}=s  \Rinv{TB_1}{CA_1} 
\Rinv{A_2C}{B_2S} \T{R}{T}\ome{B_1}{B_2} \label{commomT}
\en
\eq
\ome{A_1}{A_2} \det T=s^N r^{-2} {Q_{A_1} \over Q_{A_2}}
 (\det T) \ome{A_1}{A_2} \label{commomdetT}
\en
\eq
\ome{A_1}{A_2} \Xi=s^{-N} r^{2} {Q_{A_2} \over Q_{A_1}}
 (\Xi) \ome{A_1}{A_2} \label{commomXi}
\en
\noi where we recall that $s \equiv (c^+)^{-1} c^-$, 
cf. eq.s (\ref{Rplus}) and (\ref{Rminus}).

Notice that $\det T$ and $\Xi$ commute with all the $\om$ 
(and thus can be set to $I$) iff
all $Q_A$ are equal and for $s^N r^{-2}=1$, or
$s=r^{{2\over N}} \al$ with $\al^N =1$, 
which agrees with  the condition found in Note 3.1.5.
\sk
Using 
\eq
da = ( \cchi{A_1}{A_2} * a) ~\ome{A_1}{A_2} \label{defd22}
\en
we compute the exterior derivative
on the basis elements of $\GLqrN$, and on the $q$-determinant:
\eq
d ~\T{A}{B}=\rinv [s~\Rinv{CR}{ET} \Rinv{TE}{SB} 
\T{A}{C}-\de^S_R \T{A}{B}]
~\ome{R}{S} \label{dTAB}
\en
\eq
d~ \Xi={{s^{-N} r^2-1} \over {r-\rminus}} ~\Xi~ \tau \label{dXi}
\en
\eq
d~ \detqr T={{s^N r^{-2}-1}\over{r-\rminus}} 
~(\detqr T) \tau \label{ddetT}
\en
The reader can verify via the Leibniz rule,  and with the 
help of eq. (\ref{commomdetT}), 
that $d [(\det T)\Xi]=d [\Xi (\det T)]=0$. 
{}From (\ref{ddetT}) we find that the bi-invariant $1$-form
$\tau$ that defines the exterior differential via
$
da=\rinv [\tau a - a \tau]
$
is given by 
\eq
\tau={r -r^{1}\over s^{N}r^{-2}-1}\Xi\/ d\,{\rm{det}}T~.
\en
again, $\det T=I=\Xi$ requires 
$s^N r^{-2}=1$.
\sk
The expression of
the $\om$ in terms of a linear
combination of $\kappa (T) dT$, similar to  the classical case is:
\eqa
\ome{A}{A}&=&{r\over {s(s-r^2-r^4+s r^4)}} [(r^2-s)
~\kappa(\T{A}{B} )d\T{B}{A}
+ r^2 (s-1)~\kappa (\T{C}{B}) d\T{B}{C} \theta^{CA}+\nonumber\\
& &+(-r^2-s r^2 +s +s r^4)~\kappa(\T{C}{B}) d\T{B}{C} 
\theta^{AC}],~~~
~~~\mbox{no sum on A}~~~
\label{omAA}\\
\ome{A}{B}&=&-s^{-1} {r\over q_{BA}} \kappa (\T{B}{C} )d\T{C}{A} 
,~~~~~~~~~~~~~~~~~~~~~~~~~~~~~~~A \not= B
\label{omAB}
\ena
\noi When $s=1$, the classical limit $\ome{A}{B} 
\qronelim -\kappa (\T{A}{C})
d\T{C}{B}$  reproduces the familiar formula 
$\om=-g^{-1}dg$ for the left-invariant
$1$-forms on the group manifold.  More generally, for 
$s=r^{\alpha}, \alpha \in {\bf C}$, we have:
\eq
\ome{A}{A} \ronelim [ {{2-\alpha} \over {2 (\alpha -1)}} 
\sum_B \kappa(\T{A}{B} )d\T{B}{A}
+{\alpha \over {2 (\alpha -1)}} \sum_B \sum_{C \not= A} 
\kappa (\T{C}{B}) 
d\T{B}{C}],~~~\mbox{no sum on A},
\en
\noi which shows that the inversion formula (\ref{omAA})
diverges in the classical limit for
$s=r$.
\sk
\noi {\bf Conjugation}
\sk
Compatibility of the conjugation defined in 
(\ref{sharpconjugation})
with the differential calculus requires
$(\chi_i)^*$ to be a linear combination  
of the $\chi_i$.
From (\ref{Lstar}) and  (\ref{defchi2}), it is straightforward to find
how the $*$-conjugation acts
on the tangent vectors $\chi$. Only the conjugations i) and  iii)
are compatible
with the differential calculus:
\eqa
& & i)~~~~(\cchi{A}{B})^* = -r^{-2N} d_C \R{DA}{BC}  (\cchi{C}{D}) \nonumber\\
& & iii) ~~~(\cchi{A}{B})^*=(\cchi{B}{A})
\ena

Using the inversion formulae (\ref{omAB}) and (\ref{dTAB}), or using 
(\ref{om*chi}), one finds the
induced $*$-conjugation on the 
left-invariant $1$-forms (here $s=r^{\alpha},~ \al\in{\bf{R}}$):
\eqa
& & i)~~~~(\ome{A}{B})^*= r^{2N} \Rinv{BC}{DA} d^{-1}_{B} \ome{C}{D}~~, 
\nonumber\\
& & iii)~~~(\ome{A}{B})^*=-\ome{B}{A}~~.
\ena

\sect{Differential calculus on $\IGLqrN$}


In this section we present  a
bicovariant differential calculus on $\IGLqrN$,
based on  the following set of
functionals $f$ and elements $M$ :
\eqa
& &\ff{a_1}{a_2b_1}{b_2}= \kp (\Lp{b_1}{a_1}) \Lm{a_2}{b_2}
\nonumber\\
& &\ff{a_1}{a_2 0}{b_2}= \kp (\Lp{0}{a_1}) \Lm{a_2}{b_2}
\nonumber\\
& &\ff{0}{a_2b_1}{b_2}= \kp (\Lp{b_1}{0}) \Lm{a_2}{b_2} \equiv 0
\nonumber\\
& &\ff{0}{a_2 0}{b_2}= \kp (\Lp{0}{0}) \Lm{a_2}{b_2} \label{fin}
\ena
\eqa
& &\MM{b_1}{b_2a_1}{a_2} = \T{b_1}{a_1} \kappa (\T{a_2}{b_2})
\nonumber\\
& &\MM{b_1}{b_2 0}{a_2} = \T{b_1}{0} \kappa (\T{a_2}{b_2})
\nonumber\\
& &\MM{0}{b_2a_1}{a_2} = 0
\nonumber\\
& &\MM{0}{b_2 0}{a_2} = \T{0}{0} \kappa (\T{a_2}{b_2}) \label{Min}
\ena
The $f$ in (\ref{fin})
 are a subset of the $f$ functionals of
$\GLqrNo$, obtained by restricting  the indices
of  $\f{i}{j}$ 
to $i=ab$ and $i=0b$. The third $f$ is identically zero because
of upper triangularity of $L^+$, i.e. $\Lp{b_1}{0}=0$.

The elements $M\in \IGLqrN$ in (\ref{Min})
 are obtained with the same
restriction on the adjoint indices, and by projecting
with $P$. The effect of the projection is to replace
the coinverse in $\GLqrNo$, i.e. $\kappa_{N+1}$ ,
with the coinverse $\kappa$ of $\IGLqrN$ (see the
last of (\ref{co-igl})). The
third element in (\ref{Min}) becomes zero because of $P$.
\sk
\noi{\bf Theorem } 3.5.1 $~$ The functionals in (\ref{fin})
vanish when applied to elements of $Ker(P) \subset
GL_{q,r}(N+1)$.
\sk
{\sl Proof:} first one checks directly that the
functionals (\ref{fin}) vanish
when applied to $\T{0}{b}$.  This extends to any element
of the form $\T{0}{b} a$ ($a \in A$), i.e. to any element of
$Ker(P)$, because of the property (\ref{propf1}) and
the vanishing of the functionals in
(\ref{fzero}). 
\sk
\cvd
\sk
The theorem states that the $f$ functionals are
well defined on the quotient $\IGLqrN = \GLqrNo / Ker(P)$,
in the sense that the ``projected" functionals
\eq
\fb : \IGLqrN \rightarrow \Cb,~~
\fb (P(a)) \equiv f(a)~ , ~~~\forall a \in \GLqrNo \label{deffb}
\en
are well defined.
Indeed if $P(a)=P(b)$, then $f(a)=f(b)$ because 
$f(Ker(P))=0$. This holds for any
 functional $f$ vanishing on $Ker(P)$, not only for
the $\f{i}{j}$ functionals.
\sk
The product $fg$ of  two generic functionals 
vanishing on $KerP$ also vanishes
on $KerP$, because $KerP$ is a co-ideal (see 
Theorem  3.3.1): $fg(KerP)=(f\otimes g)\DN (KerP)=0$.
Therefore $\overline{fg}$ is well defined on $\IGLqrN$, and
\eq
\overline{fg}[P(a)]=fg(a)=(f\otimes g)\D_{N+1}(a)=
(\fb P \otimes \gb P) \DN (a)=
(\bar{f}\otimes \bar{g}) \D(P(a)) \equiv
 \bar{f} \bar{g}[P(a)] \label{fgbar}
\en
\noi (use the first of (\ref{co-igl}))
so that the product of $\bar{f}$ and $\bar{g}$ involves
the coproduct $\D$ of $\IGLqrN$. 
\sk
There is a natural way to introduce a coproduct on the $\fb$'s :
\eq
\D'\fb[P(a)\otimes P(b)]\equiv\fb[P(a)P(b)]=\fb[P(ab)]
=f(ab)=\D'_{N+1}f(a\otimes b) ~.\label{defdiDelta}
\en
It is also easy to show that 
\eq
\D'\fb^i{}_j = \fb^i{}_k\otimes \fb^k{}_j ~~\mbox{ i.e. }~ 
\fb^i{}_j[P(a)P(b)] = \fb^i{}_k[P(a)]\fb^k{}_j[P(b)] 
\label{cofb}
\en
with $i,j,k$ running over the restricted set of indices $ab, 0b$. 
Indeed due to 
\eq
\ff{0}{A_2 b_1}{B_2} \equiv 0,~~\ff{A_1}{0 B_1}{b_2} \equiv 0
\label{fzero}
\en
\noi (a consequence of upper and lower triangularity of $L^+$ and
$L^-$ respectively ), formulae (\ref{copf}) and (\ref{propf1})
involve only the $f^i{}_j$ listed in (\ref{fin}), which
 annihilate $KerP$. Then
\eq
\fb^i{}_j[P(a)P(b)]=\fb^i{}_j[P(ab)]=f^i{}_j(ab)
=f^i{}_k(a)f^k{}_j(b)=\fb^i{}_k[P(a)]\fb^k{}_j[P(b)]\label{fPaPb}
\en
and (\ref{cofb}) is proved.

With abuse of notations we
will simply write $f$ instead of $\fb$. Then the
$f$  in (\ref{fin}) will be seen as functionals on
$\IGLqrN$. Notice that with the same abuse of notations, 
the product, coproduct and antipode of the 
$f$ and $\fb$ functionals coincide.
\sk
\noi{\bf Theorem } 3.5.2 $~$ The right $A$-module ($A=\IGLqrN$)
 $\Ga$ freely generated
by $\om^i \equiv \ome{a_1}{a_2}, \ome{0}{a_2}$ is
a bicovariant bimodule over $\IGLqrN$
 with right multiplication:
\eq
\om^i a = (\f{i}{j} * a) \om^j,~~a \in \IGLqrN \label{omia}
\en
\noi where the $\f{i}{j}$ are given in (\ref{fin}), the $*$-product
is computed with the co-product $\D$ of $\IGLqrN$,  and
the left and right actions
of $\IGLqrN$ on $\Ga$ are given by
\eqa
& &\DL  (a_i \om^i) \equiv \D(a_i) I \otimes \om^i
 \label{DLin}\\
& &\DR (a_i \om^i) \equiv \D(a_i) \om^j \otimes \M{j}{i}
\label{DRin}
\ena
\noi where the $\M{j}{i}$ are given in (\ref{Min}) and from now on we
use the notation $\DL=\Delta_{\scriptstyle \Ga}$ and 
$\DR={}_{\scriptstyle \Ga}\Delta$.
\sk
{\sl Proof:} we prove the theorem
by showing  that  the
functionals $f$
and the elements $M$ listed in (\ref{fin}) and
(\ref{Min})
satisfy the properties (\ref{propf1})-(\ref{propM})
(cf. the theorem by Woronowicz discussed in the Section 2.1).
 It  is straightforward to verify directly that
 the elements $M$
in (\ref{Min}) do satisfy the properties
(\ref{copM}) and (\ref{couM}). 
We have already shown that
the functionals $f$ in (\ref{fin}) satisfy
(\ref{propf1}),
and property (\ref{propf2})
obviously also holds for this subset.

Consider now the last property (\ref{propM}).
We know that it holds for $\GLqrNo$.  Take  the free indices
$j$ and $k$  as $ab$ and $0b$, and apply the projection $P$ on
 both members of the equation.
It is an easy matter
to show that  only the $f$'s in (\ref{fin}) and the
$M$'s in (\ref{Min}) enter the sums: this
is due to the vanishing of some $P(M)$ and to (\ref{fzero}).
We still have to prove that the $*$ product
in (\ref{propM})
can be computed via the coproduct $\D$ in
$\IGLqrN$.
Consider the projection of property
 (\ref{propM}), written symbolically as:
\eq
P [ M (f \otimes id) \D_{N+1} (a)]=P [(id\otimes f)
\D_{N+1} (a) M ]~. \label{PMfa}
\en
Now apply the definition (\ref{deffb}) and the first of
(\ref{co-igl}) to rewrite 
(\ref{PMfa}) as
\eq
P(M)(\fb\otimes id)\D(P(a))=(id\otimes\fb)\D(P(a))P(M)~.
\en
This projected equation then  {\sl becomes}
 property (\ref{propM})  for the $\IGLqrN$
functionals
$f$ and adjoint elements $M$, with the
correct coproduct $\D$ of $\IGLqrN$. \cvd
\sk
Using the general formula (\ref{omia}) we can
deduce the $\om, T$ commutations for $\IGLqrN$:
\eqa
& &\ome{a_1}{a_2} \T{r}{s}=s  \Rinv{tb_1}{ca_1}
\Rinv{a_2c}{b_2s} \T{r}{t}\ome{b_1}{b_2}\\
& &\ome{a_1}{a_2} x^r= s {q_{0a_1} \over q_{0a_2}} x^r
\ome{a_1}{a_2} -(r-r^{-1}) {s r\over q_{0a_2}}
\T{r}{a_1} \ome{0}{a_2}\\
& &\ome{a_1}{a_2} u= s {q_{0a_1} \over q_{0a_2}} u~
\ome{a_1}{a_2} \\
& &\ome{a_1}{a_2} \det \T{a}{b}=s^N r^{-2} {Q_{a_1} \over Q_{a_2}}
 (\det \T{a}{b}) \ome{a_1}{a_2}\\
& &\zeta \ome{a_1}{a_2} = s^N r^{-2} {Q_{a_1} \over Q_{a_2}}
 \ome{a_1}{a_2} \zeta\\
& &\ome{0}{a_2} \T{r}{s}=s {r\over q_{0t}}
\Rinv{a_2t}{b_2s} \T{r}{t}\ome{0}{b_2}\\
& &\ome{0}{a_2} x^r= {s\over q_{0a_2}}
 x^r \ome{0}{a_2} \label{Vxcomm}\\
& &\ome{0}{a_2} u= {s\over q_{0a_2}}
 u \ome{0}{a_2}\\
& &\ome{0}{a_2} \det \T{a}{b}=s^N r^{-2} {Q_0 \over Q_{a_2}}
 (\det \T{a}{b}) \ome{0}{a_2}\\
& &\zeta \ome{0}{a_2} =s^N r^{-2} {Q_0 \over Q_{a_2}}
  \ome{0}{a_2}\zeta
\ena
\noi{\bf Note} 3.5.1 $~$ $u$ commutes with all
 $\om$ 's only if $q_{0a}=1$ (cf. Note 3.3.2) and $s=1$. This means that
$u=I$ is consistent with the differential calculus
on $IGL_{q_{0a}=1,r}(N)$ only if the additional condition
 $s=1$ is satisfied.
\sk
The 1-form $\tau \equiv \sum_a \ome{a}{a}$
is bi-invariant, as one can check by using
 (\ref{DLin})-(\ref{DRin}).
Then an exterior derivative on $\IGLqrN$ can be defined
 as in eq. (\ref{defd2}), and satisfies the Leibniz rule.
The alternative expression $da=(\chi_i * a) \om^i$
(cf. (\ref{defd22}))
continues to hold, where
\eqa
& &\cchi{a}{b}=\rinv [\ff{c}{c a}{b}-\de^a_b \epsi] \nonumber \\
& &\cchi{0}{b}=\rinv [\ff{c}{c 0}{b}] \label{chiin}
\ena
are the left-invariant vectors dual to the left-invariant
1-forms $\ome{a}{b}$ and $\ome{0}{b}$. 
They are functionals on $\IGLqrN$ and as a consequence of 
(\ref{cofb}) we have
\eqa
& &\Dp (\cchi{a}{b})=\cchi{c}{d} \otimes \ff{c}{da}{b}+
\epsi \otimes \cchi{a}{b} \label{copchi1}\\
& &\Dp (\cchi{0}{b})=\cchi{c}{d} \otimes \ff{c}{d0}{b}+
\cchi{0}{d} \otimes \ff{0}{d0}{b} +
\epsi \otimes \cchi{0}{b}\label{copchi2}
\ena
\sk
The exterior
derivative on the generators $\T{a}{b}$ is given by
formula (\ref{dTAB}) with lower case indices. For the other
generators we find:
\eqa
& & dx^a=-s {r\over q_{0s}} \T{a}{s} \ome{0}{s}+{{s-1}
\over {r-\rminus}}
x^a \tau \\
& & d u= {{s-1}\over {r-\rminus}} u \tau\\
& & d \xi={{s^{-N-1} r^2-1} \over {r-\rminus}} ~\xi~ \tau
\ena
\noi Moreover:
\eq
 d ( \det \T{a}{b})={{s^N r^{-2}-1}\over{r-\rminus}}
{}~(\det \T{a}{b}) \tau
\en
\eq
 d \zeta={{s^{-N} r^{2}-1}\over{r-\rminus}}
{}~\zeta \tau
\en
\noi ($\zeta \equiv u\xi$). Again we find that $u=I$
implies $s=1$, and $\det \T{a}{b}=\zeta=I$ requires
$s^N r^{-2}=1$.
\sk
Every element $\rho$ of $\Ga$ can be written as
$\rho=a_k db_k$ for some $a_k,b_k$ belonging to
$\IGLqrN$. In fact one has the
same formula as in (\ref{omAA}) for $\ome{m}{n}$,
where all indices now are lower case. For
$\ome{0}{n}$ we find:
\eq
\ome{0}{n}=-{q_{0n}\over {s r}} [\kappa (\T{n}{a})dx^a+
 \kappa (x^n) d u]
\en
Finally, the two properties (\ref{leftco}) and (\ref{rightco})
hold also for $\IGLqrN$, because of the bi-invariance
of $\tau=\ome{a}{a}$. 
Thus all the axioms for a bicovariant 
first order differential calculus on $\IGLqrN$
are satisfied.
\sk
The exterior product of left-invariant $1$-forms
is defined as
\eq
\om^i \we \om^j\equiv \om^i \otimes \om^j - \L{ij}{kl}
\om^k \otimes \om^l
\en
\noi where
\eq
\L{ij}{kl}=\f{i}{l} (\M{k}{j})\label{thiscaninfactbe}
\en
This $\Lambda$ tensor can in fact
be obtained from
the one of $\GLqrNo$ by restricting its indices to the
subset $ab, 0b$. This is true because
when $i,l=ab$ or $0b$ we have 
$\f{i}{l} (Ker P)=0$ so that $\f{i}{l}$ is
well defined on $\IGLqrN$, and we can write 
$\f{i}{l} (\M{k}{j})=\fb^i_{~l} [P(\M{k}{j})]$
(see discussion after
Theorem 3.5.1).
The non-vanishing components of $\Lambda$ read:
\eqa
& &\LL{a_1}{a_2}{d_1}{d_2}{c_1}{c_2}{b_1}{b_2}=
d^{f_2} d^{-1}_{c_2} \R{f_2b_1}{c_2g_1} \Rinv{c_1g_1}{e_1a_1}
    \Rinv{a_2e_1}{g_2d_1} \R{g_2d_2}{b_2f_2} \label{L1}\\
& &\LL{0}{a_2}{d_1}{d_2}{c_1}{c_2}{0}{b_2}={q_{0c_2}\over
   q_{0c_1}} \Rinv{a_2 c_1}{g_2 d_1} \R{g_2d_2}{b_2c_2}
   \label{L2}\\
& &\LL{a_1}{a_2}{0}{d_2}{c_1}{c_2}{0}{b_2}=-(r-r^{-1})
{q_{0c_2}\over
   q_{0a_2}} \de^{c_1}_{a_1}\R{a_2d_2}{b_2c_2}
   \label{L3}\\
& &\LL{a_1}{a_2}{0}{d_2}{0}{c_2}{b_1}{b_2}={q_{0a_1}\over
   q_{0a_2}} d^{f_2} d^{-1}_{c_2} \R{f_2 b_1}{c_2 a_1}
   \R{a_2d_2}{b_2f_2}
   \label{L4}\\
& &\LL{0}{a_2}{0}{d_2}{0}{c_2}{0}{b_2}={q_{0c_2}\over
   q_{0a_2}} r^{-1} \R{a_2d_2}{b_2c_2}
   \label{L5}
\ena
These components still satisfy the characteristic
equation (\ref{RRHecke}), because the $\Lambda$ tensor
of $\GLqrNo$ does satisfy this equation, and
if the free adjoint indices are taken as
$ab$, $0b$, only the components in
(\ref{L1})-(\ref{L5}) enter in (\ref{RRHecke}).
To prove this, consider $\L{ij}{kl}$
with $k,l$  of the type  $ab$
or $0b$ and observe that  it vanishes
unless also $i,j$ are of the type $ab, 0b$.
(This can be checked
directly via the formula (\ref{Lambda})).
Then equations (\ref{commom}) and (\ref{defZ}) hold also
for the $\om$'s of $\IGLqrN$.

Note that $\La^{-1}$ tensor of $\IGLqrN$ can be obtained
by specializing the indices in the $\La^{-1}$ tensor
of $\GLqrNo$ given in (\ref{RRffMMinv}), as we did
for $\La$. The reader can convince himself of this
by  i) observing that the $\Linv{ij}{kl}$ tensor of (\ref{RRffMMinv})
also vanishes when $k,l$ = $ab$ or $0b$ and $i,j$ are not
of the type $ab, 0b$; ii) considering the equation
$\Linv{ij}{rs} \L{rs}{kl}=\de^i_k \de^j_l$ for $k,l$ = $ab$ or
$0b$.
\sk
The exterior differential on $\Ga^{\we n}$
can be defined as in Section 2.2 (eq. (\ref{defd})),
and satisfies all the properties (\ref{propd1})-(\ref{propd4}).
As for $\GLqrNo$ the last two hold because of the
bi-invariance of $\tau$.
\sk
The Cartan-Maurer equations are
\eq
d\om^i=\rinv (\tau \we \om^i+\om^i \we \tau)=
-\onehalf\c{jk}{i} \om^j \we \om^k
\en
\noi with
\eqa
& &\cc{a_1}{a_2}{b_1}{b_2}{c_1}{c_2}={2\over {r^2+r^{-2}}}
[-(r-r^{-1}) \de^{b_1}_{b_2} \de^{a_1}_{c_1} \de^{c_2}_{a_2}
+\CC{a_1}{a_2}{b_1}{b_2}{c_1}{c_2}]\\
& &\cc{a_1}{a_2}{0}{b_2}{0}{c_2}={2\over {r^2+r^{-2}}}
\CC{a_1}{a_2}{0}{b_2}{0}{c_2}\\
& &\cc{0}{a_2}{b_1}{b_2}{0}{c_2}={2\over {r^2+r^{-2}}}
[-(r-r^{-1})\de^{b_1}_{b_2} \de^{c_2}_{a_2}+
\CC{0}{a_2}{b_1}{b_2}{0}{c_2}]
\ena
The structure constants $\Cb$ (appearing in the $q$-Lie
algebra of $\IGLqrN$, see later) are given by
\eqa
& &\CC{c_1}{c_2}{b_1}{b_2}{d_1}{d_2}=
\rinv [-\de^{b_1}_{b_2} \de^{c_1}_{d_1}
 \de^{d_2}_{c_2}+\LL{a}{a}{d_1}{d_2}{c_1}{c_2}{b_1}{b_2}]
 \label{C1}\\
& & \mbox{~~~~~~~~~~~~~~= structure constants of $\GLqrN$}
 \nonumber \\
& &\CC{c_1}{c_2}{0}{b_2}{0}{d_2}=-{q_{0c_2}\over q_{0c_1}}
\R{c_1d_2}{b_2c_2} \label{C2}\\
& &\CC{0}{c_2}{b_1}{b_2}{0}{d_2}=
\rinv [-\de^{b_1}_{b_2}
 \de^{d_2}_{c_2}+ d^{f_2} d^{-1}_{c_2} \R{f_2b_1}{c_2a}
 \R{ad_2}{b_2f_2}]
\label{C3}
\ena

We conclude this section by observing that
the functionals $f$ and $\chi$
in (\ref{fin}) and (\ref{chiin})
close on the algebra
(\ref{bico1}), (\ref{bico2})-(\ref{bico4}), where
the product of functionals is defined by the coproduct
$\D$ in $\IGLqrN$. This result is expected,
since the functionals in (\ref{fin}) and (\ref{chiin})
correspond to a
bicovariant differential
calculus on $\IGLqrN$.\footnote{An explicit proof is also instructive.
We first note that in
$\GLqrNo$ the subset in (\ref{fin}) and (\ref{chiin})
closes by itself on the bicovariant
algebra (\ref{bico1}), (\ref{bico2})-(\ref{bico4}).
This is due to the particular
index structure of the tensors $\Cb$ and $\Lambda$, and to
the vanishing of the $f$ components in (\ref{fzero}).
The nonvanishing components of
$\Cb$ and $\Lambda$ that enter the operatorial
bicovariance  conditions (where the free adjoint indices
are taken as $ab, 0b$), are
given in (\ref{C1})-(\ref{C3}) and (\ref{L1})-(\ref{L5}).
Finally, we know that the $f$ functionals 
vanish on $KerP$, and so do the $\chi$ functionals
(as can be deduced from their definition in terms
of the $f$ functionals, eq. (\ref{defchi2})).
{} From the discussion after Theorem 3.5.1
it follows that they are well defined on $\IGLqrN$, and
that their products involve the $\IGLqrN$ coproduct
$\D$.

Thus the  relations (\ref{bico1}), (\ref{bico2})-(\ref{bico4})
hold for the functionals (\ref{fin}) and (\ref{chiin}) 
on $\IGLqrN$. They are the bicovariance conditions 
corresponding to a consistent differential calculus on $\IGLqrN$.}

\sk
Using the values of the $\Lambda$ and $\Cb$
tensors in (\ref{L1})-(\ref{L5}) and (\ref{C1})-(\ref{C3}),
we can explicitly write the ``$q$-Lie algebra"
of $\IGLqrN$  as:
\eq
 \cchi{c_1}{c_2}\cchi{b_1}{b_2}-
\LL{a_1}{a_2}{d_1}{d_2}{c_1}{c_2}{b_1}{b_2}~
\cchi{a_1}{a_2} \cchi{d_1}{d_2}=\rinv [-\de^{b_1}_{b_2}
\de^{c_1}_{d_1} \de^{d_2}_{c_2} +
\LL{a}{a}{d_1}{d_2}{c_1}{c_2}{b_1}
{b_2}] \cchi{d_1}{d_2}\label{qLie1}
\en
\eqa
& &\cchi{c_1}{c_2} \cchi{0}{b_2} + (r-r^{-1}) {q_{0c_2}\over
q_{0a_2}} \R{a_2d_2}{b_2c_2} \cchi{c_1}{a_2} \cchi{0}{d_2}
- ~~~~~~~~~~~~~~~~~~~~~~~~~~~~~~~~~~~~~~~~~~~~~~~
{}~~~~~\nonumber\\
& &~~~~~~~~~~~~~~~~~~~~~
 -{q_{0c_2}\over q_{0c_1}} \Rinv{a_2c_1}{g_2d_1}
\R{g_2d_2}{b_2c_2} \cchi{0}{a_2} \cchi{d_1}{d_2}
=-{q_{0c_2}\over q_{0c_1}} \R{c_1d_2}{b_2c_2} ~\cchi{0}{d_2}
\label{qLie2}\\
& &\cchi{0}{c_2} \cchi{b_1}{b_2}-{q_{0a_1}\over q_{0a_2}}
d^{f_2} d^{-1}_{c_2} \R{f_2b_1}{c_2a_1} \R{a_2d_2}{b_2f_2}
\cchi{a_1}{a_2} \cchi{0}{d_2}=\nonumber\\
& &~~~~~~~~~~~~~~~~~~~~~
\rinv [-\de^{b_1}_{b_2} \de^{d_2}_{c_2} +d^{f_2} d^{-1}_{c_2}
\R{f_2b_1}{c_2a}
\R{ad_2}{b_2f_2} ] \cchi{0}{d_2} \label{qLie3}\\
& &\cchi{0}{c_2}\cchi{0}{b_2}- {q_{0c_2}\over
q_{0a_2}} r^{-1}~\R{a_2d_2}{b_2c_2}
{}~\cchi{0}{a_2} \cchi{0}{d_2}=0 \label{qLie4}
\ena
\noi where $\LL{a_1}{a_2}{d_1}{d_2}{c_1}{c_2}{b_1}{b_2}$
is the
braiding matrix
of $GL_q(n)$, given in (\ref{L1}), so that the commutations in
(\ref{qLie1}) are those of the $q$-subalgebra $GL_q(n)$.
Note that the $\rone$ limit on the right hand sides
of (\ref{qLie1})
and (\ref{qLie3}) is finite, since the terms
in square parentheses
are a (finite) series in $r-r^{-1}$ whose $0-th$
order part vanishes [see (\ref{qLieexplicit}) and (\ref{Rlam})].

\section{The universal enveloping algebra of $\IGLqrN$}

In the  previous section we have considered the Hopf 
algebra generated by the $f$ functionals in (\ref{fin}).
This, togheter with $f_0^{~00}{}_0$, is the universal enveloping algebra 
$U_{q,r}(igl(N))$ of $\IGLqrN$ [see later, after (\ref{T'close})].
We now briefly give an $L^{\pm}$ description of  $U_{q,r}(igl(N))$.
For all the details we refer to Section 4.3 where a similar construction is
performed in the orthogonal case. 

In the preceding section we have identified the
$\overline f$ functionals on $IGL_{q,r}(N)$ with the corresponding $f$ 
functionals on $GL_{q,r}(N)$, in the same perspective,
we construct $U_{q,r}(igl(N))$ as a Hopf subalgebra of $U_{q,r}(gl(N+1))$.
Let
\eq
IU\equiv [L^{+A}{}_B, L^{-a}{}_b, L^{-0}{}_{0},
\Phi, \epsi]
\label{IUgl}
\en
be the subalgebra of $U_{q,r}(gl(N+1))$ generated by
$L^{+A}{}_B, L^{-a}{}_b, L^{-0}{}_{0},
\Phi, \epsi$. ($\Phi$ is the inverse of det$L^+$det$L^-$). 

The remaining $U_{q,r}(gl(N+1))$ generators
$L^{-b}{}_0$ are the only ones that do not annihilate 
$\T{0}{a}$
(the generators of $H$) and are not included in
(\ref{IUgl}): we construct the universal enveloping algebra of $\IGLqrN$ 
as the Hopf subalgebra of $U_{q,r}(gl(N+1))$ that annihilates the ideal $H$.
\sk
Since $\Delta(IU)\subseteq IU\otimes IU$ and $\kp(IU)\subseteq IU$
(as can be immediately seen at the generators level) we have that $IU$ is
a Hopf subalgebra of $U_{q,r}(gl(N+1))$. 
Moreover one can also give the following $R$-marix formulation (cf. 
Theorem 4.4.2):
\sk
\noi {\bf Theorem 3.6.1 }$~$
The Hopf algebra $IU$ is generated by $\epsi$, $\Phi$  and
the matrix entries:
$$
L^-=\left(L^{+A}{}_{B_{{}_{}}}\right)
{}\,,~
\Lcm=\Mat2{L^{-0}{}_{0}}{0}{0}{L^{-a}{}_b}
$$
these functionals satisfy the $q$-commutation relations:
\eq
R_{12} \Lcm_2 \Lcm_1=\Lcm_1 \Lcm_2 R_{12} ~~\mbox{ or equivalently }~~
{\cal{R}}_{12}
\Lcm_2 \Lcm_1=\Lcm_1 \Lcm_2 {\cal{R}}_{12} \label{iRLcLc}
\en
\eq
R_{12} \LLp_2 \LLp_1=\LLp_1 \LLp_2 R_{12}~, \label{iRLL}
\en
\eq
{\cal{R}}_{12} \LLp_2\Lcm_1 =\Lcm_1\LLp_2  {\cal{R}}_{12}~,  
\label{iRLpLm}
\en
where
${\cal{R}}_{12}\equiv c^-[\Lcm_1(T_2)]^{-1}$
that is
$
{\cal{R}}^{ab}_{cd}={R}^{ab}_{cd}\,,~
{\cal{R}}^{AB}_{AB}={R}^{AB}_{AB}$
and otherwise  
${\cal{R}}^{AB}_{CD}
=0~.$ 
{\cvd}
Relations (\ref{iRLcLc}) and (\ref{iRLL})  explicitly read as in
(\ref{exp0})--(\ref{exp8}), just substitute $T$ with $L^{\pm}$
and ``read from right to left''; this is due to
$R_{12} L_2^{\pm}L_1^{\pm}=L^{\pm}_1L^{\pm}_2R_{12}$ while
we have $R_{12} T_1T_2=T_2T_1R_{12}$. Relations (\ref{iRLpLm}) read
\eqa
& & R^{ab}_{~ef}\Lp{f}{d}\Lm{e}{c}=\Lm{a}{e}\Lp{b}{f}R^{ef}_{~cd}
\label{RLLIGL1}\\
& &\Lp{b}{d}\Lm{0}{0}={q_{d0}\over q_{b0}}\Lm{0}{0}\Lp{b}{d}\label{RLLIGL2}\\
& &\Lp{0}{d}\Lm{a}{c}={q_{a0}\over r}R^{ef}_{~cd}\Lm{a}{e}\Lp{0}{f}\\
& &\Lp{0}{d}\Lm{0}{0}={q_{d0}\over r^2}\Lm{0}{0}\Lp{0}{d}\\
& &\Lp{0}{0}\Lm{a}{c}={q_{a0}\over q_{c0}}\Lm{a}{c}\Lp{0}{0}\label{RLLIGL5}\\
& &\Lp{0}{0}\Lm{0}{0}=\Lm{0}{0}\Lp{0}{0}\label{RLLIGL6}
\ena
\noi{\bf Note } 3.6.1 $~$ 
Apply the second espression in (\ref{iRLcLc})
to $\T{A}{B}$ to obtain the quantum Yang-Baxter equation for the matrix  
${\cal{R}}$. The need for a new $R$-matrix ${\cal{R}}$ can be seen as due to 
the impossibility of considering
$IU$ as a quotient of the algebra $U_{q,r}(GL(N+1))$. The commutation 
relation that prevents $IU$
to be a quotient of $U_{q,r}(GL(N+1))$ with respect to the ideal generated by 
the $\Lm{a}{0}$ elements is: $\Lp{0}{d}\Lm{a}{0}$=${q_{ao}\over q_{0d}} 
\Lm{a}{0}\Lp{0}{d}$+$(1-r^{-2})q_{a0}
(\Lm{a}{d}\Lp{0}{0}-\Lp{a}{d}\Lm{0}{0})$.

We stress that $IU$ is a subalgebra of $U_{q,r}(N+1)$, so that 
(\ref{iRLcLc}),
(\ref{iRLL}), (\ref{RLLIGL1})--(\ref{RLLIGL6})
hold in $GL_{q,r}(N)$ as well. On the opposite side, the $R$-matrix
of the  $IGL_{q,r}(N)$ $RTT$ relations is the same as the $R$-matrix 
of the  $GL_{q,r}(N+1)$ $RTT$ relations, but this last set of $RTT$
relations does not contain as a subset the previous one. 
\sk

We now briefly study the structure of $IU$ with respect to
$U_{q,r}(gl(N))$, that is easily seen to be a  Hopf subalgebra of $IU$. It 
is also a quotient
of $IU$ via the Hopf algebra projection [well defined only if 
$q_{a0}=const\equiv q_0~\forall a$ see (\ref{RLLIGL2}) and (\ref{RLLIGL5})]:
\[
\pi(\Lm{0}{a})=0~~,~~~~\pi(\Lm{0}{0})
=I~~,~~~~\pi(\Lpm{a}{b})=\Lpm{a}{b}~~,~~~~\pi(\Lp{a}{0})=
\Lp{a}{0}~~,~~~~\pi({\epsi'})=\epsi'~.
\]
Then the results of Theorem 3.3.6  apply to $IU$ as well, and we can write
the Hopf algebra isomorhism $IU\cong$
$B'{\mbox{$\times \!\rule{0.3pt}{1.1ex}\;\!\!\!\cdot\,$}}U_{q,r}(gl(N))$
where $B'$ is the algebra generated by $\Lp{0}{0}$ and $\Lp{0}{a}$.
Also Theorem 3.3.4 hold for  $IU$  since we 
can introduce the following 
$(${\bf Z},{\bf N}$)$ grading: the elements $\Lpm{a}{b}$ have grade $(0,0)$, 
the elements $\Lp{0}{a}$
have grade $(0,1)$, the elements $\Lpm{0}{0} $
have grade $(\pm1,0)$. This grading is compatible with the $RLL$ commutation 
relations. Notice also that 
the elements $\Lp{0}{0}$ and $\Lm{0}{0}$ are not independent (see the text 
after (\ref{epsiaepsisosp}) for a general 
discussion in the orthogonal case, the $GL_{q,r}(N)$ case is similar); 
here we give an easier argument that holds only if 
$q_{a0}=const=r~\forall a$: we fix the coefficient 
$c^-$ defined in (\ref{Rminus}) and studied after (\ref{cmeno}),
to be $c^-=(c^+)^{-1}$ (notice that the parameter 
$s=(c^+)^{-1} c^-$ entering the differential calculus is still arbitrary, 
the parameter $c^-$ is completely irrelevant).
It follows that $\Lm{0}{0}$ has a simple dependence from $\Lp{0}{0}$:
$(\Lm{0}{0})^{-1}=\Lp{0}{0}$. [Proof: $\forall {\sma{$ A,B$}}\,,$ 
$(\Lm{0}{0}{})^{-1}(\T{A}{B})=$
$\Lp{0}{0}(\T{A}{B})$, $\D'(\Lpm{0}{0})=$ $\Lpm{0}{0}\otimes \Lpm{0}{0}$].

{}From the $RLL$ relations 
a generic element of $IU$
can be written as $\eta^ia_i$ (or $a_i\eta^i$) where $a_i\in U_{q,r}(gl(N))$ 
and $\eta^i$ are ordered 
monomials in the $\Lp{0}{0}$ and $\Lp{0}{a}$  elements:
$\eta^i=(\Lp{0}{0})^{i_0}(\Lp{0}{1})^{i_{1}}...\,
(\Lp{0}{N})^{i_{N}}$. 
As in Corollary 3.1.1, we have that 
$IU$, for $q_{a0}=const~\forall a$, is a bicovariant bimodule over 
$\GLqrN$ freely generated, as a right module, 
by the elements $\eta^i$; moreover 
\eq
\Uigl ~=\sum_{(h,k)\in
({\mbox{\scriptsize \bf Z}},{\mbox{\scriptsize \bf N}})}
{}^{{\hskip -0.64 cm}\oplus} ~~\Ga'^{(h,k)}\label{grad2}
\en
where $\Ga'^{(0,0)}=U_{q,r}(gl(N))$ 
\eqa
&&\!\!\!\Ga'^{(0,1)}=\{\Lp{0}{a}\varphi^{a}\;~/~~\varphi^a\in 
U_{q,r}(gl(N))\}~~,~~~~
\Ga'^{(\pm 1,0)}=\{(\Lp{0}{0})^{\pm 1}\varphi\;~/~~\varphi\in  
U_{q,r}(gl(N))\}\nonumber\\
&&\!\!\!\Ga'^{(h,k)}=\{{(\Lp{0}{0})}^{h}
\Lp{0}{a_1}\Lp{0}{a_2}\ldots \Lp{0}{a_k}\varphi^{a_1a_2...a_k}~~/~~
\varphi^{a_1a_2...a_k}
\in 
U_{q,r}(gl(N))\}
\ena
Any submodule $\Ga'^{(h,k)}$ is a
bicovariant bimodule freely generated by the elements $\eta^i$ with degree
\sma{$(h,k)\in (\mbox{\bf{Z}},\mbox{\bf{N}})$}. 
We leave to the reader to reformulate Note 3.3.5 and Note 3.3.6 in this 
context.  
\sk
\vskip .2cm
\noi {{{\bf{Duality}} $IU\leftrightarrow \IGLqrN$}}
\sk
We now show that $IU$ is dually paired to $\IGLqrN$. This is the  
fundamental
step allowing to interpret $IU$ as the universal 
enveloping  algebra 
of $\IGLqrN$.
\sk
\noi{\bf Theorem } 3.6.2 $~~IU$ annihilates $H$.

\noi {\sl Proof :} $~~~$ This theorem has implicitly been proved in
Theorem 3.5.1 and in the comments before (\ref{fgbar}). An explicit proof
is given in  Theorem 4.4.4.
{\cvd}
In virtue of {Theorem 3.6.2} the following bracket is well defined:
\eqa \mbox{{\sl Definition}. }&\!\!\!\!\!\!\!\!&
\le ~ ,~ \re ~ :~ IU\otimes \IGLqrN \longrightarrow
\mbox{\boldmath $C$}\nonumber\\
&\!\!\!\!\!\!\!\! &\le a', P(a)\re\equiv a'(a)\label{duality}\\
&\!\!\!\!\!\!\!\! &\forall \/a'\in IU ~,~\forall \/a\in \IGLqrN\nonumber
\ena
where $P~:~GL_{q,r}(N+1)\rightarrow GL_{q,r}(N+1)/H$ $ \equiv IGL_{q,r}(N)$ 
is the canonical
projection, which is surjective.
The bracket is well
defined because two generic counterimages of $P(a)$ differ
by an addend belonging to $H$.
\sk
Since 
$IU$ is a Hopf subalgebra
of $U_{q,r}(gl(N+1))$ and $P$ is compatible with the structures and 
costructures
of $GL_{q,r}(N+1)$ and $\IGLqrN$,
the following theorem is then easily shown [cf. (\ref{fgbar}), (\ref{fPaPb})
and  Theorem 4.4.5] 
\sk
\noi{\bf Theorem } 3.6.3 $~$ The bracket (\ref{duality})
defines a pairing between $IU$ and $\IGLqrN~\mbox{:}$ 
$\forall\/ a',b'\in IU~,~\forall\/P(a),P(b)\in \IGLqrN$
\eqa
& &\le a'b' , P(a)\re = \le a'\otimes b',\Delta  
(P(a))\re\label{uuno}\nonumber\\
& &\le a',P(a)P(b)\re=\le\Delta '(a'),P(a)\otimes  
P(b)\re\label{udue}\nonumber\\
& &\le\kp(a'),P(a)\re=\le a',\kappa(P(a))\re\label{utre}\nonumber\\
& &\le I,P(a)\re=\epsi(a)~;~~\le a',P(I)\re=\epsi  
'(a')\nonumber
\ena
\cvd
We now recall that $IU$ and $\IGLqrN$, besides being dually paired, 
are bicovariant algebras with the same 
graded structure (\ref{grad1}) and (\ref{grad2}),
and can both be obtained as a cross-product cross-coproduct construction:
$\IGLqrN\cong B{\mbox{$\times \!\rule{0.3pt}{1.1ex}\;\!\!\!\cdot\,$}}
GL_{q,r}(N)$,
$IU\cong$
$B'{\mbox{$\times \!\rule{0.3pt}{1.1ex}\;\!\!\!\cdot\,$}}U_{q,r}(gl(N))$.
In particular $\IGLqrN$ and $IU$ are freely generated (as modules) by $B$ and $B'$
i.e.
by the two isomorphic sets of the ordered 
monomials in the $q$-plane plus dilatation coordinates $\Lp{0}{0},\;\Lp{0}{a}$ 
and $u,\;x^a$ respectively. We then conclude that
$IU$ is the universal enveloping algebra of $\IGLqrN$:
\eq
~\Uigl\equiv IU~.\label{UigldualIG}
\en
\sk
\vskip .2cm
\noi{\bf{Projected differential calculus}}
\sk
We have found the inhomogeneous
quantum group  $\IGLqrN$ by means of a projection from
$GL_{q,r}(N+1)$; dually, its universal
enveloping algebra is a given  Hopf subalgebra of
$U_{q,r}(gl(N+1))$.
Using the same techniques 
we conclude this section presenting another differential calculus on $\IGLqrN$,
that is obtained from the previous one considering also the dilatation 
generator
$\chi^{0}{}_{0}$. To derive this calculus one can follow the same steps of the 
Section 3.5, however the easiest way to derive it is to apply the 
results of Section  2.3.
\sk
{}From (\ref{UNO}) and (\ref{DUE}) it is immediate to see that 
$T'\equiv T\cap U_{q,r}(igl(N))$ satisfies 
\eq
\D(T')\subset T'\otimes\epsi+
U_{q,r}(igl(N))\otimes T'
\en  
\eq
[T',T']\subseteq T\cap IU=T'\label{T'close}
\en
indeed $U_{q,r}(igl(N))$ is a Hopf subalgebra of $U_{q,r}(gl(N+1))$. 
Also condition (\ref{ZERO})
is fulfilled since  $T'$ generates $U_{q,r}(igl(N))$ in the same way 
$T$ generates $U_{q,r}(gl(N+1))$ \cite{Burr}, this is a consequence of
the upper and lower triangularity of the $L^+$ and $L^-$ matrices and of the
dependence of the diagonal elements of $L^+$ from the diagonal elements of 
$L^-$.
We therefore obtain an $IGL_{q,r}(N)$ bicovariant differential
calculus with $q$-Lie algebra generators:
\eq
\cchi{a}{b}~~,~~~~\cchi{0}{b}~~,~~~~\cchi{0}{0}\label{IGLtangent}~.
\en
Since these generators  close  on the  subalgebra $T'\subset T$,
we have that the structure constants that appear in the $\IGLqrN$ Lie algebra
are a subset of the structure constants that appear in the  $GL_{q,r}(N+1)$
Lie algebra [cf. the text after  (\ref{thiscaninfactbe})].

The exterior differential reads
\eq
da= (\cchi{a}{b}*a)\omega_{a}{}^b +
(\cchi{a}{0}*a) \omega_{a}{}^0 +
(\cchi{0}{0}*a)\omega_{0}{}^{0}~~;\nonumber
\en
where $\omega_{a}{}^{b},~\omega_{a}{}^{0},$
and $\omega_{0}{}^{0}$, following Section 2.3, 
are the $1$-forms dual
to the tangent vectors (\ref{IGLtangent}).

\sect{The multiparametric quantum plane as a quantum coset space}


In this section we derive the differential calculus 
on the quantum plane
\eq
Fun_{q,r}\left( {IGL(N) \over GL(N)}\right)~,
\en
\noi i.e. the subalgebra of $\IGLqrN$
generated by the coordinates $x^a$. These coordinates 
satisfy the commutations
(\ref{RTTIGL3}):
\eq
\R{ab}{ef} x^e x^f=r x^b x^a 
\en
\noi The main difference with the more conventional approach
to the quantum plane is that now the coordinates do not trivially
commute with the $\GLqrN$ $q$-group elements, but
$q$-commute according to relations (\ref{RTTIGL2}):
\eq
\R{ab}{ef} \T{e}{c} x^f= {q_{0c} \over r} x^b \T{a}{c} 
\label{RTTIGL2bis}
\en
{}From a more mathematical viewpoint the $q$-coset space
$Fun_{q,r}({IGL(N)/GL(N)})$ is the algebra $B$ discussed in Theorem 
3.3.6 and Note 3.3.5. $B$ is generated by $u$ and $x^a$. We then study
the subalgebra of $B$ generated only by the elements $x^a$.
Expression
(\ref{matemqplane}) is the translation in Hopf algebra language of the
classical property: 
$\forall\/ g \in ISO(N)\;,~\forall\/ h\in SO(N)~~
\{gh\}=\{g\} $ ,  where $\{g\}$ is an element of the coset  
$ISO(N)/SO(N)$.
\sk

\noi {\bf Lemma } 3.7.1  $~\cchi{b}{c} (a)=0$ when a is a 
polynomial in $x^a$ and $u$ with all monomials
containing at least one $x^a$. This is easily proved
by observing that no tensor exists with the correct
index structure.  For $s=1$ we can extend this lemma
even to  $u \cdots u$, since
for example 
\eq
\cchi{b}{c} (u)={{s-1}\over {r-\rminus}} \de^b_c
\en
\noi and using the coproduct rule (\ref{copchi1})
 one finds that
$\cchi{b}{c}(u \cdots u)$
is always proportional to 
$s-1$.
\cvd
\sk
\noi {\bf Theorem} 3.7.1 $~\cchi{b}{c} * a=0$ when $a$ is
a polynomial in $x^a$ and $s=1$.
\sk
\noi {\sl Proof:} we have $\cchi{b}{c} * a=(id\otimes 
\cchi{b}{c})(a_1\otimes
a_2)=a_1 \cchi{b}{c} (a_2)$. (We use the notation
$\D (a) \equiv a_1 \otimes a_2$).  Since $a_2$ is a polynomial
in $x^a$ and $u$ (use the coproduct rule (\ref{Dx})), and
$\cchi{b}{c}$ vanishes on such a polynomial when $s=1$
(previous Lemma),  
the theorem is proved.\cvd
\sk
Because of this theorem we will henceforth set $s=1$: then
we can write the exterior derivative of an element 
of the quantum plane as
\eq
da= (\chi_s*a) V^s \label{daflat}
\en
\noi (with $\chi_s \equiv \cchi{0}{s}$, $V^s \equiv \ome{0}{s}$), 
i.e. only in terms of the ``q-vielbein" $V^s$. Notice also that
$du=0$.
\sk
The action and value of $\chi_s$ on the coordinates 
is easily computed [cf. the 
definition in  (\ref{chiin})]:
\eq
\chi_s * x^a=-{r\over q_{0s}} \T{a}{s},~~~\chi_s (x^a)=-{r\over
q_{0s}} \de^a_s
\en
\noi so that the exterior derivative of $x^a$ is:
\eq
dx^a=-{r\over q_{0s}} \T{a}{s} V^s \label{dxa}
\en
\noi and gives the relation between the $q$-vielbein $V^s$
and the differentials $dx^a$.
\sk
Using (\ref{copchi2}), the Leibniz rule for the ``$q$-partial derivatives" 
$\chi_c$ is given by :
\eq
\chi_c * (ab)=(\chi_d * a) \f{d}{c} * b + a \chi_c * b 
\label{Leibnizplane}
\en
\noi where $\f{d}{c} \equiv \ff{0}{d0}{c}$.
\sk
The $x^a$ and $V^b$ $q$-commute as (cf. (\ref{Vxcomm})):
\eq
V^a x^b= (q_{0a})^{-1} x^b V^a
\en
\noi and via eq. (\ref{dxa}) and (\ref{RTTIGL2bis})
we find the $dx^a, x^b$ commutations :
\eq
dx^a x^b=r^{-1} \Rinv{ab}{ef} x^f dx^e
\en
\noi After acting on this equation with $d$ we obtain:
\eq
dx^a \we dx^b= -r^{-1} \Rinv{ab}{ef} dx^f \we dx^e
\en
\noi which reproduce the known commutations between 
the differentials
of the quantum plane, cf. ref. \cite{WZ}, \cite{zumi}.
\sk
The commutations between the partial derivatives are given in
eq.(\ref{qLie4}). 
\sk
All the relations of this section 
are covariant under the
$\IGLqrN$ action:
\eq
x^a \longrightarrow \T{a}{b} \otimes x^b + x^a\otimes u
\en
and in particular under the $\GLqrN$ action 
$x^a \longrightarrow \T{a}{b} \otimes x^b $.
\sk
\noi{\bf Note } 3.7.1 $~$ The partial derivatives $\chi_c$,
and in general all the tangent vectors $\chi$
of this chapter  have ``flat'' indices. To compare
them with the partial derivatives discussed in
\cite{WZ}, which have "curved" indices,
we need to define the operators $\chit_s$:
\eq
\chit_s (a)\equiv - {q_{0a}\over r} (\chi_a * a) 
\kappa (\T{a}{s})
\en
\noi whose  value and action on the coordinates is
\eq
\chit_s  (x^a)=\de^a_s I
\en
\noi so that
\eq
da= \chit_s (a) ~dx^s
\en
\noi which is equation (\ref{daflat}) in ``curved" indices
[Note: ref. \cite{WZ} adopts a definition of $\chit_s$ such that
$da= dx^s~(\partial_s( a))$].
\sk
Using the Lie derivarive and the contraction operator defined
on the full inhomogeneous quantum group one can also study 
the Cartan 
calculus on the quantum plane;
this should provide an alternative 
approach to \cite{CSZ}.
\sk
The results of this section are applied to the 
multiparametric quantum plane  $IGL_{qr}(2)/GL_{qr}(2)$
at the end of the Table. The usual relations 
of the uniparametric case \cite{WZ} are recovered 
after setting $q=r$.
\vskip 1cm


\sect{{Table of $IGL_{q,r}(2)$ }}


\sk
\centerline{The quantum group $IGL_{q,r}(2)$ and its differential
calculus}
\sk
\sk
\noi {\sl Parameters:} $q (\equiv q_{12}), q_{01}, q_{02}, r~~$
 
\sk
\noi {\sl $R$ and $D$-matrices of $GL_{q,r}(2)$:}
\[
\R{ab}{cd}=\left( 
\begin{array}{cccc} 
r & 0 & 0 & 0 \\
                                 0 & {r\over q} & 0 & 0 \\
                                 0 & r-\rminus & {q\over r} & 0 \\
                                 0 & 0 & 0 & r \end{array} \right)
,~~D^a_{~b}=\left( \begin{array}{cc} r & 0  \\ 0 & r^3 \end{array} 
\right)
\]
\sk
\noi {\sl $\T{A}{B}$ (A,B=0,1,2): 
fundamental representation of $IGL_{q,r}(2)$}
\sk
\[
\T{A}{B}=\left( \begin{array}{ccc}  u & 0 & 0 \\
                                      x^1  & \al & \be \\
                                      x^2  & \ga & \de 
 \end{array} \right)
\]
\sk
\noi {\sl Determinant of  $IGL_{q,r}(2)$ and definition of $\xi$}
\sk
\[
\det \T{A}{B}=u \det \T{a}{b},  ~~\mbox{where}~
\det \T{a}{b}=\al\de-{r^2\over q} \be\ga
\]
\[
\xi \det \T{A}{B}=\det \T{A}{B} \xi=I
\]
\sk
\noi {\sl Basis elements generating $IGL_{q,r}(2)$}
\sk
\[ 
\al, \be , \ga , \de , ~x^1,~x^2, ~u, ~\xi 
\]
\sk
\noi {\sl Commutations of the basis elements }
\sk
\[
\al\be={r^2\over q}\be\al,~~\al\ga=q\ga\al,~~\be\de=q\de\be,
~~\ga\de={r^2\over q}\de\ga
\]
\[
\be\ga={q^2\over r^2}\ga\be,~~\al\de-\de\al=
{r\over q}(r-\rminus)\be\ga,
\]
\[
\al x^1={q_{01}\over r^2} x^1 \al,
 ~~~\be x^1={q_{02}\over r^2} x^1 \be,
\]
\[
\al x^2=q{q_{01}\over r^2} x^2 \al, 
~~~\be x^2=q{q_{02}\over r^2} x^2 \be,
\]
\[
\ga x^1={q_{01}\over q} x^1 
\ga-{r\over q} (r-\rminus) \al x^2, ~~~
\de x^1={q_{02}\over q} x^1
 \de-{r\over q} (r-\rminus) \be x^2,
\]
\[
\ga x^2={q_{01}\over r^2} x^2 \ga,~~~ 
\de x^2={q_{02}\over r^2} x^2 \de
\]
\[
x^1 x^2=q x^2 x^1
\]
\[
\T{a}{b} u={q_{0b}\over q_{0a}} u \T{a}{b}, 
~~~x^a u= (q_{0a})^{-1} u x^a
\]
\[
(\det \T{A}{B})\T{A}{B} ={q_{0A} q_{1A} q_{2A}
 \over q_{0B} q_{1B} q_{2B}}
 \T{A}{B} (\det \T{A}{B}), ~q_{AA} \equiv r, ~q_{AB} 
\equiv {r^2\over q_{BA}}
\]
\[
\T{A}{B} \xi={q_{0A} q_{1A} q_{2A} \over q_{0B} q_{1B} q_{2B}}
\xi \T{A}{B}
\]
\sk
\noi {\sl Conditions for centrality of 
$\det \T{A}{B}=u \det \T{a}{b}$, $\det \T{a}{b}$ and $u$}
\sk
centrality of $u \det \T{a}{b}$ ~$\Longleftrightarrow$ 
$~~q_{01} q_{02}=r^2, ~q_{01}=q$ 

centrality of $\det \T{a}{b}$ ~~~$\Longleftrightarrow$
$~~q_{01} q_{02}=r^2, ~q=r$

centrality of $u$ ~~~~~~~~~~~$\Longleftrightarrow$
$~~q_{01}= q_{02}=1$
\sk\sk

\noi {\sl Inverse of $\T{A}{B}$}

\[
\Ti{A}{B}=\left( 
\begin{array}{cc}  \det \T{a}{b} \xi & -\Ti{a}{b} x^b
                                    \det \T{a}{b} \xi  \\
                                      0  & \Ti{a}{b}
                 \end{array} \right)
\]
\[
\Ti{a}{b}= \xi u \left( \begin{array}{cc} \de & -\qm \be  \\
                                -q\ga      & \al
                 \end{array} \right)
\]
\sk

\noi {\sl Commutations of the left-invariant $1$-forms $\om$}
\sk
Notations: $\om^1 \equiv \ome{1}{1},\om^+ 
\equiv \ome{1}{2},\om^-\equiv \ome{2}{1},
\om^2 \equiv \ome{2}{2},~V^1\equiv \ome{0}{1},
V^2\equiv \ome{0}{2}$

\[\om^1 \we \om^+ + \om^+ \we \om^1 = 0\]
\[\om^1 \we \om^- + \om^- \we \om^1 = 0\]
\[\om^1 \we \om^2 + \om^2 \we \om^1 = (1-r^2)\om^+ \we \om^-
\]
\[\om^+ \we \om^- + \om^- \we \om^+ = 0\]
\[\om^2 \we \om^+ + r^2 \om^+ \we \om^2 =
 r^2 (r^2 - 1)\om^+ \we \om^1
\]
\[\om^2 \we \om^- + r^{-2} \om^- \we \om^2 =
 (1-r^2)\om^- \we \om^1\]
\[\om^2 \we \om^2 =(r^2 - 1)\om^+ \we \om^-
\]
\[\om^1 \we \om^1 = \om^+ \we \om^+ = \om^- \we \om^-=0\]
\sk
\[ \om^1 \we V^1+r^2 V^1 \we \om^1=0\]
\[ \qm {q_{02} \over q_{01}} \om^+ \we V^1 + V^1 \we \om^+=
(1-r^{-2}) \om^1 \we V^2 \]
\[ \om^- \we V^1 + {r^2 \over q}{q_{02} 
\over q_{01}} V^1 \we \om^-=0 \]
\[ \om^2 \we V^1 + V^1 \we \om^2 =(1-r^{-2}) q {q_{01}\over q_{02}}
\om^- V^2 \]
\[ \om^1 \we V^2+ V^2 \we \om^1=0\]
\[ \qm {q_{02} \over q_{01}}\om^+ \we V^2+V^2 \we \om^+=0\]
\[ {q\over r^2}{q_{01} \over q_{02}}\om^- \we V^2+V^2 \we \om^-=
(1-r^2) V^1 \we  \om^1 \]
\[ \om^2 \we  V^2+r^2 \we V^2 =(r^2-1)[(1-r^2)\om^1 \we V^2+
{r^2 \over q}{q_{02} \over q_{01}} \om^+ \we V^1] \]

\noi {\sl Cartan-Maurer equations:}
\sk
\[ d\om^1+r\om^+ \we \om^-=0 \]
\[ d\om^+ + r \om^+(-r^2 \om^1 + \om^2)=0 \]
\[ d\om^- + r (-r^2 \om^1 + \om^2)\om^-=0 \]
\[ d\om^2 - r \om^+ \we \om^-=0 \]
\[ dV^1-{q\over r}{q_{01} \over q_{02}}\om^-\we V^2-\rminus
\om^1 \we V^1 =0\]
\[ dV^2-{r\over q}{q_{02} \over q_{01}}\om^+ \we V^1-
\rminus \om^2 \we V^2-(r-\rminus) V^2 \we \om^1=0\]
\sk

\noi {\sl The q-Lie algebra:}
\sk
Notations: $\chi^1 \equiv \cchi{1}{1},\chi^+ 
\equiv \cchi{1}{2},\chi^-\equiv \cchi{2}{1},
\chi^2 \equiv \cchi{2},~P^1\equiv \cchi{0}{1},
P^2\equiv \cchi{0}{2}$

\[ \chi_1 \chi_+ - \chi_+ \chi_1 -(r^4-r^2)\chi_2 \chi_+ =
 r^3 \chi_+\]
\[ \chi_1 \chi_- - \chi_- \chi_1 +
(r^2-1)\chi_2 \chi_- = -r \chi_-\]
\[ \chi_1 \chi_2 - \chi_2 \chi_1 =0\]
\[ \chi_+ \chi_- - \chi_- \chi_+ + (1-r^2) \chi_2 \chi_1-(1-r^2) 
\chi_2 \chi_2 = r (\chi_1-\chi_2)
\]
\[ \chi_+ \chi_2-r^2 \chi_2 \chi_+ = r\chi_+\]
\[ \chi_- \chi_2-r^{-2} \chi_2 \chi_- = -\rminus \chi_-\]
\[ r^2 \chi_1 \P_1-\P_1 \chi_1+(r^2-1) \P_2 \chi_-=-r\P_1 \]
\[ q{q_{01} \over q_{02}} \chi_+ \P_1-\P_1 \chi_+-r^2 (1-r^2)
   \chi_2 \P_2=r^3 \P_2 \]
\[ \chi_- \P_1-{q\over r^2}{q_{01} \over q_{02}}\P_1 \chi_-=0 \]
\[ \chi_2 \P_1-\P_1 \chi_2=0 \]
\[ \chi_1 \P_2-\P_2 \chi_1 + (r^2-1) 
{q\over r^2}{q_{01} \over q_{02}}
    \chi_+ \P_1=0 \]
\[ \chi_+ \P_2 - \qm {q_{02} \over q_{01}}\P_2 \chi_+=0 \]
\[  {r^2\over q}{q_{02} \over q_{01}}\chi_- \P_2- \P_2 \chi_-
+(1-r^2) \chi_2 \P_1=-r \P_1 \]
\[ r^2 \chi_2 \P_2-\P_2 \chi_2=-r\P_2 \]
\sk
\[ \P_1 \P_2-{q\over r^2}{q_{01} \over q_{02}}\P_2 \P_1=0 \]
\sk

\noi {\sl The exterior derivative of the basis elements}
\sk
\[
\begin{array}{l}
 d\alpha={{s-r^2}\over{r^3-r}}\al\om^1-s 
{r \over q_{12}}\beta\om^++{{s-1}\over{r-
\rminus}}\al\om^2\\
 d\beta={{-r^2+s(1-r^2+r^4)}\over{r^3-r}}
\beta\om^1-s {q_{12}\over r} \alpha\om^-
+{{s-r^2}\over{r^3-r}}\be\om^2\\
d\gamma={{s-r^2}\over{r^3-r}}\ga\om^1-s 
{r\over q_{12}}\de\om^+
+{{s-1}\over{r-\rminus}}\ga\om^2\\
d\de={{-r^2+s(1-r^2+r^4)}\over{r^3-r}}
\de\om^1-s {q_{12}\over r} \ga\om^-+
{{s-r^2}\over{r^3-r}}\de\om^2\\
dx^1=-{s r\over q_{01}} \al V^1-{s r\over q_{02}} \be V^2+
{{s-1}\over {r-\rminus}} x^1 \tau\\
dx^2=-{s r\over q_{01}} \ga V^1-{s r\over q_{02}} \de V^2+
{{s-1}\over {r-\rminus}} x^2 \tau\\
du={{s-1}\over {r-\rminus}} u \tau\\
d\xi={{r^2 s^{-N-1} -1}\over {r-\rminus}} \xi \tau,~~d\zeta=
{{r^2 s^{-N} -1}\over {r-\rminus}} \zeta \tau\\
d(\det \T{A}{B})={{r^{-2} s^{N+1} -1}
\over {r-\rminus}} (\det \T{A}{B}) \tau,
~~d(\det \T{a}{b})=
{{r^{-2 }s^{N} -1}\over {r-\rminus}}  (\det \T{A}{B})\tau
\end{array}
\]
\noi {\sl The $\om^i$ in 
terms of the exterior derivatives on $\al,\be,\ga,
\de, x^1, x^2 ,u$:}
\[
\begin{array}{l}
 \om^1={r\over{s(-r^2-r^4+s+sr^4)}} 
[(r^2-s)(\kappa (\al)da+\kappa 
(\be)d\ga)+r^2(s-1)(\kappa (\ga)d\be + \kappa (\de) d\de)]\\
 \om^+=-{1\over s} {q_{12}\over r}[\kappa (\ga) d\al + 
\kappa (\de) d\ga]\\
\om^-=-{1\over s}{r\over {q_{12}}} [\kappa (\al) d\be + 
\kappa (\be) d\de]\\
\om^2={r\over{s(-r^2-r^4+s+sr^4)}}
[(s-r^2-sr^2+sr^4)(\kappa (\al)
d\al+\kappa (\be) d\ga)+(r^2-s)(\kappa 
(\ga) d\be+\kappa (\de) d\de)]\\
V^1=-{q_{01}\over{s r}} [\kappa (\al)dx^1+\kappa
 (\be)dx^2+\kappa (x^1)du]\\
V^2=-{q_{02}\over {s r}} [\kappa (\ga)dx^1+\kappa
 (\de)dx^2+\kappa (x^2)du]
\end{array}
\]
\sk
\noi {\sl The multiparametric quantum plane 
$Fun_{q,r}\left(IGL(2)/GL(2)\right)$}
\[
\begin{array}{l}
x^1 x^2=q x^2 x^1\\
{}\\
dx^1 x^1=r^{-2} x^1 dx^1\\
dx^1~ x^2= {q\over r^2} x^2 dx^1\\
dx^2~ x^1=(r^{-2}-1) x^2 dx^1+\qm x^1 dx^2\\
dx^2 x^2=r^{-2} x^2 dx^2\\
{}\\
dx^1 \we dx^2= -{q\over r^2}dx^2 \we dx^1
\end{array}
\]
\vfill\eject



\def\spinst#1#2{{#1\brack#2}}
\def\sk{\vskip .4cm}
\def\noi{\noindent}
\def\om{\omega}
\def\Om{\Omega}
\def\al{\alpha}
\def\la{\lambda}
\def\be{\beta}
\def\ga{\gamma}
\def\Ga{\Gamma}
\def\del{\delta}
\def\linv{{1 \over \lambda}}
\def\rinv{{1\over {r-r^{-1}}}}
\def\alb{\bar{\alpha}}
\def\beb{\bar{\beta}}
\def\gab{\bar{\gamma}}
\def\deb{\bar{\delta}}
\def\ab{\bar{a}}
\def\Ab{\bar{A}}
\def\Bb{\bar{B}}
\def\Cb{\bar{C}}
\def\Db{\bar{D}}
\def\ab{\bar{a}}
\def\cb{\bar{c}}
\def\db{\bar{d}}
\def\bb{\bar{b}}
\def\eb{\bar{e}}
\def\fb{\bar{f}}
\def\gb{\bar{g}}
\def\xih{\hat\xi}
\def\Xih{\hat\Xi}
\def\uh{\hat u}
\def\vh{\hat v}
\def\ub{\bar u}
\def\vb{\bar v}
\def\xib{\bar \xi}

\def\alp{{\alpha}^{\prime}}
\def\bep{{\beta}^{\prime}}
\def\gap{{\gamma}^{\prime}}
\def\dep{{\delta}^{\prime}}
\def\rhop{{\rho}^{\prime}}
\def\taup{{\tau}^{\prime}}
\def\rhopp{\rho ''}
\def\thetap{{\theta}^{\prime}}
\def\imezzi{{i\over 2}}
\def\unquarto{{1 \over 4}}
\def\onehalf{{1 \over 2}}
\def\unmezzo{{1 \over 2}}
\def\epsi{\varepsilon}
\def\we{\wedge}
\def\th{\theta}
\def\de{\delta}
\def\cony{i_{\de {\vec y}}}
\def\Liey{l_{\de {\vec y}}}
\def\tv{{\vec t}}
\def\Gt{{\tilde G}}
\def\deyv{\vec {\de y}}
\def\part{\partial}
\def\pdxp{{\partial \over {\partial x^+}}}
\def\pdxm{{\partial \over {\partial x^-}}}
\def\pdxi{{\partial \over {\partial x^i}}}
\def\pdy#1{{\partial \over {\partial y^{#1}}}}
\def\pdx#1{{\partial \over {\partial x^{#1}}}}
\def\pdyx#1{{\partial \over {\partial (yx)^{#1}}}}

\def\qP{q-Poincar\'e~}
\def\A#1#2{ A^{#1}_{~~~#2} }

\def\R#1#2{ R^{#1}_{~~~#2} }
\def\PA#1#2{ P^{#1}_{A~~#2} }
\def\Pa#1#2{ (P_A)^{#1}_{~~#2} }
\def\Pas#1#2{ (P_A)^{#1}_{~#2} }
\def\PI#1#2{(P_I)^{#1}_{~~~#2} }
\def\PJ#1#2{ (P_J)^{#1}_{~~~#2} }
\def\Ppp{(P_+,P_+)}
\def\Ppm{(P_+,P_-)}
\def\Pmp{(P_-,P_+)}
\def\Pmm{(P_-,P_-)}
\def\Ppo{(P_+,P_0)}
\def\Pom{(P_0,P_-)}
\def\Pop{(P_0,P_+)}
\def\Pmo{(P_-,P_0)}
\def\Poo{(P_0,P_0)}
\def\Pso{(P_{\sigma},P_0)}
\def\Pos{(P_0,P_{\sigma})}

\def\Rp#1#2{ (R^+)^{#1}_{~~~#2} }
\def\Rpinv#1#2{ [(R^+)^{-1}]^{#1}_{~~~#2} }
\def\Rm#1#2{ (R^-)^{#1}_{~~~#2} }
\def\Rinv#1#2{ (R^{-1})^{#1}_{~~~#2} }
\def\Rsecondinv#1#2{ (R^{\sim 1})^{#1}_{~~~#2} }
\def\Rinvsecondinv#1#2{ ((R^{-1})^{\sim 1})^{#1}_{~~~#2} }

\def\Rpm#1#2{(R^{\pm})^{#1}_{~~~#2} }
\def\Rpminv#1#2{((R^{\pm})^{-1})^{#1}_{~~~#2} }

\def\Rb{{\bf \mbox{\boldmath $R$}}}
\def\Rbo{{\bf \mbox{\bf $R$}}}
\def\Rbp#1#2{{ (\Rbo^+)^{#1}_{~~~#2} }}
\def\Rbm#1#2{ (\Rbo^-)^{#1}_{~~~#2} }
\def\Rbinv#1#2{ (\Rbo^{-1})^{#1}_{~~~#2} }
\def\Rbpm#1#2{(\Rbo^{\pm})^{#1}_{~~~#2} }
\def\Rbpminv#1#2{((\Rbo^{\pm})^{-1})^{#1}_{~~~#2} }

\def\RRpm{R^{\pm}}
\def\RRp{R^{+}}
\def\RRm{R^{-}}

\def\Rh{{\hat R}}
\def\Rbh{{\hat {\Rbo}}}
\def\Rhat#1#2{ \Rh^{#1}_{~~~#2} }
\def\Rbar#1#2{ {\bar R}^{#1}_{~~~#2} }
\def\L#1#2{ \La^{#1}_{~~~#2} }
\def\Linv#1#2{ \La^{-1~#1}_{~~~~~#2} }
\def\Rbhat#1#2{ \Rbh^{#1}_{~~~#2} }
\def\Rhatinv#1#2{ (\Rh^{-1})^{#1}_{~~~#2} }
\def\Rbhatinv#1#2{ (\Rbh^{-1})^{#1}_{~~~#2} }
\def\Z#1#2{ Z^{#1}_{~~~#2} }
\def\X#1#2{ X^{#1}_{~~~#2} }
\def\Rt#1{ {\hat R}_{#1} }
\def\La{\Lambda}
\def\Rha{{\hat R}}
\def\ff#1#2#3{f_{#1~~~#3}^{~#2}}
\def\MM#1#2#3{M^{#1~~~#3}_{~#2}}
\def\MMc#1#2#3{{M_{\!-}}^{#1~~~#3}_{~#2}}
\def\MMcc#1#2{M_{\!-}{}_{#1}{}^{#2}}
\def\cchi#1#2{\chi^{#1}_{~#2}}
\def\chil#1{\chi_{{}_{#1}}}
\def\ome#1#2{\om_{#1}^{~#2}}
\def\Ome#1#2{\Omega_{#1}^{~#2}}
\def\RRhat#1#2#3#4#5#6#7#8{\La^{~#2~#4}_{#1~#3}|^{#5~#7}_{~#6~#8}}
\def\RRhatinv#1#2#3#4#5#6#7#8{(\La^{-1})^
{~#2~#4}_{#1~#3}|^{#5~#7}_{~#6~#8}}
\def\LL#1#2#3#4#5#6#7#8{\La^{~#2~#4}_{#1~#3}|^{#5~#7}_{~#6~#8}}
\def\LLinv#1#2#3#4#5#6#7#8{(\La^{-1})^
{~#2~#4}_{#1~#3}|^{#5~#7}_{~#6~#8}}
\def\U#1#2#3#4#5#6#7#8{U^{~#2~#4}_{#1~#3}|^{#5~#7}_{~#6~#8}}
\def\Cb{\bf \mbox{\boldmath $C$}}
\def\CC#1#2#3#4#5#6{{\Cb}_{~#2~#4}^{#1~#3}|_{#5}^{~#6}}
\def\cc#1#2#3#4#5#6{C_{~#2~#4}^{#1~#3}|_{#5}^{~#6}}
\def\PIJ#1#2#3#4#5#6#7#8{(P_I,P_J)^{~#2~#4}_{#1~#3}|^{#5~#7}_{~#6~#8}}

\def\ZZ#1#2#3#4#5#6#7#8{Z^{~#2~#4}_{#1~#3}|^{#5~#7}_{~#6~#8}}

\def\C#1#2{ {\bf \mbox{\boldmath $C$}}_{#1}^{~~~#2} }
\def\c#1#2{ C_{#1}^{~~~#2} }
\def\q#1{   {{q^{#1} - q^{-#1}} \over {q^{\unmezzo}-q^{-\unmezzo}}}}
\def\Dmat#1#2{D^{#1}_{~#2}}
\def\Dmatinv#1#2{(D^{-1})^{#1}_{~#2}}
\def\DR{\Delta_R}
\def\DL{\Delta_L}
\def\f#1#2{ f^{#1}_{~~#2} }
\def\F#1#2{ F^{#1}_{~~#2} }
\def\T#1#2{ T^{#1}_{~~#2} }
\def\Ti#1#2{ (T^{-1})^{#1}_{~~#2} }
\def\Tp#1#2{ (T^{\prime})^{#1}_{~~#2} }
\def\Th#1#2{ {\hat T}^{#1}_{~~#2} }
\def\TP{ T^{\prime} }
\def\M#1#2{ M_{#1}^{~#2} }
\def\Mc#1#2{ {M_{\!-}}_{#1}^{~#2} }
\def\qm{q^{-1}}
\def\rminus{r^{-1}}
\def\um{u^{-1}}
\def\vm{v^{-1}}
\def\xm{x^{-}}
\def\xp{x^{+}}
\def\fm{f_-}
\def\fp{f_+}
\def\fn{f_0}
\def\D{\Delta}
\def\DN{\Delta_{N+1}}
\def\kN{\kappa_{N+1}}
\def\eN{\epsi_{N+1}}
\def\Mat#1#2#3#4#5#6#7#8#9{\left( \matrix{
     #1 & #2 & #3 \cr
     #4 & #5 & #6 \cr
     #7 & #8 & #9 \cr
   }\right) }
\def\Ap{A^{\prime}}
\def\Dp{\Delta^{\prime}}
\def\Ip{I^{\prime}}
\def\ep{\epsi^{\prime}}
\def\kp{\kappa^{\prime}}
\def\kpm{\kappa^{\prime -1}}
\def\kpsq{\kappa^{\prime 2}}
\def\km{\kappa^{-1}}
\def\gp{g^{\prime}}
\def\qone{q \rightarrow 1}
\def\rone{r \rightarrow 1}
\def\qrone{q,r \rightarrow 1}
\def\Fmn{F_{\mu\nu}}
\def\Am{A_{\mu}}
\def\An{A_{\nu}}
\def\dm{\part_{\mu}}
\def\dn{\part_{\nu}}
\def\Ana{A_{\nu]}}
\def\Bna{B_{\nu]}}
\def\Zna{Z_{\nu]}}
\def\dma{\part_{[\mu}}
\def\qsu{$[SU(2) \times U(1)]_q~$}
\def\suq{$SU_q(2)~$}
\def\su{$SU(2) \times U(1)~$}
\def\gij{g_{ij}}
\def\qL{SL_q(2,{\bf \mbox{\boldmath $C$}})}
\def\GLqrN{GL_{q,r}(N)}
\def\IGLqrN{IGL_{q,r}(N)}
\def\IGLqrtwo{IGL_{q,r}(2)}
\def\GLqrNo{GL_{q,r}(N+1)}
\def\SOqrNt{SO_{q,r}(N+2)}
\def\SpqrNt{Sp_{q,r}(N+2)}
\def\SLqrN{SL_{q,r}(N)}
\def\UglqrN{U(gl_{q,r}(N))}
\def\UglqrNo{U(gl_{q,r}(N+1))}
\def\UiglqrN{U(igl_{q,r}(N))}
\def\ISOqrN{ISO_{q,r}(N)}
\def\ISpqrN{ISp_{q,r}(N)}
\def\ISOqroN{ISO_{q,r=1}(N)}
\def\ISpqroN{ISp_{q,r=1}(N)}
\def\SqrN{S_{q,r}(N)}
\def\SqrNtwo{S_{q,r}(N+2)}
\def\USqrNtwo{U(S_{q,r}(N+2))}
\def\UISqrN{U(IS_{q,r}(N))}
\def\ISqrN{IS_{q,r}(N)}
\def\SqroNt{S_{q,r=1}(N+2)}
\def\SOqroNt{SO_{q,r=1}(N+2)}
\def\SpqroNt{Sp_{q,r=1}(N+2)}
\def\ISqroN{IS_{q,r=1}(N)}
\def\ISOqroN{ISO_{q,r=1}(N)}

\def\SOqrN{SO_{q,r}(N)}
\def\SpqrN{Sp_{q,r}(N)}
\def\SqrNt{S_{q,r}(N+2)}

\def\Tc{{\cal T}}

\def\Dtwo{\Delta_{N+2}}
\def\epsitwo{\epsi_{N+2}}
\def\kappatwo{\kappa_{N+2}}

\def\RR{R^*}
\def\rr#1{R^*_{#1}}

\def\Lpm#1#2{L^{\pm #1}_{~~~#2}}
\def\Lmp#1#2{L^{\mp#1}_{~~~#2}}
\def\LLpm{L^{\pm}}
\def\LLmp{L^{\mp}}
\def\LLp{L^{+}}
\def\LLm{L^{-}}
\def\Lp#1#2{L^{+ #1}_{~~~#2}}
\def\Lm#1#2{L^{- #1}_{~~~#2}}
\def\Dcal#1#2{{\cal D}^{#1}_{~#2}}

\def\gu{g_{U(1)}}
\def\gsu{g_{SU(2)}}
\def\tg{ {\rm tg} }
\def\Fun{$Fun(G)~$}
\def\invG{{}_{{\rm inv}}\Ga}
\def\Ginv{\Ga_{{\rm inv}}}
\def\qonelim{\stackrel{q \rightarrow 1}{\longrightarrow}}
\def\ronelim{\stackrel{r \rightarrow 1}{\longrightarrow}}
\def\limrone{\lim_{r \rightarrow 1}}
\def\ronelimeq{\stackrel{r \rightarrow 1}{=}}
\def\Pprojection{\stackrel{P}{\longrightarrow}}
\def\qronelim{\stackrel{q=r \rightarrow 1}{\longrightarrow}}
\def\viel#1#2{e^{#1}_{~~{#2}}}
\def\ra{\rightarrow}
\def\detq{{\det}}
\def\detqr{{\det}}
\def\detqrm{{\det} }
\def\detqrTAB{{\det} \T{A}{B}}
\def\detqrTab{{\det} \T{a}{b}}
\def\P{P}
\def\Qt{Q}
\def\chit{{\partial}}

\def\pp#1#2{\Pi_{#1}^{(#2)}}

\def\BCD{B_n, C_n, D_n}

\def\n2{{{N+1} \over 2}}
\def\ap{a^{\prime}}
\def\bp{b^{\prime}}
\def\cp{c^{\prime}}
\def\dpr{d^{\prime}}
\def\Dc{{\cal D}}
\def\osqrt{{1 \over \sqrt{2}}}
\def\Ntwo{{N\over 2}}

\def\bu{\bullet}
\def\ci{\circ}

\def\sma#1{\mbox{\footnotesize #1}}
\def\Q.E.D.{\rightline{$\Box$}}
\def\vt{\vartheta}

\def\cvd{{\vskip -0.49cm\rightline{$\Box\!\Box\!\Box$}}\sk}
\def\*{\star}
\chapter{Geometry of the quantum Inhomogeneous Orthogonal and Symplectic Groups
$ISO_{q,r}(N)$ and $ISp_{q,r}(N)$} 

In this chapter we study the inhomogeneous orthogonal and symplectic groups,
their universal enveloping algebras and their differential calculi. 
We tush give a detailed  analysis of the  geometry of these inhomogeneous 
groups that are canonically 
associated (via a quotient procedure) to the orthogonal and 
symplectic quantum groups studied in
\cite{FRT}.  

The method used in the previous chapter to obtain $\IGLqrN$,
is here applied to obtain $ISO_{q,r}(N)$
and $ISp_{q,r}(N)$ and to give an $R$-matrix formulation of these $q$-groups. 
This method is based on a 
projection (consistent with respect to the
Hopf structure) from the corresponding quantum groups
of higher rank $A_{n+1}, B_{n+1}, C_{n+1}, D_{n+1}$. 

In general the quantum inhomogeneous groups
we analize do contain dilatations. There exists however a  
subclass of dilatation-free cases 
for special values of the deformation parameters.
The important example of the $q$-Poincar\'e group is
contained in our construction. In particular, we
find a dilatation-free $q$-Poincar\'e group
depending on one real parameter $q$.
\sk
We next present a detailed study of the universal enveloping algebra
of the multiparametric homogeneous orthogonal and symplectic groups and 
find a suitable set of generators that can be ordered. This will clarify
the structure of their  inhomogeneous version. 
The projection procedure used to derive the $ISO_{q,r}(N)$ [$ISp_{q,r}(N)$]
$q$-groups is then used to find their universal enveloping algebras as Hopf
subalgebras of $U_{q,r}(so(N+2))$ [$U_{q,r}(sp(N+2))$]. 
An $R$-matrix formulation and the duality $ISO_{q,r}(N)\leftrightarrow
U_{q,r}(iso(N))$ [$ISp_{q,r}(N)\leftrightarrow
U_{q,r}(isp(N))$] are explicitly given. 
\sk
The quantum Lie algebras of $ISO_{q,r}(N)$ [$ISP_{q,r}(N)$] 
are subspaces (adjoint submodules)  
of
$U_{q,r}(iso(N))$ [$U_{q,r}(isp(N))$], and in the second part of the chapter 
we study these deformed Lie algebras and their associated differential calculi.
This is again done using the projection or quotient structure of 
$ISO_{q,r}(N)$ and $ISp_{q,r}(N)$.
Contrary to the $IGL_{q,r}(N)$ case, only for $r=1$ we have a quantum 
differential calculus that is a continuous deformation of the commutative one.
The necessity of taking $r=1$ is discussed.
In Section 4.5 we briefly 
introduce the multiparametric bicovariant calculus on the 
homogeneous orthogonal and symplectyc $q$-groups.
In Section 4.6 we examine
the case $r=1$. We clarify
some issues related to the classical limit and
see how in this limit some tangent vectors become
linearly dependent, thus providing the correct classical dimension
of the tangent space. A similar mechanism occurs for the
left-invariant
$1$-forms. 
In  Section 4.7 the bicovariant calculi on $ISO_{q,r=1}(N)$
[$ISp_{q,r=1}(N)$]
are studied. We first  consider  a calculus that has one more 
generator than in the classical case, this generator 
correspond to the dilatation $u$ of the quantum inhomogeneous 
group. Then we show how to restrict this calculus to 
one that has the same numbler of tangent vectors that appear in the 
classical case.
All the quantities relevant to this differential calculus are
explicitly constructed.  The results are then directly
applied to the $q$-Poincar\'e group $ISO_q(3,1)$.

\sect{$B_n,C_n,D_n$ multiparametric quantum groups}

The $B_n, C_n, D_n$ multiparametric
quantum groups are freely generated by the noncommuting
matrix elements $\T{a}{b}$ (fundamental representation) and
the identity $I$, modulo the quadratic $RTT$ and $CTT$ relations
discussed below.  The noncommutativity is controlled by the $R$
matrix:
\eq
\R{ab}{ef} \T{e}{c} \T{f}{d} = \T{b}{f} \T{a}{e} \R{ef}{cd}
\label{RTTsosp}
\en
which satisfies the quantum Yang-Baxter equation (\ref{QYB})
\eq
\R{a_1b_1}{a_2b_2} \R{a_2c_1}{a_3c_2} \R{b_2c_2}{b_3c_3}=
\R{b_1c_1}{b_2c_2} \R{a_1c_2}{a_2c_3} \R{a_2b_2}{a_3b_3}, \label{QYBE}
\en
a sufficient condition for the consistency of the
``$RTT$" relations (\ref{RTTsosp}).  The $R$-matrix components
$\R{ab}{cd}$
depend
continuously on a (in general complex)
set of  parameters $q_{ab},r$.  For $q_{ab}=r$ we
recover the uniparametric $q$-groups of ref. \cite{FRT}. Then
$q_{ab} \rightarrow 1, r \rightarrow 1$ is the classical limit for
which
$\R{ab}{cd} \rightarrow \de^a_c \de^b_d$ : the
matrix entries $\T{a}{b}$ commute and
become the usual entries of the fundamental representation. The
multiparametric $R$ matrices for the $A,B,C,D$ series can be
found in \cite{Schirrmacher}  (other ref.s on multiparametric
$q$-groups are given in \cite{multiparam1,multiparam2}). For
the $B,C,D$ case they read:
\eq
\begin{array}{ll}
\R{ab}{cd}=&\delta^a_c \delta^b_d [{r\over q_{ab}} +
(r-1) \delta^{ab}+(\rminus - 1)
\de^{a\bp}] (1-\de^{a n_2})
+\de^a_{n_2} \de^b_{n_2} \de^{n_2}_c \de^{n_2}_d  \\
&+(r-r^{-1})
[\theta^{ab} \delta^b_c \de^a_d - \epsilon_a\epsilon_c
\theta^{ac} r^{\rho_a - \rho_c}
\de^{\ap b} \de_{\cp d}]
\end{array}
\label{Rmpsosp}
\en
\noi where $\theta^{ab}=1$ for $a> b$
and $\theta^{ab}=0$ for $ a \le b$; we define
 $n_2 \equiv \n2$ and primed indices as $\ap \equiv N+1-a$. The
indices
run on $N$ values ($N$=dimension of
the fundamental representation $\T{a}{b}$),
with $N=2n+1$ for $B_n  [SO(2n+1)]$,  $N=2n$
for $C_n [Sp(2n)]$, $D_n [SO(2n)]$.
The terms with the index $n_2$ are present
only for the $B_n$ series. The $\epsilon_a$ and
$\rho_a$ vectors are given by:
\eq
\epsilon_a=
\left\{ \begin{array}{ll} +1 & \mbox{for $B_n$, $D_n$} ,\\
                            +1 & \mbox{for $C_n$ and $a < n$},\\
                              -1  & \mbox{for $C_n$ and $a > n$}.
           \end{array}
                                    \right.
\en
\eq
(\rho_1,...\rho_N)=\left\{ \begin{array}{ll}
         (\Ntwo -1, \Ntwo -2,...,
{1\over 2},0,-{1\over 2},...,-\Ntwo+1)
                   & \mbox{for $B_n$} \\
           (\Ntwo,\Ntwo -1,...1,-1,...,-\Ntwo) & \mbox{for $C_n$} \\
           (\Ntwo -1,\Ntwo -2,...,1,0,0,-1,...,-\Ntwo+1) & \mbox{for
$D_n$}
                                             \end{array}
                                    \right.
\en
Moreover the following relations reduce the number of independent
$q_{ab}$ parameters \cite{multiparam1}, \cite{Schirrmacher}:
\eq
q_{aa}=r,~~q_{ba}={r^2 \over q_{ab}}; \label{qab1}
\en
\eq
q_{ab}={r^2 \over q_{a\bp}}={r^2 \over q_{\ap b}}=q_{\ap\bp}
 \label{qab2}
\en
\noi where (\ref{qab2}) also implies $q_{a\ap}=r$. Therefore
 the $q_{ab}$ with $a < b < {N\over 2}$ give all the $q$'s.
\sk
It is useful to list the nonzero complex components
of the $R$ matrix (no sum on repeated indices):
\eqa
& &\R{aa}{aa}=r , ~~~~~~~~~~~~~~\mbox{\footnotesize
 $a \not= {n_2}$ } \nonumber\\
& &\R{a\ap}{a\ap}=r^{-1} ,  ~~~~~~~~~~\mbox{\footnotesize
 $a \not= {n_2}$ } \nonumber\\
& &\R{{n_2}{n_2}}{{n_2}{n_2}}= 1\nonumber\\
& &\R{ab}{ab}={r \over q_{ab}} ,~~~~~~~~~~~~\mbox{\footnotesize
 $a \not= b$,  $\ap \not= b$}\label{Rnonzerososp}\\
& &\R{ab}{ba}=r-r^{-1} , ~~~~~~~\mbox{\footnotesize
$a>b, \ap \not= b $}\nonumber\\
& &\R{a\ap}{\ap a}=(r-r^{-1})(1-\epsilon r^{\rho_a-\rho_{\ap}})
=(r-r^{-1})[1-C^{a'a}C_{a'a}] ,
{}~~~\mbox{\footnotesize
$a>\ap $}\nonumber\\
& &\R{a\ap}{b \bp}=-(r-r^{-1})\epsilon_a\epsilon_b
r^{\rho_a-\rho_b}=-(r-r^{-1})C^{a'a}C_{bb'} , ~~~~~~~
\mbox{\footnotesize $~~a>b ,~ \ap \not= b $}\nonumber
\ena
where $\epsilon=\epsilon_a \epsilon_{\ap}$, i.e.
 $\epsilon=1$ for $B_n$, $D_n$ and $\epsilon=-1$
for $C_n$.
\sk
\noi{\bf Note}  4.1.1 $~$  The matrix $R$ is lower triangular, that is
$\R{ab}{cd}=0$ if  [$\sma{$a=c$}$ and $\sma{$b<d$}$] or
$\sma{$a<c$}$,
and has the following properties:
\eq
R^{-1}_{q,r}=R_{q^{-1},r^{-1}}~~;~~~
(R_{\!q,r})^{ab}{}_{cd}=(R_{\!q,r})^{\cp\dpr}{}_{\ap\bp}~~;~~~
(R_{q,r})^{ab}{}_{cd}=(R_{\!p,r})^{dc}{}_{ba}
\label{Rprop1}
\en
where
$q,r$ denote  the set of parameters $q_{ab},r$, and $p_{ab}\equiv
q_{ba}$.
\sk
The inverse $R^{-1}$ is defined by
$\Rinv{ab}{cd} \R{cd}{ef}=\de^a_e \de^b_f=\R{ab}{cd}
\Rinv{cd}{ef}$.
Eq. (\ref{Rprop1}) implies
that for $|q|=|r|=1$, ${\bar R}=R^{-1}$.
\sk
\noi{\bf Note} 4.1.2 $~$ Let $R_r$ be the uniparametric $R$ matrix
for the $B, C, D$ q-groups. The multiparametric $R_{q,r}$
matrix is obtained from $R_r$ via the transformation
\cite{multiparam1,Schirrmacher}
\eq
R_{q,r}=F^{-1}R_rF^{-1}
\en
where $(F^{-1})^{ab}_{~~cd}$ is a diagonal matrix
in the index couples $ab$, $cd$:
\eq
F^{-1}\equiv diag (\sqrt{{r \over q_{11}}} ,
\sqrt{{r \over q_{12}}} , ... ~ \sqrt{{r \over q_{NN}}})
\label{effe}
\en
\noi and $ab$, $cd$ are ordered as
in the $R$ matrix.
Since $\sqrt{{r \over q_{ab}}} =(\sqrt{{r\over q_{ba} }})^{-1}$
and $q_{a\ap}=q_{b\bp}$, the non diagonal
elements of $R_{q,r}$ coincide with those of $R_r$.
The matrix $F$ satisfies $F_{12}F_{21}=1$ i.e.
$F^{ab}{}_{ef}F^{fe}{}_{dc}=\delta^a_c\delta^b_d $,
the quantum Yang-Baxter equation $F_{12}F_{13}F_{23}
=F_{23}F_{13}F_{12}$ and the relations
$(R_r)_{12}F_{13}F_{23}=F_{23}F_{13}(R_r)_{12}$.
Note that for $r=1$ the multiparametric $R$ matrix reduces to
$R=F^{-2}$. 
\sk
\noi{\bf Note} 4.1.3 $~$ Let $\Rh$ the matrix defined by
$\Rhat{ab}{cd} \equiv \R{ba}{cd}$.  Then the multiparametric
$\Rh_{q,r}$ is
obtained from $\Rh_r$ via the similarity transformation
\eq
\Rh_{q,r}=F\Rh_rF^{-1}
\en
The characteristic equation and the projector decomposition
of $\Rh_{q,r}$ are therefore the same as in the uniparametric case:
\eq
(\Rh-rI)(\Rh+r^{-1}I)(\Rh-\epsilon r^{\epsilon-N} I)=0 \label{cubic}
\en
\eq
\Rh-\Rh^{-1}=(r-r^{-1}) (I-K) \label{extrarelation}
\en
\eq
\Rh=r P_S - r^{-1} P_A+\epsilon r^{\epsilon-N}P_0  \label{RprojBCD}
\en
with 
\eq
\begin{array}{ll}
&P_S={1 \over {r+\rminus}} [\Rh+\rminus I-(\rminus+\epsilon
r^{\epsilon-N})P_0]\\
&P_A={1 \over {r+\rminus}} [-\Rh+rI-(r-\epsilon r^{\epsilon-N})P_0]\\
&P_0= Q_N(r) K\\
&Q_N(r) \equiv (C_{ab} C^{ab})^{-1}={{1-r^{-2}} \over
{(1-\epsilon r^{-N-1+\epsilon})(1+\epsilon r^{N-1-\epsilon})}}~,~~~~
K^{ab}_{~~cd}=C^{ab} C_{cd}\\
&I=P_S+P_A+P_0
\end{array}
\label{projBCD}
\en
To prove (\ref{extrarelation}) in the multiparametric case note that
$F_{12} K_{12} F_{12}^{-1}=K_{12}$.
Orthogonality  (and symplecticity) conditions can be
imposed on the elements $\T{a}{b}$, consistently
with  the $RTT$ relations (\ref{RTTsosp}):
\eq
C^{bc} \T{a}{b}  \T{d}{c}= C^{ad} I~,~~~
C_{ac} \T{a}{b}  \T{c}{d}=C_{bd} I \label{Torthogonality}
\en
\noi where the (antidiagonal) metric is :
\eq
C_{ab}=\epsilon_a r^{-\rho_a} \de_{a\bp} \label{metric}
\en
\noi and its inverse $C^{ab}$
satisfies $C^{ab} C_{bc}=\de^a_c=C_{cb} C^{ba}$.
We see
that for the orthogonal series, the matrix elements of the metric
and the inverse metric coincide,
while for the symplectic series there is a change of sign:
$C^{ab}=\epsilon C_{ab}$. Notice also the symmetries
$C_{ab}=C_{\bp\ap}$ and $C_{ba}(r)=\epsilon C_{ab}(r^{-1})$.

The consistency of (\ref{Torthogonality}) with the $RTT$ 
relations
is due to the identities:
\eq
C_{ab} \Rhat{bc}{de} = \Rhatinv{cf}{ad} C_{fe} ~,\label{crc1}
\en
\eq
\Rhat{bc}{de} C^{ea}=C^{bf} \Rhatinv{ca}{fd}  ~.\label{crc2}
\en
\noi These identities
 hold also for $\Rh \rightarrow \Rh^{-1}$ and can be proved using
the explicit expression (\ref{Rnonzerososp}) of $R$.
\sk

We also note the useful relations, easily deduced from 
(\ref{RprojBCD}):
\eq
C_{ab}\Rhat{ab}{cd}=\epsilon r^{\epsilon-N}C_{cd} ,~~~
C^{cd}\Rhat{ab}{cd}=\epsilon r^{\epsilon-N}C^{ab}  \label{CR}
\en
and, from (\ref{Rnonzerososp}),
\eq
 \R{ab}{cc'}=-(r-r^{-1})C^{ba}C_{cc'}\;,~
\R{aa'}{cd}=-(r-r^{-1})C^{a'a}C_{cd} ~~\mbox{ for } ~\sma{$a>c\;,~a\not= c'$}~.
\label{foraa'}
\en
Notice also that $\kappa^2(\T{a}{b})=D^a_e\T{e}{f}{D^{-1}}^f_b$ where 
$D^a_e=C^{as}C_{es}$ and its inverse ${D^{-1}}^f_b=C^{sf}C_{sb}$ are diagonal.
\sk
The co-structures of the $B_n,C_n,D_n$ multiparametric quantum
groups have the same form as in the uniparametric case:
the coproduct
$\D$, the counit $\epsi$ and the coinverse $\kappa$ are given by
\eqa
& & \D(\T{a}{b})=\T{a}{b} \otimes \T{b}{c}  \label{cos1sosp} \\
& & \epsi (\T{a}{b})=\delta^a_b\\
& & \kappa(\T{a}{b})=C^{ac} \T{d}{c} C_{db}
\label{cos2sosp}
\ena
\sk
\noi{\bf Note} 4.1.4 $~$
Using formula (\ref{Rmpsosp}) or (\ref{Rnonzerososp}),
we find that  the $\R{AB}{CD}$  matrix for the
$SO_{q,r}(N+2)$ and $Sp_{q,r}(N+2)$ quantum groups
can be decomposed
in terms of  $SO_{q,r}(N)$ and $Sp_{q,r}(N)$ quantities
as follows (splitting the index {\small A} as
{\small A}=$(\circ, a, \bullet)$, with $a=1,...N$):
\eq
\R{AB}{CD}=\left(  \begin{array}{cccccccccc}
   {}&\circ\circ&\circ\bullet&\bullet
          \circ&\bullet\bullet&\circ d&\bullet d
      &c \circ&c\bullet&cd\\
   \circ\circ&r&0&0&0&0&0&0&0&0\\
   \circ\bullet&0&r^{-1}&0&0&0&0&0&0&0\\
   \bullet\circ&0&f(r)&r^{-1}&0&0&0&0&0&-\epsilon C_{cd} \lambda
r^{-\rho}\\
\bullet\bullet&0&0&0&r&0&0&0&0&0\\
\circ b&0&0&0&0&{r\over q_{\circ b}} \de^b_d&0&0&0&0\\
\bullet b&0&0&0&0&0&{r\over q_{\bullet b} }
\de^b_d&0&\lambda\de^b_c&0\\
a\circ&0&0&0&0&\lambda\de^a_d&0&{r \over q_{a \circ} } \de^a_c&0&0\\
a\bullet&0&0&0&0&0&0&0&{r\over q_{a \bullet}} \de^a_c&0\\
ab&0&-C^{ba} \lambda r^{-\rho}
&0&0&0&0&0&0&\R{ab}{cd}\\
\end{array} \right) \label{Rbig}
\en
\noi where $\R{ab}{cd}$ is the $R$ matrix for  $SO_{q,r}(N)$
or $Sp_{q,r}(N)$, $C_{ab}$ is the corresponding
metric,  $\lambda \equiv r-r^{-1}$,
$\rho={{N+1-\epsilon}\over 2}~(\epsilon r^{\rho}=C_{\bullet \circ})$
and $f(r) \equiv \lambda (1-\epsilon r^{-2\rho})$.
The sign $\epsilon$ has been defined after eq. s
(\ref{Rnonzerososp}).
\sk
\subsection{Real forms: $SO_{q,r}(N,\bf R)$, $SO_{q,r}(N-1,1)$,
$SO_{q,r}(n,n)$,\\ $SO_{q,r}(n+1,n-1)$, $SO_{q,r}(n+1,n)$}

Following \cite{FRT},
a conjugation ---i.e. an algebra antiautomorphism, coalgebra automorphism
and involution, satisfying $\kappa(\kappa(T^*)^*)=T$---
can be defined
\sk
$\bullet$~~trivially as $T^*=T$. Compatibility with the
$RTT$ relations (\ref{RTT}) requires 
${\bar R}_{q,r}=R^{-1}_{q,r}=
R_{q^{-1},r^{-1}}$,
i.e. $|q|=|r|=1$. Then the $CTT$ relations are invariant under
$*$-conjugation.  The corresponding real forms are
$SO_{q,r}(n,n;\Rbo)$, $SO_{q,r}(n+1,n;\Rbo)$
(for N even and odd respectively) and $Sp_{q,r}(2n;\Rbo)$. 
A conjugation on the quantum orthogonal (symplectic) plane 
(defined respectively by 
$\PA{ab}{cd} x^c x^d=0$ and $x^ax^b=r^{-1}R^{ab}_{~cd}x^cx^d$) 
that is compatible with the natural
coaction $\delta$  of the $q$-group on the $q$-plane: $x^a\rightarrow
T^a{}_b\otimes x^b$ is given by $(x^a)^*=x^a$, indeed we
have $\delta(x^*)=T^*\otimes x^*\equiv\delta^*(x)$.

\sk
$\bullet$~~ via the metric
as
$T^{\*}=(\kappa(T))^t$ i.e. $T^{\*}=C^tTC^t$.
The condition on $R$ is
${\bar \R{ab}{cd}}=\R{dc}{ba}$, which
happens for $q_{ab} {\bar q}_{ab}=r^2, r \in $ {\bf R}.
Again the $CTT$ relations are ${\*}$-invariant.
The metric on a ``real" basis has compact signature
$(+,+,...+)$ so that the real form is  $SO_{q,r}(N;\Rbo)$.
The conjugation on the quantum orthogonal plane compatible
with the coaction  $x^a\rightarrow
T^a{}_b\otimes x^b$ is: $(x^a)^{\*}=C_{ba}x^b$, indeed 
$\delta(x^{\*})=\delta^{\*}(x)$.
\sk

We now introduce
two other conjugations that give the real forms 
$SO_{q,r}(N-1,1)$, $SO_{q,r}(n+1,n-1)$, $SO_{q,r}(n+1,n)$. 
The real form $SO_{q,r}(n+1,n-1)$ has been found in \cite{Firenze1};
the $*$-conjugation on the $T$ matrices that defines
the real form  $SO_{q,r}(2n-1,1)$ appears here for the first time.
Our study of real forms is based on the $RTT$ relations and is complementary 
to the ${\cal U}_q(\mbox{\sl g})$ $*$-structure classification of Twietmeyer 
\cite{Twietmeyer}, however, there, {\sl g} is an arbitrary 
semisimple Lie algebra. 
In \cite{prep} we explicitly construct all $*$-structures of orthogonal 
quantum groups using the $RTT$ formalism, the result agrees 
with \cite{Twietmeyer}. 
Other approaches to real forms can be found in \cite{realforms}.

In order to describe the  $SO_{q,r}(N-1,1)$, 
$SO_{q,r}(n+1,n-1)$, $SO_{q,r}(n+1,n)$ conjugations
we first notice that if we have an involution  ${\sharp}$ that is a 
Hopf algebra automorphism (algebra and 
coalgebra morphism compatible with the antipode: 
$\kappa(a^{\sharp})=[\kappa(a)]^{\sharp}\,{}$) and that commutes with a 
conjugation $*$, then the
composition of these two involutions: 
${}^{*^{\sharp}}\equiv {\sharp}{\scriptstyle{{}^{{}_{\circ}}}}*= 
*{\scriptstyle{{}^{{}_{\circ}}}}{\sharp}$ is 
again a conjugation. We now find an involution $\sharp$ that comutes with
$*$ and $\*$ as defined above.
\sk
Define the map ${\sharp}$ on the generators as:
\eq
T^{\sharp}=\Dc T{\Dc}^{-1}~~~\mbox{i.e.}~~~
(T^a{}_b)^{\sharp}=\Dc^a{}_eT^e{}_f{\Dc^f{}_b}^{-1}
\en
and extend it by linearity and multiplicativity to all $SO_{q,r}(N)$.
The entries of the $N$-dimensional $\Dc$  matrix are
\eq
\begin{array}{cc}
\Dc=
\left(  \begin{array}{cccccccc}
{1}&{}&{}&{}&{}&{}&{}&{}\\
{}&{...}&{}&{}&{}&{}&{}&{}\\
{}&{}&{1}&{}&{}&{}&{}&{}\\
{}&{}&{}&{0}&{1}&{}&{}&{}\\
{}&{}&{}&{1}&{0}&{}&{}&{}\\
{}&{}&{}&{}&{}&{1}&{}&{}\\
{}&{}&{}&{}&{}&{}&{...}&{}\\
{}&{}&{}&{}&{}&{}&{}&{1}
\end{array}\right)~,~~
&
\Dc=
\left(  \begin{array}{ccccccc}
{1}&{}&{}&{}&{}&{}&{}\\
{}&{...}&{}&{}&{}&{}&{}\\
{}&{}&{1}&{}&{}&{}&{}\\
{}&{}&{}&{-1}&{}&{}&{}\\
{}&{}&{}&{}&{1}&{}&{}\\
{}&{}&{}&{}&{}&{...}&{}\\
{}&{}&{}&{}&{}&{}&{1}
\end{array}
\right)\\
\mbox{\small for $N$ even} &
\mbox{\small for $N$ odd}
\end{array}
\en 
In the $N=2n$ case the $\Dc$ matrix
exchanges the index $n$ with the 
index $n+1$, in the $N=2n+1$ case 
$\Dc$ change the sign of the entries  of 
the $T$ matrix as 
many times as the index $n_2=(N+1)/2$ appears.
Since $\Dc^2={\bf 1}$ we immediately see that $\sharp$ is an involution.

We now prove that $\sharp$  is 
compatible with the algebra structure, 
i.e. it is compatible with the $RTT$ and $CTT$ relations;
in the $N=2n$ case this is true  if
$q_{ab}=r$ when at least one of the indices $a,b$ is equal
to $n$ or $n+1$.

Use relation
(\ref{Rnonzerososp}) and, if $N=2n$,  the above restriction on the $q_{ab}$ 
parameters to prove that
\eq
\Dc_1\Dc_2R\,\Dc_1\Dc_2=R ~.
\en
The compatibility with the $RTT$  relations is then esily 
seen to hold. 
[Hint: multiply $(R_{12}T_1 T_2)^{{\sharp}}=(T_2 T_1 R_{12})^{\sharp}$ 
by $\Dc_1\Dc_2$ from the left and 
from the right  and use $\Dc^2={\bf 1}$ to prove the equivalence with 
$R_{12}T_1T_2=T_2T_1R_{12}$]. 
Similarly the compatibility of ${\sharp}$ with the orthogonality relations 
(\ref{Torthogonality}), that we rewite in matrix notation as:
\eq
TCT^t=C\,I ~~,~~~~T^tCT=C\,I\label{Torthogonalitymat}~,
\en
is due to 
$\Dc^2={\bf 1},~\Dc^t=\Dc$ and the commutativity of the $C$ 
matrix with the $\Dc$ matrix:
\eq
\Dc C\Dc=C\label{DcCDc}~.
\en
For example we have: $(TC T^t)^{\sharp}=T^{\sharp}C(T^t)^{\sharp}=\Dc T\Dc 
C \Dc T^t \Dc=
\Dc (T C T^t) \Dc$ and using $\Dc^2={\bf 1}$ and again (\ref{DcCDc}) we
conclude that  $(T C T^t)^{\sharp}=C\,I$ is equivalent to 
$TCT^t=C\,I$. 
\sk
Next we prove that $\sharp$ is compatible with the  coalgebra
structure.
Compatibility with the coproduct is trivial,
compatibility with the antipode is easily verified:
\eq\kappa(T^{\sharp})=\kappa(\Dc T \Dc)=\Dc \kappa(T) \Dc=\Dc 
CT^tC\Dc=C\Dc T^t\Dc C=[\kappa(T)]^{\sharp}~.
\en
We now show that the two conjugations defined at the beginning of 
this subsection commutes with ${\sharp}$. For the second conjugation, 
defined by
$T^{\*}=[\kappa(T)]^t=C^tTC^t$,  we have,
since $\overline{\Dc}=\Dc$:
\eqa
(T^{\sharp})^{\*}&=&(\Dc T \Dc)^{\*}=\Dc T^{\*} \Dc
=\Dc [\kappa(T)]^t \Dc \nonumber\\
&=&\Dc C^t T C^t \Dc = C^t \Dc T \Dc C^t= 
(C^t  T  C^t)^{\sharp}= ([\kappa(T)]^t)^{\sharp}\\
&=&(T^{\*})^{\sharp}\nonumber
\ena
The two maps $\*$ and $\sharp$
not only commute when applied to the $T^a{}_b$ 
matrix entries, they also commute when applied 
to  any element of the $q$-group because they are respectively
multiplicative  and antimultiplicative. 
The proof that 
${\sharp}{\scriptstyle{{}^{{}_{\circ}}}}*= 
*{\scriptstyle{{}^{{}_{\circ}}}}{\sharp}$
for the first conjugation,  defined by $T^*=T$, is 
straighforward.

We restate the above results as a theorem:
\sk
\noi{\bf Theorem} 4.1.1 $~$ The map ${\sharp}$ is an automorphism and an 
involution
of the quantum group $SO_{q,r}(N)$; in the $N=2n$ case this holds with the
 restriction
$q_{ab}=r$ when at least one of the indices $a,b$ is equal
to $n$ or $n+1$.
The compositions ${}^{*^{\sharp}}\equiv 
{\sharp}{\scriptstyle{{}^{{}_{\circ}}}}*$ and 
${}^{{\*}^{\sharp}}\equiv {\sharp}{\scriptstyle{{}^{{}_{\circ}}}}
{\*}$
of the conjugations $*$ and $\*$ with the automorphism ${\sharp}$ is again
a conjugation. The restrictions on the parameters $r, q_{ab}$ are obtained
adding to the constraint imposed by $*$ ($\*$ respectively) the constaints 
imposed by ${\sharp}$.
\cvd

Associated to  ${}^{*^{\sharp}}$ and  ${}^{{\*}^{\sharp}}$ we
have the
conjugations 
that act on the quantum orthogonal plane 
and are compatible with the
coaction  $x^a\rightarrow
T^a{}_b\otimes x^b$. These conjugations respectively are: 
$(x^a)^{*^{\sharp}}=(x^a)^{\sharp}$ and  
$(x^a)^{{\*}^{\sharp}}=C_{ba}\Dc^b{}_e x^e$.

We  now study the real forms related to the 
  ${}^{*^{\sharp}}$ and  ${}^{{\*}^{\sharp}}$ conjugations.
\sk
$\bullet$~~ The conjugation ${}^{*^{\sharp}}$ 
for the $N$-dimensional orthogonal quantum 
groups with $N$ odd
gives the real form $SO_{q,r}(n,n+1)$.
\sk
$\bullet$~~ The conjugation ${}^{*^{\sharp}}$ for the $N$-dimensional 
orthogonal quantum 
goups with $N$ even has been studied (in the uniparametric case) in 
\cite{Firenze1}, it gives the 
real form $SO_{q,r}(n+1,n-1;\Rbo)$ and 
in particular the quantum Lorentz group $SO_{r}(3,1)$ 
with $|r|=1$. For an explicit proof see
formula (\ref{Cprimo=}).

If we require ${}^{*^{\sharp}}$ to be a conjugaton but do not require 
${\sharp}$ to be an automorphism and $*$ to be a conjugation, 
we can partially  
relax the constraints on the $r, q$ parameters.
Compatibility of ${}^{*^{\sharp}}$ with the $RTT$ relations 
is indeed easily seen to require
\eq
({\bar R})_{n \leftrightarrow n+1}=R^{-1},~~~\mbox{ i.e. }~~{\cal
D}_1 {\cal
D}_2 R_{12}
{\cal D}_1 {\cal D}_2 = {\overline R_{12}^{-1} } \label{Rprop2}
\en
\noi which implies

i) $|q_{ab}|=|r|=1$
for $a$ and $b$ both different from $n$ or $n+1$;

ii) $q_{ab}/r \in {\bf R}$
when at least one of the indices $a,b$ is equal
to $n$ or $n+1$.

\noi 
In the sequel we will denote simply by $*$ this conjugation. As
discussed  in ref.s \cite{Cas2,inson} and later in
this chapter  we will need this conjugation
to obtain the inhomogeneous Lorentz group $ISO_{q,r}(3,1;\Rbo)$. 
\sk
$\bullet$~~ The conjugation ${}^{\*^{\sharp}}$ for $N=2n+1$ gives the real 
form
$SO(2n,1)$ and has been introduced in \cite{FRT}. 
\sk
$\bullet$~~ The conjugation ${}^{\*^{\sharp}}$ for $N=2n$
as far as we know is not known 
in the literature, it gives the real form $SO_{q,r}(2n-1,1)$. 
In particular we obtain another quantum Lorentz group $SO_{r}(3,1)$,
notice that here $r\in {\bf R}$.

\sect{The quantum inhomogeneous groups
$\ISOqrN$ and $\ISpqrN$}
Following the  projection procedure  
described in Section 3.3  we here introduce the 
quantum inhomogeneous groups $\ISOqrN$ and $\ISpqrN$
and give an $R$-matrix formulation.
The $\ISOqrN$ quantum group has been independently studied 
--without an $R$-matrix
formulation--  
in the first reference of \cite{inhom}. The structure of $\ISpqrN$
cannot be derived  from $\SpqrN$ and the symplectic $q$-plane relations as
is the case  for $\ISOqrN$, however it can be easily defined via the 
projection procedure.
$\ISpqrN$ and its $R$-matrix formualtion where first introduced in
\cite{inson}.
\sk
\noi We define
 $\ISOqrN$ and $\ISpqrN$
as the quotients:
\eq
 \ISOqrN \equiv\frac{\SOqrNt}{H}~, ~~\ISpqrN\equiv\frac{\SpqrNt}{H}
\label{quotientsosp}
\en
\noi where $H$ is the Hopf ideal in
$\SOqrNt$ or $\SpqrNt$ of
all sums of monomials containing at least an element of the kind
$\T{a}{\circ}, \T{\bullet}{b}, \T{\bullet}{\circ}$.
The Hopf structure of the groups in the numerators of
(\ref{quotientsosp})
is naturally inherited by the quotient groups.

We introduce the following convenient notations:
$\Tc$ stands for
$\T{a}{\circ}$,  $\T{\bullet}{b}$ or  $\T{\bullet}{\circ}$,
$\SqrNt$ stands for either $\SOqrNt$ or
$\SpqrNt$,
and we indicate by
$\Dtwo$, $\epsitwo$ and $\kappatwo$ the corresponding
co-structures.

We denote by $P$ the canonical projection
\eq
P ~:~~  \SqrNt\longrightarrow \SqrNt/{H}
\en
\noi It is a Hopf algebra epimorphism because $H=Ker(P)$ is a Hopf
ideal. 
[The proof is as in Theorem 3.3.1, just use $\Tc$ instead of $\T{0}{b}$. 
In order to show that $\kappatwo (H)\subseteq H$,
notice that  $\kappatwo ({\cal{T}}) \propto {\cal{T}}$ and therefore,
$\forall h\in H$, 
$
\kappatwo (h)=\kappatwo(b\Tc
c)=\kappatwo(c)\kappatwo(\Tc)\kappatwo(b)
\in H$]. Then
any element of ${\SqrNt/H}$ is of the form $P(a)$ and
the Hopf algebra structure is given by:
\eq
P(a)+P(b)\equiv P(a+b) ~;~~ P(a)P(b)\equiv P(ab) ~;~~
\mu P(a)\equiv P(\mu a),~~~\mu \in \mbox{\bf C} \label{isoalgebra}
\en
\eq
\D (P(a))\equiv (P\otimes P)\Dtwo(a) ~;~~ \epsi(P(a))
\equiv\epsitwo(a) ~;~~
\kappa(P(a))\equiv P(\kappatwo(a)) \label{co-iso}~.
\en
\sk
We can also give an $R$--matrix formulation of the
 inhomogeneous $\ISOqrN$ and $\ISpqrN$
q-groups.
Indeed recall that $S_{q,r}(N+2)$ is the Hopf algebra freely
generated by the
non-commuting matrix elements
$T^A{}_B$ modulo the ideal  generated by the  $RTT$ and $CTT$
relations [$R$ matrix and metric $C$ of $S_{q,r}(N+2)$].
This can be expressed as:
\eq
S_{q,r}(N+2)\equiv \frac{<T^A{}_B>}{[RTT, CTT]}
\en
Therefore we have (recall that
$H\equiv[\T{a}{\circ},\T{\bullet}{b},\T{\bullet}{\circ}]
\equiv [\Tc]$) :
\eq
\ISqrN=  \frac{\SqrNt}{[\Tc]}
=\frac{ <T^A{}_B>/[RTT, CTT] }{[\Tc]}=
\frac{<T^A{}_B>}{[RTT, CTT, \Tc]}
\en
so that we have shown the following:
\sk
\noi{\bf Theorem} 4.2.1 $~$
The quantum inhomogeneous groups
$\ISOqrN$ and $\ISpqrN$ are freely generated by the
non-commuting matrix elements $\T{A}{B}$
 [{\small A}=$(\circ, a, \bullet)$, with $a=1,...N$)] and the
identity
$I$,
modulo the relations:
\eq
\T{a}{\circ}=\T{\bullet}{b}=\T{\bullet}{\circ}=0 , \label{Tprojectedsosp}
\en
\noi the $RTT$ relations
\eq
\R{AB}{EF} \T{E}{C} \T{F}{D} = \T{B}{F} \T{A}{E} \R{EF}{CD},
\label{RTTbigsosp}
\en
\noi and the orthogonality (symplecticity) relations
\eq
C^{BC} \T{A}{B}  \T{D}{C}= C^{AD} ,~~~
C_{AC} \T{A}{B}  \T{C}{D}=C_{BD} \label{CTTbig}
\en
The co-structures of $\ISOqrN$ and $\ISpqrN$
are simply given by:
\eq
\D (\T{A}{B})=\T{A}{C} \otimes \T{C}{B},
{}~~\kappa (\T{A}{B})=C^{AC} \T{D}{C} C_{DB} ,~~
\epsi (\T{A}{B})=\de^A_B~. \label{costructuresbigsosp}
\en
\cvd

After decomposing the indices {\small A}=$(\circ, a, \bullet)$, and
defining:
\eq
u\equiv \T{\circ}{\circ},~~v\equiv
\T{\bullet}{\bullet},~~z\equiv
\T{\circ}{\bullet},~~
x^a \equiv \T{a}{\bullet},~~y_a \equiv \T{\circ}{a} \label{names}
\en
\noi the relations (\ref{RTTbigsosp}) and (\ref{CTTbig}) become \cite{inson}:
\eqa
& &\R{ab}{ef} \T{e}{c} \T{f}{d} = \T{b}{f} \T{a}{e} \R{ef}{cd}
\label{PRTT11}\\
& &\T{a}{b} C^{bc} \T{d}{c}=C^{ad} I \label{PRTT31}\\
& &\T{a}{b} C_{ac} \T{c}{d} = C_{bd} I \label{PRTT32}
\ena
\eqa
& &\T{b}{d} x^a={r \over q_{d\bullet}} \R{ab}{ef} x^e \T{f}{d}
\label{PRTT33}\\
& &\PA{ab}{cd} x^c x^d=0 \label{PRTT13}\\
& &\T{b}{d} v={q_{b\bullet}\over q_{d\bullet}} v \T{b}{d}\\
& &x^b v=q_{b \bullet} v x^b \label{PRTT15}\\
& & uv=vu=I \label{PRTT21}\\
& &u x^b=q_{b\bullet} x^bu \label{PRTT22}\\
& &u \T{b}{d}={q_{b\bullet}\over q_{d\bullet}} \T{b}{d} u
\label{PRTT24}
\ena
\eq
y_b=-r^{\rho} \T{a}{b} C_{ac} x^c u \label{ipsilon}
\en
\eq
{(r^{-\rho}+\epsilon r^{\rho-2})}\, z=- x^b C_{ba} x^a u
\label{PRTT44}
\en
\noi where $q_{a\bullet}$ are $N$  complex parameters
related by $q_{a\bullet} = r^2 /q_{\ap\bullet}$, with $\ap = N+1-a$.
The matrix $P_A$ in eq. (\ref{PRTT13}) is the $q$-antisymmetrizer for
the
$B,C,D$ $q$-groups given by (cf.  (\ref{projBCD})):
\eq
\PA{ab}{cd}=- {1 \over {r+\rminus}}
(\Rhat{ab}{cd}-r\de^a_c \de^b_d + {r-r^{-1} \over
\epsilon r^{N-1-\epsilon} +1} C^{ab} C_{cd}). \label{PA}
\en
The last two relations (\ref{ipsilon}) - (\ref{PRTT44})
are constraints, showing that
the $\T{A}{B}$ matrix elements in eq. (\ref{RTTbigsosp})
are really a {\sl redundant} set. This redundance
is necessary if we want to express the $q$-commuations
of the $\ISOqrN$ and $\ISpqrN$ basic group elements
as $RTT=TTR$  (i.e. if we want an $R$-matrix
formulation).  Remark that,
in the $R$-matrix formulation
for $\IGLqrN$, {\sl all}  the $\T{A}{B}$
are independent. Here
we can take as independent generators the
elements
\eq
{}~~\T{a}{b} , x^a , v , u\equiv v^{-1}
\mbox{ and the identity }  I~~~~~(a=1,...N) \label{Txvu}
\en
The co-structures on the $\ISOqrN$ (or $\ISpqrN$) generators can be read from
(\ref{costructuresbigsosp})
after decomposing the indices \sma{$A=\circ,a,\bu$}:
\eqa
& &\D (\T{a}{b})=\T{a}{c} \otimes \T{c}{b}~,~~
\D (x^a)=\T{a}{c} \otimes x^c + x^a \otimes v ~,
\label{Pcoproduct2}\\
& &\D (v)=v \otimes v~,~~\D (u)=u \otimes u~,
\ena
\eqa
& &\kappa (\T{a}{b})
=C^{ac} \T{d}{c} C_{db}=\epsilon_a\epsilon_b
r^{-\rho_a+\rho_{b}}~ \T{\bp}{\ap}~,\\
& &\kappa (x^a)
=- \kappa (\T{a}{c}) x^c u~,~~\kappa (v)= u~,~~\kappa (u)= v~,
\ena
\eq
\epsi (\T{a}{b})=\de^a_b ~,~~\epsi (x^a)=0 ~,~~
\epsi (u)=\epsi (v)=\epsi(I)=1 ~.\label{cfin}
\en
\sk
In the commutative limit $q\rightarrow 1 , r\rightarrow 1$ we
recover the algebra of functions
on $ISO(N)$ (plus the dilatation $v$ that can be set to
the identity). In the $ISp(N)$ case, the $q\rightarrow 1 , r\rightarrow 1$ 
limit of 
relation (\ref{PRTT44}) 
implies $x^bC_{ba}x^a=0$ i.e. $P_0^{ab}{}_{ef}x^ex^f=0$, that with 
(\ref{PRTT13}) gives the commutativity of the
coordinates $x^a$ (both $P_A$ and $P_0$ are antisymmetrizers in the
symplectic case). We then  recover the algebra of functions on 
$ISp(N)$ plus the element $z$, 
that can be set to zero, and the dilatation $v$, that can be set to
the identity (see also Note 4.2.4).
\sk

\noi{\bf Note} {4.2.1} $~$In order to study the  $\ISOqrN$ and $\ISpqrN$ 
differential calculus and universal enveloping algebras we will use the 
definition (\ref{quotientsosp}) rather then Theorem 4.2.1 :
$\ISqrN =\frac{<T^A{}_B>}{[RTT, CTT, \Tc]}$. With abuse of
notations we therefore identify [cf. (\ref{names})]: 
$u= P(\T{\circ}{\circ}),~v=
P(\T{\bullet}{\bullet}),~z= P(\T{\circ}
{\bullet}), ~x^a = P(\T{a}{\bullet}),~y_a =
P(\T{\circ}{a}),
{}~T^a{}_b=P(T^a{}_b) ~;~~ I =P(I) $  where
$P\;:~\SqrNt\rightarrow \ISqrN\equiv\SqrNt/{H}$.
\sk
\noi{\bf Note} {4.2.2} $~$ From the commutations
(\ref{PRTT22}) - (\ref{PRTT24})
 we see that
one can set $u=I$ only when $q_{a\bullet}=1$ for all $a$.
{}From $q_{a\bullet} = r^2 /q_{\ap\bullet}$, cf. eq. (\ref{qab2}),
this implies also $r=1$.
\sk

\noi{\bf Note} {4.2..3} $~$ Eq.s (\ref{PRTT13}) are the multiparametric
orthogonal quantum plane commutations. They follow
from the   $({}^a{}_{\bullet} {}^b{}_{\bullet})$  $RTT$ components
and (\ref{PRTT44}).
\sk

\noi{\bf Note} 4.2.4 $~$ In the symplectic case eq.s (\ref{PRTT13})
alone are not sufficient to order in an arbitrary given way a monomial in
the $x$ elements; we can obtain an ordering only if we consider also the 
element $zv$ (or $z$) besides the $x$ elements. In other terms, 
the expression $x^aC_{ab}x^b$ appearing in (\ref{PRTT44})  cannot be
ordered as $\al_{ab}x^ax^b$ with $\al_{ab}\in 
{\bf \mbox{\boldmath $C$}}$ and $\al_{ab}=0$ if $a>b$.

To recover the symplectic $q$-plane commutations relations described
in \cite{FRT} one has to impose also the condition
\eq
{P_0}^{ab}_{~cd} x^c x^d=0 ~~\mbox{ i.e. }~  x^aC_{ab}x^b=0~,~~z=0
 \label{FRTplane}
\en
that arises naturally from  the characteristic equation and the
projector decomposition of the
$R$-matrix: in the symplectic case $P_0$ is an antisymmetrizer.
As a consequence the $xx$ commutations  (\ref{PRTT13}) become
\eq
\Rhat{ab}{cd} x^c x^d - r x^b x^a = 0~.
\en
Notice however that (\ref{FRTplane}) is not compatible with a
deformation of the whole symplectic group. Condition (\ref{FRTplane}) 
amounts to consider the Hopf quotient $ISp_{q,r}(N)/K$ 
where $K$ is the Hopf ideal generated by $z$. From $\D(z)=y_a\otimes
x^a+u\otimes z + z\otimes v$ 
$\in ISp_{q,r}(N)\otimes K + K\otimes ISp_{q,r}(N)$
we deduce that $y_a\otimes
x^a \in ISp_{q,r}(N)\otimes K + K\otimes ISp_{q,r}(N)$ and applying
$(m\otimes id)(id\otimes \D)$ to $y_a\otimes x^a$ we obtain that 
$x^bC_{bd}\otimes x^d$ $\in ISp_{q,r}(N)\otimes K + K\otimes ISp_{q,r}(N)$.
Projecting on $ISp_{q,r}(N)/K$ yields 
$x^bC_{bd}\otimes x^d=0$ since $K$ is the kernel of the
projection.  Now classically  $x^bC_{bd}\otimes x^d\not=0$. This
proves that the classical limit of $ISp_{q,r}(N)/K$ is not the algebra of
functions over  $ISp(N)$. 
\sk


\def\spinst#1#2{{#1\brack#2}}
\def\sk{\vskip .4cm}
\def\noi{\noindent}
\def\om{\omega}
\def\al{\alpha}
\def\be{\beta}
\def\ga{\gamma}
\def\Ga{\Gamma}
\def\del{\delta}
\def\vf{\varphi}
\def\linv{{1 \over \lambda}}
\def\rinv{{1\over {r-r^{-1}}}}
\def\alb{\bar{\alpha}}
\def\beb{\bar{\beta}}
\def\gab{\bar{\gamma}}
\def\deb{\bar{\delta}}
\def\ab{\bar{a}}
\def\Ab{\bar{A}}
\def\Bb{\bar{B}}
\def\Cb{\bar{C}}
\def\Db{\bar{D}}
\def\ab{\bar{a}}
\def\cb{\bar{c}}
\def\db{\bar{d}}
\def\bb{\bar{b}}
\def\eb{\bar{e}}
\def\fb{\bar{f}}
\def\gb{\bar{g}}
\def\xih{\hat\xi}
\def\Xih{\hat\Xi}
\def\uh{\hat u}
\def\vh{\hat v}
\def\ub{\bar u}
\def\vb{\bar v}
\def\xib{\bar \xi}
\def\ftilda{\tilde{f}}
\def\alp{{\alpha}^{\prime}}
\def\bep{{\beta}^{\prime}}
\def\gap{{\gamma}^{\prime}}
\def\dep{{\delta}^{\prime}}
\def\rhop{{\rho}^{\prime}}
\def\taup{{\tau}^{\prime}}
\def\rhopp{\rho ''}
\def\thetap{{\theta}^{\prime}}
\def\imezzi{{i\over 2}}
\def\unquarto{{1 \over 4}}
\def\onehalf{{1 \over 2}}
\def\unmezzo{{1 \over 2}}
\def\epsi{\varepsilon}
\def\we{\wedge}
\def\th{\theta}
\def\de{\delta}
\def\cony{i_{\de {\vec y}}}
\def\Liey{l_{\de {\vec y}}}
\def\tv{{\vec t}}
\def\Gt{{\tilde G}}
\def\deyv{\vec {\de y}}
\def\part{\partial}
\def\pdxp{{\partial \over {\partial x^+}}}
\def\pdxm{{\partial \over {\partial x^-}}}
\def\pdxi{{\partial \over {\partial x^i}}}
\def\pdy#1{{\partial \over {\partial y^{#1}}}}
\def\pdx#1{{\partial \over {\partial x^{#1}}}}
\def\pdyx#1{{\partial \over {\partial (yx)^{#1}}}}

\def\qP{q-Poincar\'e~}
\def\A#1#2{ A^{#1}_{~~~#2} }

\def\R#1#2{ R^{#1}_{~~~#2} }
\def\PA#1#2{ P^{#1}_{A~~#2} }

\def\Rp#1#2{ (R^+)^{#1}_{~~~#2} }
\def\Rpinv#1#2{ [(R^+)^{-1}]^{#1}_{~~~#2} }
\def\Rm#1#2{ (R^-)^{#1}_{~~~#2} }
\def\Rinv#1#2{ (R^{-1})^{#1}_{~~~#2} }
\def\Rsecondinv#1#2{ (R^{\sim 1})^{#1}_{~~~#2} }
\def\Rinvsecondinv#1#2{ ((R^{-1})^{\sim 1})^{#1}_{~~~#2} }

\def\Rpm#1#2{(R^{\pm})^{#1}_{~~~#2} }
\def\Rpminv#1#2{((R^{\pm})^{-1})^{#1}_{~~~#2} }

\def\Rb{{\bf R}}
\def\Rbo{{\bf R}}
\def\Rbp#1#2{{ (\Rbo^+)^{#1}_{~~~#2} }}
\def\Rbm#1#2{ (\Rbo^-)^{#1}_{~~~#2} }
\def\Rbinv#1#2{ (\Rbo^{-1})^{#1}_{~~~#2} }
\def\Rbpm#1#2{(\Rbo^{\pm})^{#1}_{~~~#2} }
\def\Rbpminv#1#2{((\Rbo^{\pm})^{-1})^{#1}_{~~~#2} }

\def\RRpm{R^{\pm}}
\def\RRp{R^{+}}
\def\RRm{R^{-}}

\def\Rh{{\hat R}}
\def\Rbh{{\hat {\Rbo}}}
\def\Rhat#1#2{ \Rh^{#1}_{~~~#2} }
\def\Rbar#1#2{ {\bar R}^{#1}_{~~~#2} }
\def\L#1#2{ \La^{#1}_{~~~#2} }
\def\Linv#1#2{ \La^{-1~#1}_{~~~~~#2} }
\def\Rbhat#1#2{ \Rbh^{#1}_{~~~#2} }
\def\Rhatinv#1#2{ (\Rh^{-1})^{#1}_{~~~#2} }
\def\Rbhatinv#1#2{ (\Rbh^{-1})^{#1}_{~~~#2} }
\def\Z#1#2{ Z^{#1}_{~~~#2} }
\def\Rt#1{ {\hat R}_{#1} }
\def\La{\Lambda}
\def\la{\lambda}
\def\Rha{{\hat R}}
\def\ff#1#2#3{f_{#1~~~#3}^{~#2}}
\def\MM#1#2#3{M^{#1~~~\!\!\!#3}_{~#2}}
\def\cchi#1#2{\chi^{#1}_{~#2}}
\def\chil#1{\chi_{{}_{#1}}}
\def\ome#1#2{\om_{#1}^{~#2}}
\def\RRhat#1#2#3#4#5#6#7#8{\La^{~#2~#4}_{#1~#3}|^{#5~#7}_{~#6~#8}}
\def\RRhatinv#1#2#3#4#5#6#7#8{(\La^{-1})^
{~#2~#4}_{#1~#3}|^{#5~#7}_{~#6~#8}}
\def\LL#1#2#3#4#5#6#7#8{\La^{~#2~#4}_{#1~#3}|^{#5~#7}_{~#6~#8}}
\def\LLinv#1#2#3#4#5#6#7#8{(\La^{-1})^
{~#2~#4}_{#1~#3}|^{#5~#7}_{~#6~#8}}
\def\U#1#2#3#4#5#6#7#8{U^{~#2~#4}_{#1~#3}|^{#5~#7}_{~#6~#8}}
\def\Cb{{\bf C}}
\def\CC#1#2#3#4#5#6{\Cb_{~#2~#4}^{#1~#3}|_{#5}^{~#6}}
\def\cc#1#2#3#4#5#6{C_{~#2~#4}^{#1~#3}|_{#5}^{~#6}}
\def\bu{\bullet}
\def\ci{\circ}

\def\C#1#2{ {\bf C}_{#1}^{~~~#2} }
\def\c#1#2{ C_{#1}^{~~~#2} }
\def\q#1{   {{q^{#1} - q^{-#1}} \over {q^{\unmezzo}-q^{-\unmezzo}}}}
\def\Dmat#1#2{D^{#1}_{~#2}}
\def\Dmatinv#1#2{(D^{-1})^{#1}_{~#2}}
\def\DR{\Delta_R}
\def\DL{\Delta_L}
\def\f#1#2{ f^{#1}_{~~#2} }
\def\F#1#2{ F^{#1}_{~~#2} }
\def\T#1#2{ T^{#1}_{~~#2} }
\def\t#1#2{ T^{#1}_{~~#2} }
\def\caM{{t}}
\def\M#1#2{ M_{#1}^{~\,#2} }
\def\Ti#1#2{ (T^{-1})^{#1}_{~~#2} }
\def\Tp#1#2{ (T^{\prime})^{#1}_{~~#2} }
\def\Th#1#2{ {\hat T}^{#1}_{~~#2} }
\def\TP{ T^{\prime} }
\def\qm{q^{-1}}
\def\rminus{r^{-1}}
\def\um{u^{-1}}
\def\vm{v^{-1}}
\def\xm{x^{-}}
\def\xp{x^{+}}
\def\fm{f_-}
\def\fp{f_+}
\def\fn{f_0}
\def\D{\Delta}
\def\DN{\Delta_{N+1}}
\def\kN{\kappa_{N+1}}
\def\eN{\epsi_{N+1}}
\def\Mat#1#2#3#4#5#6#7#8#9{\left( \matrix{
     #1 & #2 & #3 \cr
     #4 & #5 & #6 \cr
     #7 & #8 & #9 \cr
   }\right) }
\def\Ap{A^{\prime}}
\def\Dp{\Delta^{\prime}}
\def\Ip{I^{\prime}}
\def\ep{\epsi^{\prime}}
\def\kp{\kappa^{\prime}}
\def\kpm{\kappa^{\prime -1}}
\def\kpsq{\kappa^{\prime 2}}
\def\km{\kappa^{-1}}
\def\gp{g^{\prime}}
\def\qone{q \rightarrow 1}
\def\rone{r \rightarrow 1}
\def\qrone{q,r \rightarrow 1}
\def\Fmn{F_{\mu\nu}}
\def\Am{A_{\mu}}
\def\An{A_{\nu}}
\def\dm{\part_{\mu}}
\def\dn{\part_{\nu}}
\def\Ana{A_{\nu]}}
\def\Bna{B_{\nu]}}
\def\Zna{Z_{\nu]}}
\def\dma{\part_{[\mu}}
\def\qsu{$[SU(2) \times U(1)]_q~$}
\def\suq{$SU_q(2)~$}
\def\su{$SU(2) \times U(1)~$}
\def\gij{g_{ij}}
\def\qL{SL_q(2,{\bf C})}
\def\GLqrN{GL_{q,r}(N)}
\def\IGLqrN{IGL_{q,r}(N)}
\def\IGLqrtwo{IGL_{q,r}(2)}
\def\GLqrNo{GL_{q,r}(N+1)}
\def\SOqrNt{SO_{q,r}(N+2)}
\def\SpqrNt{Sp_{q,r}(N+2)}
\def\SLqrN{SL_{q,r}(N)}
\def\UglqrN{U(gl_{q,r}(N))}
\def\UglqrNo{U(gl_{q,r}(N+1))}
\def\UiglqrN{U(igl_{q,r}(N))}
\def\ISOqrN{ISO_{q,r}(N)}
\def\ISpqrN{ISp_{q,r}(N)}
\def\ISpqroN{ISp_{q,r=1}(N)}
\def\SqrN{S_{q,r}(N)}
\def\ISqrN{IS_{q,r}(N)}

\def\SOqrN{SO_{q,r}(N)}
\def\SpqrN{Sp_{q,r}(N)}
\def\SqrNt{S_{q,r}(N+2)}

\def\Tc{{\cal T}}
\def\Dcal#1#2{{\cal D}^{#1}_{~#2}}

\def\Dtwo{\Delta_{N+2}}
\def\epsitwo{\epsi_{N+2}}
\def\kappatwo{\kappa_{N+2}}

\def\RR{R^*}
\def\rr#1{R^*_{#1}}

\def\Lpm#1#2{L^{\pm #1}_{~~~#2}}
\def\Lmp#1#2{L^{\mp#1}_{~~~#2}}
\def\LLpm{L^{\pm}}
\def\LLmp{L^{\mp}}
\def\LLp{L^{+}}
\def\LLm{L^{-}}
\def\Lp#1#2{L^{+ #1}_{~~~#2}}
\def\Lm#1#2{L^{- #1}_{~~~#2}}
\def\Lc{{{\cal L}^+}}
\def\Ld{L^{\!+}_{\bu}}
\def\Linvpm#1#2{(L^{-1})^{\pm #1}_{~~~#2}}

\def\Li{{\pounds}}
\def\Lip#1#2{{\pounds}^{+ #1}_{~~~#2}}
\def\Lid{{\pounds}^{\!+}_{\bu}}
\def\Lipm#1#2{{\pounds}^{\pm #1}_{~~~#2}}
\def\LLipm{{\pounds}^{\pm}}
\def\LLip{{\pounds}^{+}}

\def\gu{g_{U(1)}}
\def\gsu{g_{SU(2)}}
\def\tg{ {\rm tg} }
\def\Fun{$Fun(G)~$}
\def\invG{{}_{{\rm inv}}\Ga}
\def\Ginv{\Ga_{{\rm inv}}}
\def\qonelim{\stackrel{q \rightarrow 1}{\longrightarrow}}
\def\ronelim{\stackrel{r \rightarrow 1}{\longrightarrow}}
\def\qronelim{\stackrel{q=r \rightarrow 1}{\longrightarrow}}
\def\viel#1#2{e^{#1}_{~~{#2}}}
\def\ra{\rightarrow}
\def\detq{{\det}}
\def\detqr{{\det}}
\def\detqrm{{\det} }
\def\detqrTAB{{\det} \T{A}{B}}
\def\detqrTab{{\det} \T{a}{b}}
\def\P{P}
\def\Qt{Q}
\def\chit{{\partial}}

\def\pp#1#2{\Pi_{#1}^{(#2)}}

\def\BCD{B_n, C_n, D_n}

\def\n2{{{N+1} \over 2}}
\def\ap{a^{\prime}}
\def\bp{b^{\prime}}
\def\cp{c^{\prime}}
\def\dpr{d^{\prime}}
\def\Dc{{\cal D}}
\def\osqrt{{1 \over \sqrt{2}}}
\def\Ntwo{{N\over 2}}

\def\SO{SO_{q,r}(N+2)}
\def\S{S_{q,r}(N+2)}
\def\Sp{Sp_{q,r}(N+2)}
\def\ISO{ISO_{q,r}(N)}
\def\IS{IS_{q,r}(N)}
\def\ISp{ISp_{q,r}(N)}
\def\U{U_{q,r}(so(N+2))}
\def\Up{U_{q,r}(sp(N+2))}
\def\Usosp{U_{q,r}(s(N+2))}
\def\IU{U_{q,r}(iso(N))}
\def\IUp{U_{q,r}(isp(N))}
\def\IUsosp{U_{q,r}(is(N))}

\def\ISOqrN{ISO_{q,r}(N)}
\def\ISpqrN{ISp_{q,r}(N)}
\def\ISOqroN{ISO_{q,r=1}(N)}
\def\ISpqroN{ISp_{q,r=1}(N)}
\def\SqrN{S_{q,r}(N)}
\def\SqrNtwo{S_{q,r}(N+2)}
\def\USqrNtwo{U(s_{q,r}(N+2))}
\def\UISqrN{U(is_{q,r}(N))}
\def\ISqrN{IS_{q,r}(N)}
\def\SqroNt{S_{q,r=1}(N+2)}
\def\ISqroN{IS_{q,r=1}(N)}
\def\ISOqroN{ISO_{q,r=1}(N)}

\def\SOqrN{SO_{q,r}(N)}
\def\SpqrN{Sp_{q,r}(N)}
\def\SqrNt{S_{q,r}(N+2)}

\def\H{H^\bot}
\def\le{\langle}
\def\re{\rangle}
\def\ren{\rangle_{{}_{{}_{\!N+2}}}}
\def\limrone{\lim_{r \rightarrow 1}}
\def\linv{{1 \over \lambda}}
\def\sma#1{\mbox{\footnotesize #1}}

\def\Re{R}

\sk
\noi{\bf Note} 4.2.5 $~$
We here briefly study the structure of $\ISqrN$ with respect to
$\SqrN$, that is easily seen to be a  Hopf subalgebra of $\ISqrN$. It 
is also a quotient
of $\ISqrN$ via the Hopf algebra projection [well defined only if 
$q_{a\bu}=const ~\forall a $ i.e. $ q_{a\bu}=r~\forall a$ see (\ref{PRTT24})]:
\[
\pi(x^a)=0~~,~~~~\pi(u)=I~~,~~~~\pi(\T{a}{b})=\T{a}{b}~~,~~~~\pi({I})=I~.
\]
Then the results of Theorem 3.3.6  apply to $\ISqrN$ as 
well, and we can write
the Hopf algebra isomorhism 
\eq
\ISqrN\cong B{\mbox{$\times \!\rule{0.3pt}{1.1ex}\;\!\!\!\cdot\,$}}\SqrN
\en
where $B$ is the subalgebra of $\ISqrN$ generated by $u$ and $x^{a}$
(in the orthogonal case $B$ is the 
quantum orthogonal plane with dilatation $u$).
Also Theorem 3.3.4 hold for $\ISqrN$.  
[This theorem, neglecting the graded structure, is a consequence of  
Theorem 3.3.6, an explicit proof for the $\ISqrN$ case
follows the same steps as for $IGL_{q,r}(N)$].
The $(${\bf Z},{\bf N}$)$ grading 
is introduced in the following way: 
the elements $\T{a}{b}$ have grade $(0,0)$, 
the elements $x^{a}$
have grade $(0,1)$, the elements $u$ and $v=u^{-1}$
have grade $(1,0)$ and $(-1,0)$. 
This grading is compatible with the $RTT$ commutation 
relations. 

{}For $\ISO$, the generators $u$ and $x^a$ of $B$ can be ordered using
(\ref{PRTT13}) and (\ref{PRTT22}), and the Poincar\'e series of the
subalgebra $B$ is the same as that
of the commutative algebra in the $N+1$
symbols $u$, $x^a$ \cite{FRT}.
A linear basis
of $B$ is  therefore given by the ordered monomials:
$\zeta^i=u^{i_{\ci}} (x^1)^{i_{1}}...\,
(x^N)^{i_{N}}$
with $i_{\ci}
\in{\bf \mbox{\boldmath$Z$}}$,
$i_1,... i_N\in
{\bf \mbox{\boldmath$N$}}\cup\{0\}$. In the $ISp_{q,r}(N)$ case,
if we also consider the elements $z$, a linear basis of $B\subset ISp_{q,r}(N)$
is given by the ordered monomials 
$\zeta^i=u^{i_{\ci}} (x^1)^{i_{1}}...\,
(x^N)^{i_{N}} (zv)^{i_{N+1}}$
with $i_{\ci}
\in{\bf \mbox{\boldmath$Z$}}$,
$i_1,... i_N,i_{N+1}\in
{\bf \mbox{\boldmath$N$}}\cup\{0\}$ 
($zv$ commutes with the coordinates $x^a$).

Using (\ref{ordino}), or (\ref{PRTT33}) and (\ref{PRTT24}),
a generic element of $\ISqrN$
can be written as $\zeta^ia_i$ (and also $a_i\zeta^i$)
where $a_i\in S_{q,r}(N)$.
As in Corollary 3.1.1, we have that 
$\ISqrN$, for $q_{a\bu}=r~\forall a$, is a bicovariant 
bimodule over $\SqrN$ freely generated, as a right module, 
by the elements $\zeta^i$; moreover 
\eq  
\ISqrN ~=\sum_{(h,k)\in 
({\mbox{\scriptsize \bf Z}},{\mbox{\scriptsize \bf N}})}
{}^{{\hskip -0.64 cm}\oplus} ~~\Ga^{(h,k)}\label{grad1s}
\en
where $\Ga^{(0,0)}=\SqrN$ 
\eqa
&&\!\!\!\Ga^{(0,1)}=\{x^{a}b_{a}\;~/~~b_a\in \SqrN\}~~,~~~~
\Ga^{(\pm 1,0)}=\{u^{\pm 1}b\;~/~~b\in  \SqrN\}\nonumber\\
&&\!\!\!\Ga^{(h,k)}=\{{u}^{h}
x^{a_1}x^{a_2}\ldots x^{a_k}b_{a_1a_2...a_k}~~/~~
b_{a_1a_2...a_k}
\in 
S_{q,r}(N)\}
\ena
Any submodule $\Ga^{(h,k)}$ is a
bicovariant bimodule freely generated by the ordered monomials 
$\zeta^i$ with degree
\sma{$(h,k)\in (\mbox{\bf{Z}},\mbox{\bf{N}})$}. 
We leave to the reader to 
reformulate Note 3.3.5 and Note 3.3.6 in this context.  
\sk
\noi{\bf Note} 4.2.6  $~$ 
Among  all the real forms of  $\SqrNt$ mentioned in the previous
section, only ${}^*$ and ${}^{*^{\sharp}}$  are inherited by $\ISqrN$, 
indeed only these two conjugations are compatible with
the ideal $H$: $H^*\subseteq H$ and $H^ {*^{\sharp}}\subseteq H$ 
[or, more easily, are compatible with (\ref{PRTT33})].
The 
conditions on the parameters are:
\sk
$\bullet$~~ $|q_{ab}|=|q_{a\bullet}|=|r|=1$ for $ISO_{q,r}(n,n;\Rb)$,
$ISO_{q,r}(n,n+1;\Rb)$
and $ISp_{q,r}(n;\Rb)$.
\sk
$\bullet$~~For $ISO_{q,r}(n+1,n-1;\Rb)$ : $|r|=1$;
$|q_{a\bullet}|=1$ for $a \not= n, n+1$;
$|q_{ab}|=1$ for $a$
and $b$
both different from $n$ or $n+1$;  $q_{ab}/r \in {\bf R}$
when at least one of the indices $a,b$ is equal
to $n$ or $n+1$;   $q_{a\bullet}/r \in \Rb$ for $a=n$ or $a=n+1$.
\sk
In particular,  the quantum Poincar\'e group $ISO_{q,r}(3,1;\Rb)$
is obtained by setting $|q_{1\bullet}|=|r|=1$,
$q_{2\bullet}/r \in \Rb$, $q_{12}/r  \in \Rb$.

According to  Note 4.2.2,  a dilatation-free $q$-Poincar\'e group
is found after the further restrictions
$q_{1\bullet}=q_{2\bullet}=r=1$.
The only free parameter remaining is then $q_{12} \in \Rb$.

\sect{Universal enveloping algebras  $\U$ and $\Up$}

We construct the universal enveloping algebra $\Usosp$ of $\S$ as  the
algebra of regular functionals \cite{FRT} on $\S$ (recall that $S$ stands 
for $SO$ and $Sp$).

$\Usosp$ is the algebra over $\mbox{\boldmath $C$}$ generated
by the counit $\epsi$ and by the functionals $\LLpm $
defined
by their value on the matrix elements $\T{A}{B}$  :
\eq
\Lpm{A}{B} (\T{C}{D})= \Rpm{AC}{BD}, \label{LonT}
\en
\eq
\Lpm{A}{B} (I)=\de^A_B \label{LonI}
\en
\noi with
\eq
\Rp{AC}{BD} \equiv \R{CA}{DB} \label{Rplussosp}~~;~~~
\Rm{AC}{BD} \equiv \Rinv{AC}{BD}~.
\en
To extend the definition (\ref{LonT})
to the whole algebra $\S$ we set
\eq
\Lpm{A}{B} (ab)=\Lpm{A}{C} (a) \Lpm{C}{B} (b)
{}~~~\forall a,b\in  \S~.
\label{Labsosp}
\en
{}From (\ref{LonT}),
using the upper and lower
triangularity of $R^+$ and $R^-$, we see that
$L^+$ is upper
triangular and $L^-$ is lower triangular.
\sk
The commutations between $\Lpm{A}{B}$
and $\Lpm{C}{D}$ are induced by (\ref{QYBE}) :
\eq
R_{12} \LLpm_2 \LLpm_1=\LLpm_1 \LLpm_2 R_{12} \label{RLLsosp}~,
\en
\eq
R_{12} \LLp_2 \LLm_1=\LLm_1 \LLp_2 R_{12}~, \label{RLpLmsosp}
\en
\noi where as usual the product $\LLpm_2 \LLpm_1$
is the convolution
product $\LLpm_2 \LLpm_1 \equiv (\LLpm_2 \otimes \LLpm_1)\D$.
\sk
The $\Lpm{A}{B}$ elements satisfy orthogonality conditions
analogous to (\ref{Torthogonality}):
\eqa
& &C^{AB} \Lpm{C}{B} \Lpm{D}{A} = C^{DC} \epsi\label{CLL1}\\
& &C_{AB} \Lpm{B}{C} \Lpm{A}{D} = C_{DC} \epsi \label{CLL2}
\ena
\noi as can be verified by applying them to the $q$-group
generators  and using (\ref{crc1}), (\ref{crc2}).
They provide the inverse for the matrix $\LLpm$:
\eq
[(\LLpm)^{-1}]^A{}_{\!B}=C^{DA} \Lpm{C}{D} C_{BC}
\label{Linversesosp}
\en
\sk
The co-structures of the algebra generated by the
functionals $L^{\pm}$ and $\epsi$
are defined by the duality (\ref{Labsosp}):
\eq
\D '(\Lpm{A}{B})(a \otimes b) \equiv \Lpm{A}{B}
(ab)=\Lpm{A}{G}(a) \Lpm{G}{B} (b)
\en
\eq
\ep (\Lpm{A}{B})\equiv \Lpm{A}{B} (I)
\en
\eq
\kp(\Lpm{A}{B})(a)\equiv \Lpm{A}{B} (\kappa (a))
\en
\noi so that
\eqa
\!\!\!\!& & \D ' (\Lpm{A}{B})=\Lpm{A}{G} \otimes
\Lpm{G}{B}\label{copLpm}\\
\!\!\!\!& & \epsi ' (\Lpm{A}{B})=\de^A_B \label{couLpmsosp}\\
\!\!\!\!& &\kp (\Lpm{A}{B})
= [(\LLpm)^{-1}]^A{}_{\!B}
= C^{DA} \Lpm{C}{D} C_{BC}
\label{coiLpmsosp}
\ena

{}From (\ref{coiLpmsosp}) we have that $\kp$ is an inner operation
in the algebra generated by the
functionals $\Lpm{A}{B}$ and $\epsi$, it is then easy to see that  
these
elements generate a Hopf algebra, the Hopf algebra $\Usosp$  of  regular
functionals on the quantum
group $\S$.
\sk
\noi{\bf Note} 4.3.1 $~$
{}From the $CLL$ relations
$\kappa'(\Lpm{A}{B})\Lpm{B}{C}=\Lpm{A}{B}\kappa'(\Lpm{B}{C})=\delta^A_ 
C\epsi$
we have, using upper-lower triangularity of $\LLpm$:
\eq
\Lpm{A}{A}\kappa'({\Lpm{A}{A}})=
\kappa'(\Lpm{A}{A})\Lpm{A}{A}=\epsi~~\mbox{ i.e. }~~
\Lpm{A}{A}\Lpm{A'}{A'}=\Lpm{A'}{A'}\Lpm{A}{A}=\epsi
\en
As a consequence  $\mbox{det}\LLpm\equiv \Lpm{\circ}{\circ}\Lpm{1}{1}
\Lpm{2}{2}\ldots\Lpm{N}{N}\Lpm{\bullet}{\bullet}
=\epsi$. In the $B_n$ case we also have $\Lpm{n_2}{n_2}=\epsi$.
\sk
\noi{\bf Note} 4.3.2 $~$
The $RLL$ relations imply that the subalgebra $U^0$ generated by the
elements
$\Lpm{A}{A}$ and $\epsi$ is commutative (use upper triangularity of  
$R$).
Moreover, from (\ref{copLpm}) the invertible elements $\Lpm{A}{A}$  
are
also group like, and we conclude that
$U^0$ is the group Hopf algebra of the abelian group generated by
$\Lpm{A}{A}$
and $\epsi ~.$
In the classical limit $U^0$ is a maximal commutative subgroup of
$S(N+2)$.
\sk
\noi{\bf Note} 4.3.3 $~$
When $q_{AB}=r$,
the multiparametric
$R$-matrix goes into the uniparametric $R$-matrix and we recover
the standard uniparametric orthogonal  (symplectic) quantum groups.
Then
the $\LLpm$ functionals satisfy the further relation:
\eq
\forall\, \sma{A}\;, ~~~~~~
\Lp{A}{A}\Lm{A}{A}
=\epsi \label{epsiaepsisosp}~,
\en
indeed $\Lp{A}{A}\Lm{A}{A}(a)=
\epsi(a)$ as can be easily seen when $a=\T{A}{B}$
and generalized to any $a\in \S$ using (\ref{Labsosp}).
In this case \cite{FRT}
we can avoid to realize the Hopf algebra
$U_r(s(N+2))$ as functionals on  $S_r(N+2)$ and we can define it
abstractly as
the Hopf algebra generated by  the {\sl symbols} $\LLpm$
and the unit $\epsi$
modulo relations
(\ref{RLLsosp}),(\ref{RLpLmsosp}),(\ref{CLL1}),(\ref{CLL2}), and
(\ref{epsiaepsisosp}).\sk

As discussed in \cite{FRT} in the uniparametric case, the Hopf  
algebra
$U_r(s(N+2))$
of regular functionals is a Hopf subalgebra of the
orthogonal (symplectic) Drinfeld-Jimbo universal enveloping algebra $U_h$, 
where
$r=e^h$.
In the general multiparametric case, relation (\ref{epsiaepsisosp})
does not hold any more. Here we discuss the generalization
of (\ref{epsiaepsisosp}) and the relation between $\Usosp$ and the  
multiparametric
orthogonal (symplectic) Drinfeld-Jimbo universal enveloping algebra  
$U_h^{(\cal{F})}$.
This latter is the quasitriangular Hopf algebra $U_h^{(\cal{F})}=
(U_h,\Delta^{({\cal{F}})},S,{\cal{R^{(F)}}})$ paired to the
multiparametric   $q$-group $\S$.
It is obtained from
$U_h=(U_h,\D,S,{\cal{R}})$ via a twist \cite{multiparam1}.
$U_h^{(\cal{F})}$ has the same algebra structure of $U_h$
(and the same antipode $S$),
while the coproduct $\Delta^{{\cal{(F)}}}$ and the universal element
${\cal{R^{(F)}}}$
belonging to (a completion of) $U_h\otimes U_h$
are determined by the twisting element ${\cal{F}}$
that belongs to (a completion of) a maximal commutative subalgebra of
$U_h\otimes U_h$. We have
\eq
\forall\,\phi\in U_h\,, ~~\D^{{\cal{(F)}}}(\phi)=
{\cal{F}}\D(\phi){\cal{F}}^{-1}~;~~{\cal{R^{(F)}}}=
{\cal{F}}_{21}{\cal{R}}{\cal{F}}^{-1}~;~~
{\cal{R^{(F)}}}(T\otimes T)=R_{q,r}~.
\en
The element ${\cal{F}}$ satisfies:
$(\D^{{\cal{(F)}}}\otimes  
id){\cal{F}}={\cal{F}}_{13}{\cal{F}}_{23}\,,~
(id  
\otimes\D^{{\cal{(F)}}}){\cal{F}}={\cal{F}}_{13}{\cal{F}}_{12}\,,~
{\cal{F}}_{12}{\cal{F}}_{21}=I\,,~
{\cal{F}}_{12}{\cal{F}}_{13}{\cal{F}}_{23}=
{\cal{F}}_{23}{\cal{F}}_{13}{\cal{F}}_{12}\,,$
$\, (\epsi\otimes id){\cal{F}}=(id\otimes \epsi){\cal{F}}
=\epsi\,,~(S \otimes id){\cal{F}}
=(id \otimes S ){\cal{F}}=
{\cal{F}}^{-1},\, \cdot(id\otimes S){\cal{F}}=
\cdot(S\otimes id){\cal{F}} =
\cdot(id \otimes id){\cal{F}} =\epsi\,$;
we explicitly have
\eq
{\cal{F}}(\T{A}{B}\otimes\T{C}{D})= F^{AC}_{~BD}
\en
where $F^{AC}_{~BD}$ is the diagonal matrix
\eq
F=diag (\sqrt{{q_{11}\over r}} ,
\sqrt{{ q_{12}\over r }} , ... ~ \sqrt{{ q_{NN} \over r}})
\en
It is easy to see that the definition of the $\LLpm$ functionals
given in the beginning of this section is equivalent to
the following one:
$\Lp{A}{B}(a)={\cal{R^{(F)}}} ( {a} \otimes \T{A}{B})$ and
 $\Lm{A}{B}(a)={\cal{R^{(F)}}}^{-1}(\T{A}{B}
\otimes {a})$. From
$(\D^{{\cal{(F)}}}\otimes  
id){\cal{R}}={\cal{R}}_{13}{\cal{R}}_{23}\,,~
(id  
\otimes\D^{{\cal{(F)}}}){\cal{R}}={\cal{R}}_{13}{\cal{R}}_{12}\,,~
$ we have $\D^{{\cal{(F)}}}(\Lpm{A}{B})=\Lpm{A}{C} \otimes  
\Lpm{C}{B}$
 and therefore $\D^{{\cal{(F)}}}=\D '$ on $\Usosp$. {}From $(id\otimes S)
({{\cal{R}}})=
(S\otimes id)({{\cal{R}}})={\cal{R}}^{-1}$
it is also easy to see that $S=\kp$ on $\Usosp$
and we conclude that the algebra of regular functionals $\Usosp$ is a
realization [in terms of  functionals on $\S$] of a
Hopf subalgebra of $U_h^{(\cal{F})}$ with $r=e^h$.
The generalization of (\ref{epsiaepsisosp}) lies in ${U_h^{{\cal{(F)}}}}$
and not in $\Usosp$, and it is given by
\eq
\forall\, \sma{$A$}~~~~~~~ \Lp{A}{A}\Lm{A}{A}=f_i(\T{A}{A})f^i~~  
\mbox{
where }
{\cal{F}}^4=f_i\otimes f^i~.\label{1effe1}
\en
This relation holds with $\LLpm$ considered as
abstract symbols. It can easily be checked
when $\LLpm$ are realized as functionals:
indeed $\Lp{A}{A}\Lm{A}{A}(a)=
{\cal{F}}^4(\T{A}{A}\otimes a)$ as can be seen when $a=\T{A}{B}$
[use ${\cal{F}}^2( \T{A}{A}\otimes b)={\cal{F}}(\T{A}{A}\otimes b_1)
{\cal{F}}(\T{A}{A}\otimes b_2)$] and generalized to any $a\in \S$
using ${\cal{F}}(\T{A}{A}\otimes ab)={\cal{F}}(\T{A}{A}\otimes a)
{\cal{F}}(\T{A}{A}\otimes b)$.

In order to characterize the relation between the Hopf algebra
of regular functionals $\Usosp$ and $U_h^{\cal{(F)}}$,
following \cite{FRT}, we extend the group Hopf algebra $U^0$
described in { Note 4.3.2}
to $\hat{U}^0$
by means of the elements
\footnote{
In the classical limit ${\ell^{\pm}}^A{}_A$ are the tangent vectors
to a maximal commutative subgroup of $S(N+2)$. They generate a
Cartan subalgebra of the Lie algebra $s(N+2)$.}
${\ell^{\pm}}^A{}_A=\ln \Lpm{A}{A}$. Otherwise stated this means that
in $\hat{U}^0$ we can write $\Lpm{A}{A}=
\mbox{exp}({\ell^{\pm}}^A{}_A)$ where ${\ell^{\pm}}^A{}_A\in  
\hat{U}^0$.
[Explicitly ${\ell^{\pm}}^A{}_A(\T{C}{D})=\ln
({R^{\pm}}^{AC}_{~AC})\,\delta^C_D$,
 ${\ell^{\pm}}^A{}_A(I)=0 $,
${\ell^{\pm}}^A{}_A( ab)={\ell^{\pm}}^A{}_A( a) \epsi (b) + \epsi(a)
{\ell^{\pm}}^A{}_A( b)$
and  $\kp({\ell^{\pm}}^A{}_A)=-{\ell^{\pm}}^A{}_A$ ].
It then follows that ${\cal{F}}$ belongs to (a completion of)
$\hat{U}^0\otimes\hat{U}^0$.
The corresponding extension
$\hat{U}_{q,r} (s(N+2))$ of $\Usosp$, defined as the Hopf algebra  
generated
by the {\sl symbols} $\LLpm$ and $\ell^{\pm}$ modulo relations
(\ref{RLLsosp})-(\ref{CLL2}) and (\ref{1effe1}), is isomorphic -- when  
$r=e^h$
--
to $U_h^{(\cal{F})}\;:\,~\hat{U}_{q,r}(s(N+2))
\cong U_h^{(\cal{F})}$. This relation holds because it
is the twisted version of the
known uniparametric analogue
$\hat{U}_{r}(s(N+2))
\cong U_h$ \cite{FRT,Frenkel}.
\sk

The elements $\LLpm$ [or ${1\over{r-r^{-1}}}(\Lpm{A}{B}-
\de^A_B\epsi$)]
may be seen as the quantum analogue of the tangent vectors; then
the $RLL$ relations are the quantum analogue of the Lie algebra
relations, and we can use the orthogonal (symplectic) $CLL$ conditions
to reduce the number of the $\LLpm$ generators to $(N+2)(N+1)/2$,  
(orthogonal case) or  $(N+2)(N+3)/2$ (symplectic 
case)  
i.e. the dimension of the classical group manifold.

This we proceed to do for $\U$, for $\Up$ one can proceed in a similar 
way\footnote{In the $\Up$ case relations (\ref{PLL}) [and more in general 
(\ref{PAJ'J})] do not allow to order the $\Lm{a}{\ci}$ 
(and more in general the $\Lpm{\al}{J}$) elements
because we are missing the relation with the $P_0$ projector, cf. Note 4.2.4.
However we can still order the  $\Lm{a}{\ci}$  (or $\Lpm{\al}{J}$) elements if 
we consider also  $\Lm{\bu}{\ci}$
 (or $\Lpm{J'}{J}$).   This leads to the  $(N+2)(N+3)/2$ generators of the
symplectic case. Notice that Lemma 4.3.1 and 4.3.2 still hold; 
Theorem 4.3.1 holds as well  provided that $\al$ 
can also be equal to $J'$: $J'<\alpha\leq J$.}; 
we next  study the $R\LLpm\LLpm$
commutation relations restricted to these $(N+2)(N+1)/2$ generators
and find a set of ordered monomials in the reduced $\LLpm$ that  
linearly
span all
$\hat{U}_{q,r}(so(N+2))$.

We first observe that the commutative subalgebra $\hat{U}^0$
is generated by ${(N+2)/2}$ elements ($N$ even, $N=2n$) or  
${(N+1)/2}$
elements ($N$ odd, $N=2n+1$),  for example
$\ell^{-\ci}{}_{\!\ci},$
$\ell^{-1}{}_{\!1}$ ... $\ell^{-n}{}_{\!n}$.
For the off-diagonal
$\LLpm $ elements, we can choose as free indices $(C,D)=(c,\ci)$
in  relation (\ref{CLL2}),  and using
$\Lm{\ci}{\ci}\Lm{\bu}{\bu}=\epsi$, we find:
\eq
\Lm{\bu}{c}=
-(C_{\ci \bu})^{-1} C_{ab}\Lm{b}{c}\,\Lm{a}{\ci}\Lm{\bu}{\bu}~.
\label{restrictedL1}
\en
If we choose $(C,D)=(\ci,\ci)$ we obtain
\eq
\Lm{\bu}{\ci}= - (r^{-2} C_{\bu\ci} +C_{\ci\bu})^{-1}
C_{ab}\Lm{b}{\ci}\,\Lm{a}{\ci}\Lm{\bu}{\bu}~. \label{restrictedL2}
\en
Similar results hold for $\Lp{\ci}{d}$ and $\Lp{\ci}{\bu}$.
Iterating this procedure, from
$C_{ab} \Lm{b}{c} \Lm{a}{d} = C_{dc} \epsi$
we find that $\Lm{N}{i}$ (with $i=2,...N-1$) and $\Lm{N}{1}$ are
functionally dependent on $\Lm{i}{1}$ and $\Lm{N}{N}$. Similarly
for $\Lp{1}{i}$ and $\Lp{1}{N}$. The final result is that the  
elements
$~\Lm{a}{J}$ with $\sma{$J<a<J'$}$ and
$~\Lp{a}{J}$ with $\sma{$J'<a<J$}$ --
whose number in both $\pm$ cases is ${1\over 4}N(N+2)$ for $N$ even
and ${1\over 4}(N+1)^2$ for $N$ odd --
and the elements
$\ell^{-\ci}{}_{\!\ci},$
$\ell^{-1}{}_{\!1}$... $\ell^{-n}{}_{\!n}$
generate all $\hat{U}_{q,r}(so(N+2))$.
The total number of  generators is therefore
$(N+2)(N+1)/2$.
\sk

Notice that in this derivation we have not used the $RLL$ relations
(i.e. the quantum analogue of the Lie algebra  relations)
to  further reduce the number of generators. We  therefore  expect  
that,
as in
the classical case,
monomials in the $(N+2)(N+1)/2$ generators can be ordered
(in any arbitrary way).
We begin by proving this for  polynomials  in $\Lp{A}{A}$,
$\Lp{\al}{J}$ with
$\sma{$J'<\al<J$}$,
and
for polynomials in  $\Lm{A}{A}$,
$\Lm{\al}{J}$ with $\sma{$J<\al< J'$}$ .
\sk
\noi{\bf Lemma} 4.3.1 $~$ 
Consider the $R\LLpm\LLpm$ commutation  relations
\eq
\R{AB}{EF} \Lpm{F}{D} \Lpm{E}{C} = \Lpm{A}{E} \Lpm{B}{F} \R{EF}{CD}
\label{Rell}~.
\en
{}For $\sma{$C\not= D$}$ they close respectively on the subset of the
$\Lp{\al}{J}$ with
$\sma{$J'<\al\leq J$}$
and on the subset of the $\Lm{\al}{J}$ with $\sma{$J\leq\al<J'$}$.
{}For
$\sma{$C=D$}$ they are
equivalent to the $q^{-1}$-plane commutation relations:
\eq
[P_A\sma{$(J'\!-\!J\!+\!1)$}]^{\al\be}_{~\ga\de}\Lpm{\de}{J}\Lpm{\ga}{ 
J}=0
{}~~,~~~
\label{PAJ'J}
\en where $P_A\sma{$(J'\!-\!J\!+\!1)$}$ is the antisymmetrizer in
dimension
$\sma{$J-J'+1$}$
[compare with (\ref{projBCD})].
In particular
\eq
\PA{ab}{cd} \Lm{d}{\ci} \Lm{c}{\ci}=0   \label{PLL}
\en
or equivalently
$[(P_{\!A}){}_{{}_{q^{-1}\!,r^{-1}}}]^{ab}_{~cd}\,\Lm{c}{\ci}
\Lm{d}{\ci}=0$
which coincide, for
$r\rightarrow r^{-1}$ and $q\rightarrow q^{-1}$, with the
$N$-dimensional quantum orthogonal plane relations
(\ref{PRTT13}).

\noi {\sl Proof :} $~~~$The proof is a straightforward calculation   
based
on
(\ref{foraa'}) and on upper or lower triangularity of the $R$ matrix
and of the $\LLpm$ functionals.

\cvd
\noi {\bf Lemma} 4.3.2  $~U_{q,r}(so(N))$ is a Hopf subalgebra of $\U$.

\noi {\sl Proof}:  Choosing $SO_{q,r}(N)$ indices as free indices in
(\ref{Rell})
and using upper or lower triangularity of the $\LLpm$ matrices, and
(\ref{Rnonzerososp}) or (\ref{Rbig}), we find
that only  $SO_{q,r}(N)$ indices appear in (\ref{Rell}); similarly  
for
relations
(\ref{RLpLmsosp})-(\ref{CLL2}), and for the costructures
(\ref{copLpm})-(\ref{coiLpmsosp}).
\cvd

Now we observe that in virtue of the
$R\LLp\!\LLp$  relations the $\LLp$
elements can be ordered; similarly  we can order the $\LLm$ using the
$R\LLm\!\LLm$ relations.
This statement can be proved by induction using that
$U_{q,r}(so(N))$ is a subalgebra of $\U$, and splitting the $\SO$
index in the usual way [some of the resulting formulas are given in
(\ref{1elle})-(\ref{4elle})].

It is then  straightforward to prove that the elements $\Lp{\al}{J}$  
with
$\sma{$J'<\al\leq J$}$  can
be ordered;
indeed
we can always order the $\Lp{\al}{J}~\Lp{\be}{K}$ with
$\sma{$J'<\al\leq J$}$,
$\sma{$K'<\be\leq K$}$ and $\sma{$J\not= K$}$ since
their commutation relations are a
closed subset
of (\ref{Rell}) [see { Lemma 4.3.1}].
Then there is no difficulty in ordering substrings
composed by
$\Lp{\al}{J}$ and
$\Lp{\be}{J}$ elements because
(\ref{PAJ'J}) are $q^{-1}$-plane
commutation relations, that allow for any
ordering of the  quantum plane coordinates
 \cite{FRT}.
More in general the  $\Lp{A}{A}$
and $\Lp{\al}{J}$ with
$\sma{$J'<\al<J$}$
can be ordered because of
$
\Lp{A}{A}\Lp{B}{C}$=$(q_{{}_{BA}}/q_{{}_{CA}})\Lp{B}{C}
\Lp{A}{A} \,.$
Similarly we can order the $\Lm{A}{A}$
and $\Lm{\al}{J}$ with
$\sma{$J<\al<J'$}$.
It is now easy to prove the following
\sk
\noi {\bf Theorem} 4.3.1 $~$
A set of elements spanning  $\hat{U}_{q,r}(so(N+2))$ is given by
the ordered monomials
\eq
Mon(\Lp{\al}{J};\sma{$J'<\al<J$})\;(\ell^{-\ci}{}_{\!\ci})^{p_{\ci}}
(\ell^{-1}{}_{\!1})^{p_1}\ldots(\ell^{-n}{}_{\!n})^{p_n}
\;Mon(\Lm{\al}{J};\sma{$J<\al<J'$})
\label{Starrr}
\en
where $p_{\ci}, p_1,...p_n\in {\bf \mbox{\boldmath$N$}}\cup\{0\}$,
$n=N/2$ ($N$ even), $n=(N-1)/2$ ($N$ odd) and
$Mon(\Lp{\al}{J};\sma{$J'<\al<J$})$, $ [Mon(\Lm{\al}{J};
\sma{$J<\al<J'$})]$ is
a monomial in the off-diagonal elements $\Lp{\al}{J}$ with
$\sma{$J'<\al<J$}$
[$\Lm{\al}{J}$ with $\sma{$J<\al<J'$}$] where an ordering has been  
chosen.
\cvd
\noi{\bf Note} 4.3.4 $~$  Conjecture: the above monomials are linearly  
independent
and therefore  form a basis of
$\hat{U}_{q,r}(so(N+2))\,.$
\sk

\noi {\bf Conjugation}
\sk

The canonical  $*$-conjugation on 
$\Usosp$ induced by the  $*$-conjugation on $\S$ 
is given by:
\eq
\psi^{*} (a)\equiv {\overline {\psi(\kappa^{-1}(a^*))}} 
\en
\noi where $\psi \in \Usosp$, $a \in \S$, and  the overline denotes
the usual complex conjugation.
It is not difficult 
to determine the action
on the basis elements $\Lpm{A}{B}$.
The two $\SqrNtwo$ $*$-conjugations that are compatible with $ISO_{q,r}{N}$
induce respectively the following conjugations on the
$\Lpm{A}{B}$ (we denote  ${}^{*^{\sharp}}$ simply by ${}^*$):
\eqa
& &~~~~(\Lpm{A}{B})^*=\kappa'^2(\Lpm{A}{B}) \label{conjL1}\\
& & ~~~(\Lpm{A}{B})^*=\Dcal{A}{C} \kappa'^2(\Lpm{C}{D}) \Dcal{D}{B}
\label{conjL2}
\ena
Notice that $\kappa'^2(\Lpm{A}{B})={D^{-1}}^A_E \Lpm{E}{F} D^F_B$ where
$D^a_e=C^{as}C_{es}$ and its inverse is ${D^{-1}}^f_b=C^{sf}C_{sb}$. 

\sect{Universal enveloping algebras $U_{q,r}(iso(N))$ and $U_{q,r}(isp(N))$}
Consider a generic functional $f\in\Usosp$. It is well defined on the
quotient
$\IS =\S /H$ if and only if $f(H)=0$.
It is  easy to see that the set $\H$ of all these functionals is a
subalgebra
of $\Usosp$ :
if $f(H)=0$ and $g (H)=0$ then $fg(H)=0$ because
$\Delta  (H) \subseteq
H \otimes \SqrNt + \SqrNt \otimes H.$
Moreover \cite{Sweedler} $\H$ is a Hopf subalgebra of $\Usosp$  since $H$ is 
a Hopf ideal,
cf. sections 3.5-6.
In agreement with these observations we will find the Hopf algebra
$\IU$ [dually paired to $\IS$] as a subalgebra of $\Usosp$ vanishing
on the ideal $H$.

Let
\eq
IU\equiv [L^{-A}{}_B, L^{+a}{}_b, L^{+\circ}{}_{\circ},
L^{+\bullet}{}_{\bullet}, \epsi]
\subseteq\Usosp\label{IU}
\en
be the subalgebra of $\Usosp$ generated by
$L^{-A}{}_B, L^{+a}{}_b, L^{+\circ}{}_{\circ},
L^{+\bullet}{}_{\bullet}, \epsi
.$
\sk
\noi{\bf Note} 4.4.1 $~$  These are all and only the functionals
annihilating
the generators of $H$:
$\T{a}{\circ}\;,~\T{\bullet}{b}$ and  $\T{\bullet}{\circ}\;$.
The remaining $\Usosp$ generators
$L^{+\circ}{}_b~,~L^{+a}{}_{\bullet}~,~L^{+\circ}{}_{\bullet}$ do not
annihilate
the generators of $H$ and are not included in
(\ref{IU}).
\sk
We now proceed to study  this algebra $IU$. We will show that it is a
Hopf algebra and that $IU\subseteq \H$; we will give an $R$-matrix
formulation,
and prove that $IU$ is the semidirect product of  $U_{q,r}(s(N))$ and the 
algebra $B'$ generated by the elements $\Lm{\ci}{\ci}$ and $\Lm{a}{\ci}$. 
This is the analogue of
$IS_{q,r}(N)$ being the semidirect product of $S_{q,r}(N)$
and the 
algebra $B$ generated by the elements $u$ and $x^a$,
cf. Note 4.2.5. 

We then show that $IU$ is dually paired with $\IS$. These results
lead to the conclusion that $IU$ is the
universal enveloping algebra of $IS_{q,r}(N)$.
\sk
\noi \noi{\bf Theorem} 4.4.1
$~IU$ is a Hopf subalgebra of $\Usosp$.

\noi {\sl Proof :} $~~~$
$IU$ is by definition a subalgebra. The sub-coalgebra property
$\Delta '(IU)
\subseteq IU\otimes IU$
follows immediately from the upper triangularity of $L^{+A}{}_B$:
\eq
\Delta ' (L^{+a}{}_b)=L^{+a}{}_c\otimes L^{+c}{}_b~;~
\Delta '(L^{+\circ}{}_{\circ})=L^{+\circ}{}_{\circ}
\otimes L^{+\circ}{}_{\circ}~;~
\Delta '(L^{+\bullet}{}_{\bullet})
=L^{+\bullet}{}_{\bullet}\otimes L^{+\bullet}{}_{\bullet}
\en
and the compatibility of $\Delta '$ with the product.
We conclude that $IU$ is a Hopf-subalgebra because $\kp(IU)\subseteq  
IU$
as is easily seen using (\ref{coiLpmsosp}) and
antimultiplicativity of $\kp$.
{\cvd}
We may wonder whether the $RLL$ and $CLL$ relations  of $\Usosp$ close in
$IU$.
In this case $IU$ will be given by
all and {\sl only} the polynomials in the functionals
$L^{-A}{}_B, L^{+a}{}_b, L^{+\circ}{}_{\circ},
L^{+\bullet}{}_{\bullet}, \epsi .$
This check is done by writing explicitly all $q$-commutations between
the
generators of $IU$: they do not involve the functionals
$L^{+\circ}{}_b~,~L^{+a}{}_{\bullet}~,~L^{+\circ}{}_{\bullet}$ .
Moreover one can also write them in a compact form using a new
$R$-matrix
${\cal{R}}_{12}\equiv\Lc_2({{t}}_1)$, where $\Lc$ is defined below.
Similarly
the orthogonality (symplecticity)
conditions (\ref{CLL1})-(\ref{CLL2})
do not relate elements of $IU$ with elements not belonging to $IU$.
We therefore conclude
\sk
\noi \noi{\bf Theorem} 4.4.2 $~$
The Hopf algebra $IU$ is generated by the unit $\epsi$ and
the matrix entries:
\eq
L^-=\left(L^{-A}{}_{B_{{}_{}}}\right)
{}~;~~
\Lc=\Mat{L^{+\circ}{}_{\circ}}{0}{0}{0}{L^{+a}{}_b}{0}{0}{0}
{L^{+\bullet}{}_{\bullet}}~;
\en
these functionals satisfy the $q$-commutation relations:
\eq
R_{12} \Lc_2 \Lc_1=\Lc_1 \Lc_2 R_{12} ~\mbox{ or equivalently }~
{\cal{R}}_{12}
\Lc_2 \Lc_1=\Lc_1 \Lc_2 {\cal{R}}_{12} \label{iRLcLcsosp}
\en
\eq
R_{12} \LLm_2 \LLm_1=\LLm_1 \LLm_2 R_{12}~, \label{iRLLsosp}
\en
\eq
{\cal{R}}_{12} \Lc_2 \LLm_1=\LLm_1 \Lc_2 {\cal{R}}_{12}~,
\label{iRLpLmsosp}
\en
where
$$
{\cal{R}}_{12}\equiv\Lc_2({{t}}_1) ~~~
\mbox{ that is }~~~
{\cal{R}}^{ab}_{cd}={R}^{ab}_{cd}~;~~
{\cal{R}}^{AB}_{AB}={R}^{AB}_{AB} ~\mbox{ and otherwise
}~{\cal{R}}^{AB}_{CD}
=0
$$
and the orthogonality (symplecticity) conditions :
\eq
C^{AB} \Lc^{C}{}_{B} \Lc^{D}{}_{A} = C^{DC} \epsi ~;~~
C_{AB} \Lc^{B}{}_{C} \Lc^{A}{}_{D} = C_{DC} \epsi ~;\label{iCLcLc}
\en
\eq
C^{AB} \Lm{C}{B} \Lm{D}{A} = C^{DC} \epsi ~;~~
C_{AB} \Lm{B}{C} \Lm{A}{D} = C_{DC} \epsi~, \label{iCLL}
\en
The costructures are the ones given in (\ref{copLpm})-(\ref{coiLpmsosp})
with $\LLp$ replaced by $\Lc$. 
\sk
{\cvd}
\noi{\bf Note} 4.4.2 $~$
We can consider the extension
$\hat{IU}\subset\hat{U}_{q,r}(s(N+2))$ obtained by including
the elements $\ell^{\pm A}{}_A$ (${\ell^{\pm}}^A{}_A=\ln\Lpm{A}{A}$,
see the previous section). Then
$\hat{IU}$
is generated by the symbols $\Lm{A}{B},\, \Lc^A{}_B,\, \ell^{\pm  
A}{}_A$
modulo the relations
(\ref{iRLcLcsosp})-(\ref{iCLL}) and (\ref{1effe1}) [({\ref{epsiaepsisosp}) in  
the
uniparametric case].
Equivalently,  from (\ref{restrictedL1})-(\ref{restrictedL2}),
we have that $\hat{IU}$ is generated by
$\hat{U}_{q,r}(s(N))$, the dilatation $\ell^{-\ci}{}_{\!\ci}$ and
the $N$ elements $\Lm{a}{\ci}$ (satisfying, in the orthogonal case, 
the quantum plane  relations). 
All the relations are then given by  those between the
generators of  $\hat{U}_{q,r}(s(N))$
--listed in  (\ref{RLLsosp})-(\ref{CLL2}), (\ref{1effe1}) with lower case
indices--
and by the following ones
\eqa
& &\Lm{\ci}{\ci}\Lm{a}{\ci}=q^{-1}_{\ci a}\Lm{a}{\ci}\Lm{\ci}{\ci}
\label{1elle}\\
& &P_A^{ab}{}_{\!fe}\Lm{e}{\ci}\Lm{f}{\ci}=0\label{2elle}\\
&
&\Lm{\ci}{\ci}\Lpm{b}{d}={q_{b\ci}\over{q_{d\ci}}}\Lpm{b}{d}\Lm{\ci}{\ci}
\label{3elle}\\
& &\Lm{a}{\ci}\Lpm{b}{d}={r\over {q_{d\ci}}}(R^{\pm})^{ba}{}_{\!ef}
\Lpm{e}{d}\Lm{f}{\ci}\label{4elle}
\ena
where $R^{\pm}$ is defined in (\ref{Rplussosp}).
The number of generators is $N(N-1)/2+ N + 1$ in the orthogonal case and
$N(N +1)/2+ N + 2$ in the symplectic case because we consider also
the element $\Lm{\bu}{\ci}$ so that the 
$\Lm{a}{\ci}$ elements can be ordered (cf. last footnote).
\sk
\noi{\bf Note} 4.4.3 $~$
When  $q_{a\bu}=r\;\forall a$, then   
$\Lm{\ci}{\ci}=\Lp{\bu}{\bu}\,,\,
\Lm{\bu}{\bu}=\Lp{\ci}{\ci}$ and,  in complete analogy to   
(\ref{Txvu}),
$IU$ is generated by ${U}_{q,r}(s(N))$,
$\Lm{a}{\ci},\Lm{\ci}{\ci}$ and  $\Lm{\bu}{\bu}=(\Lm{\ci}{\ci})^{-1}$.
With abuse of notations we will consider $IU$ generated by these  
elements
also
for arbitrary values of the parameters $q_{a\bu}$; this is  what  
actually
happens in $\hat{IU}$.
\sk
\noi{\bf Note} 4.4.4 $~$ From the second equation in (\ref{iRLcLcsosp}) 
applied
to ${{t}}$ we obtain the quantum Yang-Baxter equation for the matrix
${\cal{R}}$.
\sk
\sk
The results of { Note 4.2.5} holds also for $\IUsosp$ with the obvious
changes in notation. The projection $\pi$
[well defined only if  
$ q_{a\bu}=r~\forall a$ see (\ref{3elle})]
is given by: 
\[
\pi(\Lm{a}{\ci})=0~~,~~~~\pi(\Lm{\ci}{\ci})
=I~~,~~~~\pi(\Lpm{a}{b})=\Lpm{a}{b}~~,~~~~\pi({\epsi})=\epsi~.
\]
The semidirect product structure is:
\eq
\IUsosp\cong B'{\mbox{$\times \!\rule{0.3pt}{1.1ex}\;\!\!\!\cdot\,$}}
U_{q,r}(s(N))
\en
where $B'$ is the subalgebra of  $\IUsosp$ generated by $\Lm{\ci}{\ci}$ and 
$\Lm{a}{\ci}$. Moreover $B'=IU{}_{\rm{inv}}$,  the space of all right 
invariant elements of the $U_{q,r}(s(N))$--bicovariant algebra $\IUsosp$.
The ordered monomials that form a basis of $B'$
and that freely generate $IU$ as a right module, in the orthogonal case are:
\[
\eta^i=(\Lm{\ci}{\ci})^{i_{\ci}} (\Lm{1}{\ci})^{i_{1}}...\,
(\Lm{N}{\ci})^{i_{N}}~~~
\mbox{ with } i_{\ci}
\in{\bf \mbox{\boldmath$Z$}}~,~~
i_1,... i_N\in
{\bf \mbox{\boldmath$N$}}\cup\{0\}~.
\]
In the symplectic case,
if we also consider the element $z$, a linear basis of $B'$
is given by the ordered monomials 
\[
\eta^i=(\Lm{\ci}{\ci})^{i_{\ci}} (\Lm{1}{\ci})^{i_{1}}...\,
(\Lm{N}{\ci})^{i_{N}}(\Lm{\bu}{\ci}\Lm{\ci}{\ci})^{i_{N+1}}
\mbox{ with } i_{\ci}
\in{\bf \mbox{\boldmath$Z$}}~,~~
i_1,... i_N,i_{N+1}\in
{\bf \mbox{\boldmath$N$}}\cup\{0\}
\]
[$\Lm{\bu}{\ci}\Lm{\ci}{\ci}$ commutes with the elements $\Lm{a}{\ci}$.
Use
(\ref{1elle}), (\ref{2elle})  
and then (\ref{3elle}) and (\ref{4elle}),
to exlicitly  write a generic element  of $IU$ as  $\eta^ia_i$ where $a_i
\in U_{q,r}(s(N))$].

The $(${\bf Z},{\bf N}$)$ grading 
is: 
{\sl grade}$(\T{a}{b})=(0,0)$, 
{\sl grade}$(\Lm{a}{\ci})=(0,1)$, {\sl grade}$(\Lm{\ci}{\ci})=(1,0)$,
so that:
\eq
\IUsosp~=\sum_{(h,k)\in
({\mbox{\scriptsize \bf Z}},{\mbox{\scriptsize \bf N}})}
{}^{{\hskip -0.64 cm}\oplus} ~~\Ga'^{(h,k)}\label{grad2s}
\en
where $\Ga'^{(0,0)}=U_{q,r}(s(N))$ 
\eqa
&&\!\!\!\!\!\!\!\!\!\!\Ga'^{(0,1)}
=\{\Lm{a}{\ci}\varphi^{a}\;~/~~\varphi^a\in U_{q,r}(s(N))\}~~,~~~~
\Ga'^{(\pm 1,0)}=\{(\Lm{\ci}{\ci})^{\pm 1}\varphi\;~/~~\varphi\in  
U_{q,r}(s(N))\}\nonumber\\
&&\!\!\!\!\!\!\!\!\!\!\Ga'^{(h,k)}=\{{(\Lm{\ci}{\ci})}^{h}
\Lm{a_1}{\ci}\Lm{\ci}{a_2}\ldots \Lm{\ci}{a_k}\varphi^{a_1a_2...a_k}~~/~~
\varphi^{a_1a_2...a_k}
\in 
U_{q,r}(s(N))\}\nonumber
\ena
Any submodule $\Ga'^{(h,k)}$ is a
$U_{q,r}(s(N))$--bicovariant bimodule freely generated by the elements 
$\eta^i$ with degree
\sma{$(h,k)\in (\mbox{\bf{Z}},\mbox{\bf{N}})$}. 
Also the analogue of  Note 3.3.5 and Note 3.3.6 still holds for $IU$.  
\sk
\vskip .2cm
\noi {\large{{\bf{Duality}} $\IU\leftrightarrow \IS$}}
\sk
We now show that $IU$ is dually paired to $\S$. This is the
fundamental
step allowing to interpret $IU$ as the algebra of regular functionals
on $\IS$.
\sk
\noi \noi{\bf Theorem} {4.4.4} $~IU$ annihilates $H$.

\noi {\sl Proof :} $~~~$
Let ${\cal{L}}$ and $\Tc$ be generic generators of $IU$ and  $H$
respectively.
As discussed in { Note 4.4.1}, ${\cal{L}}(\Tc)=0$. A generic element
of the ideal is given by $a\Tc b$ where sum of polynomials is  
understood;
we have (using Sweedler's notation for the coproduct):
$
{\cal{L}}(a\Tc
b)={\cal{L}}_{(1)}(a){\cal{L}}_{(2)}(\Tc){\cal{L}}_{(3)}(b)=0
$
because ${\cal{L}}_{(2)}(\Tc)=0$. Indeed ${\cal{L}}_{(2)}$ is
still a generator of $IU$ since $IU$ is a sub-coalgebra of $\Usosp.$
Thus ${\cal{L}}(H)=0$. Recalling that a product of functionals
annihilating
$H$ still annihilates the co-ideal $H$, we also have $IU(H)=0$.
{\cvd}
In virtue of { Theorem 4.4.4} the following bracket is well defined:
\eqa
\!\!\!\!\!\!\!\!\!\!\!\!\!\!\!\!\!\!\mbox{{\sl Definition }$~~~~~$}
&\!\!\!\!\!\!\!\!\!\le ~ ,~ \re \; :~ &IU\otimes \IS \longrightarrow
\mbox{\boldmath $C$} \nonumber\\
&               &\le a', P(a)\re\equiv a'(a)
 ~~~\forall \/a'\in IU \,,~\forall \/a\in \S\label{dualityiso}
\ena
where $P~:~\S\rightarrow \S/H\equiv \IS$ is the canonical
projection, which is surjective.
The bracket is well
defined because two generic counterimages of $P(a)$ differ
by an addend belonging to $H$.

Note that when we use the bracket $\le ~,~\re$, $a'$ is seen as an
element
of $IU$ , while in the expression $a'(a)$, $a'$ is seen
as an element of $\Usosp$ (vanishing on $H$).
\sk
\noi {\bf Theorem} 4.4.5 $~~~$ The bracket (\ref{dualityiso})
defines a pairing between $IU$ and $\IS~\mbox{:}$ 
$\forall\/ a',b'\in IU~,~\forall\/P(a),P(b)\in \IS$
\eqa
& &\le a'b' , P(a)\re = \le a'\otimes b',\Delta
(P(a))\re\label{uunososp}\\
& &\le a',P(a)P(b)\re=\le\Delta '(a'),P(a)\otimes
P(b)\re\label{uduesosp}\\
& &\le\kp(a'),P(a)\re=\le a',\kappa(P(a))\re\label{utresosp}\\
& &\le I,P(a)\re=\epsi (P(a))~~;~~~\le a',P(I)\re=\epsi
'(a')\label{uquattro}
\ena

\noi{\sl Proof : } The proof is easy since $IU$ is a Hopf subalgebra
of $\Usosp$ and $P$ is compatible with the structures and costructures
of $\S$ and $\IS$. Indeed we have
\[
\le a',P(a)P(b)\re=\le a',P(ab)\re=a'(ab)=\Delta '(a')(a\otimes b)
=\le\Delta ' (a'), P(a)\otimes P(b)\re
\]
\[
\le a'b',P(a)\re=a'b'(a)=(a'\otimes b')\Dtwo (a)=
\le a'\otimes b',(P\otimes P)\Dtwo(a)\re=\le a'\otimes b',\Delta
(P(a))
\re
\]
\[
\le \kp(a'),P(a)\re=\kp(a')(a)=a'(\kappa_{N+2}(a))
=\le a',P(\kappa_{N+2}(a))\re=\le a' ,\kappa(P(a))\re
\]
\sk
\cvd

We now recall that $IU$ and $\IS$, besides being dually paired, 
are bicovariant algebras with the same 
graded structure (\ref{grad1s}) and (\ref{grad2s}),
and can both be obtained as a cross-product cross-coproduct construction:
$\IS\cong B{\mbox{$\times \!\rule{0.3pt}{1.1ex}\;\!\!\!\cdot\,$}}
S_{q,r}(N)$,
$IU\cong$
$B'{\mbox{$\times \!\rule{0.3pt}{1.1ex}\;\!\!\!\cdot\,$}}U_{q,r}(s(N))$.
In particular $\IS$ and $IU$ are freely generated (as modules) by $B$ and $B'$
i.e.
by the two isomorphic sets of the 
monomials in the $q$-plane plus dilatation coordinates $\Lm{\ci}{\ci},\;
\Lm{a}{\ci}$ 
and $u,\;x^a$ respectively. We then conclude that
$IU$ is the universal enveloping algebra of $\IS$:
\eq
~\IUsosp\equiv IU~.\label{IUsospdualIS}
\en
\sk
\noi{\bf Note} 4.4.5 $~$Given a $*$-structure on $\IS$, the duality
$\IS\leftrightarrow \IUsosp$
induces a $*$-structure on $\IUsosp$. If in particular the
$*$-conjugation
on $\IS$ is
found by
projecting a $*$-conjugation on $\S$, then the induced $*$ on $\IU$
is simply the
restriction to $\IUsosp$ of
the $*$ on $\Usosp$. This is the case for the $*$-structures that lead
to the real
forms
$IS_{q,r}(N, \mbox{\boldmath $R$})$ and $ISO_{q,r}(n+1, n-1)$ and in
particular to
the quantum Poincar\'{e} group.



\sect{Bicovariant calculus on $SO_{q,r}(N+2)$ and \mbox{$Sp_{q,r}(N+2)$} }

The bicovariant differential calculus on the multiparametric 
$q$-groups of the $B,C,D$ series can be formulated following
Section 2.2. We list here some formulae and comments that do not appear
in that section.
\sk
The commutations between the generators $\T{R}{S}$ and the  $1$-forms 
$\ome{A_1}{A_2}$ are explicitly given by
\eq
\ome{A_1}{A_2} \T{R}{S}= \Rinv{TB_1}{CA_1}
\Rinv{A_2C}{B_2S} \T{R}{T}\ome{B_1}{B_2} \label{commomTsosp}
\en
\noi Using (\ref{defd2}) we compute the exterior derivative
on the basis elements of $\SqrNtwo$:
\eq
d ~\T{A}{B}=\rinv [\Rinv{CR}{ET} \Rinv{TE}{SB}
\T{A}{C}-\de^R_S \T{A}{B}]
{}~\ome{R}{S} \equiv\T{A}{C}\X{CR}{BS}
 \ome{R}{S}\label{dTABsosp}
\en
where we have
\eq
\X{A_1B_1}{A_2B_2}\equiv \rinv [\Rinv{A_1 B_1}{ET} \Rinv{TE}{B_2A_2}-
\de^{B_1}_{B_2} \de^{A_1}_{A_2} ]
=zK^{A_1B_1}_{~~A_2B_2} -\Rhatinv{A_1B_1}{A_2B_2}
\label{X}
\en
with  $z \equiv \epsilon r^{N-\epsilon},
K^{A_1B_1}_{~~A_2B_2}=C^{A_1B_1}C_{A_2B_2} ~.$ [From
(\ref{extrarelation}) and (\ref{CR}), the second equality in (\ref{X}) is 
easily
proven.] Notice also that from (\ref{dTABsosp}) and (\ref{defchi}), 
$\X{A_1B_1}{A_2B_2}$ is the fundamental representation of the $q$-Lie
algebra generators $\cchi{B_1}{B_2}$: 
\[
\X{A_1B_1}{A_2B_2}=\cchi{B_1}{B_2}(\T{A_1}{A_2})
\]
Every element  $\rho$ of $\Gamma$, which by definition is
written in a unique way as
$\rho=a^{A_1}{}_{A_2}\ome{A_1}{A_2}$,  can also be written as
\eq
\rho=\sum_k a_k db_k  \label{propd}
\en
\noi for some $a_k,b_k$ belonging to $A$. This can
be proven directly by inverting
the relation (\ref{dTABsosp}). The
result  is an expression of
 the $\om$ in terms of a linear
combination of $\kappa (T) dT$, as in the classical case:
\eq
\ome{A_1}{A_2}=Y_{A_1B_1}^{~~~A_2B_2} \kappa(\T{B_1}{C}) d\T{C}{B_2}
\label{omABsosp}
\en
\noi where $Y$ satisfies
$
\X{A_1B_1}{A_2B_2}Y_{B_1C_1}^{~~~B_2C_2}=
\de^{A_1}_{C_1}\de^{C_2}_{A_2 
} ~,~
Y_{A_1B_1}^{~~~A_2B_2} \X{B_1
C_1}{B_2C_2}=\de^{C_1}_{A_1}\de^{A_2}_{C_2}
$
and is given explicitly by
\eq
Y_{A_1B_1}^{~~~A_2B_2}=
\alpha[(z-\la) C_{A_1B_1}C^{A_2B_2} + C_{A_1D}
\R{DA_2}{CB_1} C^{CB_2}
 - {\la\over z(z-z^{-1})}
D^{A_2}_{~A_1} (D^{-1})^{B_2}_{~B_1}] \label{Y}
\en
\noi with $\alpha={1 \over z(z-z^{-1}-\la)}$ and
$D^{E}_{~C} \equiv C^{EF} C_{CF}$.
The $r=1$ limit of  (\ref{dTABsosp}) is discussed in the next section.
\sk
The braiding matrix $\Lambda$ that defines the exterior product of forms is
given by (\ref{RRffMM}). It satisfies the characteristic equation:
\eq
\begin{array}{rcl}
& & (\La+r^2 I)~(\La+r^{-2}I)~(\La + \epsilon r^{\epsilon + 1-N} I)
(\La + \epsilon r^{-\epsilon - 1+N} I) \times\\
& &~~~~~~~(\La - \epsilon r^{-\epsilon + 1+N} I)
{}~(\La - \epsilon r^{\epsilon - 1-N} I)~
(\La-I)=0
\end{array} \label{Laeigen}
\en
\noi due to the characteristic equation  (\ref{cubic}).
For simplicity
we will at times use the adjoint indices $i,j,k,...$
with ${}^i={}_A^{~B},~{}_i={}^A_{~B}$.
Define
\eq
\PIJ{a_1}{a_2}{d_1}{d_2}{c_1}{c_2}{b_1}{b_2}
\equiv  d^{f_2} d^{-1}_{c_2}
\Rhat{b_1f_2}{c_2g_1} \PI{c_1g_1}{a_1e_1}
\Rhatinv{a_2e_1}{d_1g_2} \PJ{d_2g_2}{b_2f_2} \label{PIJ}
\en
\noi where  $P_I=P_S,P_A,P_0$ are given in
(\ref{projBCD}) and $d^{f_2}\equiv D^{f_2}_{~f_2}\,,\;
d^{-1}_{c_2}\equiv (D^{-1})^{c_2}_{~c_2}$.
The $(P_I,P_J)$ are themselves
projectors, i.e.:
\eq
(P_I,P_J) (P_K,P_L) =\de_{IK} \de_{JL} (P_I,P_J) \label{projIJ}
\en
Moreover
\eq
(I,I)=I \label{IIeqI}
\en
{}From $\om^i \we \om^j\equiv \om^i \otimes \om^j - \Lambda^{ij}_{~kl}
\om^k \otimes \om^l$ we find
\eq
\om^i \we \om^j = - \Z{ij}{kl} \om^k \we \om^l \label{commomsosp}
\en
\noi with
\eq
Z=(P_S,P_S)+(P_A,P_A)+(P_0,P_0)-I
\en
\noi see ref. \cite{Monteiro}.
The inverse of $\Lambda$ always
exists, and is given by
\eq
\begin{array}{rl}
\!\!\!\RRhatinv{A_1}{A_2}{D_1}{D_2}{B_1}{B_2}{C_1}{C_2}=&
\ff{D_1}{D_2B_1}{B_2}
(\T{A_2}{C_2} \km (\T{C_1}{A_1})) = \\
=&\R{F_1B_1}{A_1G_1} \Rinv{A_2G_1}{E_2D_1} \Rinv{D_2E_2}{G_2
C_2} \R{G_2C_1}{B_2F_1} (d^{-1})^{C_1} d_{F_1}
\end{array}
\label{Lambdainv}
\en
\noi Note that for $r=1$, $\La^2 = I$ and $(\La +I)(\La-I)=0$
replaces the
seventh-order spectral equation (\ref{Laeigen}). In this special
case, one
finds
the simple formula:
\eq
\om^i \we \om^j= - \L{ij}{kl} \om^k \we \om^l
\label{omcomrone}~~~~\mbox{ i.e.}
{}~~Z=\Lambda~.
\en
\sk
The {\bf \mbox{\boldmath $q$}
-Cartan-Maurer equations} are given by:
\eq
d\ome{C_1}{C_2}= \rinv (\ome{B}{B} \we \ome{C_1}{C_2} +
 \ome{C_1}{C_2} \we
\ome{B}{B}) \equiv
-\onehalf \cc{A_1}{A_2}{B_1}{B_2}{C_1}{C_2}
{}~\ome{A_1}{A_2} \we \ome{B_1}{B_2} \label{CartanMaurersosp}
\en
\noi with:
\eq
\cc{A_1}{A_2}{B_1}{B_2}{C_1}{C_2}= -{2\over (r-r^{-1})}
[ \ZZ{B}{B}{C_1}{C_2}{A_1}{A_2}{B_1}{B_2} + \de^{A_1}_{C_1}
\de^{C_2}_{A_2} \de^{B_1}_{B_2} ]\label{cc}
\en
To derive this formula we have used the flip operator $Z$
on $\ome{B}{B} \we \ome{C_1}{C_2}$.
\sk

Finally, we recall that the $\chi$ operators close
on the {\bf $q$-Lie algebra} :
\eq
\chi_i \chi_j - \L{kl}{ij} \chi_k \chi_l = \C{ij}{k} \chi_k
\label{bico1sosp}
\en
\noi where the $q$-structure constants are given by
\eq
{\C{jk}{i}=\chi_k(\M{j}{i})} ~~\mbox{ explicitly : }~~
{\CC{A_1}{A_2}{B_1}{B_2}{C_1}{C_2}} =\rinv [- \de^{B_1}_{B_2}
\de^{A_1}_{C_1} \de^{C_2}_{A_2} +
\LL{B}{B}{C_1}{C_2}{A_1}{A_2}{B_1}{B_2}]. \label{CCsosp}
\en
The $C$ structure constants
appearing in the Cartan-Maurer equations are  related
to the $\Cb$ constants of the $q$-Lie algebra by:
\eq
\C{jk}{i}=\onehalf [\c{jk}{i}-\L{rs}{jk} \c{rs}{i}]~.
\en
\noi In the particular case $\Lambda^2=I$ (i.e. for $r=1$) it is not
difficult to
see that in fact
$C= \Cb$, and that
the $q$-structure constants are $\Lambda$-antisymmetric:
\eq
\C{jk}{i}= - \L{rs}{jk} \C{rs}{i}. \label{antisymC}
\en
\sk
\noi{\bf Note} 4.5.1 $~$ The formulae characterizing the bicovariant
calculus have
been written in the
 basis $\{\cchi{A}{B}\},~ \{\ome{C}{D}\}$  because of the
particularly simple
expression of the
$\ff{A}{BC}{D}$ and $\cchi{A}{B}$ functionals
in terms of $L^{\pm A}{}_B$, see (\ref{defff}) and (\ref{defchi2}).
Obviously the calculus is independent from the basis chosen. If we
consider the
linear transformation
\[
\om^i\rightarrow \om '^i=X^i{}_j\om^j
\]
(where we use adjoint indices ${~}^i={}_{A_1}{}^{A_2} ,
{}_j={}^{B_1}{}_{B_2}$), from the exterior differential
\eq
da=(\chi_i*a)\om^i=(\chi'_i*a)\om '^i
\en
we find
\[
\chi_i\rightarrow\chi'_i=\chi_j\,(X^{-1}){}^j{}_i~,
\]
and from the coproduct rule (\ref{copf}) of the $\chi_i$ we find
$f^{i}{}_{j}\rightarrow f'^i{}_j=X^i{}_lf^l{}_m\,(X^{-1}){}^m{}_j ;$
while from
(\ref{adjoint})
we have $M_i{}^j
\rightarrow M'_i{}^j=(X^{-1}){}^l{}_iM_l{}^mX^j{}_m.$

A useful change of basis is obtained via the following
transformation:
\eq
\begin{array}{rcl}
&
&\ome{A_1}{A_2}\rightarrow\vt^{A_1}_{~A_2}=
\X{A_1B_1}{A_2B_2}\ome{B_1}{B_2}\\
& &\cchi{A_1}{A_2}\rightarrow
\psi_{A_1}^{~A_2}=\cchi{B_1}{B_2}Y_{B_1A_1}^{~~B_2A_2}
\end{array}
\en
where $X$ and its (second) inverse $Y$ are 
defined in (\ref{X}) and
(\ref{Y}).
Using (\ref{omABsosp}) it is immediate to see that
\eq
\vt^{A_1}_{~A_2}=\kappa(\T{A_1}{C})d\T{C}{A_2}~.
\en
We also have:
\eq
\psi_{A_1}^{~A_2}({T}^{B_1}_{~B_2})=
\psi_{A_1}^{~A_2}({\widetilde
T}^{B_1}_{~B_2})=\de_{A_1}^{B_1}\de^{A_2}_{B_2}~~
\mbox{ where }~~ {\widetilde T}^{B_1}_{~B_2}\equiv
\T{B_1}{B_2}-\de^{B_1}_{B_2}I~.\label{dualchix}
\en
Formula (\ref{dualchix}) follows from $\psi_{A_1}^{~A_2}(I)=0$ and:
\eq
\begin{array}{rcl}
\vt^{A_1}_{~A_2}&=&\kappa(\T{A_1}{C})d\T{C}{A_2}=\kappa(\T{A_1}{C})
(\psi_{B_1}^{~B_2}*\T{C}{A_2})\vt^{B_1}_{~B_2}\\
&=&\kappa(\T{A_1}{C})\T{C}{D}
\psi_{B_1}^{~B_2}(\T{D}{A_2})\vt^{B_1}_{~B_2}=\psi_{B_1}^{~B_2}
(\T{A_1}{A_2})
\vt^{B_1}_{~B_2}~.
\end{array}
\en
The analogue of the coordinates ${\widetilde T}^{B_1}_{~B_2}$ in the
old basis
is given by
\eq
x_{B_1}^{~B_2}\equiv Y_{B_1C_1}^{~~B_2C_2}
{\widetilde T}^{C_1}_{~C_2}~~,~~~\cchi{A_1}{A_2}
(x_{B_1}^{~B_2})=\de^{A_1}_{B_1}\de^{B_2}_{A_2}~.
\en
The set of coordinates $x_{B_1}^{~B_2}$ and ${\widetilde T}^{C_1}_{~C_2}$
span the space $X$ described in (\ref{axri}) and dual to the $q$-Lie algebra 
of $\S$. 
\sk
\noi {\bf Conjugation}
\sk
From the $*$-structures (\ref{conjL1}), (\ref{conjL2})  and the definition  
(\ref{defchi2}) it is straightforward to find
how the $*$-conjugation acts
on the tangent vectors $\chi$.
Both conjugations (\ref{conjL1})
and (\ref{conjL2})
are  compatibile
with the differential calculus. Indeed they 
 respectively yield [use (\ref{defchi2}),
(\ref{coiLpmsosp}),
(\ref{RLLsosp}), (\ref{crc1}), (\ref{crc2}) and ({\ref{CR}) with
\sma{$N\rightarrow N+2$} since we have capital indices]:
\eq
(\cchi{A}{B})^*=- r^{-N-1} \cchi{C}{D} \Dcal{F}{~B}
 \Dcal{A}{~G} \R{EG}{FC} D^D_{~E}~~~
\sma{\mbox{ for } $SO_{q,r}(n+2,n ;\Rbo)~~;~~~  2n+2=N+2$} \label{conjchiAB}
\en
\[
(\cchi{A}{B})^*=- \epsilon r^{\epsilon-(N+2)} \cchi{C}{D}
\R{EA}{BC} D^D_{~E}~~~
\sma{\mbox{ for } $SO_{q,r}(n+1,n+1 ;\Rbo)~\mbox{ or }~ 
Sp_{q,r}(n+1,n+1 ;\Rbo)$} 
\]
\noi with $D^E_{~C} \equiv C^{EF} C_{CF}$.  
As for the $L$ matrices 
(and similarly to the $T$ matrices) we have 
$\kappa^2(\cchi{A}{B})={D^{-1}}^A_E\cchi{E}{F}D^F_B$.
\sk
In a basis $\{\cchi{A}{B}\}$ relation (\ref{om*chi}) 
reads
\eq
(\cchi{A}{B})^*= V^{A~~Q}_{~BP}\cchi{P}{Q}~~~\mbox{ if and only if }~~~
\ome{C}{D}=-\ome{E}{F}\,{\overline {V^{E~~D}_{~FC}}}
\label{explicom*}
\en
where $V$ is a matrix with complex entries and 
${\overline V}$ is its complex conjugate . Using this espression 
[or the inversion formulae (\ref{omABsosp})] one finds
the induced conjugation on the
left invariant 1-forms (use $\overline{D^A_{~B}}=D^{-1\, A}{}_B$):
\eq
(\ome{A}{B})^*=  r^{N+1}\Dcal{F}{D} \Dcal{C}{G} {D^{-1}}^B{}_{E}
{R^{-1}}^{EG}_{~FA}\ome{C}{D}~~~~
\sma{\mbox{ for } $SO_{q,r}(n+2,n ;\Rbo)~~;~~~  2n+2=N+2$}
\label{conjomAB}
\en
\[
(\ome{A}{B})^*= \epsilon r^{N+2-\epsilon}  {D^{-1}}^B{}_{E}
{R^{-1}}^{EC}_{~DA}\ome{C}{D}~~~\sma{\mbox{ for } $SO_{q,r}(n+1,n+1 ;\Rbo)~
\mbox{ or }~ Sp_{q,r}(n+1,n+1 ;\Rbo)$}.
\]

\sect{Differential calculus on $SO_{q,r=1}(N+2)$ 
and \\
{}\mbox{$Sp_{q,r=1}(N+2)$}}

As discussed in Section 4.2,  we have obtained the quantum
inhomogeneous groups $\ISqrN$ via the projection
{\vskip -0.5cm}
\eq
P~:~~~\SqrNtwo -\!\!\!-\!\!\!\!\!\longrightarrow {\SqrNtwo \over
H}=\ISqrN
\label{Pproj}
\en
{\vskip -0.1cm}
\noi with $H$=Hopf ideal in $\SqrNtwo$ defined after (\ref{quotientsosp}).  As a
consequence, the universal enveloping algebra
$\UISqrN$ is a Hopf subalgebra of $\USqrNtwo$
and contains all the functionals that annihilate $H=Ker(P)$.

Let us consider now the $\chi$ functionals in the differential
calculus
on $\SqrNtwo$. Decomposing the indices we find:
{\vskip -0.4cm}
\eqa
& & \cchi{a}{b}=\rinv [\ff{c}{c a}{b}-\de^a_b \epsi]~~~~~~~+\rinv
\ff{\bu}{\bu
a}{b} \label{chiN2}\\
& & \cchi{a}{\ci}=\rinv \ff{c}{c a}{\ci} ~~~~~~~~~ ~~~~~~~ +\rinv
\ff{\bu}{\bu
a}{\ci} \\
& &\cchi{\ci}{b}=~~~~~~~~~~~~~~~~~~~~~~~~~~~~~~~~~~ + \rinv
 [\ff{c}{c \ci}{b}+\ff{\bu}{\bu \ci}{b}]  \\
& & \cchi{a}{\bu}=~~~~~~~~~~~~~~~~~~~~~~~~~~~~~~~~~~+ \rinv
\ff{\bu}{\bu a}{\bu}\\
& & \cchi{\bu}{b}=\rinv \ff{\bu}{\bu\bu}{b} ~~~~~~~~~~~~~\\
& &\cchi{\ci}{\ci}=\rinv [\ff{\ci}{\ci\ci}{\ci}-\epsi]~~~~~~~~~~ +
\rinv
[\ff{c}{c\ci}{\ci}
+\ff{\bu}{\bu \ci}{\ci} ] \\
& &\cchi{\ci}{\bu}=~~~~~~~~~~~~~~~~~~~~~~~~~~~~~~~~~~ + \rinv
\ff{\bu}{\bu \ci}{\bu}\\
& & \cchi{\bu}{\ci}=\rinv \ff{\bu}{\bu \bu}{\ci}~~~~~~~~~~~~~~~\\
& & \cchi{\bu}{\bu} = \rinv [\ff{\bu}{\bu\bu}{\bu} -
\epsi]~~~~~\label{chiN2end}\\
& &~~~~~~~~\underbrace{~~~~~~~~~~~~~~~~~~~~~~~~~~~~~~~~~}_
{\hbox{terms annihilating $H$}} \nonumber
\ena
\noi where we have indicated the terms that do and do not
annihilate the Hopf ideal $H$, i.e. that belong or do
not belong to $\IUsosp$.
We see that only the functionals
$\cchi{\bu}{b}$,  $\cchi{\bu}{\ci} $ and $\cchi{\bu}{\bu}$ do
annihilate
$H$, and therefore belong to $\UISqrN$. The resulting bicovariant
differential calculus, see Chapter 5,  contains dilatations and
translations,
but does not contain the tangent vectors of $\SqrN$, i.e.
the functionals $\cchi{a}{b}$.  Indeed these contain
 $\ff{\bu}{\bu a}{b}$, in general not vanishing on $H$. We can,
however, try to find restrictions on the parameters
$q,r$ such that  $\ff{\bu}{\bu a}{b} (H)=0$. As we will see, this
happens
 for $r=1$. For this reason we consider in the following
the particular multiparametric deformations called ``minimal
deformations" or twistings, corresponding to $r=1$.
\sk
We first examine what happens to the
bicovariant calculus on $\SqrNtwo$ in the $r=1$
limit\footnote{By $\lim_{\rone}a$,
where the generic element $a\in\SqrNtwo$  is a polynomial in the
matrix
elements $\T{A}{B}$ with complex coefficients $f(r)$ depending on
$r$,
we
understand the
element of $S_{q,r=1}(N+2)$ with coefficients given by
$\lim_{\rone}f(r)$. The
expression
$\lim_{\rone}\phi=\varphi$, where $\phi\in \USqrNtwo$  and
$\varphi\in
U(S_{q,r=1}(N+2))$ means that
$\lim_{\rone}\phi(a)=\varphi(\lim_{\rone}a)$ for any $a\in \SqrNtwo$
such
that
$\lim_{\rone}a$ exists.
Finally, the left invariant $1$-forms $\omega^i$ are symbols, and
therefore
$\lim_{\rone}a_i\omega^i\equiv(\lim_{\rone}a_i)\omega^i$.}.
The $R$ matrix is given by, cf. (\ref{Rnonzerososp}):
\eqa
& & \R{AB}{AB}=q^{-1}_{AB} + O(\la) \\
& & \R{AB}{BA}=\la~~~~~~~~~~~~~~~~~~~~~~~~
\mbox{\footnotesize $A>B, A' \not= B$}\\
& & \R{AA'}{A'A}=\la~ (1-\epsilon r^{\rho_A - \rho_{A'}})~~~~~
\mbox{\footnotesize  $A>A'$} \\
& & \R{AA'}{BB'}=-\la \epsilon_A \epsilon_B + O(\la^2)~~~~
\mbox{\footnotesize  $A>B, A' \not= B$}
\ena
\noi where $O(\la^n)$ indicates an infinitesimal of order
$\geq\la^n$;
the $q_{AB}$ parameters satisfy:
\eq
q_{AB}=q^{-1}_{AB'} = q^{-1}_{A'B} = q^{-1}_{BA}~~;~~~~ q_{AA} =
q_{AA'} = 1
\en
\noi up to order $O(\la)$.  Note that the components $\R{AA'}{A'A}$
are of order $O(\la^2)$ for the orthogonal case ($\epsilon=1$) and
of order $O(\la)$ for the symplectic case ($\epsilon=-1$).
The $RTT$ relations simply become:
\eq
\T{B_1}{A_1} \T{B_2}{A_2}={q_{B_1B_2} \over q_{A_1A_2}}\T{B_2}{A_2}
\T{B_1}{A_1}~.\label{4.6.16}
\en
For $r=1$ the metric is $C_{AB}=\epsilon_A\delta_{AB'}$
and therefore we have
$C_{AB}=\epsilon C_{BA}$.
Using the
definition
(\ref{LonT}), it is easy to see that
\eqa
{}\!\!\!\!\!{}\!\!\!\!\!&{}\!\!\!\!\!{}\!\!\!\!\! &
\Lpm{A}{A} (\T{C}{D})=\de^C_D q_{AC} +O(\la)\\
{}\!\!\!\!\!{}\!\!\!\!\!&{}\!\!\!\!\!{}\!\!\!\!\! &
\Lpm{A}{B} (\T{B}{A})=\pm \la {}~~~~~~~\mbox{\footnotesize $A\not=B,
A'\not= B$;~$A<B$ for $L^+$, $A>B$ for $L^-$}\label{pmlanda1}\\
{}\!\!\!\!\!{}\!\!\!\!\!&{}\!\!\!\!\!{}\!\!\!\!\! &
\Lpm{A}{A'} (\T{A'}{A})=\pm \la~[1-\epsilon r^{\pm
(\rho_A-\rho_{A'})}]
{}~~~\mbox{\footnotesize $A<A'$ for $L^+$, $A>A'$ for $L^-$}\label{spmeno1}\\
{}\!\!\!\!\!{}\!\!\!\!\!&{}\!\!\!\!\!{}\!\!\!\!\! &
\Lpm{A}{B} (\T{A'}{B'}) = \mp \la \epsilon_A \epsilon_B + O(\la^2)
{}~~\mbox{\footnotesize $A\not=B,
A'\not= B$;~~$A<B$ for $L^+$, $A>B$ for $L^-$}\label{pmlanda2}
\ena
\noi all other $L^{\pm}(T)$ vanishing. Relations (\ref{pmlanda1})
and (\ref{pmlanda2}) imply that for any generator $\T{C}{D}$ we have
$
\Lpm{A}{B}(\T{C}{D}) = -\epsilon_A \epsilon_B \Lpm{B'}{A'}(\T{C}{D})
+
O(\la^2)~
\mbox{ with }\mbox{\footnotesize ~~ $A\not=B$, $A\not= B'$}~.$\\
In general, since
$$\Delta (\Lpm{A}{A})=\Lpm{A}{A} \otimes \Lpm{A}{A} ~~;~~~
\Delta(\Lpm{A}{B})=\Lpm{A}{A} \otimes \Lpm{A}{B}
+ \Lpm{A}{B} \otimes \Lpm{B}{B} + O(\la{}^2), ~~~\mbox{\footnotesize
$A\not=
B$}
$$
\noi we find that
\eqa
& &\Lpm{A}{A}=O(1)\\
& &\Lpm{A}{B}= O(\la),~~~~~\mbox{\footnotesize $A\not= B$,~$A\not=
B'$}\\
& &\Lpm{A}{A'}=O(\la^2) \mbox{~for $SO_q$},~~~~O(\la) \mbox{~for
$Sp_q$}
\ena
where, by definition, $\phi =O(\la{}^n)$ ($\phi$ being a functional)
means that
for any element $a\in \SqrNtwo$ with well-defined classical limit, we
have $\phi(a)=O(\la{}^n)$.

\noi Moreover the following relations hold:
\eqa
& &\Lpm{A}{A} = \Lmp{A}{A} + O(\la)~,\\
& &\kappa (\Lpm{A}{B})=\epsilon_A\epsilon_B\Lpm{B'}{A'}+O(\la) ~
\mbox{ and therefore,}~ \kappa^2=id + O(\la) ~.
\ena
Similarly one can prove the relations involving the
$f$ functionals (no sum on repeated indices):
\eqa
& &\ff{A}{AA}{A}=\epsi+O(\la)\label{fpropepsi}\\
& &\ff{A}{BA}{B}=O(1)~~~~ \mbox{ and
}~\ff{A}{BA}{B}=\ff{B'}{A'B'}{A'}+O(\la)
\label{relationf}\\
& &\ff{C}{CA}{A}=O(\la^2)~~~~\mbox{\footnotesize $C\not=
A$}\label{relationfCCAA}\\
& &\ff{C}{CA}{B}=O(\la^2)~~~~\mbox{\footnotesize [$A<B,C\not= B$] or
[$A>B,C\not= A$]}\label{6.30}
\ena
[hint: check (\ref{relationf})-(\ref{6.30}) first on the generators,
then use the coproduct in (\ref{copf})].
\noi From the last relation we deduce
\eqa
& &\cchi{A}{B}= \linv \ff{B}{BA}{B}+O(\la),~~~~~\mbox{\footnotesize
$A<B$}\\
& &\cchi{A}{B}= \linv \ff{A}{AA}{B}+O(\la),~~~~~\mbox{\footnotesize
$A>B$}
\ena
\noi and from (\ref{fpropepsi}) and (\ref{relationfCCAA}) one has
\eq
\cchi{A}{A}=\linv [\ff{A}{AA}{A}-\epsi]
\en
\noi Next one can verify that
\eq
\left. \begin{array}{ll} &\cchi{A}{B}(\T{B}{A})=-q_{BA}+O(\la)\\
                    & \cchi{A}{B}(\T{A'}{B'})=\epsilon_A
\epsilon_B+O(\la)\\
                    & \cchi{A}{B}(\T{C}{D})=0~ \mbox{ otherwise}
         \end{array}
                                   {}~~~~~~~~ \right\}
{}~\mbox{\footnotesize $A\not= B$, $A\not= B'$} \label{relationschiT}
\en
\eq
\forall
{}~\T{C}{D}\;,~~~\cchi{A}{A}(\T{C}{D})
=-\cchi{A'}{A'}(\T{C}{D})+O(\la) . 
\label{relationchiAAT}
\en
Eq.s (\ref{relationschiT}) yield the relation between $\chi$
functionals:

\eq
\forall ~\T{C}{D}\;,~~~\cchi{B'}{A'}(\T{C}{D})=-{\epsilon_A
\epsilon_B \over
q_{BA}}
\cchi{A}{B}(\T{C}{D})+O(\la),~
\sma{$A \not= B$, $A \not= B'.$}
\label{relationchiABT}
\en

It is not difficult to prove that the coproduct rule in (\ref{copf})
is compatible with
(\ref{relationchiABT}) and (\ref{relationchiAAT})
 making them valid on arbitrary polynomials in the $\T{A}{B}$
elements:
\eq
\cchi{B'}{A'}=-{\epsilon_A \epsilon_B \over q_{BA}}
\cchi{A}{B}+O(\la),~
\sma{$A \not= B$, $A \not= B'$}~~~;~~~~
\cchi{A}{A}=-\cchi{A'}{A'}+O(\la)~.\label{relationchiAA}
\en
Finally:
\eq
\cchi{A}{A'}=O(\la) \mbox{~for $SO_q$}~,~~~~O(1) \mbox{~for $Sp_q$}~,
 ~~~\sma{$A \not= A'$}. \label{relationchiAAp}
\en
Summarizing, in the $r \rightarrow 1$ limit, only the following
$\chi$
functionals survive:
\eqa
& &\cchi{A}{A}\equiv\limrone \linv~
[\ff{A}{AA}{A}-\epsi]\label{limuno}\\
& &\cchi{A}{B}\equiv\limrone \linv
{}~\ff{A}{AA}{B},~~~\sma{$A>B,A\not=B'$}\\
& &\cchi{A}{B}\equiv\limrone \linv ~
\ff{B}{BA}{B},~~~\sma{$A<B,A\not=B'$}\\
& &\cchi{A}{A'}\equiv\limrone \linv ~\sum_C \ff{C}{CA}{A'}=0
{}~~\sma{for
$SO_q$},
{}~~\not= 0~~\sma{for  $Sp_q$}\label{limquattro}
\ena
\noi Notice that  (\ref{relationchiAA}) and
 (\ref{relationchiAAp}) are all contained in the formula:
\eq
\cchi{B'}{A'}=-{\epsilon_A \epsilon_B \over q_{BA}}
\cchi{A}{B}+O(\la)\label{relationcchi}
\en
thus in the $\rone$ limit there are $\sma{$(N+2)(N+1)/2$}$ tangent
vectors for
$SO_q(N+2)$
and $\sma{$(N+2)(N+3)/2$}$ tangent vectors for $Sp_q(N+2)$, exactly
as in the
classical case.
\sk
The $r=1$ limit of  (\ref{dTABsosp}) reads:
\eq
d \T{A}{B}= -\sum_{C}\T{A}{C}  q_{CB} (\ome{B}{C}
- \epsilon_B\epsilon_C q_{BC}  \ome{C'}{B'})~, \label{dTABrone}
\en
and therefore, for $r=1$,  $\om$ appears only in the combination
\eq
\Ome{A}{B} \equiv \ome{A}{B} - \epsilon_A\epsilon_B q_{AB}
\ome{B'}{A'}~,
\label{defOme}
\en
Only $\sma{$(N+2)(N+1)/2$}$ [$\sma{$(N+2)(N+3)/2$}$ for $Sp_q(N+2)$]
of the $\sma{$(N+2)^2$}$ one forms $\Ome{A}{B}$
are linearly independent because [compare with (\ref{relationcchi})]:
\eq
\Ome{B'}{A'}=-{\epsilon_A\epsilon_B\over
q_{AB}}\Ome{A}{B}\label{relationOme}~.
\en
In the sequel, instead of considering the left module of $1$-forms
freely
generated by $\ome{A}{B}$, we consider the submodule $\Gamma$
freely generated by $\Ome{A}{B}$ with \sma{$A'<B$} for $SO_q$ and
 \sma{$A'\leq B$} for $Sp_q$.
In fact only this submodule will be relevant for the $r=1$
differential
calculus.
As in the classical case
 \footnote{
To make closer contact with the classical case one may define:
\[
\Omega^{AB} \equiv\Ome{C}{B} C^{CA}
=\epsilon_A\Ome{A'}{B}~~;~~~
\chil{AB} \equiv C_{AC} \cchi{C}{B} =\epsilon_A \cchi{A'}{B}
{}~,\nonumber
\]
and  retrieve  the more familiar  $q$-antisymmetry:
\[
\Omega^{AB}=-\epsilon q_{BA} \Omega^{BA}~~;~~~
\chil{AB}= - \epsilon q_{AB} \chil{BA}~.\nonumber
\] } , in order to simplify notations in sums
we often use $\cchi{A}{B}$ and $\Ome{A}{B}$
without the restriction \sma{$A' \leq B$} see for ex. (\ref{defd3})
below.
The bimodule structure on $\Gamma$, see {Theorem 4.5.1}, is given
by
the $\rone$ limit of the $\f{i}{j}$ functionals. These are diagonal
in the
$i,j$
indices [i.e. they vanish for $i \not= j$, see
(\ref{fpropepsi})-(\ref{6.30})] and still satisfy the property
(\ref{propf1}).
We have:
\eq
\begin{array}{rcl}
\Ome{A}{B}\,a &=&
(\ome{A}{B}-\epsilon_A\epsilon_B
q_{AB}\ome{B'}{A'})a\\
&=&(\ff{A}{BC}{D}*a)\ome{C}{D}-
\epsilon_A\epsilon_B q_{AB}(\ff{B'}{A'D'}{C'}*a)\ome{D'}{C'}\\
&=&(\ff{A}{BA}{B}*a)\Ome{A}{B}~
\label{bimOme}
\end{array}
\en
where in the last equality we have  used (\ref{relationf}) and no sum
is
understood.
We see that the bimodule structure is very simple since it does not
mix
different
$\Omega$'s. Moreover, relation (\ref{dTABrone}) is invertible and
yields:
\eq
\Ome{A}{B}=-q_{AB} \kappa (\T{B}{C}) d \T{C}{A} \label{Ome}~;
\en
in the limit $q_{AB}=1$, the $ \Ome{A}{B}$ are to be identified with
the
classical
$1$-forms, and indeed for $q_{AB}=1$ eq. (\ref{Ome}) reproduces
the correct classical limit $\Omega=-g^{-1}dg$ for the left-invariant
$1$-forms on the group manifold.
\sk
The bimodule commutation rule (\ref{bimOme}) yields a formula
similar to (\ref{commomTsosp}), replacing the values of the $R$ matrix
for $r=1$ we find the commutations:
\eq
 \Ome{A_1}{A_2} \T{R}{S}={q_{A_2 S} \over q_{A_1 S}} \T{R}{S}
 \Ome{A_1}{A_2} \label{commomTRone}
\en

For $r=1$ the coproduct on the $\chi$ functionals
reads
\eq
\Delta'(\cchi{A}{B})=\cchi{A}{B}\otimes\ff{A}{BA}{B} +
\epsi\otimes \cchi{A}{B} ~~~~~~~\mbox{ no sum on repeated
indices.}\label{copchiAB}
\en
We then consider the $r=1$ limit of (\ref{defdiff}): $da=(\chi_i*a)\om^i$ 
and therefore obtain  the following definition of the exterior differential:
\eq
da\equiv{1\over 2}(\cchi{A}{B}*a)\Ome{A}{B}=
\sum_{A'\leq B}(\cchi{A}{B}*a)\Ome{A}{B},~~~~~~\forall \/a\in A~,
\label{defd3}
\en
where in the second expression we have
used the basis of linear independent tangent
vectors  $\{\cchi{A}{B}\}_{A'\leq B}$
and dual $1$-forms $\{\Ome{A}{B}\}_{A'\leq B}$
(notice that in the $SO_q$ case we have
$\sma{$A'< B$}$ because $\cchi{A}{A}=\Ome{A}{A}=0$). The Leibniz
 rule is satisfied
for $d$ defined in (\ref{defd3}) because of (\ref{copchiAB}) and
(\ref{bimOme}).
Moreover any $\rho=a^A{}_B\Ome{A}{B}\in\Gamma$
can be written as $\rho=\sum_k a_kdb_k$, [use (\ref{Ome})].

We now introduce a left and a right action on the bimodule $\Gamma$
of
$1$-forms:
\eqa
&&\DL (a\Ome{A}{B})\equiv\D(a)I\otimes\Ome{A}{B}~,\label{defDL}\\
&&\DR (a\Ome{A}{B})\equiv\D(a)(\Ome{C}{D}\otimes
\MM{C}{DA}{B})~.\label{defDR3}
\ena
where $\MM{C}{DA}{B}=\T{C}{A}\kappa(\T{B}{D})$. [Using (\ref{defOme})
one can check that this is the $r=1$ limit of
$\DL(a \ome{A_1}{A_2})
=$ $\D(a)\,(I\otimes \ome{A_1}{A_2})$ and
$\DR(a \ome{A_1}{A_2})=$ $\D(a)\,(\ome{B_1}{B_2}\otimes
\MM{B_1}{B_2A_1}{A_2})$]. Relation (\ref{defDR3}) is well defined i.e.
$ \DR(\Ome{B'}{A'})=\DR(-{\epsilon_A\epsilon_B\over
q_{AB}}\Ome{A}{B})$ because
$\epsilon_F\epsilon_Eq_{FE}\MM{F'}{E'A}{B}$ 
$=
\epsilon_A\epsilon_Bq_{AB}\MM{E}{FB'}{A'}$.
Since in the $r=1$ case the bicovariant bimodule conditions
(\ref{propf1}),
(\ref{copM}) and (\ref{propM}) are still satisfied, it is easy to
deduce
 that $\DL$ and $\DR$ give
a bicovariant bimodule structure to $\Gamma$.
\sk
The differential (\ref{defd3}) gives a bicovariant differential
calculus if it
is
compatible with $\DL$ and $\DR$, i.e. if:
\eqa
&&\DL (adb)=\D(a)(id\otimes d)\D(b)~,\label{propdaRunoa}\\
&&\DR (adb)=\D(a)(d\otimes id)\D(b)~.\label{propdaRuno}
\ena
The proof of the compatibility of $d$ with $\DL$ is straightforward,
just use
(\ref{defd3}) and the coassociativity of the coproduct $\D$. In order
to prove
(\ref{propdaRuno}) we recall, from Proposition 2.3.1,
that  in the case $r\not= 1$ property (\ref{propdaRuno}) holds if and only if
\eq
b_1\om^j\otimes b_2 \chi_i
(b_3)\M{j}{i}=b_1\om^j\otimes\chi_j(b_2)b_3\label{OmechiM}
\en
and this last relation is equivalent to
\eq
b_1\chi_i( b_2)\M{j}{i}=\chi_j(b_1)b_2~,
{}~~\mbox{ i.e. }~~~\chi_i*b=(b*\chi_j)\kappa(\M{i}{j})
\label{leftrightchi}
\en
as one can verify by applying $m(\kappa \otimes id) \DL \otimes id $
($m$ denotes multiplication) to (\ref{OmechiM}), and using the linear
independence of the $\om^i$. Now  formula (\ref{leftrightchi}) holds
also in the limit $r=1$.  Indeed if we
consider $b$ to be a polynomial in the $\T{A}{B}$ with well behaved
coefficents
in the
$\rone$ limit, then
$\lim_{\rone}[b_1\chi_i(b_2)\M{j}{i}]=\lim_{\rone}[\chi_j(b_1)b_2]$
i.e.
$b_1[\lim_{\rone}\chi_i(b_2)]\M{j}{i}=
[\lim_{\rone}\chi_j(b_1)]b_2$ so that relation (\ref{leftrightchi})
remains valid  for $r=1$, cf .(\ref{limuno})-(\ref{limquattro}).
At this point one can prove (\ref{propdaRuno}) in the $r=1$ case
simply by substituting $\Omega$ to $\om$ in (\ref{p1}), (\ref{p2})
and (\ref{OmechiM}). Since (\ref{leftrightchi}) holds for $r=1$, then
also (\ref{OmechiM}) holds in this limit and the theorem is proved.
\sk
\cvd
We conclude that (\ref{defd3}) defines
a bicovariant differential calculus on $S_q(N+2)$.
\sk
\sk
\noi{\bf Note} 4.6.1 $~$  We have found the $S_{q,r=1}(N)$ differential
calculus studying
the $r=1$ limit of the $\chi$ functionals and of the 
bicovariant bimodule of $1$-forms
[see (\ref{defDL}), (\ref{defDR3}), (\ref{propdaRunoa}),
(\ref{propdaRuno})]. 
This has given a comprehensive analysis of the 
$r\rightarrow 1$ limit. The classical limit $r\rightarrow 1$, 
$q\rightarrow 1$ and the 
classical differential calculus are  now easily
recovered. 
From a slightly different perspective, since the
calculus can be defined from the $q$-Lie algebra alone, we could just have 
studied only the limit of the $q$-Lie algebra. It is immediate to see 
that (\ref{ZERO}), (\ref{UNO}), (\ref{DUE}) or (\ref{adbico})
still hold. This is another proof of the bicovariance of the 
$S_{q,r=1}(N)$ calculus. 

We have chosen to study the $\rone$ 
limit of the quantum Lie algebra in the $\chi$ basis because this gives 
the classical Lie algebra. 
Another possiblity is to perform the limit in the $\psi$
basis. In this case the $\psi$ are linearly independent also
when $r=1$, see (\ref{dualchix}) and the classical differential calculus is
contained in this calculus.
\sk
\noi{\bf Note} 4.6.2 $~$
In (\ref{leftrightchi}) the sum on
the indices $j=(\mbox{\sma{$C,D$}})$ can be restricted to \sma{$C'
\leq D$},
thus using the basis $\{\cchi{C}{D}\}_{C'\leq D}$,
provided             
one replaces $M$ by  
\eq                   
\begin{array}{lll}  
\MMc{C}{DA}{B} \equiv\MM{C}{DA}{B}-\epsilon_A \epsilon_B
q_{AB}            
\MM{C}{DB'}{A'} ~& &~\sma{$\mbox{ for } ~C'
\not= D ,~A'\not= B$} \\
\MMc{C}{C' A}{B} \equiv 0 ~&
\!\!\!\!\!\!\!\!\!\!\!\!\!\!\!\!\!\!\!
\!\!\!\!\!\!\!\!\!\!\!\!\!\!\!\! 
\!\!\!\!\!\!\!\!\!\!\!\!\!\!\!\!\!\!\!
,~~\MMc{C}{D A}{A'} \equiv 0 &~\sma{$\mbox{ for } SO_q$} \\
\MMc{C}{C' A}{B}\equiv\MM{C}{C' A}{B} ~&
\!\!\!\!\!\!\!\!\!\!\!\!\!\!\!\!\!\!
\!\!\!\!\!\!\!\!\!\!\!\!\!\!\!\!\! 
\!\!\!\!\!\!\!\!\!\!\!\!\!\!\!\!\!\!\!
,~~\MMc{C}{D A}{A'}\equiv\MM{C}{D A}{A'}
 &~\sma{$\mbox{ for } Sp_q$}\label{Mantisymm}
\end{array}      
\en              
This is easily seen from
(\ref{relationcchi}). 
We have thus obtained the fundamental relation (\ref{leftrightrans}):
$(\cchi{B_1}{B_2}*b)\MMc{A_1}{A_2B_1}{B_2}=
(b*\cchi{A_1}{A_2})$, with \sma{$ A_1'\leq A_2\;,\,~{{B'_1}\leq {B_2}}$}
for the $IS_{q,r=1}(N)$ differential calculus.
Notice that 
$\MMc{C}{DA}{B} = \MM{C}{DA}{B}-\epsilon_C \epsilon_D
q_{DC}            
\MM{D'}{C'A}{B} $, this equality is due to the $RTT$ 
relations (\ref{4.6.16}).
We can also write
$\DR (a\Ome{A}{B})=\sum_{C'\leq D}\D(a)(\Ome{C}{D}\otimes
\MMc{C}{DA}{B})~$  cf.(\ref{defDR3}), thus using the basis
$\{\Ome{C}{D}\}_{C'\leq D}$. According to the general theory
the
elements
$\MMc{C}{DA}{B}$ with \sma{$C' \leq D,~A' \leq B$} are
the
adjoint representation for the differential calculus on
$S_{q,r=1}(N+2)$.
Since the calculus is bicovariant [cf.(\ref{propdaRunoa}),
(\ref{propdaRuno})]
we  know a priori that
the $\MMc{C}{DA}{B}$ with \sma{$C' \leq D,~A' \leq B$} satisfy the
properties
(\ref{copM}) and
(\ref{propM}).\footnote{A direct proof in the \sma{$ SO_q$} 
case is also instructive.
We call $P_{\!-}$ the ``q-antisymmetric" projector
defined by:
\[P_{\!-}{}^{A~~~D}_{~BC}\equiv
\onehalf(\de^A_C\de^D_B-q_{BA}\de_C^{B'}\de_{A'}^D)=
\onehalf(\de^A_C\de^D_B-q_{CD}\de_C^{B'}\de_{A'}^D)~.
\]
Then one easily shows that
$~P_{\!-}{}^{A~~~D}_{~BC}=-q_{BA}
P_{\!-}{}^{B'~~~C}_{~A'D}~
,~~ P_{\!-}{}^{A~~~D}_{~BC}=-q_{CD}
P_{\!-}{}^{A~~~C'}_{~BD'}$ and
\[
\begin{array}{c}
\Omega^jP_{\!-}{}_j{}^i=\Omega^i ~,~~~P_{\!-}{}_i{}^j\chi_j=
\chi_i{}~,~~~ P_{\!-}{}_k{}^if^k{}_j=f^i{}_kP_{\!-}{}_j{}^k=
P_{\!-}{}_k{}^if^k{}_nP_{\!-}{}_j{}^n~,\\
\MMcc{i}{j}=2P_{\!-}{}_i{}^lM_{l}{}^{j}=
2M_{i}{}^{l}P_{\!-}{}_l{}^j~,~{}~~
\MMcc{i}{j}=P_{\!-}{}_i{}^l\MMcc{l}{j}=
\MMcc{i}{l}P_{\!-}{}_l{}^j=
2P_{\!-}{}_i{}^{\alpha}\MMcc{\alpha}{j}
=2\MMcc{i}{\beta}P_{\!-}{}_{\beta}{}^{j}~,\\
\end{array}
\]
where greek letters $\alpha, \beta$ represent adjoint indices
\sma{$(A_1,A_2),~(B_1,B_2)$} with the restriction
\sma{$A_1' < A_2,~B_1' < B_2$}. It is then straightforward to show
that
$\D(\MMcc{i}{j})=\MMcc{i}{\alpha}\otimes\MMcc{\alpha}{j}$ and
$\epsi(\MMcc{\alpha}{\beta})=\de_{\alpha}^{\beta}.$ Applying
$P_{\!-}$ to
(\ref{propM}) and
using $f^i{}_j=0$ unless $i=j$ cf. (\ref{fpropepsi})-(\ref{6.30}) 
one also proves
$\MMcc{\alpha}{j} (a * \f{\alpha}{k})=
(\f{j}{\beta} * a) \MMcc{k}{\beta}$. 
These
formulae
hold in particular if all indices are greek, thus proving
(\ref{copM}) and (\ref{propM}) for $SO_{q,r=1}(N+2)$.}
\sk
\sk
It is useful to express the bicovariant
algebra
(\ref{bico1sosp}), (\ref{bico2})-(\ref{bico4})  in the $\rone$ limit.
Due to the $R$ matrix being diagonal for $r=1$, the $\Lambda$ 
tensor
$
\LL{A_1}{A_2}{D_1}{D_2}{C_1}{C_2}{B_1}{B_2}
\equiv \ff{A_1}{A_2B_1}{B_2} (\MM{C_1}{C_2D_1}{D_2})
$
takes the simple form:
\eq
\LL{A_1}{A_2}{B_1}{B_2}{B_1}{B_2}{A_1}{A_2}=q_{A_1B_2}q_{A_2B_1}
q_{B_1A_1}q_{B_2A_2}~,
{}~~~~0 \mbox{ otherwise}
\en
Therefore (\ref{bico2})-(\ref{bico4}) read (no sum on repeated
indices):
\eq
\f{i}{i}\f{j}{j}=\f{j}{j}\f{i}{i}
\en
\eq
\C{jk}{i} \f{j}{j} \f{k}{k} + \f{i}{j} \chi_k= \L{kj}{jk} \chi_k
\f{i}{j} + \C{jk}{i} \f{i}{i} \label{1bico3}
\en
\eq
\chi_k  \f{i}{i}=\L{ik}{ki} \f{i}{i} \chi_k~.  \label{1bico4}
\en
Explicitly the $q$-Lie algebra (\ref{bico1sosp}) reads:
\eqa
& &\cchi{C_{1}}{C_2}  \cchi{B_{1}}{B_2} - q_{B_1C_2} q_{C_1B_1}
q_{B_2C_1} q_{C_2B_2}
{}~\cchi{B_{1}}{B_2} \cchi{C_{1}}{C_2}=~~~~~~~~~\nonumber\\
& &~~~~-q_{B_1C_2} q_{C_2B_2} q_{B_2B_1} \de^{C_1}_{B_2}
{}~\cchi{B_{1}}
{C_2}+
q_{C_1B_1} q_{B_2B_1} C_{B_2C_2} ~\cchi{B_{1}}{C_1'}+\nonumber\\
& &~~~~+q_{C_2B_2} q_{B_1C_2} C^{C_1B_1}~ \cchi{B_2'} {C_2}-
q_{B_2C_1} \de^{B_1}_{C_2}~\cchi{B_2'}{C_1'}~. \label{Lierone}
\ena
The Cartan-Maurer equations are obtained by differentiating
(\ref{Ome}):
\eq
d\Ome{A}{B}=q_{AB} q_{BC} q_{CA} C_{CD} ~\Ome{C}{B} \we \Ome{A}{D}
\label{CMrone}
\en

The commutations between $\Omega$ 's are easy to find using
(\ref{omcomrone}):
\eq
\Ome{A_1}{A_2} \we \Ome{D_1}{D_2} = - q_{A_1 D_2} q_{D_1A_1}
q_{A_2D_1} q_{D_2A_2} \Ome{D_1}{D_2} \we \Ome{A_1}{A_2}
\label{commOm}
\en
\sk

Finally, we turn to the $*$-conjugations given by equations
(\ref{conjchiAB}) and
(\ref{conjomAB}).
Their $\rone$ limit yields, for (\ref{conjL2})
\eq
(\Ome{A}{B})^*= q_{BA} \Dcal{C}{A} \Ome{C}{D} \Dcal{B}{D}~~~;~~~~
(\cchi{A}{B})^*=-q_{CD} \Dcal{A}{C} \cchi{C}{D} \Dcal{D}{B}=
-{\overline{q_{BA}}}\Dcal{A}{E}\cchi{E}{F}\Dcal{F}{B} ~~,
\label{conjchiABo}
\en
while for the conjugation (\ref{conjL1})
\eq
(\Ome{A}{B})^*= \epsilon q_{BA} \Ome{A}{B}~~~;~~~~
(\cchi{A}{B})^*=
-\epsilon {{q_{AB}}}\cchi{A}{B} ~~.
\en
This shows that we have a bicovariant $*$-differential calculus.


\sect{Differential calculus on $ISO_{q,r=1}(N)$ and \\
$ISp_{q,r=1}(N)$}


We have found the inhomogeneous
quantum group  $\IS$ by means of a projection from
$\S$; dually, its universal
enveloping algebra is a given  Hopf subalgebra of
$U_{q,r}(s(N+2))$.
Using the same techniques and the results of Section  2.3
we here derive the differential calculus on $IS_{q,r=1}(N)$.
\sk
{}From (\ref{UNO}), (\ref{DUE}) and (\ref{adbico}) 
it is immediate to see that 
$T'\equiv T\cap U_{q,r=1}(is(N))$ satisfies 
\eq
\D(T')\subset T'\otimes\epsi+
U_{q,r=1}(is(N))\otimes T'\label{UNOsosp}
\en  
\eq
[T',T']\subseteq T\cap U_{q,r=1}(is(N))=T'\label{T'closesosp}
\en
\eq
\forall\psi\in {\IUsosp}~, ~~~~~~~~~~~~~
a\!d_{\psi}T'\subseteq T'~ \label{adbicososp}
\en
indeed $U_{q,r}(is(N))$ is a Hopf subalgebra of $U_{q,r}(s(N+2))$. 
Also condition (\ref{ZERO})
is fulfilled since  $T'$ generates $U_{q,r}(is(N))$ in the same way 
$T$ generates $U_{q,r}(s(N+2))$ \cite{Burr}, this is a consequence of
the upper and lower triangularity of the $L^+$ and $L^-$ matrices and of the
dependence of the diagonal elements of $L^+$ from the diagonal elements of 
$L^-$; this is true for $r\not=1$ and therefore also for $r=1$.
{}From this last statement, (\ref{UNOsosp}) and 
(\ref{T'closesosp})--or just from (\ref{UNOsosp}) and (\ref{adbicososp})--
we obtain that $T'$ generates  an $IS_{q,r=1}(N)$ bicovariant differential
calculus.
\sk
We reconsider now, in the $\rone$ limit,  the functionals given
in eq.s (\ref{chiN2})-(\ref{chiN2end}). We list below the functionals
among
these that belong to $T'$:
\eqa
& & \cchi{a}{b}=\rinv [\ff{c}{c a}{b}-\de^a_b \epsi] \nonumber \\
& & \cchi{a}{\ci}=\rinv \ff{c}{c a}{\ci} \nonumber \\
& & \cchi{\bu}{b}=\rinv \ff{\bu}{\bu\bu}{b} \nonumber \\
& &~~\nonumber \\
& &\cchi{\ci}{\ci}=\rinv [\ff{\ci}{\ci\ci}{\ci}-\epsi]  \nonumber \\
& & \cchi{\bu}{\bu} = \rinv [\ff{\bu}{\bu\bu}{\bu} - \epsi]\nonumber\\
& & \cchi{\bu}{\ci}=\rinv \ff{\bu}{\bu\bu}{\ci}
\label{chiinisso}
\ena
\noi Note that in the $\rone$ limit
$\cchi{\bu}{\ci}$ vanishes  for
$\SOqroNt$, and does not vanish in the case $\SpqroNt$.

For $r=1$ the $\chi$'s in (\ref{chiinisso}) are not linearly independent, cf.
relation (\ref{relationcchi}) of previous section, and we have:
\eq
\cchi{b'}{a'}=-q_{ab} \cchi{a}{b},~~\cchi{b'}{\ci}=
-{1\over q_{b \bu}} \cchi{\bu}{b},~~\cchi{\ci}{\ci}=-\cchi{\bu}{\bu}
\en
A  basis of tangent vectors for $T'$, in the orthogonal case, is 
therefore given by
\eq
\cchi{a}{b}=\limrone \linv~ [\ff{c}{c a}{b}-\de^a_b \epsi]~,~~~
\mbox{with  \sma{$a+b> N+1~~~$ i.e. $~~a'< b$}};\label{ISOtang1}
\en
\eq
\cchi{\bu}{b}= \limrone \linv \ff{\bu}{\bu\bu}{b}~~;~~~
\cchi{\bu}{\bu} = \limrone \linv~ [\ff{\bu}{\bu\bu}{\bu} - \epsi]~,
\label{ISOtang2}
\en
The
$q$-Lie algebra commutations are a subset  of (\ref{Lierone})
obtained specializing the capital indices of (\ref{Lierone}) 
to the indices
${}^a{}_b\,,~{}^{\bu}{}_b$ and ${}^{\bu}{}_{\bu}$.
We have the $SO_{q,r=1}(N)$ $q$-Lie algebra that  
reads as in eq. (\ref{Lierone}) with lower case indices; the remaining
commutations are:
\eqa
& &\cchi{c_1}{c_2} \chi_{b_2} - {q_{c_1\bu} \over q_{c_2\bu}}
q_{b_2c_1} q_{c_2b_2} \chi_{b_2} \cchi{c_1}{c_2}=
{q_{c_1\bu} \over q_{c_2\bu}}[C_{b_2c_2} \chi_{c_1'}-\de^{c_1}_{b_2}
q_{c_2c_1} \chi_{c_2}]~, \label{Commchi0}\\
& &\chi_{c_2} \chi_{b_2} - {q_{b_2\bu} \over q_{c_2\bu}}
q_{c_2b_2} \chi_{b_2} \chi_{c_2}
=0~,\label{Commchi}\\
& &\cchi{c_1}{c_2} \cchi{\bu}{\bu} -\cchi{\bu}{\bu}
\cchi{c_1}{c_2}=0~~~~,~~~~
\chi_{c_2}\cchi{\bu}{\bu}-\cchi{\bu}{\bu}\chi_{c_2}=-\chi_{c_2}~
\ena
\noi where we have defined
$\chi_a \equiv  \cchi{\bu}{a}~.$
The exterior differential reads, $\forall a\in ISO_{q,r=1}(N)$
\eq
da=\sum_{a'<b}(\cchi{a}{b}*a)\Omega_{a}{}^b +
(\cchi{\bu}{b}*a)\Omega_{\bu}{}^b  +
(\cchi{\bu}{\bu}*a)\Omega_{\bu}{}^{\bu}~~\label{provd}
\en
where $\Omega_{a}{}^{b},~\Omega_{\bu}{}^{b},$
and $\Omega_{\bu}{}^{\bu}$ are the $1$-forms dual
to the tangent vectors (\ref{ISOtang1}) and (\ref{ISOtang2}).
As discussed in \cite{ISODUAL}, these $1$-forms can be seen as the
projection of the $S_{q,r=1}(N+2)$ $1$-forms : $P(\Omega_A{}^B)
=-q_{AB}P[\kappa(\T{B}{C})]d P(\T{C}{A})$.
\sk
The adjoint representation, defined by  (\ref{adbico}):
$ad_{\psi}=M_i{}^j(\psi)\chi_j$,
is given by the elements
$P(\MMc{C}{DA}{B})$
$\in \ISOqroN$  with \sma{$C' \leq D,~A' \leq B$}   obtained by
projecting with
$P$
those of $\SOqroNt$.

\noi{\sl Proof }:$~$ In $SO_{q,r}(N+2)$ we have 
$a\!d_{\psi}\cchi{C}{D}=\MMc{C}{DA}{B}(\psi)\cchi{C}{D}$
with \sma{$C' \leq D,~A' \leq B$}, since
$\MMc{C}{DA}{B}$ is the adjoint representation of the $SO_{q,r}(N+2)$ 
calculus (see Note 4.6.2).
Now $\MMc{C}{DA}{B}(\psi)=\psi(\MMc{C}{DA}{B})=\le\psi\,,\,P(\MMc{C}{DA}{B})
\re$
where the last bracket is the duality bracket between $ISO_{q,r}(N)$ and 
$U_{q,r}(iso(N))$ [cf. (\ref{dualityiso})].
We then obtain:
$$a\!d_{\psi}\cchi{C}{D}=\le\psi\,,\,P(\MMc{C}{DA}{B})\re\,\cchi{C}{D}~~~
\mbox{ with } \sma{$C' \leq D,~A' \leq B$}~,$$
this is the defining formula for the adjoint representation associated to
the quantum Lie algebra $T'$.
The nonvanishing elements are:
\eqa
& &P(\MMc{b_1}{b_2a_1}{a_2}) = \T{b_1}{a_1} \kappa (\T{a_2}{b_2})
- q_{b_2b_1} \T{b_2'}{a_1} \kappa (\T{a_2}{b_1'})\nonumber\\
& &P(\MMc{b_1}{b_2 \bu}{a_2}) = x^{b_1} \kappa (\T{a_2}{b_2})
- q_{b_2b_1} x^{b_2'} \kappa (\T{a_2}{b_1'})
\nonumber\\
& &P(\MMc{\bu}{b_2 \bu}{a_2}) = v  \kappa (\T{a_2}{b_2})\nonumber\\
& &P(\MMc{\bu}{\bu \bu}{a_2}) = v  \kappa (x^{a_2})\nonumber \\
& &P(\MMc{\bu}{\bu \bu}{\bu}) = I 
\label{Minsospext}
\ena
We will later use the 
relation between left invariant and right invariant vectorfields; 
in our case (\ref{leftrightrans})
reads:
\eq
(\cchi{A_1}{A_2}*b)\,P(\MMc{B_1}{B_2A_1}{A_2})=b*\cchi{B_1}{B_2}
~~~\mbox{ with \sma{$ A_1'<A_2,~B_1'<B_2$} }\label{perisoris}~.
\en 
\sk
The $ISp_{q,r=1}(N)$ differential calculus has the same structure as
the  $ISO_{q,r=1}(N)$ one, provided one considers \sma{$a'\leq b$} in 
(\ref{ISOtang1}) and (\ref{provd}), and includes the extra generator 
$\cchi{\bu}{\ci}$
in (\ref{ISOtang2}) and his dual form $\Omega_{\ci}{}^{\bu}$ 
in the definition of the exterior differential. The adjoint representation 
is  obtained by projecting with $P$ the adjoint representation of 
the $Sp_{q,r}(N+2)$ differential calculus.
\sk
We now show that it is possible to exclude
the generator  $\cchi{\bu}{\bu}$
(and $\cchi{\ci}{\ci}$)  and obtain a
dilatation-free
bicovariant differential calculus on $ISO_{q,r=1}(N)$. 

We study the $ISO_{q,r=1}(N)$ subspace {\sl{g}} linearly spanned
by the functionals $\cchi{a}{b}\,,~\chi_b$:
\eq
\mbox{ {\sl{g}} }\equiv {\mbox{\sl{span}}}
\{\cchi{a}{b}\,,~\chi_b\}\label{restrictedlie}~~.
\en
The space {\sl{g}} is
our candidate 
$q$-Lie algebra. A basis of {\sl{g}} is  
\mbox{$\{\chi_{\al}\} = \{\cchi{a}{b} \/ (a'<b), \cchi{\bu}{b}\}$.}
In the sequel  greek  letters will denote adjoint indices
$\al= (a_1,a_2)$ with $a'_1 < a_2 \mbox{, and } \al=(\bu,a_2)$.
The coproduct on the elements $\chi_{\al}$ reads
$\D'\chi_{\al}=\chi_{\al}\otimes f^{\al}{}_{\al} +\epsi\otimes\chi_{\al}$;
this shows that {\sl{g}} satisfies condition (\ref{UNO}).
We also have $\D'
f^{\al}{}_{\al}=f^{\al}{}_{\al}\otimes f^{\al}{}_{\al}$.
To prove that {\sl{g}}  defines a bicovariant differential calculus we can
proceed as in Section 3.5. We here give an alternative proof based on 
the results of Section 2.3. Recalling 
Theorem 2.3.1, {\sl{g}}  defines a bicovariant differential calculus if
there exists a set of elements $M_i{}^j\in ISO_{q,r=1}(N)$ that satysfy 
(\ref{UNO}) and(\ref{leftrightrans}): $(\chi_j*b)\M{i}{j}=b*\chi_i$. 
It is immediate to verify that the subset of (\ref{Minsospext}) given by
\eqa
& &P(\MMc{b_1}{b_2a_1}{a_2}) = \T{b_1}{a_1} \kappa (\T{a_2}{b_2})
- q_{b_2b_1} \T{b_2'}{a_1} \kappa (\T{a_2}{b_1'})\nonumber\\
& &P(\MMc{b_1}{b_2 \bu}{a_2}) = x^{b_1} \kappa (\T{a_2}{b_2})
- q_{b_2b_1} x^{b_2'} \kappa (\T{a_2}{b_1'})
\nonumber\\
& &P(\MMc{\bu}{b_2 \bu}{a_2}) = v  \kappa (\T{a_2}{b_2})
\label{Minsosp}
\ena
satisfies
\eq
(\chi_{\be}*a)M_{\al}{}^{\be}=a*\chi_{\al}
\en
indeed $P(\MMc{b_1}{b_2\bu}{\bu})=P(\MMc{\bu}{b_2\bu}{\bu})=0$ and therefore
(\ref{perisoris})
closes also on the subset of $\chi$ and $M_-$ with greek indices.
We have therefore shown:
\sk

\noi{\bf Theorem }4.7.1 $~$ 
{\sl g} is a quantum Lie algebra and  defines a bicovariant differential
calculus on $ISO_{q,r=1}(N)$ that has the same dimension as in the 
commutative case.
\sk
\cvd
We now analize this differential calculus.
The exterior derivative is
\eq
da=(\chi_{\al} * a) \Omega^{\al}\label{isp}
\en
The left $\ISOqroN$--module 
$\Ga$ freely generated
by the $1$-forms $\Omega^{\al}$ dual to the tangent vectors $\chi_{\al}$
is a bicovariant bimodule over $\ISOqroN$
with the right multiplication (no sum on repeated indices):
\eq
\Omega^{\al} a = (\f{\al}{\al} * a) \Omega^{\al}~,
{}~~~~a \in \ISOqroN \label{omiasosp}
\en
\noi and with
the left and right actions
of $\ISOqroN$ on $\Ga$  given by:
\eqa
& &\DL (a_{\al} \Omega^{\al}) \equiv \D(a_{\al}) I \otimes
\Omega^{\al}
 \label{DLinsosp}\\
& &\DR (a_{\al} \Omega^{\al}) \equiv \D(a_{\al}) \Omega^{\be} \otimes
P(\Mc{{\be}}{{\al}})
\label{DRinsosp}~.
\ena
\sk

Using the general formula (\ref{omiasosp}) we can
deduce the $\Omega, T$ commutations:
\eqa
& &\Ome{a_1}{a_2} \T{r}{s}={q_{a_2 s} \over q_{a_1 s}} \T{r}{s}
 \Ome{a_1}{a_2}\\
& &\Ome{a_1}{a_2} x^r= {q_{a_2 \bu} \over q_{a_1 \bu}} x^r
\Ome{a_1}{a_2} \\
& &\Ome{a_1}{a_2} u= {q_{a_1 \bu} \over q_{a_2 \bu}} u~
\Ome{a_1}{a_2} \\
& &\Ome{\bu}{a_2} \T{r}{s}=q_{s\bu} q_{a_2 s} \T{r}{s}
\Ome{\bu}{a_2}\label{VTcomm}\\
& &\Ome{\bu}{a_2} x^r=q_{a_2\bu} x^r \Ome{\bu}{a_2} \label{Vxcommsosp}\\
& &\Ome{\bu}{a_2} u={1\over q_{a_2\bu}} u~\Ome{\bu}{a_2}
\ena
\noi {\bf Note} 4.7.1  $u$ commutes with all
 $\Omega$ 's only if $q_{a\bu}=1$ (cf.  Note 4.2.2). This means
that
$u=I$ is consistent with the differential calculus
on $ISO_{q_{ab},r=1,q_{a\bu}=1}(N)$.
\sk
The exterior
derivative on the generators $\T{A}{B}$ is given by:
\eqa
& & d \T{a}{b}= -\sum_{c}\T{a}{c} q_{cb} \Ome{b}{c}\nonumber\\
& & dx^a=-\sum_c \T{a}{c} q_{c\bu} V^c \label{dTABiso}\\
& & du=dv=0\nonumber
\ena
\noi  where we have defined $V^a \equiv \Ome{\bu}{a}$.
Again, for $q_{a\bu}=1$,  $u=v=I$ is a consistent  choice.
\sk
Inverting (\ref{dTABiso}) yields:
\eqa
& & \Ome{a}{b}=-q_{ab} \kappa (\T{b}{c}) ~d \T{c}{a}\label{Omein}\\
& & V^b=-{1\over q_{b\bu}}
\kappa (\T{b}{c}) ~ dx^c \label{Vin}
\ena
\sk
The exterior product of the left-invariant $1$-forms
is defined as
\eq
\Om^{\al} \we \Om^{\be}\equiv \Om^{\al} \otimes \Om^{\be} -
\L{\al\be}{\ga\de}
\Om^{\ga} \otimes \Om^{\de}\label{26}
\en
\noi where
\eq
\L{\al\be}{\ga\de}\equiv\le\f{\al}{\de}\,,\,P(\Mc{\ga}{\be})\re=
\f{\al}{\de}(\Mc{\ga}{\be}) \label{27}
\en
[cf. (\ref{dualityiso})]; so that this $\Lambda$ tensor is 
obtained from
the one of $\SOqroNt$ by restricting its indices to the
subset $ab,  \bu b$.
We therefore just
specialize the indices
in equation (\ref{commOm}) to deduce
the $q$-commutations for the $1$-forms $\Omega$ and $V$:
\eq
\Ome{a_1}{a_2} \we \Ome{d_1}{d_2} = - q_{a_1 d_2} q_{d_1a_1}
q_{a_2d_1} q_{d_2a_2} \Ome{d_1}{d_2} \we \Ome{a_1}{a_2}
\en
\eq
\Ome{a_1}{a_2} \we V^{d_2} = - {q_{a_2\bu} \over q_{a_1\bu}}
q_{a_1 d_2} q_{d_2a_2} V^{d_2} \we \Ome{a_1}{a_2}
\en
\eq
V^{a_2} \we V^{d_2}=-{q_{a_2\bu} \over q_{d_2\bu}} q_{d_2a_2}
V^{d_2} \we V^{a_2}
\en
\sk
The Cartan-Maurer equations
\eq
d\Om^{\al}=
-\onehalf\c{\be\ga}{\al} \Om^{\be} \we \Om^{\ga}\label{CMin0}
\en
\noi can be explicitly written for the $\Omega$ and $V$
by differentiating eq.s (\ref{Omein}) and (\ref{Vin}) [or
again secializing the indices in (\ref{CMrone})]:
\eqa
& &d\Ome{a}{b}= q_{ab} q_{bc} q_{ca} ~\Ome{c}{b} \we \Ome{a}{c}
\label{CMin1}\\
& &d V^b={q_{a\bu} \over q_{b\bu}} q_{ba} ~\Ome{a}{b} \we V^a
\label{CMin2}
\ena
where the $1$-forms
$\Ome{a}{b}$ with $a'>b$ are given by  $\Ome{a}{b}=-q_{ab}
\Ome{b'}{a'}$;
i.e. we consider (as it is usually done in the classical limit),
the $1$-forms $\Ome{a}{b}$ to be
``$q$-antisymmetric''
$\Ome{a}{b}=-q_{ab} \Ome{b'}{a'}$,
cf. eq. (\ref{relationOme}).
\sk
The $*$-conjugation on the $\chi$ functionals and on the $1$-forms
$\Omega$ can be deduced from (\ref{conjchiABo}):
\eq
(\cchi{a}{b})^*=-q_{cd} \Dcal{a}{c} \cchi{c}{d}
\Dcal{d}{b},~~(\chi_b)^*=
-(q_{d\bu})^{-1} \chi_d \Dcal{d}{b}=-{\overline{q_{b\bu}}} \chi_d \Dcal{d}{b}
\en
\eq
(\Ome{a}{b})^*=q_{ba} \Dcal{c}{a} \Ome{c}{d} \Dcal{b}{d},~~(V^b)^*=
q_{b\bu} V^d \Dcal{b}{d}
\en
\sk
\noi{\bf Note} 4.7.2  $~$ As discussed at the end of
Section 4.2,
a  $q$-Poincar\'e group without dilatations (i.e. $u=I$) has only
one free real parameter $q_{12}$, which is the real parameter
related to the $q$-Lorentz subalgebra. Then the formulas of this
section
can be specialized to describe a bicovariant calculus on the
dilatation-free
$ISO_{q,r=1}(3,1)$ provided  $q_{a\bu}=1$ and  $q_{12} \in \Rb$. It
is
however possible to have a bicovariant calculus without the
dilatation
generator
$\cchi{\bu}{\bu}$ even on $ISO_{q,r=1}(3,1)$ with $u \not= I$.
The $q$-Poincar\'e
algebra presented in \cite{Cas2} corresponds to the case
$q \equiv q_{1\bu}$, $q_{2\bu}=q_{12}=1$, for which the Lorentz
subalgebra is undeformed and the $q$-Poincar\'e group contains
$u\not= I$.
 The
possibility
of having a dilatation-free $q$-Lie algebra describing a bicovariant
calculus
on
a $q$-group containing dilatations $u$ was already observed in the
case of $IGL$ $q$-groups (see Section 3.5).
\sk

\noi{\bf Note} 4.7.3 $~$ We here study a differential calculus on 
$ISp_{q,r=1}(N)$ that has the same number of tangent vectors as in the 
classical case. Following the same arguments given after (\ref{restrictedlie}) 
we have an $ISp_{q,r=1}(N)$ differential calculus with quantum Lie algebra
generators $\cchi{a}{b}$ with \sma{$a'\leq b$},
$\chi_b$ and $\cchi{\bu}{\ci}$. To further restrict the quantum 
Lie algebra to the one spanned by the basis  $\{\cchi{a}{b}\,
\mbox{\sma{$(a'\leq b)$}}~, \chi_b\}$, observe that
from (\ref{spmeno1}),
\eq
\cchi{\bu}{\ci}=\limrone \linv~ \kappa'(\Lp{\bu}{\bu})\Lm{\bu}{\ci}
\en
is different from zero only on monomials that contain the element $z$.
However in the $\qrone$ limit, as noticed after (\ref{cfin}), 
$z$ can be set to zero since there is no more any constraint between $z$ and 
the generators $\T{a}{b}, x^a, u, v=u^{-1}$. Then $\cchi{\bu}{\ci}$ is zero
as well and we have an $[N(N+1)/2 + N]$--dimensional bicovariant 
differential calculus 
on the twisted inhomogeneous symplectic group generated by 
$\T{a}{b}, x^a, u, v=u^{-1}$. 
The adjoint representation is given by the elements
$P(\Mc{\al}{\be}) \in ISp_{q,r}(N)$ obtained by
projecting with
$P$
those of  $Sp_{q,r=1}(N+2)$.
The explicit formulae carachterizing this differential calculus
are as in (\ref{isp})--(\ref{CMin0}), where now 
greek  letters denote adjoint indices
$\al= (a_1,a_2)$ with $a'_1 \leq a_2 \mbox{, and } \al=(\bu,a_2)$.


\def\spinst#1#2{{#1\brack#2}}
\def\sk{\vskip .4cm}
\def\noi{\noindent}
\def\om{\omega}
\def\Om{\Omega}
\def\al{\alpha}
\def\la{\lambda}
\def\be{\beta}
\def\ga{\gamma}
\def\Ga{\Gamma}
\def\del{\delta}
\def\linv{{1 \over \lambda}}
\def\rinv{{1\over {r-r^{-1}}}}
\def\alb{\bar{\alpha}}
\def\beb{\bar{\beta}}
\def\gab{\bar{\gamma}}
\def\deb{\bar{\delta}}
\def\ab{\bar{a}}
\def\Ab{\bar{A}}
\def\Bb{\bar{B}}
\def\Cb{\bar{C}}
\def\Db{\bar{D}}
\def\ab{\bar{a}}
\def\cb{\bar{c}}
\def\db{\bar{d}}
\def\bb{\bar{b}}
\def\eb{\bar{e}}
\def\fb{\bar{f}}
\def\gb{\bar{g}}
\def\xih{\hat\xi}
\def\Xih{\hat\Xi}
\def\uh{\hat u}
\def\vh{\hat v}
\def\ub{\bar u}
\def\vb{\bar v}
\def\xib{\bar \xi}

\def\alp{{\alpha}^{\prime}}
\def\bep{{\beta}^{\prime}}
\def\gap{{\gamma}^{\prime}}
\def\dep{{\delta}^{\prime}}
\def\rhop{{\rho}^{\prime}}
\def\taup{{\tau}^{\prime}}
\def\rhopp{\rho ''}
\def\thetap{{\theta}^{\prime}}
\def\imezzi{{i\over 2}}
\def\unquarto{{1 \over 4}}
\def\onehalf{{1 \over 2}}
\def\unmezzo{{1 \over 2}}
\def\epsi{\varepsilon}
\def\we{\wedge}
\def\th{\theta}
\def\de{\delta}
\def\cony{i_{\de {\vec y}}}
\def\Liey{l_{\de {\vec y}}}
\def\tv{{\vec t}}
\def\Gt{{\tilde G}}
\def\deyv{\vec {\de y}}
\def\part{\partial}
\def\pdxp{{\partial \over {\partial x^+}}}
\def\pdxm{{\partial \over {\partial x^-}}}
\def\pdxi{{\partial \over {\partial x^i}}}
\def\pdy#1{{\partial \over {\partial y^{#1}}}}
\def\pdx#1{{\partial \over {\partial x^{#1}}}}
\def\pdyx#1{{\partial \over {\partial (yx)^{#1}}}}

\def\qP{q-Poincar\'e~}
\def\A#1#2{ A^{#1}_{~~~#2} }

\def\R#1#2{ R^{#1}_{~~~#2} }
\def\PA#1#2{ P^{#1}_{A~~#2} }
\def\PS#1#2{ P^{#1}_{S~~#2} }
\def\Pa#1#2{ (P_A)^{#1}_{~~#2} }
\def\Pas#1#2{ (P_A)^{#1}_{~#2} }
\def\PI#1#2{(P_I)^{#1}_{~~~#2} }
\def\PJ#1#2{ (P_J)^{#1}_{~~~#2} }
\def\Ppp{(P_+,P_+)}
\def\Ppm{(P_+,P_-)}
\def\Pmp{(P_-,P_+)}
\def\Pmm{(P_-,P_-)}
\def\Ppo{(P_+,P_0)}
\def\Pom{(P_0,P_-)}
\def\Pop{(P_0,P_+)}
\def\Pmo{(P_-,P_0)}
\def\Poo{(P_0,P_0)}
\def\Pso{(P_{\sigma},P_0)}
\def\Pos{(P_0,P_{\sigma})}

\def\Rp#1#2{ (R^+)^{#1}_{~~~#2} }
\def\Rpinv#1#2{ [(R^+)^{-1}]^{#1}_{~~~#2} }
\def\Rm#1#2{ (R^-)^{#1}_{~~~#2} }
\def\Rinv#1#2{ (R^{-1})^{#1}_{~~~#2} }
\def\Rsecondinv#1#2{ (R^{\sim 1})^{#1}_{~~~#2} }
\def\Rinvsecondinv#1#2{ ((R^{-1})^{\sim 1})^{#1}_{~~~#2} }

\def\Rpm#1#2{(R^{\pm})^{#1}_{~~~#2} }
\def\Rpminv#1#2{((R^{\pm})^{-1})^{#1}_{~~~#2} }

\def\Rb{{\bf \mbox{\boldmath $R$}}}
\def\Rbo{{\bf \mbox{\boldmath $R$}}}
\def\Rbp#1#2{{ (\Rbo^+)^{#1}_{~~~#2} }}
\def\Rbm#1#2{ (\Rbo^-)^{#1}_{~~~#2} }
\def\Rbinv#1#2{ (\Rbo^{-1})^{#1}_{~~~#2} }
\def\Rbpm#1#2{(\Rbo^{\pm})^{#1}_{~~~#2} }
\def\Rbpminv#1#2{((\Rbo^{\pm})^{-1})^{#1}_{~~~#2} }

\def\RRpm{R^{\pm}}
\def\RRp{R^{+}}
\def\RRm{R^{-}}

\def\Rh{{\hat R}}
\def\Rbh{{\hat {\Rbo}}}
\def\Rhat#1#2{ \Rh^{#1}_{~~~#2} }
\def\Rbar#1#2{ {\bar R}^{#1}_{~~~#2} }
\def\L#1#2{ \La^{#1}_{~~~#2} }
\def\Linv#1#2{ \La^{-1~#1}_{~~~~~#2} }
\def\Rbhat#1#2{ \Rbh^{#1}_{~~~#2} }
\def\Rhatinv#1#2{ (\Rh^{-1})^{#1}_{~~~#2} }
\def\Rbhatinv#1#2{ (\Rbh^{-1})^{#1}_{~~~#2} }
\def\Z#1#2{ Z^{#1}_{~~~#2} }
\def\X#1#2{ X^{#1}_{~~~#2} }
\def\Rt#1{ {\hat R}_{#1} }
\def\La{\Lambda}
\def\Rha{{\hat R}}
\def\ff#1#2#3{f_{#1~~~#3}^{~#2}}
\def\MM#1#2#3{M^{#1~~~#3}_{~#2}}
\def\MMc#1#2#3{{M_{\!-}}^{#1~~~#3}_{~#2}}
\def\MMcc#1#2{M_{\!-}{}_{#1}{}^{#2}}
\def\cchi#1#2{\chi^{#1}_{~#2}}
\def\chil#1{\chi_{{}_{#1}}}
\def\ome#1#2{\om_{#1}^{~#2}}
\def\Ome#1#2{\Omega_{#1}^{~#2}}
\def\RRhat#1#2#3#4#5#6#7#8{\La^{~#2~#4}_{#1~#3}|^{#5~#7}_{~#6~#8}}
\def\RRhatinv#1#2#3#4#5#6#7#8{(\La^{-1})^
{~#2~#4}_{#1~#3}|^{#5~#7}_{~#6~#8}}
\def\LL#1#2#3#4#5#6#7#8{\La^{~#2~#4}_{#1~#3}|^{#5~#7}_{~#6~#8}}
\def\LLinv#1#2#3#4#5#6#7#8{(\La^{-1})^
{~#2~#4}_{#1~#3}|^{#5~#7}_{~#6~#8}}
\def\U#1#2#3#4#5#6#7#8{U^{~#2~#4}_{#1~#3}|^{#5~#7}_{~#6~#8}}
\def\Cb{\bf \mbox{\boldmath $C$}}
\def\CC#1#2#3#4#5#6{{\Cb}_{~#2~#4}^{#1~#3}|_{#5}^{~#6}}
\def\cc#1#2#3#4#5#6{C_{~#2~#4}^{#1~#3}|_{#5}^{~#6}}
\def\PIJ#1#2#3#4#5#6#7#8{(P_I,P_J)^{~#2~#4}_{#1~#3}|^{#5~#7}_{~#6~#8}}

\def\ZZ#1#2#3#4#5#6#7#8{Z^{~#2~#4}_{#1~#3}|^{#5~#7}_{~#6~#8}}

\def\C#1#2{ {\bf \mbox{\boldmath $C$}}_{#1}^{~~~#2} }
\def\c#1#2{ C_{#1}^{~~~#2} }
\def\q#1{   {{q^{#1} - q^{-#1}} \over {q^{\unmezzo}-q^{-\unmezzo}}}}
\def\Dmat#1#2{D^{#1}_{~#2}}
\def\Dmatinv#1#2{(D^{-1})^{#1}_{~#2}}
\def\DR{\Delta_R}
\def\DL{\Delta_L}
\def\f#1#2{ f^{#1}_{~~#2} }
\def\F#1#2{ F^{#1}_{~~#2} }
\def\T#1#2{ T^{#1}_{~~#2} }
\def\Ti#1#2{ (T^{-1})^{#1}_{~~#2} }
\def\Tp#1#2{ (T^{\prime})^{#1}_{~~#2} }
\def\Th#1#2{ {\hat T}^{#1}_{~~#2} }
\def\TP{ T^{\prime} }
\def\M#1#2{ M_{#1}^{~#2} }
\def\Mc#1#2{ {M_{\!-}}_{#1}^{~#2} }
\def\qm{q^{-1}}
\def\rminus{r^{-1}}
\def\um{u^{-1}}
\def\vm{v^{-1}}
\def\xm{x^{-}}
\def\xp{x^{+}}
\def\fm{f_-}
\def\fp{f_+}
\def\fn{f_0}
\def\D{\Delta}
\def\DN{\Delta_{N+1}}
\def\kN{\kappa_{N+1}}
\def\eN{\epsi_{N+1}}
\def\Mat#1#2#3#4#5#6#7#8#9{\left( \matrix{
     #1 & #2 & #3 \cr
     #4 & #5 & #6 \cr
     #7 & #8 & #9 \cr
   }\right) }
\def\Ap{A^{\prime}}
\def\Dp{\Delta^{\prime}}
\def\Ip{I^{\prime}}
\def\ep{\epsi^{\prime}}
\def\kp{{\kappa^{\prime}}}
\def\kpm{\kappa^{\prime -1}}
\def\kpsq{\kappa^{\prime 2}}
\def\km{\kappa^{-1}}
\def\gp{g^{\prime}}
\def\qone{q \rightarrow 1}
\def\rone{r \rightarrow 1}
\def\qrone{q,r \rightarrow 1}
\def\Fmn{F_{\mu\nu}}
\def\Am{A_{\mu}}
\def\An{A_{\nu}}
\def\dm{\part_{\mu}}
\def\dn{\part_{\nu}}
\def\Ana{A_{\nu]}}
\def\Bna{B_{\nu]}}
\def\Zna{Z_{\nu]}}
\def\dma{\part_{[\mu}}
\def\qsu{$[SU(2) \times U(1)]_q~$}
\def\suq{$SU_q(2)~$}
\def\su{$SU(2) \times U(1)~$}
\def\gij{g_{ij}}
\def\qL{SL_q(2,{\bf \mbox{\boldmath $C$}})}
\def\GLqrN{GL_{q,r}(N)}
\def\IGLqrN{IGL_{q,r}(N)}
\def\IGLqrtwo{IGL_{q,r}(2)}
\def\GLqrNo{GL_{q,r}(N+1)}
\def\SOqrNt{SO_{q,r}(N+2)}
\def\SpqrNt{Sp_{q,r}(N+2)}
\def\SLqrN{SL_{q,r}(N)}
\def\UglqrN{U(gl_{q,r}(N))}
\def\UglqrNo{U(gl_{q,r}(N+1))}
\def\UiglqrN{U(igl_{q,r}(N))}
\def\ISOqrN{ISO_{q,r}(N)}
\def\ISpqrN{ISp_{q,r}(N)}
\def\ISOqroN{ISO_{q,r=1}(N)}
\def\ISpqroN{ISp_{q,r=1}(N)}
\def\SqrN{S_{q,r}(N)}
\def\SqrNtwo{S_{q,r}(N+2)}
\def\USqrNtwo{U(S_{q,r}(N+2))}
\def\UISqrN{U(IS_{q,r}(N))}
\def\ISqrN{IS_{q,r}(N)}
\def\SqroNt{S_{q,r=1}(N+2)}
\def\SOqroNt{SO_{q,r=1}(N+2)}
\def\SpqroNt{Sp_{q,r=1}(N+2)}
\def\ISqroN{IS_{q,r=1}(N)}
\def\ISOqroN{ISO_{q,r=1}(N)}

\def\SOqrN{SO_{q,r}(N)}
\def\SpqrN{Sp_{q,r}(N)}
\def\SqrNt{S_{q,r}(N+2)}

\def\Tc{{\cal T}}

\def\Dtwo{\Delta_{N+2}}
\def\epsitwo{\epsi_{N+2}}
\def\kappatwo{\kappa_{N+2}}

\def\RR{R^*}
\def\rr#1{R^*_{#1}}

\def\Lpm#1#2{L^{\pm #1}_{~~~#2}}
\def\Lmp#1#2{L^{\mp#1}_{~~~#2}}
\def\LLpm{L^{\pm}}
\def\LLmp{L^{\mp}}
\def\LLp{L^{+}}
\def\LLm{L^{-}}
\def\Lp#1#2{L^{+ #1}_{~~~#2}}
\def\Lm#1#2{L^{- #1}_{~~~#2}}
\def\Dcal#1#2{{\cal D}^{#1}_{~#2}}

\def\gu{g_{U(1)}}
\def\gsu{g_{SU(2)}}
\def\tg{ {\rm tg} }
\def\Fun{$Fun(G)~$}
\def\invG{{}_{{\rm inv}}\Ga}
\def\Ginv{\Ga_{{\rm inv}}}
\def\qonelim{\stackrel{q \rightarrow 1}{\longrightarrow}}
\def\ronelim{\stackrel{r \rightarrow 1}{\longrightarrow}}
\def\limrone{\lim_{r \rightarrow 1}}
\def\ronelimeq{\stackrel{r \rightarrow 1}{=}}
\def\Pprojection{\stackrel{P}{\longrightarrow}}
\def\qronelim{\stackrel{q=r \rightarrow 1}{\longrightarrow}}
\def\viel#1#2{e^{#1}_{~~{#2}}}
\def\ra{\rightarrow}
\def\detq{{\det}}
\def\detqr{{\det}}
\def\detqrm{{\det} }
\def\detqrTAB{{\det} \T{A}{B}}
\def\detqrTab{{\det} \T{a}{b}}
\def\P{P}
\def\Qt{Q}
\def\chit{{\partial}}

\def\pp#1#2{\Pi_{#1}^{(#2)}}

\def\BCD{B_n, C_n, D_n}

\def\n2{{{N+1} \over 2}}
\def\ap{a^{\prime}}
\def\bp{b^{\prime}}
\def\cp{c^{\prime}}
\def\dpr{d^{\prime}}
\def\Dc{{\cal D}}
\def\osqrt{{1 \over \sqrt{2}}}
\def\Ntwo{{N\over 2}}

\def\bu{\bullet}
\def\ci{\circ}

\def\sma#1{\mbox{\footnotesize #1}}
\def\Q.E.D.{\rightline{$\Box$}}
\def\vt{\vartheta}

\def\Nmezzi{{N\over 2}}
\def\fmi#1#2{f^{#1}_{~#2}}
\def\LM#1#2{\Lambda^{#1}_{~~~#2}}
\def\Acal{{\cal A}}
\def\parleft{{\stackrel{\leftarrow}{\chit}}}
\def\lati{{\tilde \la}}
\def\muti{{\tilde \mu}}
\def\le{\langle}
\def\re{\rangle}


\chapter{\bf Geometry of the 
quantum orthogonal plane}

We present here a bicovariant calculus on the full multiparametric
$ISO_{q,r}(N)$
 {\sl without the restriction}  $r=1$. This calculus, however,  is
trivial on
the $SO_{q,r}(N)$
quantum subgroup: it can really be seen as a non-trivial calculus
only on the
coset $Fun_{q,r}[ISO(N)/SO(N)]$, i.e. on the quantum orthogonal  
plane.
We therefore call this calculus on the quantum plane
$ISO_{q,r}(N)$--bicovariant. We find that in the $r\not= 1$ case this
$ISO_{q,r}(N)$--bicovariant calculus necessarily contains  
dilatations.

If we break $ISO_{q,r}(N)$ bicovariance and require     
right covariance under $ISO_{q,r}(N)$ and
left covariance only under $SO_{q,r}(N)$,
i.e.  compatibility of the exterior differential
on the quantum plane with
the right $ISO_{q,r}(N)$-coaction and the
left $SO_{q,r}(N)$-coaction,
the calculus can be expressed in terms of coordinates $x$,
differentials $dx$ and partial derivatives $\partial$,
{\sl without the need of dilatations}.  In this case
$q$-commutations between $x$, $dx$ and $\partial$
close by themselves, and in
fact generalize to the multiparametric case  the known results of
ref.s
\cite{Manin,qplane,Fiore1}.  Here these results emerge from the
broader setting of the bicovariant calculus on $ISO_{q,r}(N)$.

The two $*$-conjugations of the previous sections, consistent with
the $q$-group structure, lead to a $ISO_{q,r}(n+1,n-1)$, and a
$ISO_{q,r}(n,n)$ or $ISO_{q,r}(n,n+1)$
bicovariant calculus on the quantum orthogonal plane
respectively with $(n+1,n-1)$, $(n,n)$ or $(n,n+1)$ signature.
We will be concerned with the conjugation that gives the 
 $ISO_{q,r}(n-1,n+1)$ calculus. [To retrieve the other conjugations,
both for $N$=even and $N$=odd, just
take $\Dcal{A}{B}=\delta^A_B$ in the formulae where $\Dcal{A}{B}$ appears].

Using this conjugation  one can define real coordinates $X$ and
hermitian partial derivative operators $P$ i.e. momenta.
This is achieved by a canonical procedure, using the compatibility of 
the $*$-structure with the
bicovariant calculus on $ISO_{q,r}(N)$, i.e. the property that $*$ 
is a linear operation on the $q$-Lie algebra.
The $q$-commutations of the momenta $P$ with the coordinates $X$
define a  deformed version of the Heisenberg $X$, $P$ commutation relations
(with no extra operator as in ref.s
\cite{Lorek}).  
In the same spirit as in ref.s \cite{Lorek} it would be 
interesting to investigate 
the Hilbert space representations of this deformed phase-space algebra.
\sk
In Section 5.1 we present the 
$ISO_{q,r}(N)$ bicovariant differential calculus with $r\not= 1$,
then, in Section 5.2 we restrict this calculus to the 
quantum orthogonal plane. We find that
in order to obtain a space of $1$-forms that
has the same dimension as in classical case we have to break
$ISO_{q,r}(N)$-bicovariance. 
This naturally leads to a  right $ISO_{q,r}(N)$-covariant
and  $SO_{q,r}(N)$-bicovariant calculus.
The commutation relations 
caracterizing this calculus are explicitly given in the tables at the end 
of the chapter.

\section{Bicovariant  calculus on $ISO_{q,r}(N)$ with $r\not=1$}
\sk
In this section we study, with projection techniques, a differential 
calculus on $ISO_{q,r}(N)$ with $r\not= 1$,
a similar calculus exists also for $ISp_{q,r}(N)$;
for physical reasons we here treat in detail the orthogonal case.

As discussed  at the beginning of Section 4.6, in the $r\not= 1$ case,
the quantum tangent space $T'\equiv T\cup  \IUsosp$   contains dilatations and
translations, but does not 
contain the tangent vectors of $\SqrN$, i.e.
the functionals $\cchi{a}{b}$. However $T'$ defines a bicovariant differential 
calculus on $ISO_{q,r}(N)$ or $ISp_{q,r}(N)$ because conditions (\ref{UNO}) 
and (\ref{adbico})
are satisfied. The proof is as in (\ref{UNOsosp}) and (\ref{adbicososp}).
\sk
The $q$-Lie algebra in the orthogonal case is explicitly given by
\eqa
&&\cchi{\bu}{\ci}\cchi{\bu}{b}-(q_{\bu b})^{-2}
\cchi{\bu}{b}\cchi{\bu}{\ci}=0\label{qlieprim}\\
&&\cchi{\bu}{c}\cchi{\bu}{\bu}-r^{-2}\cchi{\bu}{\bu}
\cchi{\bu}{c}=-r^{-1}\cchi{\bu}{c}\label{qlieMin}\\
&&\cchi{\bu}{\ci}\cchi{\bu}{\bu}-r^{-4}\cchi{\bu}{\bu}
\cchi{\bu}{\ci}={-(1+r^2)\over {r^{3}}}\cchi{\bu}{\ci}\\
&& q_{\bu a} \PA{ab}{cd}\cchi{\bu}{b} \cchi{\bu}{a}=0
\label{quattro}
\ena
Relation (\ref{quattro}) is equivalent to 
$q_{b \bu } \PA{ab}{cd}\cchi{\bu}{b} \cchi{\bu}{a}=0$ and 
$q_{a \bu } [(P_{\!A}){}_{{}_{q^{-1}\!,r^{-1}}}]^{ab}_{~cd}\cchi{\bu}{a} 
\cchi{\bu}{b}=0$.
A combination of (\ref{qlieprim})-(\ref{quattro}) yields:
\eq
\cchi{\bu}{\ci} + \lambda \cchi{\bu}{\ci} \cchi{\bu}{\bu} = \lambda
{-r^{{N\over 2}}\over {r^2 + r^N}}
{1\over q_{d\bu}} \cchi{\bu}{b} C^{db} \cchi{\bu}{d}  \label{chibuci}
\en
Notice the similar structure of eq.s (\ref{PRTT44}),  
(\ref{restrictedL2})
and
(\ref{chibuci}).
\sk
Following the same arguments as in (\ref{Minsospext}),
the adjoint representation is given by the elements
\eq
P({\MM{\bu}{B\bu}{D}}) =
P(\T{\bu}{\bu}\/\kappa_{N+2}(\T{D}{B}))=vP(\kappa_{N+2}(\T{D}{B}))
\label{adMink}
\en
that explicitly read
\eq
\begin{array}{lll}
P(\MM{\bu}{\ci\bu}{\ci})=v^2     &P(\MM{\bu}{\ci\bu}{d})=0
&P(\MM{\bu}{\ci\bu}{\bu})=0     \\
P(\MM{\bu}{b \bu}{\ci})=vr^{-{N\over 2}}x^eC_{eb}      &P(\MM{\bu}{b
\bu}{d})=v\kappa(\T{d}{b})    &P(\MM{\bu}{b \bu}{\bu})=0     \\
P(\MM{\bu}{\bu\bu}{\ci})=-{1\over{r^N(r^{N\over 2}+r^{-{N\over
2}+2})}}x^eC_{ef}x^f
&P(\MM{\bu}{\bu\bu}{d})=v\kappa(x^d)     &P(\MM{\bu}{\bu\bu}{\bu})=I
{}~~~
\label{cccinque}
\end{array}
\en
The differential related to this calculus is given by
\eq
\forall a\in \ISO ~~~~~~da=(\cchi{\bu}{b}*a)\om_{\bu}{}^b
+(\cchi{\bu}{\bu}*a)\om_{\bu}{}^{\bu} +
(\cchi{\bu}{\ci}*a)\om_{\bu}{}^{\ci}\label{dMink}
\en
where $\om_{\bu}{}^{b},~\om_{\bu}{}^{\bu}$ and $\om_{\bu}{}^{\ci}$
are the
$1$-forms dual
to the tangent vectors $\cchi{\bu}{b}$,  $\cchi{\bu}{\ci} $ and
$\cchi{\bu}{\bu}$.
The left and right actions $\DL~:~\Ga\rightarrow ISO_{q,r}(N)\otimes \Ga$
and $\DR~:~\Ga\rightarrow \Ga\otimes ISO_{q,r}(N)$ are defined by:
\eq
\DL(\ome{\bu}{A})=I\otimes \ome{\bu}{A}~~,~~~~\DR(\ome{\bu}{A})=
\ome{\bu}{B}\otimes P({\MM{\bu}{B\bu}{A}})\label{coadjqort}
\en
\sk
We now explicitly give the relation characterizing this differential 
calculus. These formulae will be needed in the next chapter.

To simplify notations, we write the composite indices
as follows:
\eq
{}_{\bu} {}^a \rightarrow {}^a,~{}_{\bu} {}^{\bu} \rightarrow
{}^{\bu},~{}_{\bu} {}^{\ci} \rightarrow {}^{\ci};~~~{}^{\bu} {}_a
\rightarrow
{}_a,~{}^{\bu} {}_{\bu} \rightarrow {}_{\bu},~{}^{\bu} {}_{\ci}
\rightarrow
{}_{\ci}
\en
Similarly we'll write $q_b$ instead of  $q_{ b \bu}  $.
The explicit expression for the tangent vectors then reads:
\eqa
& &\chi_b=\rinv \fmi{\bu}{b}   \nonumber \\
& &\chi_{\ci}=\rinv \fmi{\bu}{\ci} \nonumber \\
& &  \chi_{\bu}=\rinv[\fmi{\bu}{\bu}-\epsi]    \label{chiinso}
\ena
and  their coproduct is given by
\eqa
& &\D (\chi_b)=\chi_{\bu}  \otimes \fmi{\bu}{b}+\chi_c \otimes
\fmi{c}{b} +
 \epsi \otimes \chi_b  \label{copchi1so}\\
& &\D (\chi_{\bu}) = \chi_{\bu} \otimes \fmi{\bu}{\bu} + \epsi
\otimes
\chi_{\bu}
\label{copchi2so}\\
& &\D (\chi_{\ci})=\chi_{\ci} \otimes \fmi{\ci}{\ci} + \chi_{\bu}
\otimes
 \fmi{\bu}{\ci} + \chi_{c}\otimes \fmi{c}{\ci}
+\epsi\otimes\chi_{\ci}
 \label{copchi3}
\ena
\sk
\noi Using the general formula (\ref{omiasosp}) we can
deduce the $\om, T$ commutations for $\ISOqrN$:
\eqa
& &\om^{b} \T{c}{d}= {q_{f} \over r} \Rinv{bf}{ed} \T{c}{f}
\om^{e}\label{VTcommso}\\
& &\om^{b} x^c={q_{b}\over r^2}  x^c \om^{b}
+ \lambda r^{\Nmezzi -1 } q_{d} C^{bd} \T{c}{d}
\om^{\ci}\label{Vxcommso}\\
& &\om^{b} u={r^2\over q_{b}} u~\om^{b}\\
& &\om^{b} v={q_{b}\over r^2} v~\om^{b}\\
& &\om^{\bu} \T{c}{d}= \T{c}{d}
\om^{\bu}\label{tauTcomm}\\
& &\om^{\bu} x^c={1\over r^2}  x^c \om^{\bu}
- \lambda {q_{b}\over r} \T{c}{b} \om^{b}   \label{tauxcomm}\\
& &\om^{\bu} u=r^2 u \om^{\bu}  \\
& &\om^{\bu} v=r^{-2} v \om^{\bu}  \\
& &\om^{\ci} \T{c}{d}= {q^2_{d} r^{-2} \T{c}{d} \om^\ci}
\label{omciTcomm}\\
& &\om^{\ci} x^c=x^c \om^{\ci}   \label{omcixcomm}\\
& &\om^{\ci} u= u \om^{\ci}  \\
& &\om^{\ci} v= v \om^{\ci}  \label{omvcomm}
\ena
The 1-form $\tau \equiv \om^{\bu} \equiv \ome{\bu}{\bu}$
is bi-invariant, and one can check that
$\forall\,a\in A~,~~da=\lam [\tau a - a \tau].$
The exterior
derivative on the generators of $\ISOqrN$  reads:
\eqa
& & d\T{c}{d}=0 \label{trivial} \\
& & dx^c=-q_{b} r^{-1} \T{c}{b} \om^b - r^{-1}  x^c \om^{\bu}
\label{donx}\\
& & du = ru\om^{\bu} \\
& & dv = -r^{-1} v \om^{\bu} \label{donv}\\
& & dz=  - q_{b} r^{-1} y_b \om^b -r(1-r^N) u\om^{\ci} - r^{-1} z
\om^{\bu}
\label{donz}
\ena
\noi where we have included the exterior derivative on $z$ for
convenience.  Note that the calculus is trivial on the $\SOqrN$
subgroup of $\ISOqrN$, as is evident from (\ref{trivial}). Thus
effectively
we are discussing a bicovariant calculus on the orthogonal
$q$-plane generated by the coordinates $x^a$ and the
``dilatations" $u,v$.
\sk
Every element $\rho$ of $\Ga$ can be written as
$\rho=\sum_k a_k db_k$ for some $a_k,b_k$ belonging to
$\ISOqrN$. Indeed inverting  the  relations (\ref{donx})-(\ref{donz})
yields:
\eqa
& & \om^a=-{r\over q_{a}} \kappa (\T{a}{c})[ dx^c - x^c u dv]=
{r\over q_{a}} [d\kappa(x^a)]v=r^{-1}vd\kappa(x^a)
\label{omdT1}\\
& & \om^{\bu} = -r u dv=r^{-1} v du \\
& & \om^{\ci}=-{{vdz + r^{-N}zdv + r^{-\Nmezzi} C_{ab} x^a dx^b}\over
{r(1-r^N)}}
\label{omdT3}
\ena
The exterior product of left-invariant $1$-forms
is as usual defined by
\eq
\om^i \we \om^j\equiv \om^i \otimes \om^j - \L{ij}{kl}
\om^k \otimes \om^l \label{omcomgen}
\en
\noi where
\eq
\L{ij}{kl}=\f{i}{l} (\M{k}{j})
\en
As in  (\ref{27}) this $\Lambda$ tensor is
obtained from
the one of $\SOqrNt$ by restricting its indices to the
subset  {\small $\bu b, \bu\bu, \bu\ci$}.
The non-vanishing components of $\Lambda$ read:
\[
\begin{array}{lll}
\LM{ad}{cb}={q_a \over q_{c}} r^{-1} \R{ad}{bc} &
 \LM{\bu\ci}{cb}= - {r^{-\Nmezzi-1}\over q_{c}} \lambda C_{bc} &
 \LM{\bu d}{c\bu}=r^{-2} \de^d_c \\
\LM{a\ci}{c \ci}=r^{-1} \lambda \de^a_c &
\LM{\ci d}{c \ci}=({r\over q_c})^2 \de^d_c &
\LM{a\bu}{\bu b}=\de^a_b \\
\LM{\bu d}{\bu b}=r^{-1} \lambda \de^d_b &
\LM{a\ci}{\ci b}=r^{-4} (q_a)^2 \de^a_b &
\LM{\bu\bu}{\bu\bu}=1 \\
\LM{a d}{\bu\ci}=-q_ar^{-\Nmezzi-1}\la C^{da} &
\LM{\bu\ci}{\bu\ci}=\la r^{-1}(1-r^{-N}) &
\LM{\ci \bu}{\bu \ci}=1 \\
\LM{\bu\ci}{\ci\bu}=r^{-4}&
\LM{\ci\ci}{\ci\ci}=1& \end{array}
\]
{}From (\ref{omcomgen}) it is not difficult to deduce the
commutations
between the $\om$'s:
\eqa
& & {1\over q_c} \PS{ab}{cd} ~ \om^d \we \om^c=0\label{PSom}\\
& &  \om^a \we \om^{\bu} = -r^2 \om^{\bu} \we \om^{a} \\
& & \om^a \we \om^{\ci} = -r^{-4} (q_a )^2 \om^{\ci} \we \om^{a} \\
& & \om^{\bu} \we \om^{\bu} = \om^{\ci} \we \om^{\ci} = 0 \\
& & \om^{\bu} \we \om^{\ci}=-r^{-4} \om^{\ci} \we \om^{\bu}
+{\lambda r^{-\Nmezzi -1} \over  q_a (1-r^{-N})}C_{ba} ~\om^a \we
\om^b
\label{omomcomm}
\ena
Notice that the dimension of the space of  2-forms generated by
$\om^a \we \om^b$ is
larger than in the commutative case since $P_S$ project into an
$N(N+1)/2-1$
(and not into an $N(N+1)/2$) dimensional space.
This is not surprising since the exterior algebra of
homogeneous orthogonal quantum groups is known to be larger than its
classical counterpart.
\sk
The Cartan-Maurer equations
\eq
d\om^i=\rinv (\tau \we \om^i+\om^i \we \tau)
\en
can be explicitly found after use of the commutations  (\ref{PSom})-
(\ref{omomcomm}):
\eqa
& & d\om^a=r^{-1} \om^a \we \om^{\bu} \\
& & d\om^{\bu}=0 \\
& & d\om^{\ci}=-r(1+r^2) \om^{\bu} \we \om^{\ci} +{r^3 \over
r^{\Nmezzi} -
r^{-\Nmezzi}}  {C_{ba} \over q_a} \om^a \we \om^b
\ena
Finally, the nonvanishing structure constants $\Cb$, 
given by $\C{jk}{i}=\chi_k(\M{j}{i})$,  read:
\[
\begin{array}{lll}
\C{ab}{\ci}=- q_a^{-1}  r^{-\Nmezzi - 1} C_{ba} &
\C{a\bu}{c} = -r^{-1} \de^c_a &
\C{\bu b}{c}=r^{-1} \de^c_b \\
\C{\ci\bu}{\ci}=-r^{-3} (1+r^2)&
\C{\bu\ci}{\ci}=r^{-1} (1-r^{-N})& \end{array}
\]
\noi  These structure constants
 can be obtained from those of $\SOqrNt$ by
specializing indices, for the same reason
as for the $\Lambda$ components.
\sk
The $*$-conjugation on the $\chi$ functionals 
can be deduced from (\ref{conjchiAB}) [use
$(q_{f})^{-1} \Dcal{f}{b}$ $={\overline{q}_{b}}\Dcal{f}{b}$]
\eqa
& &\!\!\!\!\!\!(\chi_b)^*=-r^{-N} 
\Dcal{f}{b}{1\over q_f}D^d_f \chi_d=-r^{-N} {\overline q_b}
\Dcal{f}{b}D^d_f\chi_d=-r^{-N} {\overline q_b}
D^f_b\Dcal{d}{f}\chi_d\label{*46}\\
& &\!\!\!\!\!\!(\chi_{\bu})^*=- \chi_{\bu}\\
& &\!\!\!\!\!\!(\chi_{\ci})^*=-r^{-2N-2}\chi_{\ci}\label{*48}
\ena
whereas the conjugation on the $\om$ $1$-forms can be deduced 
from (\ref{om*chi}) and (\ref{*46})- (\ref{*48}) or directly from
their expression in terms of $dx,du,dv$ differentials
(\ref{omdT1})-(\ref{omdT3}) remembering
that $(da)^*=d(a^*)$:
\eqa
& &(\om^a)^*={\overline q_a}^{-1} r^{N} (D^{-1})^a_{~b}\Dcal{b}{c}
={\overline q_a}^{-1} r^{N} \Dcal{a}{b} (D^{-1})^b_{~c}
\om^c \\
& &(\om^{\bu})^*=\om^{\bu} \\
& &(\om^{\ci})^*=r^{2N+2} \om^{\ci}
\ena


\section{Calculus on the multiparametric orthogonal quantum plane}


In this section we concentrate on the
orthogonal quantum plane
\eq
M\equiv Fun_{q,r}\left( {ISO(N) \over SO(N)}\right)~,
\en
\noi i.e. the $ISO_{q,r}(N)$ subalgebra
generated by the coordinates $x^a$ and the dilatations
$u,v$. This is the algebra we called $B$ in the study of the cross-product
cross-coproduct construction 
$\ISOqrN\cong B{\mbox{$\times \!\rule{0.3pt}{1.1ex}\;\!\!\!\cdot\,$}}\SOqrN$
of Section 4.2. 

We study  the action of the exterior differential $d$ on
$M$  and the corresponding space $\Ga_M$ of $1$-forms. $\Ga_M$ is
the sub-bimodule of $\Ga$ formed by all the elements $adb$ or
$(da')b'$
where $a, b, a', b'$ are polynomials in $x^a,u$ and $v$
[of course $adb=d(ab)-(da)b$].

We will see that a generic element $\rho$ of $\Ga_M$ cannot be
generated,
as a left module, only by the differentials $dx, dv$, i.e.
it cannot be written as $\rho=a_i dx^i+ a dv$. We need also to
introduce
the  differential $dz$ (or equivalently $dL\equiv d(x^aC_{ab}x^b)$).
Thus the basis of differentials  is
given by $dx^a,dv,dz$ and corresponds to the
intrinsic basis of {\sl independent} $1$-forms $\om^a, \om^{\bu}$ and
$\om^{\ci}$.
Note that $du$ can be expressed in terms of $dv$
since $du=-u(dv)u= -r^2u^2dv= -r^{-2} (dv)u^2$ [ see (\ref{udv})
below].

In Subsection 5.2.1 we consistently impose an extra conditon
in order to relate $dz$ to $dx$ and $dv$. This is done 
in two different ways:
checking explicilty the consistency of the extra condition as in 
\cite{WZ} \cite{qplane}
and also deriving it using $\ISO$ symmetry principles.
\sk
{\bf Commutations}
\sk
The commutations between the coordinates $x^a , u$ and $v$ have been
given
in Section 4.2.
The commutations between coordinates and differentials are found
by expressing the differentials in terms of the $1$-forms $\om$ as in
(\ref{donx})-(\ref{donz}), and using then the $x,u,v$ commutations
with the $\om$'s given
in (\ref{VTcommso})-(\ref{omvcomm}). The resulting $q$-commutations
between $x$ and $dx$ are found to be:
\eq
(r^{-2}P_S-P_A)(x\otimes dx)=(P_S+P_A)(dx\otimes x)
\en
where we have used the tensor
notation
$A^{ab}_{~~cd} x^c dx^d \equiv A (x\otimes dx)$
etc. The remaining commutations are
given in formulae (\ref{xdu})--(\ref{zdzcomm}) in Table 1.

Let us consider the above formula, giving the $x,dx$
commutations. If we
multiply it by $P_0$ we find $0=0$. Thus from this equation we have
no
information on $P_0 ( x \otimes dx) $. Applying instead the
projectors
$P_S$ and $P_A$ yields
\eq
P_S (x \otimes dx)= r^2 P_S (dx \otimes x)~;~~~
P_A (x \otimes dx)= - P_A (dx \otimes x)  \label{PSPAxdx}
\en
which does not  allow  to express $x^a dx^b$ only in terms of linear
combinations of  $(dx)x$ since no linear combination of $P_S$ and
$P_A$
is invertible. The space of $1$-forms has therefore one more
dimension
than
his classical analogue because we are missing a condition involving
the one dimensional projector
$P_0^{ab}{}_{\!\!\!ef}=(C^{lm}C_{~lm})^{-1}C^{ab}C_{~ef}$,
see (\ref{projBCD}).

However,  if we consider the $1$-form  $dL\equiv d(x^eC_{ef}x^f)$
-- an exterior derivative of {\sl polynomials} in the basic elements  
--
we can  write  the
commutations
between
the $x$ and $dx$ elements  as follows:
\eqa
dx\otimes x
&=&-(P_S+P_A+P_0)x\otimes dx +(P_S+P_A)d(x\otimes x)+P_0d(x\otimes x)
\nonumber\\
&=&P_S dx\otimes x +P_A dx\otimes x - P_0x\otimes dx +P_0d(x\otimes
x)
\nonumber\\
&=&(r^{-2}P_S-P_A-P_0)x\otimes dx + P_0d(x\otimes x) \label{xdxcomm2}
\ena
where we have used the Leibniz rule, the commutations
(\ref{PSPAxdx})
and  $P_S+P_A+P_0=I$. Equivalently we have
\eq
dx \otimes x=(r^{-2}P_S-P_A-P_0) x\otimes dx
- C {r^{\Nmezzi-2} (1-r^2)\over 1-r^N} (v dz + zdv) \label{xdxcomm3}
\en
\noi involving the $dv$ and $dz$ differentials.
\sk
The presence of $dz$ can also be explained within the general theory
by
recalling
that $\Ga$ is a free right module [see paragraph following
(\ref{eta})].
A basis of right invariant $1$-forms is given by (\ref{eta}):
$\eta^A \equiv \kappa^{-1}(\M{B}{A})\om^B$,
we explicitly have:
\eqa
\eta^{a}&=&-r^{-1} dx^a\,u=-r^{-1} d\T{a}{\bu}\,\kappa(\T{\bu}{\bu})\\
\eta^{\bu}&=&-r^{-1} dv\,u =-r^{-1}  d\T{\bu}{\bu}\,\kappa(\T{\bu}{\bu})
\label{right1} \\
\eta^{\ci} &=& {r^{\Nmezzi-1}\over{(1-r^N)(1+r^{N-2})}}[dx^eC_{ef}x^f  -
                   r^{N-2} x^eC_{ef}dx^f] u^2
\label{righteta}\\ 
        &=&{-r^{N-1}\over {r^N-1}}[dz\,u +dy_b\,\kappa(x^b)+du\,\kappa(z)]
         = {-r^{N-1}\over {r^N-1}} d\T{\ci}{B}\,\kappa(\T{B}{\bu})
\ena
To derive the expressions for $\eta^{\ci}$ use:
$y_b=-r^{-\Nmezzi}ux^eC_{ef}\T{f}{b};$
$dy_b\,\kappa(x^b)=r^{-\Nmezzi}du\,xxu +r^{-\Nmezzi}dx\,xu^2;$
$dz=d({-r^{-\Nmezzi}\over{1+r^{2-N}}}uxx)$
$={-r^{-\Nmezzi}\over{1+r^{2-N}}}(du\,xx+
dx\,xu+xdx\,u)$;
$\kappa(z)=\kappa(\T{\ci}{\bu})=r^{-N}z$; $uxx=r^2xxu$;
$udx\,x=dx\,xu$,
where $xx\equiv L\equiv x^eC_{ef}x^f$.

The  $1$-forms (\ref{right1})-(\ref{righteta})
in $\Ga$ do not contain any $\T{a}{b}$ element
and therefore belong to $\Ga_M$ as well; they are linearly
independent
and freely generate $\Ga_M$ as a right module because they freely
generate
the full $\Ga$ as a right module. The extra $1$-form $\eta^{\ci}$ (or $dz$)
is therefore a natural consequence of the right module structure of
$\Ga$.
\sk
In summary: either
$dL$
or $dz$ or $\eta^{\ci}$ are necessary in order to close the commutation
algebra
between
coordinates and differentials. Thus the commutations involving $z$
and
$dz$ appear in Table 1.
\sk
We have seen that $dv\,u$; $dx^a\,u$ and $\eta^{\ci}$ freely generate
$\Ga_M$ as
a right module; recalling that $\Ga$ is also a free left module, we
have the :
\sk
\noi {\bf Proposition 5.1} The $M$-bimodule $\Ga_M$, as a left module
(or as a right module), is freely generated by the differentials
$dx$, $dL$ (or $dz$)
and $dv$.
\noi{\sl Proof:}  to show that
$a_idx^i+adL+a_{\bu}dv=0$ $\Rightarrow a_i=0, a=0, a_{\bu}=0$
express $dx^i, dL, dv$ in terms of $\om^a,\om^{\ci}\om^{\bu}$,
see (\ref{donx})-(\ref{donz}), and recall that $\Ga$ is a free left
module.
\sk
\noi {\bf Note } 5.2.1 $~$
{}From (\ref{xdu}), (\ref{xdz}) and the commutations of $L$ with $x$
and
$u$
we have $x^cdL=dL\,x^c$,
$udL=dL\,u$ and $vdL=dL\,v$.
These relations and (\ref{xdxcomm}) show that inside $\Ga_M$ there is
the
smaller bimodule generated by the differentials $dx^a$ and $dL$.
\sk
We now  examine the space of 2-forms.
By simply applying the exterior derivative $d$ to the relations
(\ref{xdxcomm})-(\ref{zdzcomm})  we deduce
the commutations between the differentials given in Table 1.
As with the $\om^a$'s in eq. (\ref{PSom}), the relations in
(\ref{dxdxcomm})
are not
sufficient to order the differentials $dx^a$.
\sk
{\bf $ISO_{q,r}(N)$ - coactions}
\sk
All the relations we have been deriving have
many symmetry properties because they are covariant,
under  the actions on $M$ and $\Ga$, of the full  $ISO_{q,r}(N)$
$q$-group. In fact we have the following three $ISO_{q,r}(N)$
actions:

{\bf 1)} the coproduct of $ISO_{q,r}(N)$ can be seen as a left-coaction
$\D :~M\rightarrow ISO_{q,r}(N)\otimes M$:
\eq
\D x^a= \T{a}{b}\otimes x^b + x^a\otimes v~,~~~~
\D(u)=u\otimes u~,~~~~\D(v)=v\otimes v
\label{1)} 
\en

{\bf 2)} the left coaction
$\DL:\Ga\rightarrow ISO_{q,r}(N)\otimes \Ga$,
when restricted to $\Ga_M$ gives
\eq
\left.\DL\right|_{{}_{\Ga_M}}:\Ga_M\rightarrow ISO_{q,r}(N)\otimes
\Ga_M\label{2)}
\en
and defines a left coaction of $ISO_{q,r}(N)$ on $\Ga_M$ compatible  
with
the bimodule structure of $\Ga_M$ and the exterior differential:
$\left.\DL\right|_{{}_{\Ga_{\!M}}}(adb)
=\D(a)(id\otimes d)\D(b)$.

{\bf 3)} the right coaction $\DR:\Ga\rightarrow \Ga\otimes ISO_{q,r}(N)$,
does not become a right coaction of $ISO_{q,r}(N)$ on
$\Ga_M$; however we have
\eq
\left.\DR\right|_{{}_{\Ga_{\!M}}}:\Ga_M\rightarrow \Ga\otimes  
M\,\subset
\Ga\otimes ISO_{q,r}(N)\label{3)}
\en
this map is obviously well defined and satisfies
$\left.\DR\right|_{{}_{\Ga_{\!M}}} (a db)=
\D(a)(d\otimes id)\D(b)$
$\forall a, b \in M$ since $M\subset ISO_{q,r}(N)$.

We call this calculus $ISO_{q,r}(N)$-bicovariant
because $\left.\DL\right|_{{}_{\Ga_{\!M}}}$ and
$\left.\DR\right|_{{}_{\Ga_{\!M}}}$ are compatible with the
bimodule structure
of $\Ga_M$ and with the exterior differential.
\sk
\subsection{$\ISO$-covariant and $SO_{q,r}(N)$-bicovariant
calculus}
\sk
{\bf Commutations}
\sk
Since the $\PA{ab}{cd} x^c x^d=0 $ commutation relations allow for an
ordering of the coordinates (moreover the Poincar\'e series of the
polynomials on the quantum orthogonal plane is the same as the
classical
one), it is tempting to impose extra conditions on the differential
algebra
of the $q$-Minkowski plane, so that the space of $1$-forms has the  
same
dimension as in the classical case. We require that the
commutation relations between $x$ and $dx$ close on the algebra
generated
by $x$ and $dx$:
\eq
dx^a x^b= \alpha^{ab}{}_{\!\!ef}x^edx^f \label{dxxalfa}
\en
where $\alpha$ is an unknown matrix whose entries are complex
numbers.
Any matrix can be expanded as $\al=aP_S + bP_A + cP_0$ with
$a,b,c=const$.
{} From (\ref{PSPAxdx}) we have $\alpha=r^{-2}P_S-P_A+cP_0$;
 therefore condition (\ref{dxxalfa})
is equivalent to
\eq
P_0(dx\otimes x)=
 c P_0(x\otimes dx)\label{P0dxx}
\en
and supplements eq.s  (\ref{PSPAxdx}). Taking its exterior derivative
leads to a supplementary condition on the $dx,dx$ products (for $c
\not= -1$):
\eq
P_0 (dx \we dx)=0 \label{P0dxdx}~.
\en
{}From  (\ref{dxdxcomm}) and (\ref{P0dxdx}) it follows that
$dx \we dx = (P_S + P_A + P_0) (dx \we dx )= P_A (dx \we dx)$, or
[see the definition of $P_A$ in (\ref{projBCD})] :
\eq
dx\we dx=-r\Rha ~dx\we dx\label{dxdxcomm2}~.
\en
\noi which allows the ordering of $dx,dx$ products.

Using (\ref{xdxcomm2}), (\ref{P0dxx}) and (\ref{RprojBCD}), we find
\eqa
dx\otimes x
& =& (r^{-2}P_S-P_A) (x\otimes dx) + P_0 (dx \otimes x)\\
& =& (r^{-2}P_S-P_A) (x\otimes dx) + c P_0 (x \otimes dx)\\
&=& (r^{-2}P_S-P_A+r^{N-2}P_0) (x\otimes dx) + (c-r^{N-2}) P_0 (x
\otimes
dx)\\
&=& r^{-1}{\hat{R}}^{-1}(x\otimes dx)
+ (c-r^{N-2})P_0(x\otimes dx)
\label{xdxP0}~.
\ena
The consistency of the
commutation relations (\ref{dxdxcomm2}) and (\ref{xdxP0}) with the
associativity condition on the triple $dx^i\,dx^j\,x^k$ fixes
 $c=r^{N-2}$ i.e.:
\eq
P_0(dx\otimes x)=
r^{N-2} P_0(x\otimes dx)\label{P0dxx2}~;
\en
the  $x,dx$ commutations (\ref{xdxP0}) then become:
\eq
x\otimes dx = r \Rha (dx \otimes x) \label{xdxcomm4}
\en
\noi and reproduce (in the uniparametric case) the
known $x,dx$ commutations of the quantum orthogonal plane
\cite{Fiore1}.
\sk
{\bf Coactions}
\sk
This calculus is no more bicovariant under the $\ISOqrN$ action,
\eq
x^a \longrightarrow \T{a}{b} \otimes x^b~+x^a\otimes v~,~~~~ u
\longrightarrow u
\otimes u~,~~~~v \longrightarrow  v \otimes v \label{ISOcov}
\en
but we are left with bicovariance under the $SO_{q,r}(N)$ action
\eq
x^a \longrightarrow \T{a}{b} \otimes x^b~. \label{SOcov}
\en
In other words, $\delta_L : \Ga_M'\rightarrow SO_{q,r}\otimes \Ga_M'$
defined by $\delta_L(adb)=\delta(c)(id\otimes d)\delta(b)$
with $\delta(x^a)=\T{a}{b}\otimes x^b$ is a left coaction of
$SO_{q,r}(N)$ on the bimodule $\Ga_M'$ where $\Ga_M'$ is $\Ga_M$
with the extra condition (\ref{P0dxx}) [cf. (\ref{2)})].
Similarly, the map $\delta_R(adb)=\delta(a)(d\otimes id)\delta(b)$
is well defined [cf. (\ref{3)})].

Left covariance under (\ref{ISOcov}) is broken only by (\ref{P0dxx}).
Indeed, while relations (\ref{PSPAxdx}) are left and right
$ISO_{q,r}(N)$-covariant,
the extra condition (\ref{P0dxx}) is {\sl not}
left $ISO_{q,r}(N)$-covariant :
$\DL[P_0(dx\otimes x)- c P_0(x\otimes dx)]\not=0, \forall c$. It is
right $ISO_{q,r}(N)$-covariant,
$\DR[P_0(dx\otimes x)- c P_0(x\otimes dx)]=0$,
only for $c=r^{N-2}$, as can be seen
using
$\T{b}{d}dx^a=d(\T{b}{d}x^a)={r\over {q_{d}}}R^{ab}_{~ef}dx^e\,
\T{f}{d}$
and (\ref{CR}).
Therefore the choice $c=r^{N-2}$ preserves the right  
coaction $\DR$ i.e. the  right $\ISO$-covariance. 
\sk
\noi {\bf Note }5.2.2 $~$  We can reformulate the quotient
procedure $\Ga_M\rightarrow \Ga_M'$
in a more abstact setting by considering that $\Ga_M$ is a
subbimodule of
the {\sl bicovariant} bimodule $\Ga$.
In  (\ref{righteta}) we have expressed the
$x^eC_{ef}dx^f\leftrightarrow
dx^eC_{ef}x^f$ commutation via the right invariant $1$-form $\eta^{\ci}$.
A condition on $\Ga$ (and therefore on $\Ga_M$) that preserves
{\sl right $ISO_{q,r}(N)$ covariance}, i.e. compatible with $\DR$,
 is: $\eta^{\ci}$ linearly dependent
from the remaining right invariant $1$-forms: $dv\,u$ and $dx^a\,u$.
It is
easily seen that since $\eta^{\ci}$ is quadratic in the basis elements
$x^a$ the
only possible linear condition is $\eta^{\ci}=0$, and this gives exactly
(\ref{P0dxx2}).
The $M$-bimodule $\Gamma_M'$ is therefore
generated by the differentials $dx^b$ and
$dv$. Since left $ISO_{q,r}(N)$ covariance, contrary to right
$ISO_{q,r}(N)$ covariance, is broken, the relations between the left
invariants $1$-forms is nonlinear. Explicitly we have
\eq
\om^{\ci}=-{q_a \over r^2} vy_a ~\om^a+{r^{\Nmezzi-2} \over r^N +
r^2}
C_{ab} x^a x^b \om^{\bu}
\en
\noi [express $dz$ in terms of $dx^i$,$dv$ in (\ref{omdT3}) and use
the expansion of $dx^b$ and $dv$ on $\om^a$ and $\om^{\bu}$ as given
in (\ref{donx}),(\ref{donv})].
\sk
{\bf Partial derivatives}
\sk

The tangent vectors $\chi$ in (\ref{chiinso})
and the corresponding  vector fields $\chi *$
 have ``flat" indices. To compare $\chi *$
with partial derivative operators with ``curved" indices,
we need to define the operators $\stackrel{\leftarrow}{\chit}$ 
such that
\eq
d\,a= \stackrel{\leftarrow}{\chit_c}\!(a) ~dx^c + \parleft_{\bu} (a)
dv \equiv \parleft_{C} (a) dx^C \label{dacoord1}
\en
\noi [{\small $C=(c,\bu)$}, $dx^C = (dx^c,dv)$].
The action of $\stackrel{\leftarrow}{\chit_C}$ on the coordinates
$x^C = (x^c,v)$ is given by
\eq
\stackrel{\leftarrow}{\chit_C}\!(x^A)=\de^A_C I~,
\en
{}From the Leibniz rule $d(ab)=(da)b+a(db)$,
 using (\ref{dacoord1}) and the
fact that  $dx^C = (dx^c,dv)$  is a basis for 1-forms, we
find
\eqa
& &\parleft_c (ax^b)=a\delta^b_c + \parleft_d (a) r^{-1}
(\Rha^{-1})^{db}_{~~ec}
x^e-(1-r^2) \parleft_{\bu} \de^b_c v \label{Leibpart2}\\
& &\parleft_{\bu} (ax^b)=q_b^{-1} \parleft_{\bu} (a) x^b
\ena
\eqa
& &\parleft_c(av)=r^{-2} q_c \parleft_c(a) v\\
& &\parleft_{\bu} (av)= r^{-2} \parleft_{\bu} (a) v+a
\ena
Note the dilatation operator $\parleft_{\bu}$ appearing on the
right-hand side of (\ref{Leibpart2}).

{}From $d^2(a)=0=d(\parleft_C (a) dx^C)=\parleft_B (\parleft_C (a))
dx^B \we dx^C$ and the $q$-commutations of the differentials
(\ref{dxdxcomm})-(\ref{dzdzcomm})
one finds the commutations between the ``curved" partial
derivatives:
\eqa
& & (P_A)^{ab}_{~~cd} \parleft_a \parleft_b=0 \label{parleftcomm} \\
& & \parleft_b \parleft_{\bu} - { q_b\over r^2 } \parleft_{\bu}
\parleft_{b}=0 
\ena
\sk
We can also define the partial derivatives $\partial_C$ so that
\cite{WZ}, \cite{qplane},
\eq
d\,a=dx^C\;\chit_C(a)~;\label{ioio}
\en
again the action  of $\chit_{{}_{\scriptstyle C}}$ on the coordinates is
\eq
\chit_C (x^A)=\de^A_C I~.
\en
We now give an explicit relation between the
$\partial_C$ and the $q$-Lie algebra generators $\chi_{{}_{\scriptstyle C}}$
(a similar expression holds also for the $\parleft$ derivatives).
{}From (\ref{etabbb}) we have:
\eq
d\,a=-\eta^C(a*\kp(\chi_{{}_C}))\label{etakchi}
\en
where {\small $C=(c,\bu)$} because we have set $\eta^{\ci}=0$.
Relations (\ref{ioio}) and (\ref{etakchi})
give, $\forall a\in \ISO$  :
\eq
\chit_c(a)=r^{-1}u(a*\kp(\chi_c))~~~,~~~~~~
\chit_{\bu}(a)=r^{-1}u(a*\kp(\chi_{\bu}))
\label{pechi}
\en
The commutations between the partial derivatives can be derived as 
done above for $\parleft$, or via (\ref{pechi}) and the $q$-Lie algebra 
(\ref{qlieprim})-(\ref{quattro}). They are given in Table 2.
Similarly we can introduce the right invariant 
vectorfields 
\eq
\tilde{h}_C\equiv \tilde{h}_{\kp(\chi_{{}_C})}\equiv 
[\kp(\chi_{{}_C})\otimes id] \D 
\en
and use their Leibniz rule 
[it follows from  $\D(\kp(\chi_{{}_C}))=\kp(\chi_{{}_C})\otimes\epsi +
\kp(\f{D}{C})\otimes \kp(\chi_{{}_D})\,$]:
\eq
\tilde{h}_C(ab)=\tilde{h}_C(a)b+\kp(\f{D}{C})(a_1)\,a_2\tilde{h}_D(b)
\en
to derive the $\partial, x,u$ commutations. 
For example we have $\tilde{h}_ax^b=rv\delta^b_a+(r/q_b) R^{lb}_{~ac}x^c
\tilde{h}_l$ $+r\lambda\tilde{h}_{\bu}$ that 
togheter with $\partial_C=r^{-1}u\tilde{h}_C$
gives 
\[\chit_ax^b=\delta^i_j[I+(r^2-1)v\partial_{\bu}]
+rR^{lb}_{~ac}x^c\chit_l~.
\]
Similarly for the
the other relations, see Table 2.
\sk

{\bf Conjugation}
\sk
The commutations in Table 2 are consistent under the conjugation
(already defined for $x^a$ and $dx^a$)
\eq
(x^a)^*=\Dcal{a}{b}x^b,~(dx^a)^*=\Dcal{a}{b}dx^b,~(\partial_a)^*=
-r^Nd^{-1}_b \Dcal{b}{a}  \partial_b
\en
\eq
v^*=v,~ (dv)^*=dv,~(\partial_{\bu})^*=u- \partial_{\bu}
\en
\noi where we have used the notation $D^a_a=d^a\,,~ D^{-1}{}^a_a=d^{-1}_a$
($D^a_b=C^{ae}C_{be}$ is diagonal).
This consistency can be checked directly by taking the $*$-conjugates of
the relations in Table 2, and by using the identity (\ref{Rprop2})
and:
\eq
{\bar C}=C^T;~ [Q_N(r)]^*=Q_N(r);
\en
\eq
\begin{array}{ll}
{}~d^c d^{-1}_h R^{cg}_{~~ha} (R^{-1})^{ea}_{~~cd}=
\de^e_h \de^g_d; &~ R^{ab}_{~~cd} d^a d^b=R^{ab}_{~~cd} d^c d^d
\nonumber\\
d^c R^{cg}_{~~hc}=r^{N-1} \de^g_h; &~
R^{ab}_{~~cd} d^{-1}_a d^{-1}_b=R^{ab}_{~~cd} d^{-1}_c d^{-1}_d
\nonumber
\end{array}
\en
\eq
\overline{q_a}={1 \over q_a} \mbox{for $a \not=n,n+1$},~~\overline{q_n}=
{1\over q_{n+1}}
\en
We now {\sl derive}
the conjugation on the partial derivatives from the differential calculus on
$\ISO$. This is achieved by studying the conjugation on the right invariant 
vectorfields $\tilde{h}$.
\sk
For a general Hopf algebra, with tangent vectors $\chi_i$, we  deduce 
the conjugation on $\tilde{h}$ 
from the commutation relations between $\tilde{h}$ and a generic element
of the Hopf algebra:
\eq
\tilde{h}_jb=\tilde{h}_j(b)+\le\kp(\f{s}{j})\,,\,b_1\re
\,b_2\tilde{h}_s~=\le \kp(\chi_j)\,,\,b_1\re\, b_2+\le\kp(\f{s}{j})\,,\,b_1\re
\,b_2\tilde{h}_s~
\en 
We multiply this expression by $\le {\kp}^2(\f{j}{i})\,,\,b_0\re$ [where
we have used the notation $(id\otimes \D)\D(b)=b_0\otimes b_1\otimes b_2$]
to obtain 
\eq
\le\kp^2(\f{j}{i})\,,\,b_1\re\,\tilde{h}_jb_2+
\le\kp^2(\chi_i)\,,\,b_1\re b_2=b\tilde{h}_i
\en
Now, using $\le\psi\,,\,b\re=\overline{\le[\kp(\psi)]^*,b^*\re}$
and then applying $*$ we obtain (here $a=b^*$)
\eq
\tilde{h}_j^*a=\le [\kp^3(\chi_j)]^*\,,
\,a_1\re\, a_2+\le[\kp^3(\f{s}{j})]^*\,,\,a_1\re
\,a_2\tilde{h}_s~.
\en
This last relation implies
\eq
\tilde{h}_i^*\equiv [\tilde{h}_{\kp(\chi_i)}]^*
=\tilde{h}_{[\kp^3(\chi_i)]^*}\label{genfor}
\en
notice that $*\sma{$\circ$}\kp^2$ is a  well defined conjugation since 
$(*\sma{$\circ$}\kp^2)^2=id$ .
\sk
We now apply  formula (\ref{genfor}), valid for a generic Hopf algebra, 
to the  $*$-conjugation
and the right invariant vectorfields  of this section; we have:
\eqa
& &[\tilde{h}_{\kp(\chi_a)}]^*=-r^N\overline{q}_ad^{-1}_a\Dcal{b}{a}
\tilde{h}_{\kp(\chi_b)}\\
& &[\tilde{h}_{\kp(\chi_{\bu})}]^*=-
\tilde{h}_{\kp(\chi_{\bu})}~.
\ena
From these last relations and (\ref{pechi}) we then finally deduce
$(\partial_a)^*=
-d^{-1}_b \Dcal{b}{a} r^{N} \partial_b$ and 
$(\partial_{\bu})^*=u- \partial_{\bu}~.$
  
\subsection{The reduced 
$x^a,dx^a, \partial_a$ algebra and the quantum \\
Minkowski phase-space.}
\sk
Note that the algebra in Table 2 actually contains a subalgebra
generated only by $x^a,dx^a, \partial_a$, indeed $\partial_{\bu}$
vanishes when acting on monomials containing only the coordinates
$x^b$,
as can be seen from (\ref{Leibpart5}). 
This calculus is $\ISO$-right covariant because it can also be obtained 
imposing the conditions $\eta^{\bu}=0$ and $\chi_{\bu}=0$ that are compatible 
with the right coaction $\DR$ and the bimodule structure given by the $f^i_j$
functionals.
\sk
Table 3 contains the multiparametric orthogonal quantum plane algebra
of coordinates, differentials and partial derivatives, together with
a consistent conjugation for any even dimension. We emphasize here
that this conjugation {\sl does not } require an additional scaling
operator
as in ref. \cite{Lorek}. Thus the algebra in Table 3 can be taken
as
starting point for a deformed Heisenberg algebra (i.e. a deformed
phase-space) .
\sk
{\bf Real coordinates and hermitean momenta}
\sk
We note that the transformation
\eqa
& &X^a={1\over \sqrt{2}}(x^a+x^{a'})\,,~~~\sma{$a \leq n$}\label{54ii}\\
& &X^{n+1}={i\over \sqrt{2}}(x^n - x^{n+1})\\
& &X^a={1\over \sqrt{2}}(x^a-x^{a'})\,,~~~\sma{$a > n+1$}\label{56ii}
\ena
\noi defines real coordinates $X^a$.
 On this basis the metric becomes $C'=(M^{-1})^T C M^{-1}$
(where $M$ is the transformation matrix $X=Mx$):
\eq C'=
{\small
{1\over 2}  \left (\matrix{
r^{\Nmezzi-1}+r^{-\Nmezzi+1}& 0 &{}&{}&0&
-(r^{\Nmezzi-1}-r^{-\Nmezzi+1})\cr
0 & r^{\Nmezzi-2}+r^{-\Nmezzi+2}&{}&{}&
-(r^{\Nmezzi-2}-r^{-\Nmezzi+2})&0 \cr
 \vdots & \vdots  &{}&{}&\vdots & \vdots\cr
0&0&2&0&0&0\cr
0&0&0&2&0&0\cr
\vdots & \vdots &{}&{} & \vdots & \vdots\cr
0 & r^{\Nmezzi-2}-r^{-\Nmezzi+2} &{}&{}&
-(r^{\Nmezzi-2}+r^{-\Nmezzi+2})&0 \cr
r^{\Nmezzi-1}-r^{-\Nmezzi+1}& 0 &{}&{}&0&
-(r^{\Nmezzi-1}+r^{-\Nmezzi+1})\cr
}\right ) }
\label{Cprimo=}
\en
and reduces for $r\rightarrow 1$ to the usual $SO(n+1,n-1)$
diagonal metric with $n+1$  plus signs and $n-1$ minus signs.
Notice that the diagonal elements of $C'$ are real while the off diagonal ones 
are pure imaginary, moreover $C'$ is hermitian (and can therefore be 
diagonalized via a unitary matrix).
\sk
As for the coordinates $X$, it is possible to define antihermitian $\chi$ and
$\partial$, and real
$\om$ and $dx$. 
To  define hermitian momenta 
we first notice that the partial derivatives
\eq
\tilde\partial_a\equiv r^{N\over 2} d_a^{-{1\over 2}}\partial_a
\en
behave,  under the hermitian conjugation $*$, similarly to the 
coordinates $x^a$: 
\eqa
& &(\tilde{\partial}_a)^* =-\tilde{\partial}_a~~~~~\sma{$a\not= n,\,n+1$}\\
& &(\tilde{\partial}_n)^* =-\tilde{\partial}_{n+1}~.
\ena
As in (\ref{54ii})--(\ref{56ii}) we then define:
\eqa
& &P_a={-i\hbar\over \sqrt{2}}(\tilde\partial_a+\tilde\partial_{a'})
\,,~~~\sma{$a \leq n$}\\
& &P_{n+1}={-\hbar\over \sqrt{2}}(\tilde\partial_n-\tilde\partial_{n+1})\\
& &P_a={-i\hbar\over \sqrt{2}}(\tilde\partial_a-\tilde\partial_{a'})
\,,~~~\sma{$a > n+1$}
\ena

It is easy to see that the $P_a$ are hermitian: ${P_a}{}^*=P_a$ and that
in the classical limit are the momenta conjugated to the coordinates $X^a$:
$P_a(X^b)=-i\hbar\delta_a^b$. In the $r\not=1$ case we explicitly have
(use $d_{a'}=d_a^{-1}\,,~ d_n=d_{n+1}=1$):
\[
\begin{array}{l}
P_a(X^a)={-i\over 2}r^{N\over 2}
\hbar(d_a^{{1\over 2}}+d_a^{-{1\over 2}})\\[1em]
P_a(X^{a'})=\epsilon_a{-i\over 2}r^{N\over 2}
\hbar(d_a^{{1\over 2}}-d_a^{-{1\over 2}})
~~\sma{ where $\epsilon_a=1$ if $a < n$ and
$\epsilon_a=-1$ if $a > n+1$}
\end{array}
\]
while the other entries of the  $P_a(X^b)$ matrix are zero.
\sk
By defining the transformation matrix $N_{a}^{~b}$ as:
\eq
P_a \equiv -i\hbar N_a^{~b} \partial_b
\en
\noi  we find the deformed canonical commutation relations:
\eq
P_a X^b -r S^{bc}_{~~ad} X^d P_c = - i \hbar E_a^b I \label{PXcomm}
\en
\noi where
\eq
S^{bc}_{~~ad}=N_a^{~e} M^b_{~f} \Rhat{fh}{eg} (M^{-1})^g_{~d}
(N^{-1})_h^{~c},
{}~~~E_a^b \equiv {i\over \hbar} P_a(X^b) = N_a^{~c} M^b_{~c}
\en
Similarly one finds all the remaining commutations of the 
$P$, $X$ and $dX$ algebra.
Notice that no unitary operator appears on the right-hand
side of  (\ref{PXcomm}). Our conjugation is consistent with
(\ref{PXcomm})
without the need of the extra operator of  ref.  \cite{Lorek} .
\sk
For $n=2$ the results of this section immediately yield the
bicovariant calculus on the quantum Minkowski
space, i.e. on the
multiparametric orthogonal quantum plane
$Fun_{q,r}(ISO(3,1)/SO(3,1))$.
\sk
\noi{\bf Note }5.2.3 $~$  In the $\rone$ limit the reduced differential
calculus on the coordinates $x^a$ coincides with the $r=1$ bicovariant
differential calculus on $ISO_{q,r=1}(N)$ of Section 4.7 [see (\ref{isp})]. 
This is so 
because 
the $ISO_{q,r=1}(N)$ bicovariant differential can be
written $da=(\cchi{a}{b}*a)\omega_a{}^b+
(\cchi{\bu}{c}*a)\omega_{\bu}{}^c=$ 
$-\eta_a{}^b(a*\kp(\cchi{a}{b}))-
\eta_{\bu}{}^c(a*\kp(\cchi{\bu}{c}))$. 
Similarly to Theorem 3.7.1, we have that
$(a*\kp(\cchi{a}{b}))=0$ when $a$ is
a polynomial in $x^a$. Then the exterior differential
on such polynomials reads
$da=-\eta_{\bu}{}^c(a*\kp(\cchi{\bu}{c}))$ as in the reduced differential
calculus on the coordinates $x^a$.
\sk
\sk
\section {Table 1: the $ISO_{q,r}(N)$-bicovariant
algebra}
\sk

\eqa
& &\PA{ab}{cd} x^c x^d=0 \\
& &x^b v=q_{b} v x^b~;~~~x^b u=q^{-1}_{b} u x^b \\
& &z=-{1\over {(r^{-\Nmezzi}+ r^{\Nmezzi-2})}} x^b C_{ba} x^a u \\
& &zv=r^2vz~;~~~zu=r^{-2}uz\\
& &q_a x^a z = z x^a
\ena

\eqa
& &
(x \otimes dx)=  (r^2 P_S-P_A-P_0 )(dx \otimes x) +P_0 d(x\otimes x)
 \label{xdxcomm}\\
& & x^c du= {1\over q_c} (du) x^c - {\la \over r} (dx^c) u ;~~
x^c dv= q_c (dv)x^c +\la r (dx^c) v \label{xdu}\\
& & x^c dz={1\over q_c} (dz) x^c\label{xdz}\\
& & \nonumber\\
& & u dx^c = {q_c \over r^2} (dx^c) u\label{udx}\\
& & u du=r^{-2} (du) u;~~u dv = r^{-2} (dv) u\label{udv}\\
& & u dz=(dz) u\label{udz}\\
& &\nonumber\\
& & v dx^c = {r^2 \over q_c } (dx^c) v\label{vdx}\\
& & v du=r^{2} (du) v;~~v dv = r^{2} (dv) v\label{vdu}\\
& & v dz=(dz) v\\
& & \nonumber\\
& & z dx^c = q_c (dx^c) z \label{zdx}\\
& & z du =  r^{-2} (du) z + (r^{-2} - 1) (dz)u
\label{zdu}\\
& & z dv = r^2 (dv) z  + (r^2-1) (dz)v \label{zdv}\\
& & z dz=r^{-2} (dz) z \label{zdzcomm}
\ena

\eqa
& & P_S (dx \we dx) =0 \label{dxdxcomm}\\
& & dx^c \we du=-{r^2\over q_c} du \we dx^c;~~dx^c \we dv=-{q_c\over
r^2} dv
\we dx^c\\
& & dx^c \we dz=-{1\over q_c} dz \we dx^c \\
& & du \we du=dv \we dv = 0\\
& & du \we dv = - r^{-2} dv \we du \\
& & dz \we du = -du \we dz;~~dz \we dv = - dv \we dz\\
& & dz \we dz = 0 \label{dzdzcomm}
\ena

\section {Table 2: the $ISO_{q,r}(N)$-covariant  
$x^a,v,\partial_a,\partial_{\bu},dx^a,dv$ algebra}

\eqa
& &\PA{ab}{cd} x^c x^d=0 \\
& &x^b v=q_{b} v x^b
\ena

\eqa
& & x \otimes dx=  r\Rha (dx \otimes x)  \\
& & x^c dv= q_c (dv)x^c +\la r (dx^c) v \\
& & v dx^c = {r^2 \over q_c } (dx^c) v
\ena

\eqa
& & dx \we dx =-r {\hat R} dx \we dx \\
& & dx^c \we dv=-{q_c\over  r^2} dv \we dx^c\\
& & dv \we dv = 0
\ena

\eqa
& &\partial_c x^b= r \Rha^{be}_{~~cd} x^d \partial_e+
\delta^b_c [ I +(r^2-1) v \partial_{\bu}] \label{Leibpart4}\\
& &\partial_{\bu} x^b=q_b x^b \partial_{\bu}  \label{Leibpart5}\\
& &\partial_{\bu} v= r^{2} v \partial_{\bu} +I
\ena

\eqa
& & (P_A)^{ab}_{~~cd} \partial_b \partial_a=0 \label{partialcomm2}\\
& & \partial_b \partial_{\bu} - {q_b \over r^2} \partial_{\bu}
\partial_{b}=0
\ena

{\bf Conjugation:}
\eq
(x^a)^*=\Dcal{a}{b}x^b,~(dx^a)^*=\Dcal{a}{b}dx^b,~(\partial_a)^*=
-r^Nd^{-1}_a \Dcal{b}{a}  \partial_b
\en
\eq
v^*=v,~ (dv)^*=dv,~(\partial_{\bu})^*=u-\partial_{\bu}
\en
\eq
r^*=r^{-1},~~q_a^*={1 \over q_a} \mbox{for $a
\not=n,n+1$},~~q_n^*={1\over q_{n+1}}
\en

\sk\sk
\section{Table 3: the reduced $ISO_{q,r}(N)$-covariant
$x^a,\partial_a,dx^a$
algebra}

\eq
\PA{ab}{cd} x^c x^d=0
\en
\eq
 x \otimes dx=  r\Rha (dx \otimes x)
\en
\eq
dx \we dx =-r {\hat R} (dx \we dx)
\en
\eq
\partial_c x^b= r \Rha^{be}_{~~cd} x^d \partial_e+
\delta^b_c I
\en
\eq
 \PA{ab}{cd} \partial_b \partial_a=0
\en
\sk
{\bf Conjugation:}
\eq
(x^a)^*=\Dcal{a}{b}x^b,~(dx^a)^*=\Dcal{a}{b}dx^b,~(\partial_a)^*=
-r^Nd^{-1}_a \Dcal{b}{a}  \partial_b
\en

\chapter{Appendix}

\app{The Hopf algebra axioms}

A Hopf algebra over the field $K$ is a
unital algebra over $K$ endowed
with the  linear maps:
\eq
\D~: ~~A \rightarrow A\otimes A,~~
\epsi~:~~ A\rightarrow K,~~
\kappa~:~~A\rightarrow A
\en
\noi satisfying the following
properties $\forall a,b \in A$:
\eq
(\D \otimes id)\D(a) = (id \otimes \D)\D(a) \label{axiom1}
\en
\eq
(\epsi \otimes id)\D(a)=(id \otimes \epsi)\D(a) =a \label{axiom2}
\en
\eq
m(\kappa\otimes id)\D(a)=m(id \otimes \kappa)\D(a)
=\epsi(a)I \label{kappadelta}
\en
\eq
\D(ab)=\D(a)\D(b)~;~~\D(I)=I\otimes I
\en
\eq
\epsi(ab)=\epsi(a)\epsi(b)~;~~\epsi(I)=1
\en
\noi where $m$ is the multiplication map $m(a\otimes b)
= ab$.
{} From  these axioms we deduce:
\eq
\kappa(ab)=\kappa(b)\kappa(a)~;~~\D[\kappa(a)]=\tau(\kappa\otimes
\kappa)\D(a)~;~~\epsi[\kappa(a)]=\epsi(a)~;~~\kappa(I)=I
\en
where $\tau(a \otimes b)=b\otimes a$ is the twist map.

\sk

\sect{The derivation of two equations}

In this Appendix we derive the two equations (\ref{bic4}) and 
(\ref{bic3}). Consider the exterior derivative of eq. (\ref{omb}):
\eq
d(\om^i a)=d[ (\f{i}{j} * a) \om^j]. \label{oma}
\en
The left-hand side is equal to:
\eqa
\lefteqn{
d(\om^ia)=}\nonumber\\
& &=d\om^i \we a - \om^i \we da = -\C{jk}{i} \om^j \otimes \om^k a - 
\om^i \we (\chi_j * a)\om^j=\nonumber\\
& &=-\C{jk}{i}\om^j \otimes \om^k a - (\f{i}{s} * \chi_j * a)\om^s \we 
\om^j=\nonumber\\
& &=-\C{jk}{i} (\f{j}{p}* \f{k}{q} * a) \om^p \otimes \om^q - (\f{i}{s}  
*\chi_j * a) (\om^s \otimes \om^j - \La^{sj}_{~pq} \om^p \otimes \om^q) =
\nonumber\\
& &=[(-\C{jk}{i} \f{j}{p}  \f{k}{q} - \f{i}{p}  \chi_q + \La^{sj}_{~pq} 
\f{i}{s} \chi_j)*a](\om^p\otimes \om^q) \label{A1}
\ena
The right-hand side reads:
\eqa
\lefteqn{d[(\f{i}{j} * a)\om^j] =}\nonumber\\
& &=d(\f{i}{j} * a)\om^j + (\f{i}{j} *a) d\om^j=\nonumber\\
& &=(\chi_k* \f{i}{j} * a)\om^k \we \om^j - 
(\f{i}{j} * a) \C{pq}{j} \om^p 
\otimes \om^q =\nonumber\\
& &=(\chi_k *\f{i}{j} * a)(\om^k \otimes \om^j - 
\La^{kj}_{~pq} \om^p \otimes \om^q)-(\f{i}{j} * a) \C{pq}{j} \om^p 
\otimes \om^q =\nonumber\\
& &= [(\chi_p \f{i}{q} - \La^{kj}_{~pq} \chi_k \f{i}{j} - \C{pq}{j}
\f{i}{j})*a](\om^p \otimes \om^q), \label{A2}
\ena
\noi so that we deduce the equation
\eqa
& &-\C{jk}{i} \f{j}{p}  \f{k}{q} - \f{i}{p}  \chi_q + \La^{sj}_{~pq} 
\f{i}{s} \chi_j=\nonumber\\
& &=\chi_p \f{i}{q} - \La^{kj}{~pq} \chi_k \f{i}{j} - \C{pq}{j}
\f{i}{j}. \label{A3}
\ena

We now need two lemmas.
\sk
{\sl Lemma 1}
\sk
\eq
\f{n}{l} * a \theta = (\f{n}{r} *a)(\f{r}{l} * \theta),~~~a\in A,~\theta
\in \Ga^{\otimes n}. \label{A6}
\en

\noi{\sl Proof}
\eqa
\lefteqn{\f{n}{l} * a\theta=}\nonumber\\
& &(id \otimes \f{n}{l})\D (a) \DR (\theta)= a_1 \theta_1 \f{n}{l}
(a_2 \theta_2)=\nonumber\\
& &a_1 \theta_1 \f{n}{r}(a_2) \f{r}{l}(\theta_2)=
a_1\f{n}{r} (a_2)\theta_1 \f{r}{l} (\theta_2)=\nonumber\\
& &(\f{n}{r} * a) \theta_1 \f{r}{l} (\theta_2) = (\f{n}{r} * a)
(\f{r}{l} *\theta). \label{A7}
\ena
\sk
{\sl Lemma 2}
\eq
\f{r}{l} * \om^j = \La^{rj}_{~kl} \om^k. \label{A8}
\en

\noi{\sl Proof}
\eqa
\lefteqn{\f{r}{l} * \om^j=}\nonumber\\
& &(id \otimes \f{r}{l})\DR (\om^j)= (id\otimes \f{r}{l})[\om^k \otimes 
\M{k}{j}]=\nonumber\\
& &=\om^k \f{r}{l} (\M{k}{j})=\La^{rj}_{~kl}. 
\ena
\sk

Consider now eq. (\ref{dhthe}) with $h=\f{n}{l}$:
\eq
d(\f{n}{l} *a)=\f{n}{l} *da. \label{A9}
\en
\noi The first member is equal to $(\chi_k*\f{n}{l}*a)\om^k$, while the 
second member is:
\eqa
\f{n}{l} * da &=& \f{n}{l} * [(\chi_j * a)\om^j]=(\f{n}{r} * 
\chi_j *a)(\f{r}{l} * \om^j) \nonumber\\
&=&(\f{n}{r} * \chi_j * a)( \La^{rj}_{~kl} \om^k)
\ena
\noi We have used here the two lemmas (\ref{A6}) and (\ref{A8}). 
Therefore
the following equation holds:
\eq
\chi_k * \f{n}{l}= \La^{rj}_{~kl} ~\f{n}{r} * \chi_j,  \label{A10}
\en
\noi which is just eq. (\ref{bic4}). Equation (\ref{bic3}) is obtained 
simply by subtracting (\ref{A10}) from eq. (\ref{A3}).
\sk

\sect{Two theorems on $i_t$ and $\ell_t$}
\sk

\noi {\bf Theorem} 2.4.11 $~$
\[
\ell_{t_{i}}= i_{t_{i}} d +d i_{t_{i}} 
\]
that is
\begin{eqnarray}
\label{54}
\forall \;
a_{ i_{1} \ldots i_{n} }
\omega^{i_{1}} \wedge \ldots \omega^{i_{n}} \; \in \Gamma^{\wedge n },
 & & \nonumber \\
\ell_{t_{i}} 
( a_{ i_{1} \ldots i_{n} }
\omega^{i_{1}} \wedge \ldots \omega^{i_{n}} )
&=&  
( i_{t_{i}} d +d i_{t_{i}} ) 
( a_{ i_{1} \ldots i_{n} }
\omega^{i_{1}} \wedge \ldots \omega^{i_{n}} ).
\end{eqnarray}

We will show this theorem by induction on the integer $n$. To do this, 
we need the following:
\sk
\noi {\sl Lemma}

If $n=1$, the theorem is true, i.e.
\begin{equation}
\label{540}
\ell_{t_{i}} ( b_{k} \omega^{k} ) = 
( i_{t_{i}} d +d i_{t_{i}} ) 
( b_{k} \omega^{k} ). 
\end{equation} 
First we show that:
\begin{equation}
\label{541}
\ell_{t_{i}} ( \omega^{k} ) = 
( i_{t_{i}} d +d i_{t_{i}} ) 
\omega^{k}. 
\end{equation} 
We already know that $\ell_{t_i}(\om^k)=\om^j\C{ji}{k}$. The 
right-hand side of (\ref{541}) yields:
\[\begin{array}{lcl}
( i_{t_{i}} d +d i_{t_{i}} )  ( \omega^{k} ) & = & 
 i_{t_{i}} d\omega^{k}  +d (i_{t_{i}} \omega^{k} ) = \\ 
&=&-{1\over 2}C_{nj}{}^{k}i_{t_i}(\omega^n\wedge\omega^j)=\\
& = & 
-{1\over 2}C_{nj}{}^{k} \left(
f^n{}_{i} * \omega^j  - \omega^n \delta_{i}^{j} \right) = \\
& = & 
-{1\over 2}C_{nj}{}^{k} \left[ \left( id \otimes
f^n{}_{i}  \right) \left( \omega^{\ell} \otimes M{_\ell}^{j}   \right)  
- \delta_{i}^{j} \omega^n \right] = \\
& = & 
-{1\over 2}C_{nj}{}^{k} \left[ 
\omega^{\ell} \La^{nj}_{~\ell i}
- \delta_{i}^{j} \omega^n \right] = \\
& = & 
+ {1\over 2}C_{nj}{}^{k} \left[ 
\delta_{\ell}^{n} \delta_{i}^{j} 
-\La^{nj}_{~\ell i}
\right] 
\omega^{\ell} =  \\
& = & \C{\ell i}{k} \omega^{\ell} 
\end{array}
\]
and (\ref{541}) is thus proved.

The right-hand side of (\ref{540}) gives:
\[
\begin{array}{lcl}
( i_{t_{i}} d +d i_{t_{i}} ) ( b_{k} \omega^{k} ) & = &
i_{t_{i}} \left( d b_{k} \wedge \omega^{k} + b_{k} d \omega^{k} \right) +
d \left( b_k i_{t_{i}} (\omega^k) \right) = \\
& = &
i_{t_{j}} ( d b_{k} ) f{^j}_i * \omega^{k} -
(d b_{k}) i_{t_{i}} (\omega^{k}) + \\
& &+ b_{k} i_{t_{i}} (d \omega^{k}) +
(d b_{k}) i_{t_{i}} (\omega^{k}) = \\
& = &
i_{t_{j}} ( ( \chi_n  * b_{k} ) \omega^n ) f{^j}_i * \omega^{k} +
 b_{k} i_{t_{i}} (d \omega^{k}) = \\
& = &
( \chi_n  * b_{k} ) \delta_j^n f{^j}_i * \omega^{k} +
b_{k} ( i_{t_{i}} d + d i_{t_{i}} )
\omega^{k} = \\
& = &
( \chi_n  * b_{k} ) f{^n}_i * \omega^{k} +
b_{k}  \ell_{t_{i}} (\omega^k ) = \\
& = &
\ell_{t_{n}} (b_{k} )
f{^n}_i * \omega^{k} +
b_{k}  \ell_{t_{i}} (\omega^k ) = \\
& = &
\ell_{t_{i}} (b_k \omega^k ),
\end{array}
\]
\noi and the lemma is proved. We 
now finally prove the theorem.\\
Let us suppose it to be true for a ($n-1$)-form:
\begin{equation}
\ell_{t_{a}} 
( a_{ i_2 \ldots i_n}
\omega^{i_{2} } \wedge  \ldots \omega^{i_{n}} ) =
( i_{t_{i}} d + d i_{t_{i}} ) 
( a_{ i_2 \ldots i_n}
\omega^{i_{2} } \wedge  \ldots \omega^{i_{n}} ).\nonumber 
\end{equation}
Then it holds also for an $n$-form. Indeed,
the left-hand side of (\ref{54}) yields
\[
\begin{array}{l}
\ell_{t_{i}} 
( a_{ i_1 \ldots i_n}
\omega^{i_{1} } \wedge  \ldots \omega^{i_{n}} ) = \\
= \ell_{t_{j}} 
( a_{ i_1 \ldots i_n}
\omega^{i_{1}})  \wedge f{^j}_i * 
( \omega^{i_{2}} \wedge  \ldots \omega^{i_{n}} ) + 
 a_{ i_1 \ldots i_n}
\omega^{i_{1}}  \wedge \ell_{t_{i}} 
( \omega^{i_{2}} \wedge  \ldots \omega^{i_{n}} ) 
\end{array}
\]
\noi while the right-hand side of (\ref{54}) is given by :\\
\begin{eqnarray*}
\lefteqn{( i_{t_{i}} d +d i_{t_{i}} ) 
( a_{ i_{1} \ldots i_{n} }
\omega^{i_{1}} \wedge \ldots \omega^{i_{n}} ) = } \\
 = &&
i_{t_{i}} [ d
( a_{ i_{1} \ldots i_{n} } \omega^{i_1}) \wedge 
\omega^{i_{2}} \wedge \ldots \omega^{i_{n}} - 
( a_{ i_{1} \ldots i_{n} } \omega^{i_1}) \wedge 
d ( \omega^{i_{2}} \wedge \ldots \omega^{i_{n}} )] + \\
& &
d [ i_{t_{j}} 
( a_{ i_{1} \ldots i_{n} } \omega^{i_1}) 
f{^j}_i * 
(\omega^{i_{2}} \wedge \ldots \omega^{i_{n}}) - 
( a_{ i_{1} \ldots i_{n} } \omega^{i_1}) \wedge 
i_{t_{i}} 
(\omega^{i_{2}} \wedge \ldots \omega^{i_{n}}) ] = \\ 
=& &
 i_{t_{j}} ( d 
a_{ i_{1} \ldots i_{n} } \omega^{i_1}) \wedge 
f{^j}_i * 
(\omega^{i_{2}} \wedge \ldots \omega^{i_{n}}) +
d ( a_{ i_{1} \ldots i_{n} } \omega^{i_1}) \wedge 
i_{t_{i}} 
(\omega^{i_{2}} \wedge \ldots \omega^{i_{n}}) + \\ 
& &
- i_{t_{j}} 
(a_{ i_{1} \ldots i_{n} } \omega^{i_1})
f{^j}_i * 
d (\omega^{i_{2}} \wedge \ldots \omega^{i_{n}}) +
 a_{ i_{1} \ldots i_{n} } \omega^{i_1} \wedge 
i_{t_{i}} 
d (\omega^{i_{2}} \wedge \ldots \omega^{i_{n}}) + \\ 
& &
+ d i_{t_{j}} 
(a_{ i_{1} \ldots i_{n} } \omega^{i_1})
f{^j}_i * 
(\omega^{i_{2}} \wedge \ldots \omega^{i_{n}}) +
i_{t_{j}} ( a_{ i_{1} \ldots i_{n} } \omega^{i_1}) \wedge 
f{^j}_i *
d (\omega^{i_{2}} \wedge \ldots \omega^{i_{n}}) + \\ 
& &
- d(a_{ i_{1} \ldots i_{n} } \omega^{i_1})
\wedge i_{t_{i}}  
(\omega^{i_{2}} \wedge \ldots \omega^{i_{n}}) +
 a_{ i_{1} \ldots i_{n} } \omega^{i_1} \wedge 
d i_{t_{i}} 
(\omega^{i_{2}} \wedge \ldots \omega^{i_{n}}) = \\ 
 =& &
[ ( i_{t_{j}} d + d i_{t_{j}} )
( a_{ i_{1} \ldots i_{n} } \omega^{i_1} ) ] \wedge
f{^j}_i   * 
(\omega^{i_{2}} \wedge \ldots \omega^{i_{n}}) +\\
&&+a_{ i_{1} \ldots i_{n} } \omega^{i_1} 
( i_{t_{i}} d + d i_{t_{i}} )
(\omega^{i_{2}} \wedge \ldots \omega^{i_{n}}) = \\
=& &
\ell_{t_j} 
( a_{ i_{1} \ldots i_{n} } \omega^{i_1} )\wedge
f{^j}_i *
(\omega^{i_{2}} \wedge \ldots \omega^{i_{n}}) +
 a_{ i_{1} \ldots i_{n} }\omega^{i_1} \wedge
\ell_{t_i} 
(\omega^{i_{2}} \wedge \ldots \omega^{i_{n}})
\end{eqnarray*}
and the theorem is proved.

\cvd

\noi{\bf Lemma }
$\left[i_{t_i},\ell_{t_k} \right]=
i_{t_i}{\scriptstyle{{}^{{}_{\circ}}}}\ell_{t_k} -\La_{~ik}^{ef}
\ell_{t_e}{\scriptstyle{{}^{{}_{\circ}}}}i_{t_f}$
\sk
\noi{\sl Proof } By definition we have
\eq
\left[i_{t_i}, \ell_{t_k}, \right](\vart)\equiv
\ell_{\kappa'{({\chi_k}_{1})*}}{\scriptstyle{{}^{{}_{\circ}}}}
i_{t_i}{\scriptstyle{{}^{{}_{\circ}}}}\ell_{{\chi_k}_2 *}\,(\vart)=
\kappa'(\chi_j)*(i_{t_i}(f^j{}_k*\vart))+
i_{t_i}( \ell_{t_k}(\vart))
\en
we therefore have to prove that
\eq
\kappa'(\chi_j)*(i_{t_i}(f^j{}_k*\vart))=
-\La_{~ik}^{ef}\ell_{t_e}{\scriptstyle{{}^{{}_{\circ}}}} i_{t_f}\label{lemma}
\en
First notice that on a generic covariant tensor $\tau\in\Ga^{\otimes}$
\eq
i_{t_i}(f^j{}_k*\tau)=\La_{~ik}^{ef}f^j{}_e*i_{t_f}(\tau)\label{f*i}
\en
as can be easily proved by induction with $\tau=\tau'\otimes\om^{\ell}$, 
$\tau'\in\Ga^{\otimes}$.
To complete the proof recall that 
$\kappa(\chi_j)f^j{}_e=-\chi_e\;$ [apply $m(\kappa'\otimes id)$
to $\D'(\chi_e)$]. \cvd
\sk

\noi{\bf Theorem} $~$ 
$\left[i_{t_i},\ell_{t_k} \right] = i_{[t_i,t_k]}$=$ \C{ik}{l}i_{t_l}$ 
\sk
\noi {\sl Proof} 

\noi The proof is by induction, it holds on $1$-forms, 
let assume it holds for  a generic $\vart$ form of order $n$, we prove
it holds also for the generic $n+1$ form $\om^a\wedge\vart$.
Use the previous lemma to rewrite  $\left[i_{t_i},\ell_{t_k} \right]$ as 
(1)+(2):
\eqa
\mbox{(1)}&&i_{t_i}\ell_{t_j}(\om^a\wedge\vart)=
i_{t_i}[(\om^b\C{bp}{a}\wedge f^p{}_j*\vart+\om^a\wedge\ell_{t_j}(\vart)]
\nonumber\\
&&~~=[\C{bp}{a}f^b{}_if^p{}_j+f^a{}_i\chi_j]*\vart-\om^b\wedge
\C{bp}{a}i_{t_i}(f^p{}_j*\vart)-\om^a \wedge i_{t_i}\ell_{t_j}(\vart)
\nonumber\\[1em]
\mbox{(2)}&&-\La^{kl}_{~ij}\ell_{t_k}i_{t_l}(\om^a\wedge\vart)=
-\La^{kl}_{~ij}\ell_{t_k}[f^a{}_l*\vart-\om^a i_{t_l}(\vart)]\nonumber\\
&&~~=[-\La^{kl}_{~ij}\chi_kf^a{}_l]*\vart +
\La^{kl}_{~ij}\om^b\wedge\C{bp}{a}f^p{}_k*i_{t_l}(\vart) -
\om^a\wedge (-1)\La^{kl}_{~ij}\ell_{t_k}i_{t_l}(\vart)\nonumber\\[1em]
\mbox{(3)}&&\C{ij}{k}i_{t_k}(\om^a\wedge\vart)=\C{ij}{k}f^a{}_k*\vart
-\C{ij}{k}\om^a\wedge i_{t_k}(\vart)\nonumber\\
\ena
We then have
\eqa
\mbox{(1)$+$(2)$-$(3) }&=&[\C{bp}{a}f^b{}_if^p{}_j+f^a{}_i\chi_j
-\La^{kl}_{~ij}\chi_kf^a{}_l-\C{ij}{k}f^a{}_k]*\vart\\
&&-\om^b\wedge\C{bp}{a}[i_{t_i}(f^p{}_j*\vart) 
-\Lambda^{kl}_{~ij}f^p{}_k*i_{t_l}(\vart)]\\
&&-\om^a\wedge[i_{t_i}\ell_{t_j}-\Lambda^{kl}_{~ij}
\ell_{t_k}i_{t_l}-\C{ij}{k}i_{t_k}](\vart)\\
&=&0
\ena
Where the first addend is zero because of (\ref{bico3}), the second is zero 
because of (\ref{f*i}), the third  because of the  inductive hypothesis.
\cvd

\vfill\eject
\end{document}